\documentclass[10pt,num,fleqn]{elsarticle}

\newcommand{\wrt}{\textit{w}.\textit{r}.\textit{t}. }
\newcommand{\ie}{\textit{i}.\textit{e}., }
\newcommand{\eg}{\textit{e}.\textit{g}.\ }
\newcommand{\etc}{\textit{etc}.}

\usepackage[a4paper, margin=0.75in]{geometry}  % Adjust margin size

\usepackage{soul} 

\usepackage[utf8]{inputenc}
\usepackage{tikz}
\usetikzlibrary{shapes.geometric, arrows, positioning}
\usetikzlibrary {shapes.misc}
\usepackage{lscape}
\usepackage{comment}
\usepackage{mathrsfs}

\usepackage{rotating}
\usepackage{tabularx}
\usepackage{subcaption}
\usepackage{blindtext}

\usepackage{multicol}

\usepackage{listings}
\usepackage{booktabs}
\usepackage{multirow}

\usetikzlibrary {arrows.meta,bending, positioning}

\usepackage[rightcaption]{sidecap}

\usepackage{graphicx} %package to manage images

\usepackage{textcomp}

\usepackage{underscore} %https://chatgpt.com/share/85479189-2a22-4a4a-bfd3-22e8de690eba
\usepackage{amssymb} %% The amssymb package provides various useful mathematical symbols
\usepackage{hyperref} % For clickable references

\usepackage{orcidlink}
\usepackage{enumitem}
\usepackage{amsmath}
\usepackage{cleveref} % For handling multiple references better %https://chatgpt.com/share/472aa5a2-a470-47dd-b2c5-aba4e99a429c
\usepackage{caption} % This is important as it gives a minimum space to the captions above

\usepackage{float} %https://www.youtube.com/watch?v=k2ntGEGelvA&t=54s

\usepackage{MnSymbol,bbding,pifont} %https://www.physicsread.com/latex-star-symbol/
\usepackage{longtable} %https://stackoverflow.com/questions/2896833/how-to-stretch-a-table-over-multiple-pages
\usepackage{mathtools} %https://tex.stackexchange.com/questions/103988/rightarrow-with-text-above-it
\usepackage{array} % for vertical centering with m{} column specifier

\newcommand{\instbox}[1]{\tikz[baseline=(X.base)] \node[draw, rectangle, rounded corners, inner sep=2pt](X) {#1};} %✅ 2. Define a custom macro for the boxed label

\usepackage{pdflscape} % To use resize box for table resizing
\usepackage{pbox} % For ensuring neat text wrapping with a table's cell
\usepackage{algorithm}
\usepackage{algpseudocode}

\usepackage{algorithmicx}
\algdef{SE}[SEARCH]{Search}{EndSearch}[1]{\textbf{Search} #1:}{\algorithmicend\ \textbf{Search}}

\algrenewcommand{\algorithmicdo}{}

\usepackage{makecell} % seems necessary to break multicolumn cells and wrap them

\usepackage{xr}
\usepackage{fancyhdr}

\fancypagestyle{myfancy}{
  \fancyhf{}
  \fancyhead[L]{\scriptsize\textit{S. Banerjee et al.}}
  \fancyhead[R]{\scriptsize\textit{Preprint available at https://arxiv.org/abs/2509.13227}}
  \fancyfoot[C]{\small \thepage}
}

\fancypagestyle{titlepagefancy}{
  \fancyhf{}
  \fancyfoot[L]{\small\textit{Pre-print version II of: \href{https://arxiv.org/abs/2509.13227}{arXiv:2509.13227}}}
  \fancyfoot[C]{\small\thepage}
  \fancyfoot[R]{\small\textit{Date: \today}}
}

\makeatletter
\let\ps@plain\ps@myfancy            % all normal pages use myfancy
\let\ps@pprintTitle\ps@titlepagefancy % title page uses no-header style
\makeatother

\begin{document}

\begin{frontmatter}



\title{Rich Vehicle Routing Problem in Disaster Management enabling Temporally-causal Transhipments across Multi-Modal Transportation Networks}


\author[1]{\href{https://www.linkedin.com/in/santanu-banerjee-093929150/}{Santanu Banerjee}~\orcidlink{0000-0001-9861-7030}\corref{cor1}}\ead{santanu@kgpian.iitkgp.ac.in}
\author[1]{Goutam Sen}\ead{gsen@iem.iitkgp.ac.in}
\author[1]{Siddhartha Mukhopadhyay~\orcidlink{0009-0005-5764-7299}}\ead{sid.mpadhyay@kgpian.iitkgp.ac.in}

\cortext[cor1]{Corresponding author}

\affiliation[1]{organization={Department of Industrial and Systems Engineering (ISE), Indian Institute of Technology (IIT) Kharagpur},%Department and Organization
            %addressline={}, 
            city={Kharagpur},
            postcode={721302}, 
            state={West Bengal},
            country={India}}

\begin{abstract}
A rich vehicle routing problem is considered allowing multiple trips of heterogeneous vehicles stationed at geographically distributed vehicle depots having access to different modes of transportation. The problem arises from the real world requirement of optimizing the disaster response time by minimizing the makespan of vehicular routes. Multiple diversely-functional vertices are considered, including Transhipment Ports as inter-modal resource transfer stations. Both simultaneous and split pickup and delivery is considered, for multiple cargo-types, along with Vehicle-Cargo and Transhipment Port-Cargo compatibilities. The superiority of the proposed cascaded minimization approach is demonstrated over the existing makespan minimization approaches through our developed Mixed-Integer Linear Programming formulation. To solve the problem quickly for practical implementation in a Disaster Management-specific Decision Support System, an extensive Heuristic Algorithm is devised which utilizes Decision Tree based structuring of possible routes; the Decision Tree approach helps to inherently capture the compatibility issues, while also explore the solution space through stochastic weights. Preferential generation of small route elements are performed, which are integrated into route clusters; we consider multiple different logical integration approaches, as well as shuffling the logics to simultaneously produce multiple independent solutions. Finally perturbation of the different solutions are done to find better neighbouring solutions. The computational performance of the PSR-GIP Heuristic, on our created novel datasets, indicate that it is able to give good solutions swiftly for practical problems involving large integer instances which the MILP is unable to solve.
\end{abstract}



\begin{highlights}
\item Multi-Trip rich-VRP with multiple Depots, Warehouses, Relief Centres and Nodes
\item Ports allow delayed inter-modal resource transfer across Transportation Networks
\item Diverse Pickup and Delivery Cargos having Vehicle and Port compatibilities
\item MILP and Heuristic developed for Cascaded Makespan Minimization Objective
\item Applications in post-disaster Relief-and-Rescue and in pre-disaster Evacuations
\end{highlights}

\footnotesize

\end{frontmatter}

% \small

% Refer to this for understanding EJOR Format and Supplementary Materials: https://www.sciencedirect.com/science/article/pii/S0377221723002333

% Refer to this for general Instructions: https://www.elsevier.com/wps/find/journaldescription.cws_home/505543?generatepdf=true

% Refer to this for Latex Instructions: https://www.elsevier.com/authors/policies-and-guidelines/latex-instructions

% Refer to this for Policies: https://www.elsevier.com/authors/policies-and-guidelines/documents/elsdoc-1.pdf

\small
\section{Introduction} \label{Introduction}

Route planning is an essential part of relief and rescue missions, and optimization of vehicle routes helps in ensuring effective disaster response. One of the major needs in emergency/disaster missions is to fully satisfy all delivery resources requested (say relief kits) from certain geographical locations and pick up stranded individuals (including stray animals, cattle and livestock) stuck at majorly affected places. In many cases, this pickup and delivery must be necessarily performed together, allowing the entire pickup and/or delivery requirement to be satisfied in one go and without load splitting. In other cases, load splitting could be preferred allowing an overall better approach to satisfy the requirements raised at the different geographical locations affected by the disaster. Furthermore, in the case of a preemptive approach towards disaster preparedness and mitigation of associated risks of an incoming disaster, similar pickup and delivery requirements arise from the geographies where some natural disaster is being forecasted.

Often, such emergency missions encompass broad geographies involving multiple modes of transportation. Therefore, transhipment considerations must be taken into account while planning for these missions across transhipment points that are connected to more than one transportation system acting as temporary transfer hubs. These transhipment points are many a time associated with complex compatibilities, allowing only certain resources to pass through, while the transfer of some other type of resource may be restricted.

Furthermore, the locations of resource storage (like water bottles, plastic sheets, medicines, relief kits, \etc) vehicles and final locations for 
Generally the delivery resources in such cases include water bottles, plastic sheets, medicines, relief kits \etc which may need to be sourced from different locations once the exact demand is known. Similarly, the rescued individuals may be sent to hospitals or shelters depending on their condition (post disaster) or based on the predicted intensity/severity of an incoming disaster (pre-disaster phase). These necessitate the development and usage of formulations and heuristics which consider different functionalities of vertex (in the graph, generally used as depots/nodes). Conventional vehicle routing problems (VRPs) and their subsequent enhancements have generally considered a single node element to have all the three features of stationing vehicles, storing resources, as well as the final delivery location of collected pickups. This may never be the case in real-life disaster-specific VRPs as disasters are generally a one-time occurrence and stationing and maintaining large amounts of resources at one single place only for providing relief during the aftermath of disasters or bringing valuable lives at one safe place during the preparation phase of an incoming disaster may not always be a good approach and may be practically infeasible due to the costs involved. Instead, a decentralized approach with robust and swift transfer capability is required; in this case, existing warehouses (WHs) could reserve some specific amount of different types of resources to be used during emergencies allowing for reduced maintenance costs since the WHs would be able to function normally as well. Similarly, Vehicle Depots (VDs) could reserve a limited number of (rotating) vehicles to be used during emergencies while able to use a wider fleet for other different operations allowing for better maintenance and regular upkeep of working vehicles. In any case, we believe a decentralised distributed network with regular operational points (like vehicle depots, warehouses, relief centres) would help create a wider, deeper and leaner infrastructure network hat will be able to leap into action before, during and after a disaster has struck.

Keeping these considerations in mind, we conduct this study, reviewing the literature (\autoref{Previous Works and State of Art}), defining the detailed featured of the novel problem identified (\autoref{Problem Description}), developing a exact Mixed-Integer Linear Program (MILP) formulation (\autoref{Mathematical Formulation}) along with a heuristic algorithm (\autoref{Heuristic Development}), and finally perform computational analysis on small and large datasets developed by us (\autoref{Discussion and Computational Study}, along with some instances from the literature). We consider a multi-trip multi-modal capacity-and-compatibility constrained VRP with a network consisting of multiple types of vertices each with their own unique function. Our problem is extremely flexible and can be fit in cases where two functionalities need to overlap (like in the traditional sense of multi-depot VRPs), by creating two (or more) separate functionally different vertices in the same geographical location. A list of abbreviations used in this study is available in \ref{List of Abbreviations}.

\section{Previous Works and State of Art} \label{Previous Works and State of Art}
We conduct a survey of existing literature related to our problem; we start with specific keyword searches related to our problem's features and tabulate the other features considered by the authors. To keep this section brief, we present our extensive literature review of individual papers in Table S1 (presenting this single table by breaking its horizontal extensions into A, B, C and D), and mainly discuss the review papers available in the literature in the subsequent text. We also refer the readers to \cite{KRAMER2019162} for an interesting overview of associated problems in the rich-VRP domain.

\subsection{Related works focusing on Rich VRP}
In \cite{Simeonova2020}, multiple optimization objectives are discussed multi-objectives and pareto-optimalities. They review papers discussing Genetic Algorithms, Adaptive Large Neighbourhood Search, Iterated Local Search. The network features encompassed within their papers include VRPs considering closed and open tours, time windows and multi-modal transportation, while the considered vehicle features of reviewed papers include homogeneous and heterogeneous fleets performing simultaneous and/or split tasks.

The concept of Road Capacity Reduction due to disaster or damaged roads, congestion or traffic, is discussed in \cite{EKSIOGLU20091472}, \cite{CORONAGUTIERREZ2022108054} and \cite{VIDAL2020401}. The authors in \cite{EKSIOGLU20091472} discuss about the multi-objective transportation network design proposed in \cite{CURRENT19934}, where the authors discuss about the multi-criteria transhipment problem with multiple objectives proposed in \cite{OGRYCZAK198953} for the facility location problem. We don't consider dynamic cases allowing real-time road condition alterations due to ongoing disaster effects; in our case, we consider that the shortest distance (updated with the predictions of available roads throughout the disaster) across nodes will be used to develop the static travel time matrix between geographical points which are used throughout our problem. These shortest distances could be pre-processes with the predictive disaster data and made to avoid disaster-prone areas. Therefore, implementing such a problem as ours within disaster management systems requires this pre-requisite of intelligent pre-processing for the case when our formulation and algorithm is applied in the pre- or in- disaster phase; for post-disaster cases, our study contributions would be more robust without the need for any predictive disaster-impact on the transportation network, however in this case, the transportation network needs to be updated to the post-disaster situation where some transportation links could be broken. 

We believe that the optimization sense in the context of VRP plays a crucial role, especially while differentiating between cost, distance, or time minimization, and whether these are being optimized as a minimization problem of their sum or minimization of the highest solution characteristics for among the vehicles. Since among the resources of time, cost and distance, only time is possible to be parallelized, the cascaded minimization of vehicle travel time should be included among the prominent objectives which should be considered or integration within the VRP taxonomy aspects as presented in the Fig. 6 in \cite{EKSIOGLU20091472}.

The survey papers of \cite{CORONAGUTIERREZ2022108054} and \cite{VIDAL2020401} both review papers that have considered constraints on the number of trips of a single vehicle, VRP with Time Windows, VRP with homogenous and heterogeneous fleet with or without limited number of capacitated vehicles, vehicles performing split simultaneous delivery, as well as the single or multiple depot problem. Some extended problem features of constraining the maximum route time, constraining the maximum number of customers per route, and discussions about having separate vehicle (capacitated) compartments for different types of cargo (like for tankers carrying multiple types of fluids at different pickup stages) are also highlighted in \cite{CORONAGUTIERREZ2022108054}; the authors also dwelve into the optimization sense, discussing about objective functions like cumulative objectives (arc summation till last node), weighted distance with load quantity, as well as priority indexes. Additionally, \cite{VIDAL2020401} describes literature which have considered many-to-many transportation, papers allowing the multi-visit of nodes as necessary, and studies on multi-modal networks; many papers in their study include real-world problem features of profitability, service quality, equity, consistency, simplicity, reliability, and externalities.

The authors in \cite{10.1145/2666003} perform a survey on the rich VRP literature available. They mention rich-VRPs as emerging from VRPs applied in the real world, and especially now more than ever due to the advent of advanced computing infrastructure. They discuss and review papers having considered cross-docking, multi-trips, vehicle site dependance, vehicle road dependance, intra route replenishments, and multiple time windows among other features.

% \subsection{Humanitarian Logistics Literature}
% For the humanitarian logistics papers, our reviews may be found below:

Based on these review papers, and the individual studies, we conclude that the problem considered by is novel and has never been considered before. Specifically, as per our knowledge, no one in literature has considered transhipment with time consideration as we are considering, where time causality comes into the picture; therefore this should be the first study considering the development of exact and heuristic algorithms with causality considerations. Furthermore, we are considering multi-modal aspects which further complicate the problem. Another novety lies in our cascaded makespan minimization approach, which we propose to be superior (more correct, allows deeper optimization not reachable by the single objective function of makespan) than conventional makespan minimization. We next describe the problem in detail, develop the formulation and heuristic algorithm, and prove our claim through computational results.

\subsection{Applications of VRP in disaster relief}
VRP plays an important role in disaster management; optimized routing of vehicles helps reach victims swifter allowing faster rescues, as well as enable holistic distribution of relief materials during large scale natural disasters. A multi-depot VRP is discussed in \cite{VIEIRA2021102193} for the case of drought relief to solve emergency water trucking; the authors hybridize metaheuristics and demonstrate applications within a case study in Brazil. The authors develop a two-phase multi-period location-routing-inventory model in \cite{VAHDANI2018290} focusing on resilient distribution of relief during the aftermath of an earthquake developing the model for uncertain conditions using a robust optimization approach; they test two meta-heuristic algorithms for this problem.
The study in \cite{XU2024110506} discuss about the problem of pre-disaster evacuation network design; the authors optimize the total evacuation time, investment cost, and network congestion degree focusing on equilibrium flow.

In \cite{MOMENI2023104027}, the authors develop an uncertain multi-objective MILP for a VRP with trucks and drones for wild-fire detection; the drones being able to use the trucks as launchpads to reach diverse heights and monitor hard-to-reach areas thereby minimizing the maximum regret (this objective may be further improved by using a cascaded minimization of the dynamic maximum regret). They mention coordinated system of trucks and drones will be effective in reducing the makespan (based on their numerical calculations), citing \cite{Wang2017} where the authors minimize the makespan of a similar VRP with drones. A latest review paper on this truck-drone routing problem \cite{LUO2025116074} also discusses about the frequently used objective of minimizing \ie makespan, but fail to find any indication of cascading the same in the literature; their described modes of extended collaboration is highly applicable in disaster situations especially to survey damage and detect (/predict potentially) stranded victims.

The paper \cite{TUKENMEZ2024104776} discussed how drone monitoring can prevent forest fires by swift detection and mitigation; they introduce the degree of sensitivity of fire-prone areas for the first time and develop matheuristics to provide swift solutions. They refer to \cite{9837784}, where the authors developed ant colony optimization heuristics for multi-depot drone location-routing to detect post-earthquake damage; impacted regions were divided into grids and allocated importance values based on multiple factors including the number of buildings damaged, time of day, \etc
Unmanned aerial VRP for forest-fire surveillance is also discussed in \cite{SADATI2025111087} where the author develops an MILP for the covering VRP and use Bender's decomposition solution strategy.

In \cite{ALDAHLAWI2024104638}, the authors focus on large-scale evacuation planning while conducting a literature review of human behavior modelling and route optimization, revealing agent based simulations are commonly used in such scenarios.
The authors in \cite{callen2025digital} discuss how direct digital payments could be a game-changer as an aid-delivery mechanism during humanitarian crises.
The paper \cite{vandenberg2024cfr} discusses about community-volunteers as first responders for cardiac arrest, which can be extended to optimizing volunteer response during localized-emergencies.
A multi-criteria decision making model is developed in \cite{ABDULVAHITOGLU2024110639} to select the best facility locations for disaster-responder units, and for task optimization; the authors highlight the importance of strategic positioning in minimizing response times and maximizing operational effectiveness.

\section{Problem Description} \label{Problem Description}
The real-world disaster management problem (discussed in \autoref{Introduction}) broadly considers the transfer of different types of delivery resources (relief materials) and pickups (rescue) to and from disaster-affected geographies or locations where a disaster is predicted to strike. In both cases, the entire operation completion time of the emergency mission is of the essence, and our aim in this study is to minimize this total emergency operation time. This corresponds to the optimization of the vehicle route with the highest route duration. We don't find any relevant literature that even comes close to considering all the problem features considered by us at once. Our implementation of transhipment ports (TPs) is novel and very appropriate for modelling real-world scenarios involving multi-modal fleet coordination.

To solve this problem at hand, we consider a rich-VRP with specific arcs (Travel Time edge network pre-processed with Traffic information and reduced network due to disaster) allowing multiple heterogeneous vehicles present in various Vehicle Depots to travel on their specific (modal) network layers for performing pickup of relief supply from various capacity constrained Warehouses to be provided to multiple Nodes as a simultaneous delivery-and-pickup allowing rescue of stranded victims/livestock to be sent to strategically located capacity constrained Relief Centers (RCs). Transhipment Ports are provided for temporary holding and cross-docking of cargo during the entire operation. The problem features are highlighted here:

\begin{enumerate}[label=(\Roman*)]
  \item Travel Time Network across modes considering vehicle-mode compatibility: travel time between same two points could be different among Vehicle Types (VTs) as they might have access to different types of roads/rails \etc,
  
  \item Multiple Types of Pickup as well as Delivery Loads (both Weight and Volume-constrained) are considered, while also taking into account their respective loading/unloading times specific to each VT
  
  \item Functional categorization of the graph vertices (namely: VDs, TPs, WHs, Node Locations, RCs) which are also considered multiple in number, thereby enhancing the model functionality and evolving the fundamental meaning of Depots. We consider multiple capacity-constrained:
    \begin{enumerate}[label=(\alph*)]
    \item WHs containing diverse Delivery-Load Types (DCTs) as relief materials,
    \item RCs as rescue destinations for diverse PickUp-Load Types (PCTs),
    \end{enumerate}

  \item Heterogeneous capacity (volume and weight separately) constrained VTs are present in various numbers across respective VDs.
  
  \item The vehicles can have multiple trips, i.e. visit vertices multiple times (unless not allowed as for the Simultaneous Pickup-and-Delivery Nodes).

  \item  VTs are also associated with Open or Closed routes; open routes indicate the final return-journey from a vertex to the vehicle's respective VD is not considered during trip-time calculation,
  
  \item Allow flexible Hierarchical Transportation across various vehicular networks through introduction of TPs ensuring transhipment happens in a temporally-causal manner.
  
  \item Compatibility of Cargo-Types with TPs are considered such that incompatible loads may be present in the vehicles at the TPs but cannot be transhipped.
  
  \item Waiting times are considered separately for each vehicle for each of its unique visits to a TP.
  
  \item Vehicle to Load preferences are considered such that a VT is allowed to carry only specific CTs. This also ensures the compatibility in cases when certain vehicles are meant to swiftly carry humans and therefore are not allowed to evacuate cattle which must be done by another vehicle,
  
  \item Allow Simultaneous Delivery \& Pickup (this ensures no detachment of loads during rescue missions) at Simultaneous Nodes and Split Delivery and PickUp at Split Nodes.
  
  \item The MILP Formulation minimizes the total operation time (and is flexible to allow choosing of other objectives). A multi-step approach for solving our Cascaded Makespan Minimization objective is shown (where the objective function during every step is updated to be include association with only the vehicles which were not having the largest routes during any of the optimization steps performed before).
\end{enumerate}

Some parts of our problem's features have been considered in different papers in the literature, however we find no single paper coming close to considering all the features of our problem. We understand this could be due to formulation complexity and logical complications during the heuristic development, both of which we have been able to address. Together, the WHs are RCs are termed as the Primary Resource Vertices (PRVs) as they have capacity constrained DCTs and a threshold on the acceptability of PCTs.

\subsection{Illustrative Example} \label{An Example Case Study}
\begin{landscape}
% \begin{figure}
%     \centering
%     \includegraphics[width=1\linewidth]{Picture_Instance.pdf}
%     \caption{\textbf{Case Study Instance}}
%     \label{fig:Picture_Instance}
% \end{figure}
\begin{figure}
    \centering
    \includegraphics[width=1.035\linewidth]{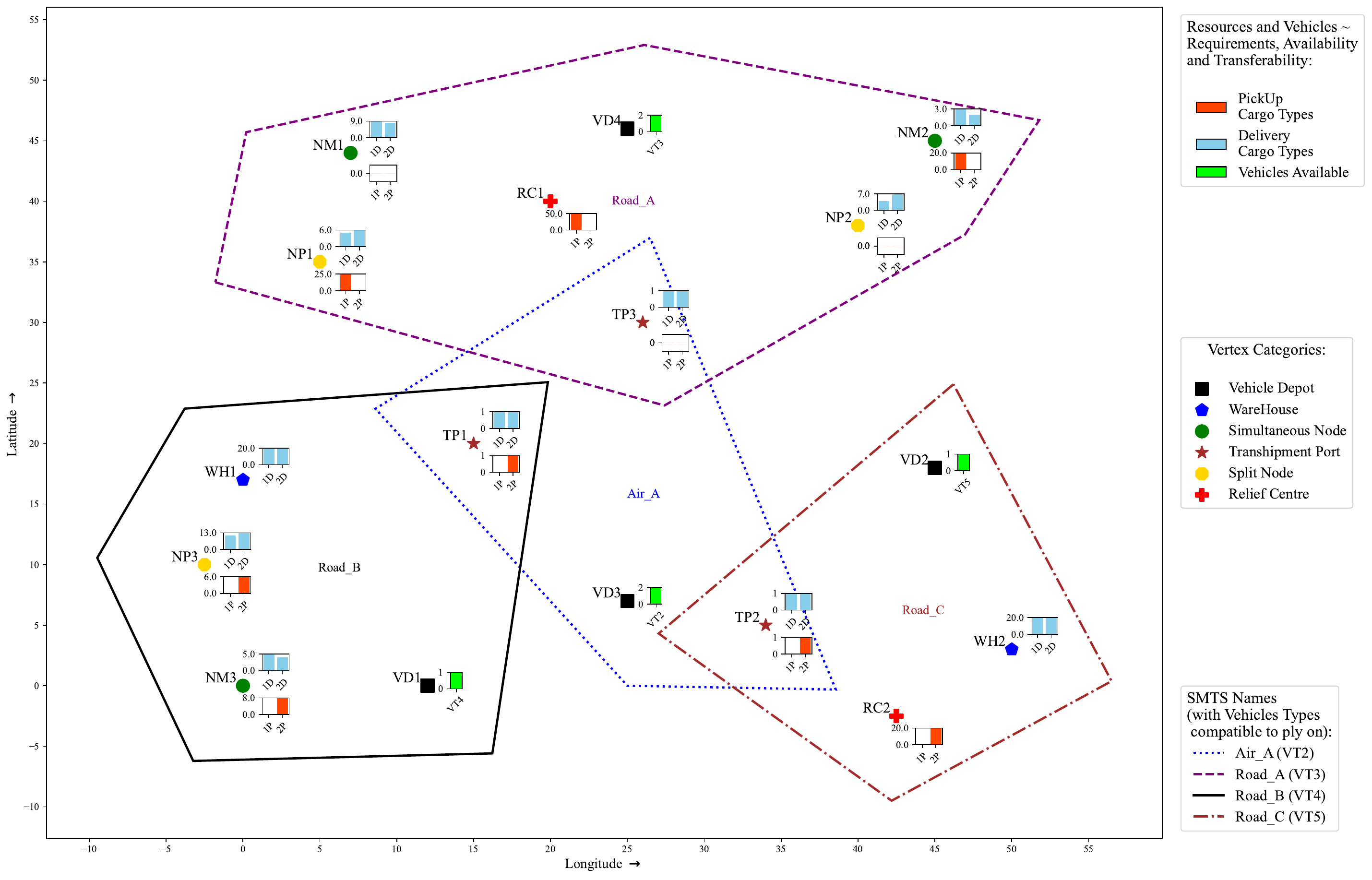} % 0.975 fits
    \caption{\textbf{Illustrative Example} - This MMTN consists of four distinct SMTS, interconnected by 3 TPs having their respective CT-specific transhipment preferences}
    \label{fig:Picture_Instance}
\end{figure}
\end{landscape}

An example problem is shown in \autoref{fig:Picture_Instance} where four different Single-Mode Transportation Segments (SMTSs) form a Multi-Modal Transportation Network (MMTN). Each SMTS has a Vehicle Depot consisting of vehicles of the specific type which can ply on the respective SMTS. Three Transhipment Ports can be observed with one of them allowing the transhipment of only Delivery type resources. The problem is ideated for the case of an archipelago type system having their individual road networks and connected together by an air transportation segment. The full details of this instance is available in the shared datasets.
The bar charts alongside Nodes (Simultaneous or Split Nodes) indicate the request for the Delivery resources (quantity depicted in blue) and the requirement to pickup stranded victims (quantity depicted in orange). The bar-charts alongside Warehouses depict the amount of Delivery resource of each type that is available to be leveraged from the respective Warehouse (and similarly bar-charts alongside each Relief Centre depicts the amount of corresponding PickUp type the specific Relief Centre can accept within itself). For Vehicle Depots, the bar-charts depict the number of vehicles of each type present at the specific depot. For Transhipment Ports, the bar-charts represent the compatibility of resource transfer across that Port (individually for each resource type); $1$ indicates the transfer of that specific resource type is allowed, whereas $0$ indicates transfer not allowed (the incompatible resource can still be within the vehicle, but cannot leave the vehicle).

\section{Mathematical Formulation} \label{Mathematical Formulation}

To solve the problem having the features discussed in \autoref{Problem Description} at hand, a rich vehicle route planning problem (VRP) is considered, having various types of vertices of a graph; these diverse vertex types correspond to different functionalities (Table S2 in Supplementary Material). We initially allow all network linkages, i.e. all possible connection combinations between all vertices clustered within sets mentioned in \autoref{longTab: Sets and Parameters} is allowed. It needs to be evident that Nodes ($N$) are the controlling (\ie demanding) elements of the formulation and a problem without any Nodes will not require any vehicle to ply. Keeping this in mind we remove linkages connecting Vehicle Depots, linkages connecting Warehouses to Vehicle Depots, as well as linkages connecting Vehicle Depot to Relief Centres; since these would never be used for our considered problem.

We consider multiple types of loads being carried by compatible vehicles. Each vehicle can have multiple trips, which is formulated here using a layered approach such that each vehicle starts and ends at its respective vehicle depot while travelling through other vertices as many times as required and allowed, to perform necessary tasks of relief-and-rescue. The  total no. of vehicles allowed to cater to each Simultaneous Node is one. Each individual vehicle ($v$) may only be accurately identified by its corresponding VD, VT, and a unique number from among the number of vehicles of that specific type available at that depot.

% Due to the complexity of the problem proposed, we present the full exact formulation in the supplementary (\autoref{Detailed Mathematical Formulation (Expanded)}) due to its large length; readers are highly recommended to visit this section. We proceed by presenting the results of the case study as found using our exact formulation.

% \section{Detailed Mathematical Formulation (Expanded)} \label{Detailed Mathematical Formulation (Expanded)}

\footnotesize

% \begin{longtable}{|p{0.12\linewidth}|p{0.88\linewidth}|}
\begin{longtable}{|m{0.08\linewidth}|m{0.88\linewidth}|}
\caption{\textbf{Sets and Parameters}\label{longTab: Sets and Parameters}}\\
 % \hline
 % \multicolumn{2}{| c |}{Begin of Table}\\
 \hline
 \textbf{Notation} & \textbf{Description}\\
 \hline
 \endfirsthead

 \hline
 \multicolumn{2}{|c|}{Continuation of \autoref{longTab: Sets and Parameters} \textbf{Sets and Parameters}}\\
 \hline
 \textbf{Notation} & \textbf{Description}\\
 \hline
 \endhead

 % \hline
 % \endfoot

 % \hline
 % \multicolumn{2}{| c |}{End of Table}\\
 % \hline\hline
 % \endlastfoot

 \hline
\hline $H$ & Set of Vehicle Depots (each Vehicle Depot is accessible by all modes of transport which its hosted Vehicles can ply on)\\
\hline $W$ & Set of Warehouses\\
\hline $N^M$ & Set of Simultaneous Nodes (allowing only a single vehicle to perform the pickup and delivery operation at one go)\\
\hline $N^P$ & Set of Split Nodes (allowing any number of vehicles to visit at different times to satisfy portions of the total requirement, which includes pickup or delivery or both types of operations)\\
\hline $N$ & Set of all Nodes, $\ie N^M \cup N^P$\\
\hline $S$ & Set of Transhipment Ports\\
\hline $R$ & Set of Relief Centers\\
\hline $V_0$ & Vertex set consisting of $W \cup N \cup S \cup R$\\

\hline
\hline $W_k$ & Set of all Warehouses that are accessible by (at least one among) the mode(s) of transport on which VT $k$ plies on\\
\hline $N^M_k$ & Set of all Simultaneous Nodes that are accessible by (at least one among) the mode(s) of transport on which VT $k$ plies on\\
\hline $N^P_k$ & Set of all Split Nodes that are accessible by (at least one among) the mode(s) of transport on which VT $k$ plies\\
\hline $N_k$ & Set of all Nodes that are accessible by (at least one among) the mode(s) of transport on which VT $k$ plies on, $\ie N^M_k \cup N^P_k$\\
\hline $S_k$ & Set of all Transhipment Ports that are accessible by (at least one among) the mode(s) of transport on which VT $k$ plies on\\
\hline $R_k$ & Set of all Relief Centres that are accessible by (at least one among) the mode(s) of transport on which VT $k$ plies on\\
\hline $V_k$ & Vertex set consisting of $W_k \cup N_k \cup S_k \cup R_k$\\

%\hline $V$ & Entire vertex set consisting of $VD \cup W \cup N \cup TP \cup RC$\\
%\hline $E$ & Set of all Edges joining all two-vertex permutations from the vertex set $V$. However, a few permutations are not considered since they will never be used during the routing. All self-linkages alongwith the Edges $W\rightarrow VD$, $VD \rightarrow VD$ and $VD\rightarrow RC$ are hence removed from the network.\\
%\hline $G$ & Network graph consisting of $V$ vertices and $E$ edges\\

\hline
\hline $v$ & Refers to each individual vehicle which is referenced as the combination of a vehicle depot $h$, a vehicle type $k$ within that vehicle depot, and the vehicle number $u$ among those vehicle types at the concerned depot, where $h \in H, k \in K_h, u \in G_{h,k}, v=(h,k,u)$.
The Vehicle Depot refers to the originating vertex of the considered vehicle; the Vehicle Type should be among one of the categories of vehicles that are present initially at the concerned vehicle depot, and the Vehicle Number may refer to sequential integers starting from $1$ till a maximum value equal to the number of vehicles of the concerned type present at that depot.\\
\hline $L_v$ & Refers to an array of all the levels $l$ available to vehicle $v$. The levels are indicated by integer numbers starting from $1$ (which is the minimum value, \ie $l_v^{min}$) being the base level, till a maximum value $l_v^{max}$.\\

\hline
\hline $K_i$ & $K_i$ represents the set of all VTs, which can access any of the modes of transport that connect to the Vertex $i$, where $i \in V_0$. In other words, it is the set of all Vehicle Types that can access the Node $i$ via the multi-modal transportation network.\\
\hline $K_h$ & $K_h$ indicates the Vehicle Types present at some specific Vehicle Depot which is referenced through $h$, $h \subseteq H$. Each $H$ has some specific types of vehicles available for performing tasks.\\
\hline $G_{h,k}$ & $G$ is referenced using the two parameters ($h$,$k$), which indicates the vehicles of type $k \in K_h$ available at a specific vehicle depot $h$. It is the entire list of numbers starting from $1$ to the maximum number of vehicles of type $k$ available at $h$.\\
\hline $A_c$ & $A_c$ represents the set of all VTs that are compatible to carry the load type $c$. If $c$ is a set of load types instead of an individual load type (\eg as in Eqns. \ref{eq:40}, \ref{eq:41}, \ref{eq:42}, \ref{eq:43} and \ref{eq:0r}), then $A_c$ would represent the union of all the vehicle types that are compatible to carry any of the load types represented within $c$.\\
\hline $O$ & Set of vehicle types which shall not go back to their starting depots (\ie Open Trips)\\

\hline
\hline $D$ & The set of all Delivery cargos\\
\hline $P$ & The set of all Pickup cargos\\
\hline $C_k$ & Compatible Cargo: $C_k$ refers to the set of all load types that are allowed to be carried by a vehicle type $k$\\
\hline $Q^c_i$ & Represents quantity of Cargo Type $c$ at PickUp and/or Delivery Vertex $i$ (availability of each DCT in $W$, or, capacity for each PCT at $R$, or, requirement of each compatible cargo at $N$)\\
\hline $B_i$ & Set of delivery and pickup load types that are allowed to be transhipped at transhipment port $i$, where $i \in S$. In other words, loading or unloading of all load-types within the set $B_i$ is allowed at $i$. It is assumed that the rest of the incompatible cargo would not be taken out of vehicles during transhipment at the concerned incompatible ports (as the formulation does not allow this).\\

\hline
\hline $U^{k,c}$ & Represents the loading/unloading time of cargo type $c$ into/from vehicle type $k$. ($U^{k,c} \geq 0, ~ \forall k, \forall c$)\\
\hline $T_{i,j}^k$ & $T_{i,j}^k$ Represents the shortest time taken to travel from vertex $i$ to vertex $j$ by vehicle type $k$ (allowing the vehicle to change modes of transport on which it is able to travel, if possible, enabling swift robust transportation)\\
\hline $M$ & A very large positive real number\\

\hline
\hline $E_c$ & Volume of a unit of load-type $c$\\
\hline $E^k$ & Volume of a vehicle of type $k$\\
\hline $F_c$ & Weight of a unit of load-type $c$\\
\hline $F^k$ & Weight of a vehicle of type $k$\\
\hline

\end{longtable}

%\begin{longtable}[c]{| c | c |}
% \begin{longtable}{|p{1.75cm}|p{10.1cm}|}
\begin{longtable}{|m{0.08\linewidth}|m{0.88\linewidth}|}
\caption{\textbf{Variable Definitions}\label{longTab: Variable Definitions}}\\

 % \hline
 % \multicolumn{2}{| c |}{Begin of Table}\\
 \hline
 \textbf{Notation} & \textbf{Description}\\
 \hline
 \endfirsthead

 \hline
 \multicolumn{2}{|c|}{Continuation of \autoref{longTab: Variable Definitions} \textbf{Variable Definitions}}\\
 \hline
 \textbf{Notation} & \textbf{Description}\\
 \hline
 \endhead

 % \hline
 % \endfoot

 % \hline
 % \multicolumn{2}{| c |}{End of Table}\\
 % \hline\hline
 % \endlastfoot

\hline
\hline $z_x$ & Continuous variable equal to the largest tour time among all vehicles essentially representing the entire Emergency Operation Time \\
\hline $z_s$ & Continuous variable equal to the MakeSpan of the Entire Emergency Operation Time, \ie  the sum of the routes times of all vehicles\\
\hline $z_{s_{OT}}$ & Intermediate continuous variable used to populate $z_s$; this takes the sum of the route-times of all vehicles which need not go back to their respective starting depots essentially performing open trips (like most trains) \\
\hline $z_{s_{CT}}$ & Intermediate continuous variable used to populate $z_s$; this takes the sum of the route-times of all vehicles which must go back to their respective starting depots essentially performing closed trips (like most cars and buses in general) \\

\hline
\hline $x_{i,j}^{v,l}$ & Binary decision variable referring to a vehicle $v$ travelling from vertex $i$ to vertex $j$ on its level $l$\\
\hline $x_i^{v,l,m}$ & Binary decision variable referring to a vehicle $v$ at vertex $i$ travelling from level $l$ to level $m$\\

\hline
\hline $y_{i,j}^{v,l,c}$ & Continuous-positive variable referring to the amount of cargo load of category $c$ being taken by vehicle $v$ traveling from vertex $i$ to vertex $j$ both on its level $l$, the carried cargo being compatible with the vehicle type. Ideally this should be an integer variable since we are considering indivisible units of commodities being carried, however we relax this consideration.\\
\hline $y_i^{v,l,m,c}$ & Continuous-positive variable referring to the amount of cargo load of category $c$ being taken by vehicle $v$ at vertex $i$ traveling from level $l$ to level $m$, the carried cargo being compatible with the vehicle type\\

\hline
\hline $b_i^c$ & Intermediate continuous variable representing the amount of load transfer of type $c$ at only the lowermost/first level/layer, for all vehicles, at Node $i$\\
\hline $q_i^c$ & Intermediate continuous variable representing the amount of load transfer of type $c$ at layers apart from the lowermost, for all vehicles, at Node $i$\\

\hline
\hline $r_i^{v,l,c}$ & Residue refers to the amount of actual transhipment operation of each cargo-type at each $TP$. It is a Real variable allowed to take negative or positive values depending on whether the compatible cargo-type was deposited or collected by the vehicle during any of its layer-wise visits. \\
& In this formulation, we calculate the Residue similar to deposit at $TP$s, and therefore the variable takes negative values whenever the resources are collected for subsequent transhipment.\\
\hline $n_i^{v,l,c}$ & This is an intermediate binary variable used during calculations of modulus of the Residue. This is constrained to take a value of $1$ when the corresponding residue is positive, and $0$ when the residue is negative.\\
\hline $e_i^{v,l,c}$ & This continuous-positive variable stores the modulus or absolute value of the Residues.\\

\hline
\hline $a_i^{v,l}$ & Continuous-positive variable, refers to the arrival time of vehicle $v$ at vertex $i$ for a specific visit \ie at layer $l$\\
\hline $d_i^{v,l}$ & Continuous-positive variable, refers to the departure time of vehicle $v$ at vertex $i$ after a specific visit \ie at layer $l$\\
\hline $t_{i,j}^{v,l}$ & Continuous-positive variable, refers to the travel time of vehicle $v$ when journeying from vertex $i$ to $j$ on it's layer $l$\\
\hline $t_i^{v,l,m}$ & Continuous-positive variable, refers to the cascaded (or, continually added) travel time of the vehicle $v$ at vertex $i$ when journeying from layer $l$ to layer $m$\\
\hline $w_i^{v,l}$ & Continuous-positive variable, refers to the non-working or waiting time of the vehicle $v$ at vertex $i$ during a specific visit \ie at layer $l$. The waiting time would be non-zero only during the collection at any $S$, when the resource to be collected is yet to arrive at that Transhipment Port, \ie the cumulative addition of the working-durations (working durations refer to the loading/unloading done by the vehicle $v$ at $i$ during its visit at level $l$) is less than the required resource arrival time which is still on its way onboard another vehicle.\\

\hline
\hline $o_{i,(\hat{v},\hat{l},\hat{s})}^{c,(v,l)}$ & The residues obtained needs to be populated over load-type and Transhipment Port specific space-time matrices. This continuous variable $o$ represents the cells in each of these matrices. For each $i \in S$ and $c \in B_i$, matrices are created with the space axis varying over $(v,l)$ and the time axis being varied over $(\hat{v},\hat{l},\hat{s})$, where $l$ refers to a layer of vehicle $v$, and $\hat{l}$ refers to a layer of vehicle $\hat{v}$. Here $\hat{s} \in \{Arr,Dep\}$, where $Arr$ refers to arrival and $Dep$ refers to departure.
These space-time matrices are used to ensure that the temporal accumulation of the residues always remain positive.\\
\hline $g_{i,(\hat{v},\hat{l},\hat{s})}^{(v,l,s)}$ & This intermediate binary variable compares the arrival or departure times of the subscripted vehicle $\hat{v}$ \wrt the super-scripted vehicle $v$ during their visits on layers $\hat{l}$ and $l$ respectively, only at the Transhipment Ports \ie $i \in S$. Here $s,\hat{s} \in \{Arr,Dep\}$ representing arrival or departure for the vehicle $v$ and $\hat{v}$ respectively. The variable takes value $1$ when:
    \begin{minipage}[t]{1\textwidth}
    \begin{itemize}
        \item for $s=\hat{s}=Arr$: if $a_i^{\hat{v},\hat{l}}$ $\leq$ $a_i^{v,l}$, else takes value $0$
        \item for $s=Arr$ and $\hat{s}=Dep$: if $d_i^{\hat{v},\hat{l}}$ $\leq$ $a_i^{v,l}$, else takes value $0$
        \item for $s=Dep$ and $\hat{s}=Arr$: if $a_i^{\hat{v},\hat{l}}$ $\geq$ $d_i^{v,l}$, else takes value $0$
        \item for $s=\hat{s}=Dep$: if $d_i^{\hat{v},\hat{l}}$ $\geq$ $d_i^{v,l}$, else takes value $0$
    \end{itemize}
    \end{minipage}\\
\hline
 \end{longtable}

\small

% \subsection{Objective Function} \label{Objective Function}

% Depending on the choice of the minimization approach, the objective function may be any one among the equations in \ref{eq:1}. If the optimization approach is set to be done using the Min-Max approach which minimizes the largest route-time among all vehicles then Eq. \ref{eq:1a} is used. Otherwise Eq. \ref{eq:1b} is used to minimise the MakeSpan, \ie the sum of the route-times of all vehicles. We suggest combining the cascaded optimization approach with the min-max objective function; this is discussed in detail after the formulation description.

% \begin{subequations} \label{eq:1}
% \begin{equation} \label{eq:1a}
% Minimize \quad z_x
% \end{equation}
% \begin{equation} \label{eq:1b}
% Minimize \quad z_s
% \end{equation}
% \end{subequations}

The Min-Max objective in Eq. \ref{eq:1} minimizes the largest route-time among all vehicles (\ie the Makespan). Makespan is defined as the largest time taken among all simultaneous tasks, or more precisely the time difference between the start and finish of the entire operation (if all functions don't start at the same time), in our case it corresponds to vehicle trip durations. Generally, we have found that researchers have focused on minimizing the makespan for various problems but overlook the subsequent steps for further internal optimization; this is a general problem in the OR research field and is not limited to VRP. No study has focused on the time minimization aspect as done by us; the proposition of cascaded makespan minimization as a replacement to traditional makespan minimization is novel and an improvement over traditional makespan minimization (\cite{VENKATANARASIMHA201363}) which is simply equivalent to performing a single cascade of our optimization process. We combine the cascaded optimization approach with the min-max objective function, to solve this MILP using a multi-step minimization approach; this Cascaded Makespan Minimization is discussed in detail after the formulation description. The variables used in the formulation are described in \autoref{longTab: Variable Definitions}.

{\noindent \textbf{Objective Function}:}
\begin{equation} \label{eq:1}
Minimize \quad z_x
\end{equation}

% \subsection{Routing Constraints} \label{Routing Constraints}
{\noindent \textbf{Routing Constraints}:}

\begin{equation} \label{eq:2}
\sum_{j \in W_k,N_k,S_k} x_{h,j}^{v,1}\leq 1, \quad
\forall h \in H,
\forall k \in K_h,
\forall u \in G_{h,k},
v=(h,k,u)
\end{equation}

\begin{equation} \label{eq:3}
\sum_{l \in L_v}~
\sum_{j \in N_k,S_k,R_k}
x_{j,h}^{v,l} = \sum_{j \in W_k,N_k,S_k} x_{h,j}^{v,1}, \quad
\forall h \in H,
\forall k \in K_h,
\forall u \in G_{h,k},
v=(h,k,u),
\end{equation}

\begin{subequations} \label{eq:4}
    \begin{equation} \label{eq:4a}
    x^{v,l+1,l}_i +
    \sum_{\substack{j \in V_k,h \\ i \neq j \\ if~i \in R \Rightarrow j\neq h}}
    x_{j,i}^{v,l}
    =
    x^{v,l,l+1}_i +
    \sum_{\substack{j \in V_k,h \\ i \neq j \\ if~i \in W \Rightarrow j \neq h}}
    x_{i,j}^{v,l},
    \quad
    \forall h \in H,
    \forall k \in K_h,\\
    \forall u \in G_{h,k},
    v=(h,k,u),
    l=1,
    \forall i \in V_k,
    \end{equation}    
    
    \begin{multline} \label{eq:4b}
    x^{v,l-1,l}_i +
    x^{v,l+1,l}_i +
    \sum_{\substack{j \in V_k \\ i \neq j}}
    x_{j,i}^{v,l}
    =
    x^{v,l,l-1}_i +
    x^{v,l,l+1}_i +
    \sum_{\substack{j \in V_k,h \\ i \neq j \\ if~i \in W \Rightarrow j \neq h}}
    x_{i,j}^{v,l}, \quad
    \forall h \in H,
    \forall k \in K_h,
    \forall u \in G_{h,k},
    v=(h,k,u),
    \forall i \in V_k,\\
    \forall l \in L_v \smallsetminus{\{1, l_v^{max}\}},
    \end{multline}
    
    \begin{equation} \label{eq:4c}
    x^{v,l-1,l}_i +
    \sum_{\substack{j \in V_k \\ i \neq j}}
    x_{j,i}^{v,l}
    =
    x^{v,l,l-1}_i +
    \sum_{\substack{j \in V_k,h \\ i \neq j \\ if~i \in W \Rightarrow j \neq h}}
    x_{i,j}^{v,l}, \quad
    \forall h \in H,
    \forall k \in K_h,
    \forall u \in G_{h,k},
    v=(h,k,u),
    \forall i \in V_k,
    if~l_v^{max}>1,\\
    l=l_v^{max},
    \end{equation}
\end{subequations}

\begin{equation} \label{eq:7}
\sum_{\substack{h \in H\\
                k \in K_h \cap K_i\\
                u \in G_{h,k}\\
                v=(h,k,u)}}
~\sum_{ \substack{j \in V_k,h \\ i \neq j}} 
~\sum_{l \in L_v}
x^{v,l}_{i,j} \leq 1, \quad \forall i \in N^M,
\end{equation}

% \subsection{Flow Constraints} \label{Flow Constraints}
{\noindent \textbf{Flow Constraints (at Vehicle Depots)}:}

% \subsubsection{For Vehicle Depots}

\begin{equation} \label{eq:8}
y^{v,1,c}_{h,j} \leq 0, \quad
\forall h \in H,
\forall k \in K_h,
\forall u \in G_{h,k},
v=(h,k,u),
\forall c \in C_k,
\forall j \in W_k \cup S_k \cup N_k,
\end{equation}

\begin{equation} \label{eq:9}
y^{v,l,c}_{i,h} \leq 0, \quad
\forall h \in H,
\forall k \in K_h,
\forall u \in G_{h,k},
v=(h,k,u),
\forall l \in L_v,
\forall c \in C_k,
\forall i \in S_k \cup N_k \cup R_k,
\end{equation}

%However, for generalizing the formulation, we retain these terms.

% \subsubsection{For Warehouses}
{\noindent \textbf{Flow Constraints (at Warehouses)}:}

\begin{subequations} \label{eq:10}
\begin{equation} \label{eq:10a}
y^{v,l,l+1,c}_i +
\sum_{\substack{j \in V_k \\ i \neq j}}  y_{i,j}^{v,l,c}
\geq y^{v,l+1,l,c}_i +
\sum_{\substack{j \in V_k,h \\ i \neq j}}  y_{j,i}^{v,l,c}, \quad
\forall h \in H,
\forall k \in K_h,
\forall u \in G_{h,k},
v=(h,k,u),
l=1,
\forall c \in D \cap C_k,\\
\forall i \in W_k,
\end{equation}

\begin{multline} \label{eq:10b}
y^{v,l,l+1,c}_i +
y^{v,l,l-1,c}_i +
\sum_{\substack{j \in V_k \\ i \neq j}}  y_{i,j}^{v,l,c}
\geq
y^{v,l+1,l,c}_i +
y^{v,l-1,l,c}_i +
\sum_{\substack{j \in V_k \\ i \neq j}} 
y_{j,i}^{v,l,c}, \quad
\forall h \in H,
\forall k \in K_h,
\forall u \in G_{h,k},
v=(h,k,u),\\
\forall l \in L_v \setminus \{1,l_v^{max}\},
\forall c \in D \cap C_k,
\forall i \in W_k,
\end{multline}

\begin{multline} \label{eq:10c}
y^{v,l,l-1,c}_i +
\sum_{\substack{j \in V_k \\ i \neq j}}  y_{i,j}^{v,l,c}
\geq
y^{v,l-1,l,c}_i +
\sum_{\substack{j \in V_k \\ i \neq j}}  y_{j,i}^{v,l,c}, \quad
\forall h \in H,
\forall k \in K_h,
\forall u \in G_{h,k},
v=(h,k,u),
if~l_v^{max}>1,
l=l_v^{max},\\
\forall c \in D \cap C_k,
\forall i \in W_k,
\end{multline}
\end{subequations}

\begin{subequations} \label{eq:11}
\begin{equation} \label{eq:11a}
y^{v,l,l+1,c}_i +
\sum_{\substack{j \in V_k \\ i \neq j}}  y_{i,j}^{v,l,c}
=
y^{v,l+1,l,c}_i +
\sum_{\substack{j \in V_k,h \\ i \neq j}}  y_{j,i}^{v,l,c}, \quad
\forall h \in H,
\forall k \in K_h,
\forall u \in G_{h,k},
v=(h,k,u),
l=1,
\forall c \in P \cap C_k,\\
\forall i \in W_k,
\end{equation}

\begin{multline} \label{eq:11b}
y^{v,l,l+1,c}_i +
y^{v,l,l-1,c}_i +
\sum_{\substack{j \in V_k \\ i \neq j}}  y_{i,j}^{v,l,c}
=
y^{v,l+1,l,c}_i +
y^{v,l-1,l,c}_i +
\sum_{\substack{j \in V_k \\ i \neq j}}  y_{j,i}^{v,l,c}, \quad
\forall h \in H,
\forall k \in K_h,
\forall u \in G_{h,k},
v=(h,k,u),\\
\forall l \in L_v \setminus \{1,l_v^{max}\},
\forall c \in P \cap C_k,
\forall i \in W_k,
\end{multline}

\begin{multline} \label{eq:11c}
y^{v,l,l-1,c}_i +
\sum_{\substack{j \in V_k \\ i \neq j}}  y_{i,j}^{v,l,c}
=
y^{v,l-1,l,c}_i +
\sum_{\substack{j \in V_k \\ i \neq j}}  y_{j,i}^{v,l,c}, \quad
\forall h \in H,
\forall k \in K_h,
\forall u \in G_{h,k},
v=(h,k,u),
if~l_v^{max}>1,
l=l_v^{max},\\
\forall c \in P \cap C_k,
\forall i \in W_k,
\end{multline}
\end{subequations}

{\noindent \textbf{Capacity Constraints (for Warehouses)}:}

\begin{subequations} \label{eq:12}

\begin{equation} \label{eq:12a}
\sum_{\substack{h \in H\\
                k \in K_h \cap A_c \cap K_i\\
                u \in G_{h,k}\\
                v=(h,k,u)}}
~\sum_{l = 1}~
\left(
\sum_{ \substack{j \in V_k \\ i \neq j}} y^{v,l,c}_{i,j}
-
\sum_{ \substack{j \in V_k,h \\ i \neq j}} y^{v,l,c}_{j,i}
\right)
= b_i^c, \quad
\forall i \in W,
\forall c \in D,
\end{equation}

\begin{equation} \label{eq:12b}
\sum_{\substack{h \in H\\
                k \in K_h \cap A_c \cap K_i\\
                u \in G_{h,k}\\
                v=(h,k,u)}}
~\sum_{l \in L_v \setminus{\{1\}}}
~\sum_{ \substack{j \in V_k \\ i \neq j}}
\Bigl( y^{v,l,c}_{i,j}
-
y^{v,l,c}_{j,i} \Bigr)
= q_i^c,\quad
\forall i \in W,
\forall c \in D,
\end{equation}

\begin{equation} \label{eq:12c}
b^c_i + q^c_i \leq Q_i^c,
\quad \forall i \in W, \forall c \in D,
\end{equation}
\end{subequations}

{\noindent \textbf{Flow Constraints (at Relief Centres)}:}
% \subsubsection{For Relief Centres}

\begin{subequations} \label{eq:13}
    \begin{equation} \label{eq:13a}
    y^{v,l,l+1,c}_i +
    \sum_{\substack{j \in V_k,h \\ i \neq j}}
    y_{i,j}^{v,l,c} \leq
    y^{v,l+1,l,c}_i +
    \sum_{\substack{j \in V_k \\ i \neq j}}
    y_{j,i}^{v,l,c}, \quad
    \forall h \in H,
    \forall k \in K_h,
    \forall u \in G_{h,k},
    v=(h,k,u),
    l=1,
    \forall c \in P \cap C_k,\\
    \forall i \in R_k,
    \end{equation}
    
    \begin{multline} \label{eq:13b}
    y^{v,l,l+1,c}_i + y^{v,l,l-1,c}_i +
    \sum_{\substack{j \in V_k,h \\ i \neq j}}  y_{i,j}^{v,l,c} \leq
    y^{v,l+1,l,c}_i + y^{v,l-1,l,c}_i +
    \sum_{\substack{j \in V_k \\ i \neq j}}  y_{j,i}^{v,l,c}, \quad
    \forall h \in H,
    \forall k \in K_h,
    \forall u \in G_{h,k},
    v=(h,k,u),\\
    \forall l \in L_v \setminus \{1,l_v^{max}\},
    \forall c \in P \cap C_k,
    \forall i \in R_k,
    \end{multline}
    
    \begin{multline} \label{eq:13c}
    y^{v,l,l-1,c}_i +
    \sum_{\substack{j \in V_k,h \\ i \neq j}}
    y_{i,j}^{v,l,c} \leq
    y^{v,l-1,l,c}_i +
    \sum_{\substack{j \in V_k \\ i \neq j}}  y_{j,i}^{v,l,c}, \quad
    \forall h \in H,
    \forall k \in K_h,
    \forall u \in G_{h,k},
    v=(h,k,u),
    if~l_v^{max}>1, 
    l=l_v^{max},\\
    \forall c \in P \cap C_k,
    \forall i \in R_k,
    \end{multline}
\end{subequations}

\begin{subequations} \label{eq:14}
    \begin{equation} \label{eq:14a}
    y^{v,l,l+1,c}_i + 
    \sum_{\substack{j \in V_k,h \\ i \neq j}}  y_{i,j}^{v,l,c} =
    y^{v,l+1,l,c}_i +
    \sum_{\substack{j \in V_k \\ i \neq j}}  y_{j,i}^{v,l,c}, \quad
    \forall h \in H,
    \forall k \in K_h,
    \forall u \in G_{h,k},
    v=(h,k,u),
    l=1,
    \forall c \in D \cap C_k,\\
    \forall i \in R_k,
    \end{equation}

    \begin{multline} \label{eq:14b}
    y^{v,l,l+1,c}_i +
    y^{v,l,l-1,c}_i +
    \sum_{\substack{j \in V_k,h \\ i \neq j}}  y_{i,j}^{v,l,c} =
    y^{v,l+1,l,c}_i +
    y^{v,l-1,l,c}_i +
    \sum_{\substack{j \in V_k \\ i \neq j}}  y_{j,i}^{v,l,c}, \quad
    \forall h \in H,
    \forall k \in K_h,
    \forall u \in G_{h,k},
    v=(h,k,u),\\
    \forall l \in L_v \setminus \{1,l_v^{max}\}, 
    \forall c \in D \cap C_k,
    \forall i \in R_k,
    \end{multline}
    
    \begin{multline} \label{eq:14c}
    y^{v,l,l-1,c}_i + 
    \sum_{\substack{j \in V_k,h \\ i \neq j}}  y_{i,j}^{v,l,c} =
    y^{v,l-1,l,c}_i +
    \sum_{\substack{j \in V_k \\ i \neq j}}  y_{j,i}^{v,l,c}, \quad
    \forall h \in H,
    \forall k \in K_h,
    \forall u \in G_{h,k},
    v=(h,k,u),
    if~l_v^{max}>1,    
    l=l_v^{max},\\
    \forall c \in D \cap C_k,
    \forall i \in R_k,
    \end{multline}
\end{subequations}

{\noindent \textbf{Capacity Constraints (for Relief Centers)}:}

\begin{equation} \label{eq:15}
\sum_{\substack{h \in H\\
                k \in K_h \cap A_c \cap K_i\\
                u \in G_{h,k}\\
                v=(h,k,u)}}
~\sum_{l \in L_v}~
\left(
\sum_{ \substack{j \in V_k \\ i \neq j}}y^{v,l,c}_{j,i}
-
\sum_{ \substack{j \in V_k,h \\ i \neq j}}y^{v,l,c}_{i,j}
\right)
\leq Q_i^c, \quad
\forall i \in R,
\forall c \in P,
\end{equation}

% \subsubsection{Flow Constraints at Nodes}
{\noindent \textbf{Flow Constraints (at Nodes)}:}

\begin{subequations} \label{eq:16}
\begin{equation} \label{eq:16a}
y^{v,l,l+1,c}_i +
\sum_{\substack{j \in V_k,h \\ i \neq j}}  y_{i,j}^{v,l,c} \geq
y^{v,l+1,l,c}_i +
\sum_{\substack{j \in V_k,h \\ i \neq j}}  y_{j,i}^{v,l,c}, \quad
\forall h \in H,
\forall k \in K_h,
\forall u \in G_{h,k},
v=(h,k,u),
l=1,
\forall c \in P \cap C_k,\\
\forall i \in N_k,
\end{equation}

\begin{multline} \label{eq:16b}
y^{v,l,l+1,c}_i +
y^{v,l,l-1,c}_i +
\sum_{\substack{j \in V_k,h \\ i \neq j}}  y_{i,j}^{v,l,c} \geq
y^{v,l+1,l,c}_i +
y^{v,l-1,l,c}_i +
\sum_{\substack{j \in V_k \\ i \neq j}}  y_{j,i}^{v,l,c}, \quad
\forall h \in H,
\forall k \in K_h,
\forall u \in G_{h,k},
v=(h,k,u),\\
\forall l \in L_v \setminus \{1,l_v^{max}\}, 
\forall c \in P \cap C_k,
\forall i \in N_k,
\end{multline}

\begin{multline} \label{eq:16c}
y^{v,l,l-1,c}_i +
\sum_{\substack{j \in V_k,h \\ i \neq j}}  y_{i,j}^{v,l,c} \geq
y^{v,l-1,l,c}_i +
\sum_{\substack{j \in V_k \\ i \neq j}}  y_{j,i}^{v,l,c}, \quad
\forall h \in H,
\forall k \in K_h,
\forall u \in G_{h,k},
v=(h,k,u),
if~l_v^{max}>1,
l=l_v^{max},\\
\forall c \in P \cap C_k, 
\forall i \in N_k,
\end{multline}
\end{subequations}

\begin{subequations} \label{eq:17}
\begin{equation} \label{eq:17a}
y^{v,l,l+1,c}_i +
\sum_{\substack{j \in V_k,h \\ i \neq j}}  y_{i,j}^{v,l,c} \leq
y^{v,l+1,l,c}_i +
\sum_{\substack{j \in V_k,h \\ i \neq j}}  y_{j,i}^{v,l,c}, \quad
\forall h \in H,
\forall k \in K_h,
\forall u \in G_{h,k},
v=(h,k,u),
l=1,
\forall c \in D \cap C_k,\\
\forall i \in N_k,
\end{equation}

\begin{multline} \label{eq:17b}
y^{v,l,l+1,c}_i +
y^{v,l,l-1,c}_i +
\sum_{\substack{j \in V_k,h \\ i \neq j}}  y_{i,j}^{v,l,c} \leq
y^{v,l+1,l,c}_i +
y^{v,l-1,l,c}_i +
\sum_{\substack{j \in V_k \\ i \neq j}}  y_{j,i}^{v,l,c}, \quad
\forall h \in H,
\forall k \in K_h,
\forall u \in G_{h,k},
v=(h,k,u),\\
\forall l \in L_v \setminus \{1,l_v^{max}\}, 
\forall c \in D \cap C_k,
\forall i \in N_k,
\end{multline}

\begin{multline} \label{eq:17c}
y^{v,l,l-1,c}_i +
\sum_{\substack{j \in V_k,h \\ i \neq j}}  y_{i,j}^{v,l,c} \leq
y^{v,l-1,l,c}_i +
\sum_{\substack{j \in V_k \\ i \neq j}}  y_{j,i}^{v,l,c}, \quad
\forall h \in H,
\forall k \in K_h,
\forall u \in G_{h,k},
v=(h,k,u),
if~l_v^{max}>1,
l=l_v^{max},\\
\forall c \in D \cap C_k, 
\forall i \in N_k,
\end{multline}
\end{subequations}

{\noindent \textbf{Resource Constraints (for Nodes)}:}

% \begin{multicols}{2}
\begin{subequations} \label{eq:18}
    \begin{equation} \label{eq:18a}
    \sum_{\substack{h \in H\\
                    k \in K_h \cap A_c \cap K_i\\
                    u \in G_{h,k}\\
                    v=(h,k,u)}}
    ~\sum_{l = 1}
    ~\sum_{ \substack{j \in V_k,h \\ i \neq j}}
    \Bigl(
    y^{v,l,c}_{i,j}
    -
    y^{v,l,c}_{j,i}
    \Bigr)
    = b_i^c,\quad
    \forall i \in N,
    \forall c \in P,
    \end{equation}
    
    \begin{equation} \label{eq:18b}
    \sum_{\substack{h \in H\\
                k \in K_h \cap A_c \cap K_i\\
                u \in G_{h,k}\\
                v=(h,k,u)}}
    ~\sum_{l \in L_v \setminus{\{1\}}}
    \left(
    \sum_{ \substack{j \in V_k,h \\ i \neq j}} y^{v,l,c}_{i,j}
    -
    \sum_{ \substack{j \in V_k \\ i \neq j}} y^{v,l,c}_{j,i}
    \right)
    = q_i^c,
    \quad
    \forall i \in N,
    \forall c \in P,
    \end{equation}

    \begin{equation} \label{eq:18c}
    b_i^c + q_i^c
    \geq Q_i^c,
    \quad \forall i \in N,
    \forall c \in P,
    \end{equation}
    
\end{subequations}

\begin{subequations} \label{eq:19}
    \begin{equation} \label{eq:19a}
    \sum_{\substack{h \in H\\
                k \in K_h \cap A_c \cap K_i\\
                u \in G_{h,k}\\
                v=(h,k,u)}}
    ~\sum_{ \substack{j \in V_k,h \\ i \neq j}}
    \Bigl(
    y^{v,1,c}_{j,i}
    -
    y^{v,1,c}_{i,j}
    \Bigr)
    = b_i^c, \quad
    \forall i \in N,
    \forall c \in D,
    \end{equation}
    
    \begin{equation} \label{eq:19b}
    \sum_{\substack{h \in H\\
            k \in K_h \cap A_c \cap K_i\\
            u \in G_{h,k}\\
            v=(h,k,u)}}
    ~\sum_{l \in L_v \setminus{\{1\}}}
    \left(
    \sum_{ \substack{j \in V_k \\ i \neq j}} y^{v,l,c}_{j,i}
    -
    \sum_{ \substack{j \in V_k,h \\ i \neq j}} y^{v,l,c}_{i,j}
    \right)
    = q_i^c,
    \quad
    \forall i \in N,
    \forall c \in D,
    \end{equation}

    \begin{equation} \label{eq:19c}
    b_i^c + q_i^c
    \geq Q_i^c,
    \quad \forall i \in N,
    \forall c \in D,
    \end{equation}
\end{subequations}
% \end{multicols}

% \subsubsection{Flow Constraints at Transhipment Ports:}
{\noindent \textbf{Flow Constraints (at Transhipment Ports)}:}

\begin{subequations} \label{eq:20}

\begin{multline} \label{eq:20a}
y^{v,l+1,l,c}_i -
y^{v,l,l+1,c}_i +
\sum_{\substack{j \in V_k,h \\ i \neq j}} \Bigl( y_{j,i}^{v,l,c} - y_{i,j}^{v,l,c} \Bigr) =
r_i^{v,l,c},
\quad
\forall h \in H,
\forall k \in K_h,
\forall u \in G_{h,k},
v=(h,k,u),
l=1,
\forall i \in S_k,\\
\forall c \in C_k \cap B_i, 
\end{multline}

\begin{multline} \label{eq:20b}
y^{v,l+1,l,c}_i -
y^{v,l,l+1,c}_i +
y^{v,l-1,l,c}_i -
y^{v,l,l-1,c}_i +
\sum_{\substack{j \in V_k \\ i \neq j}} y_{j,i}^{v,l,c}
-
\sum_{\substack{j \in V_k,h \\ i \neq j}} y_{i,j}^{v,l,c}
=
r_i^{v,l,c}, \quad
\forall h \in H,
\forall k \in K_h,
\forall u \in G_{h,k},
v=(h,k,u),\\
\forall i \in S_k,
\forall c \in C_k \cap B_i,
\forall l \in L_v \setminus \{1,l_v^{max}\},
\end{multline}

\begin{multline} \label{eq:20c}
y^{v,l-1,l,c}_i -
y^{v,l,l-1,c}_i +
\sum_{\substack{j \in V_k \\ i \neq j}} y_{j,i}^{v,l,c}
-
\sum_{\substack{j \in V_k,h \\ i \neq j}} y_{i,j}^{v,l,c}
=
r_i^{v,l,c},\quad
\forall h \in H,
\forall k \in K_h,
\forall u \in G_{h,k},
v=(h,k,u),
if~l_v^{max}>1,
l=l_v^{max},\\
\forall i \in S_k,
\forall c \in C_k \cap B_i,
\end{multline}

\end{subequations}

\begin{subequations} \label{eq:21}

\begin{multline} \label{eq:21a}
y^{v,l+1,l,c}_i -
y^{v,l,l+1,c}_i +
\sum_{\substack{j \in V_k,h \\ i \neq j}}  \Bigl( y_{j,i}^{v,l,c} - y_{i,j}^{v,l,c} \Bigr) = 0, \quad
\forall h \in H,
\forall k \in K_h,
\forall u \in G_{h,k},
v=(h,k,u),
l=1,
\forall i \in S_k,\\
\forall c \in C_k \setminus B_i,
\end{multline}

\begin{multline} \label{eq:21b}
y^{v,l+1,l,c}_i -
y^{v,l,l+1,c}_i +
y^{v,l-1,l,c}_i -
y^{v,l,l-1,c}_i +
\sum_{\substack{j \in V_k \\ i \neq j}} y_{j,i}^{v,l,c}
-
\sum_{\substack{j \in V_k,h \\ i \neq j}} y_{i,j}^{v,l,c}
= 0, \quad
\forall h \in H,
\forall k \in K_h,
\forall u \in G_{h,k},
v=(h,k,u),\\
\forall l \in L_v \setminus \{1,l_v^{max}\},
\forall i \in S_k,
\forall c \in C_k \setminus B_i,
\end{multline}

\begin{multline} \label{eq:21c}
y^{v,l-1,l,c}_i -
y^{v,l,l-1,c}_i +
\sum_{\substack{j \in V_k \\ i \neq j}} y_{j,i}^{v,l,c}
-
\sum_{\substack{j \in V_k,h \\ i \neq j}} y_{i,j}^{v,l,c}
= 0,\quad
\forall h \in H,
\forall k \in K_h,
\forall u \in G_{h,k},
v=(h,k,u),
if~l_v^{max}>1,
l=l_v^{max},\\
\forall i \in S_k,
\forall c \in C_k \setminus B_i,
\end{multline}

\end{subequations}

% \subsection{Time constraints}

{\noindent \textbf{Vehicle Start-Time Constraints}:}

\begin{equation} \label{eq:22}
d^{v,1}_h \leq 0, \quad
\forall h \in H,
\forall k \in K_h,
\forall u \in G_{h,k},
v=(h,k,u),
\end{equation}

% \subsubsection{Calculating the cascaded travel times}
{\noindent \textbf{Calculating the cascaded travel times (on the same layer)}:}

\begin{subequations} \label{eq:23}
\begin{multline} \label{eq:23a}
t^{v,1}_{i,j} \leq M \cdot x_{i,j}^{v,1}, \quad
\forall h \in H,
\forall k \in K_h,
\forall u \in G_{h,k},
v=(h,k,u),
\forall i,j \in V_k \cup h,
i \neq j,
(if~i \in W \Rightarrow j \neq h),\\
(if~i=h \Rightarrow j \notin R),
\end{multline}

\begin{multline} \label{eq:23b}
t^{v,l}_{i,j} \leq M \cdot x_{i,j}^{v,l}, \quad
\forall h \in H,
\forall k \in K_h,
\forall u \in G_{h,k},
v=(h,k,u),
\forall l \in L_v \smallsetminus \{{1}\},
\forall i \in V_k,
\forall j \in V_k \cup h,
i \neq j,\\
(if~i \in W \Rightarrow j \neq h),
\end{multline}
\end{subequations}

\begin{subequations} \label{eq:24}
\begin{multline} \label{eq:24a}
t^{v,1}_{i,j}
\geq -M \cdot \Bigl( 1-x_{i,j}^{v,1} \Bigr) +
d_i^{v,1} +
T_{i,j}^k \cdot x_{i,j}^{v,1},
\quad
\forall h \in H,
\forall k \in K_h,
\forall u \in G_{h,k},
v=(h,k,u),
\forall i,j \in V_k \cup h,
i \neq j,\\
(if~i \in W \Rightarrow j \neq h),
(if~i=h \Rightarrow j \notin R),
\end{multline}

\begin{multline} \label{eq:24b}
t^{v,l}_{i,j}
\geq -M \cdot \Bigl( 1-x_{i,j}^{v,l} \Bigr) +
d_i^{v,l} +
T_{i,j}^k \cdot x_{i,j}^{v,l},
\quad
\forall h \in H,
\forall k \in K_h,
\forall u \in G_{h,k},
v=(h,k,u),
\forall l \in L_v \smallsetminus \{{1}\},
\forall i \in V_k, \forall j \in V_k \cup h,\\
i \neq j,
(if~i \in W \Rightarrow j \neq h),
\end{multline}
\end{subequations}

\begin{subequations} \label{eq:25}
\begin{multline} \label{eq:25a}
t^{v,1}_{i,j} \leq
M \cdot \Bigl( 1-x_{i,j}^{v,1} \Bigr) +
d_i^{v,1} + T_{i,j}^k \cdot x_{i,j}^{v,1},
\quad
\forall h \in H,
\forall k \in K_h,
\forall u \in G_{h,k},
v=(h,k,u),
\forall i,j \in V_k \cup h,
i \neq j,\\
(if~i \in W \Rightarrow j \neq h),
(if~i=h \Rightarrow j \notin R),
\end{multline}

\begin{multline} \label{eq:25b}
t^{v,l}_{i,j} \leq
M \cdot \Bigl( 1-x_{i,j}^{v,l} \Bigr) +
d_i^{v,l} + T_{i,j}^k \cdot x_{i,j}^{v,l},
\quad
\forall h \in H,
\forall k \in K_h,
\forall u \in G_{h,k},
v=(h,k,u),
\forall l \in L_v \smallsetminus \{{1}\},
\forall i \in V_k,
\forall j \in V_k \cup h,\\
i \neq j,
(if~i \in W \Rightarrow j \neq h),
\end{multline}
\end{subequations}

{\noindent \textbf{Calculating the cascaded travel times (between layers)}:}

\begin{equation} \label{eq:26}
t^{v,l,m}_i \leq
M \cdot x^{v,l,m}_i,\quad
\forall h \in H,
\forall k \in K_h,
\forall u \in G_{h,k},
v=(h,k,u),%\\
\forall l,m \in L_v,|l-m|=1,
\forall i \in V_k,
\end{equation}

\begin{equation} \label{eq:27}
t^{v,l,m}_i \geq
-M \cdot \Bigl(1-x^{v,l,m}_i \Bigr) +
d^{v,l}_i,\quad
\forall h \in H,
\forall k \in K_h,
\forall u \in G_{h,k},
v=(h,k,u),\\
\forall l,m \in L_v,|l-m|=1,
\forall i \in V_k,
\end{equation}

\begin{equation} \label{eq:28}
t^{v,l,m}_i \leq
M \cdot \Bigl( 1-x^{v,l,m}_i \Bigr) +
d^{v,l}_i,\quad
\forall h \in H,
\forall k \in K_h,
\forall u \in G_{h,k},
v=(h,k,u),\\
\forall l,m \in L_v,|l-m|=1,
\forall i \in V_k,
\end{equation}

{\noindent \textbf{Calculating Arrival Times}:}
% \subsubsection{Calculating arrival times}

\begin{subequations} \label{eq:29}
\begin{equation} \label{eq:29a}
a^{v,l}_i =
t_i^{v,l+1,l} +
\sum_{\substack{j \in V_k,h
\\ i \neq j
\\ if~i \in R \Rightarrow j \neq h
\\if~i=h \Rightarrow j \notin W}} t_{j,i}^{v,l},
\quad
\forall h \in H,
\forall k \in K_h,
\forall u \in G_{h,k},
v=(h,k,u),
l=1,
\forall i \in V_k,
\end{equation}

\begin{equation} \label{eq:29b}
a^{v,l}_i =
t_i^{v,l+1,l} +
t_i^{v,l-1,l} +
\sum_{\substack{j \in V_k
\\ i \neq j
\\ if~i = h \Rightarrow j \notin W}}t_{j,i}^{v,l},
\quad
\forall h \in H,
\forall k \in K_h,
\forall u \in G_{h,k},
v=(h,k,u),
\forall i \in V_k,
\forall l \in L_v \setminus{\{1,l_v^{max}\}},
\end{equation}

\begin{equation} \label{eq:29c}
a^{v,l}_i =
t_i^{v,l-1,l} +
\sum_{\substack{j \in V_k
\\ i \neq j
\\ if~i = h \Rightarrow j \notin W}}t_{j,i}^{v,l},
\quad
\forall h \in H,
\forall k \in K_h,
\forall u \in G_{h,k},
v=(h,k,u),
if~l_v^{max}>1,
l=l_v^{max},
\forall i \in V_k,
\end{equation}

\end{subequations}

% ////////// Commented Equation ////////// \newline
\begin{equation} \label{eq:29.5}
a^{v,l}_h =
\sum_{j \in N_k,S_k,R_k}
t_{j,h}^{v,l},
\quad
\forall h \in H,
\forall k \in K_h,
\forall u \in G_{h,k},
v=(h,k,u),
\forall l \in L_v,
\end{equation}
% ////////// Commented Equation ////////// \newline

{\noindent \textbf{Calculating modulus of Residues}:}
% \subsubsection{Calculating modulus of the residue}

\begin{equation} \label{eq:30}
r^{v,l,c}_i \geq -M \cdot \Bigl( 1-n^{v,l,c}_i \Bigr),
\quad
\forall h \in H,
\forall k \in K_h,
\forall u \in G_{h,k},
v=(h,k,u),
\forall l \in L_v,
\forall i \in S_k,
\forall c \in B_i \cap C_k,
\end{equation}

\begin{equation} \label{eq:31}
r^{v,l,c}_i \leq M \cdot n^{v,l,c}_i,
\quad
\forall h \in H,
\forall k \in K_h,
\forall u \in G_{h,k},
v=(h,k,u),
\forall l \in L_v,
\forall i \in S_k,
\forall c \in B_i \cap C_k,
\end{equation}

\begin{equation} \label{eq:32}
e^{v,l,c}_i \geq r^{v,l,c}_i,
\quad
\forall h \in H,
\forall k \in K_h,
\forall u \in G_{h,k},
v=(h,k,u),
\forall l \in L_v,
\forall i \in S_k,
\forall c \in B_i \cap C_k,
\end{equation}

\begin{equation} \label{eq:33}
e^{v,l,c}_i \geq -r^{v,l,c}_i,
\quad
\forall h \in H,
\forall k \in K_h,
\forall u \in G_{h,k},
v=(h,k,u),
\forall l \in L_v,
\forall i \in S_k,
\forall c \in B_i \cap C_k,
\end{equation}

\begin{equation} \label{eq:34}
e^{v,l,c}_i \leq r^{v,l,c}_i + M \cdot \Bigl(1-n^{v,l,c}_i\Bigr),
\quad
\forall h \in H,
\forall k \in K_h,
\forall u \in G_{h,k},
v=(h,k,u),
\forall l \in L_v,
\forall i \in S_k,
\forall c \in B_i \cap C_k,
\end{equation}

\begin{equation} \label{eq:35}
e^{v,l,c}_i \leq -r^{v,l,c}_i + M \cdot n^{v,l,c}_i,
\quad
\forall h \in H,
\forall k \in K_h,
\forall u \in G_{h,k},
v=(h,k,u),
\forall l \in L_v,
\forall i \in S_k,
\forall c \in B_i \cap C_k,
\end{equation}

{\noindent \textbf{Calculating Departure Times (at Transhipment Ports)}:}

% \subsubsection{Calculating departure times}

% At Transhipment Ports:
\begin{equation} \label{eq:36}
d^{v,l}_i =
a_i^{v,l} +
w_i^{v,l} +
\sum_{c \in C_k \cap B_i}
\Bigl( e^{v,l,c}_i \cdot U^{k,c} \Bigr),
\quad
\forall h \in H,
\forall k \in K_h,
\forall u \in G_{h,k},
v=(h,k,u),
\forall l \in L_v,
\forall i \in S_k,
\end{equation}

{\noindent \textbf{Calculating Departure Times (at Nodes)}:}
% At Nodes:
\begin{subequations} \label{eq:37}
\begin{multline} \label{eq:37a}
d^{v,l}_i =
a_i^{v,l} +
\sum_{c \in C_k \cap P}
\Biggl( U^{k,c} \cdot \biggl(
y^{v,l,l+1,c}_i
- y^{v,l+1,l,c}_i +
\sum_{\substack{j \in V_k,h \\ i \neq j}}
\Bigl( y^{v,l,c}_{i,j}-y^{v,l,c}_{j,i} \Bigr)
\biggr)\Biggr)
+ \sum_{c \in C_k \cap D}
\Biggl( U^{k,c} \cdot \biggl(
y^{v,l+1,l,c}_i - y^{v,l,l+1,c}_i
+ \sum_{\substack{j \in V_k,h \\ i \neq j}}
\Bigl( y^{v,l,c}_{j,i}\\
-y^{v,l,c}_{i,j} \Bigr)
\biggr)\Biggr),\quad
\forall h \in H,
\forall k \in K_h,
\forall u \in G_{h,k},
v=(h,k,u),
l=1,
\forall i \in N_k,
\end{multline}

\begin{multline} \label{eq:37b}
d^{v,l}_i =
a_i^{v,l} +
\sum_{c \in C_k \cap P}
\Biggl( U^{k,c} \cdot \biggl(
y^{v,l,l+1,c}_i
+ y^{v,l,l-1,c}_i
- y^{v,l+1,l,c}_i
- y^{v,l-1,l,c}_i
+ \sum_{\substack{j \in V_k,h \\ i \neq j}} y^{v,l,c}_{i,j}
- \sum_{\substack{j \in V_k \\ i \neq j}}y^{v,l,c}_{j,i} 
\biggr)\Biggr) +
\sum_{c \in C_k \cap D}
\Biggl( U^{k,c} \cdot \biggl(
y^{v,l+1,l,c}_i \\
+ y^{v,l-1,l,c}_i
- y^{v,l,l+1,c}_i
- y^{v,l,l-1,c}_i +
\sum_{\substack{j \in V_k \\ i \neq j}} y^{v,l,c}_{j,i} -
\sum_{\substack{j \in V_k,h \\ i \neq j}} y^{v,l,c}_{i,j}
\biggr)\Biggr),
\quad
\forall h \in H,
\forall k \in K_h,
\forall u \in G_{h,k},
v=(h,k,u),\\
\forall l \in L_v \setminus \{1,l_v^{max}\},
\forall i \in N_k,
\end{multline}

\begin{multline} \label{eq:37c}
d^{v,l}_i =
a_i^{v,l} +
\sum_{c \in C_k \cap P}
\Biggl( U^{k,c} \cdot \biggl(
y^{v,l,l-1,c}_i
- y^{v,l-1,l,c}_i +
\sum_{\substack{j \in V_k,h \\ i \neq j}} y^{v,l,c}_{i,j}
- \sum_{\substack{j \in V_k \\ i \neq j}}y^{v,l,c}_{j,i} 
\biggr)\Biggr) +
\sum_{c \in C_k \cap D}
\Biggl( U^{k,c} \cdot \biggl(
y^{v,l-1,l,c}_i
- y^{v,l,l-1,c}_i +
\sum_{\substack{j \in V_k \\ i \neq j}} y^{v,l,c}_{j,i} \\
-\sum_{\substack{j \in V_k,h \\ i \neq j}} y^{v,l,c}_{i,j}
\biggr)\Biggr),
\quad
\forall h \in H,
\forall k \in K_h,
\forall u \in G_{h,k},
v=(h,k,u),
if~l_v^{max}>1,
l=l_v^{max},
\forall i \in N_k,
\end{multline}
\end{subequations}

{\noindent \textbf{Calculating Departure Times (at Warehouses)}:}
% At Warehouses:
\begin{subequations} \label{eq:38}
\begin{multline} \label{eq:38a}
d^{v,l}_i =
a_i^{v,l} +
\sum_{c \in C_k \cap D}
\Biggl( U^{k,c} \cdot \biggl(
y^{v,l,l+1,c}_i
- y^{v,l+1,l,c}_i +
\sum_{\substack{j \in V_k \\ i \neq j}} y^{v,l,c}_{i,j}
- \sum_{\substack{j \in V_k,h \\ i \neq j}} y^{v,l,c}_{j,i}
\biggr)\Biggr),
\quad
\forall h \in H,
\forall k \in K_h,
\forall u \in G_{h,k},
v=(h,k,u),\\
l=1,
\forall i \in W_k,
\end{multline}

\begin{multline} \label{eq:38b}
d^{v,l}_i =
a_i^{v,l} +
\sum_{c \in C_k \cap D}
\Biggl( U^{k,c} \cdot \biggl(
y^{v,l,l+1,c}_i
+ y^{v,l,l-1,c}_i
- y^{v,l+1,l,c}_i
- y^{v,l-1,l,c}_i
+ \sum_{\substack{j \in V_k \\ i \neq j}}
\Bigl(
y^{v,l,c}_{i,j} - y^{v,l,c}_{j,i}
\Bigr)
\biggr)\Biggr),
\quad
\forall h \in H,
\forall k \in K_h,\\
\forall u \in G_{h,k},
v=(h,k,u),
\forall l \in L_v \setminus \{1,l_v^{max}\},
\forall i \in W_k,
\end{multline}

\begin{multline} \label{eq:38c}
d^{v,l}_i =
a_i^{v,l} +
\sum_{c \in C_k \cap D}
\Biggl( U^{k,c} \cdot \biggl(
y^{v,l,l-1,c}_i
- y^{v,l-1,l,c}_i
+ \sum_{\substack{j \in V_k \\ i \neq j}}
\Bigl(
y^{v,l,c}_{i,j} - y^{v,l,c}_{j,i} 
\Bigr)
\biggr)\Biggr),
\quad
\forall h \in H,
\forall k \in K_h,
\forall u \in G_{h,k},
v=(h,k,u),\\
if~l_v^{max}>1,
l=l_v^{max},
\forall i \in W_k,
\end{multline}
\end{subequations}

{\noindent \textbf{Calculating Departure Times (at Relief Centers)}:}
% At Relief Centres:
\begin{subequations} \label{eq:39}

\begin{multline} \label{eq:39a}
d^{v,l}_i =
a_i^{v,l} +
\sum_{c \in C_k \cap P}
\Biggl( U^{k,c} \cdot \biggl(
y^{v,l+1,l,c}_i - y^{v,l,l+1,c}_i
+ \sum_{\substack{j \in V_k \\ i \neq j}} y^{v,l,c}_{j,i}
- \sum_{\substack{j \in V_k,h \\ i \neq j}} y^{v,l,c}_{i,j}
\biggr)\Biggr),
\quad
\forall h \in H,
\forall k \in K_h,
\forall u \in G_{h,k},
v=(h,k,u),\\
l=1,
\forall i \in R_k,
\end{multline}

\begin{multline} \label{eq:39b}
d^{v,l}_i =
a_i^{v,l} +
\sum_{c \in C_k \cap P}
\Biggl( U^{k,c} \cdot \biggl(
y^{v,l+1,l,c}_i
+ y^{v,l-1,l,c}_i
- y^{v,l,l+1,c}_i
- y^{v,l,l-1,c}_i +
\sum_{\substack{j \in V_k \\ i \neq j}} y^{v,l,c}_{j,i} -
\sum_{\substack{j \in V_k,h \\ i \neq j}} y^{v,l,c}_{i,j}
\biggr)\Biggr),
\quad
\forall h \in H,
\forall k \in K_h,\\
\forall u \in G_{h,k},
v=(h,k,u),
\forall l \in L_v \setminus \{1,l_v^{max}\},
\forall i \in R_k,
\end{multline}

\begin{multline} \label{eq:39c}
d^{v,l}_i =
a_i^{v,l} +
\sum_{c \in C_k \cap P}
\Biggl( U^{k,c} \cdot \biggl(
y^{v,l-1,l,c}_i
- y^{v,l,l-1,c}_i +
\sum_{\substack{j \in V_k \\ i \neq j}} y^{v,l,c}_{j,i} -
\sum_{\substack{j \in V_k,h \\ i \neq j}} y^{v,l,c}_{i,j}
\biggr)\Biggr),
\quad
\forall h \in H,
\forall k \in K_h,
\forall u \in G_{h,k},
v=(h,k,u),\\
if~l_v^{max}>1,
l=l_v^{max},
\forall i \in R_k,
\end{multline}
\end{subequations}

{\noindent \textbf{Comparing Arrival and Departure Times of Vehicles at Transhipment Ports (Temporal Residue Constraints)}:}
% \subsection{Temporal Residue Constraints} \label{sec:3.4}
% The residues need to be arranged \wrt time so that we may ensure that at any point in time, the sum of residues of a specific load-type at a specific $TP$ is never negative.

% \subsubsection{Comparing Arrival and Departure times of vehicles at Transhipment Ports} \label{sec:3.4.1}

% Populating the intermediate binary variable $g$ when $s=Dep$ and $\hat{s}=Arr$:
    \begin{subequations} \label{eq:40}
        \begin{multline} \label{eq:40a}
        a_i^{\hat{v},\hat{l}} - d^{v,l}_i \geq
        -M \cdot \Bigl( 1 - g_{i,(\hat{v},\hat{l},\hat{s})}^{(v,l,s)} \Bigr),
        \quad
        \forall i \in S,
        \forall h \in H,
        \forall k \in K_h \cap A_{B_i} \cap K_i,
        \forall u \in G_{h,k}, v=(h,k,u),
        \forall \hat{h} \in H,\\
        \forall \hat{k} \in K_{\hat{h}} \cap A_{B_i} \cap K_i,
        \forall \hat{u} \in G_{\hat{h},\hat{k}},
        \hat{v}=(\hat{h},\hat{k},\hat{u}),
        \forall l \in L_v,
        \forall \hat{l} \in L_{\hat{v}},
        \forall s \in \{Dep\},
        \forall \hat{s} \in \{Arr\},
        \end{multline}
        
        \begin{multline} \label{eq:40b}
        d^{v,l}_i - a_i^{\hat{v},\hat{l}} \geq
        -M \cdot g_{i,(\hat{v},\hat{l},\hat{s})}^{(v,l,s)},
        \quad
        \forall i \in S,
        \forall h \in H,
        \forall k \in K_h \cap A_{B_i} \cap K_i,
        \forall u \in G_{h,k}, v=(h,k,u),
        \forall \hat{h} \in H,
        \forall \hat{k} \in K_{\hat{h}} \cap A_{B_i} \cap K_i,\\
        \forall \hat{u} \in G_{\hat{h},\hat{k}},
        \hat{v}=(\hat{h},\hat{k},\hat{u}),
        \forall l \in L_v,
        \forall \hat{l} \in L_{\hat{v}},
        \forall s \in \{Dep\},
        \forall \hat{s} \in \{Arr\},
        \end{multline}
    \end{subequations}

% Populating the intermediate binary variable $g$ when $s=Dep$ and $\hat{s}=Dep$:
    \begin{subequations} \label{eq:41}
        \begin{multline} \label{eq:41a}
        d_i^{\hat{v},\hat{l}} - d^{v,l}_i \geq
        -M \cdot \Bigl( 1 - g_{i,(\hat{v},\hat{l},\hat{s})}^{(v,l,s)} \Bigr),
        \quad
        \forall i \in S,
        \forall h \in H,
        \forall k \in K_h \cap A_{B_i} \cap K_i,
        \forall u \in G_{h,k}, v=(h,k,u),
        \forall \hat{h} \in H,\\
        \forall \hat{k} \in K_{\hat{h}} \cap A_{B_i} \cap K_i,
        \forall \hat{u} \in G_{\hat{h},\hat{k}},
        \hat{v}=(\hat{h},\hat{k},\hat{u}),
        \forall l \in L_v,
        \forall \hat{l} \in L_{\hat{v}},
        \forall s \in \{Dep\},
        \forall \hat{s} \in \{Dep\},
        \end{multline}
        \begin{multline} \label{eq:41b}
        d^{v,l}_i - d_i^{\hat{v},\hat{l}} \geq
        -M \cdot g_{i,(\hat{v},\hat{l},\hat{s})}^{(v,l,s)},
        \quad
        \forall i \in S,
        \forall h \in H,
        \forall k \in K_h \cap A_{B_i} \cap K_i,
        \forall u \in G_{h,k}, v=(h,k,u),
        \forall \hat{h} \in H,
        \forall \hat{k} \in K_{\hat{h}} \cap A_{B_i} \cap K_i,\\
        \forall \hat{u} \in G_{\hat{h},\hat{k}},
        \hat{v}=(\hat{h},\hat{k},\hat{u}),
        \forall l \in L_v,
        \forall \hat{l} \in L_{\hat{v}},
        \forall s \in \{Dep\},
        \forall \hat{s} \in \{Dep\},
        \end{multline}
    \end{subequations}

% Populating the intermediate binary variable $g$ when $s=Arr$ and $\hat{s}=Arr$:
    \begin{subequations} \label{eq:42}
        \begin{multline} \label{eq:42a}
        a_i^{\hat{v},\hat{l}} - a^{v,l}_i \leq
        M \cdot \Bigl( 1 - g_{i,(\hat{v},\hat{l},\hat{s})}^{(v,l,s)} \Bigr),
        \quad
        \forall i \in S,
        \forall h \in H,
        \forall k \in K_h \cap A_{B_i} \cap K_i,
        \forall u \in G_{h,k}, v=(h,k,u),
        \forall \hat{h} \in H,\\
        \forall \hat{k} \in K_{\hat{h}} \cap A_{B_i} \cap K_i,
        \forall \hat{u} \in G_{\hat{h},\hat{k}},
        \hat{v}=(\hat{h},\hat{k},\hat{u}),
        \forall l \in L_v,
        \forall \hat{l} \in L_{\hat{v}},
        \forall s \in \{Arr\},
        \forall \hat{s} \in \{Arr\},
        \end{multline}
        
        \begin{multline} \label{eq:42b}
        a^{v,l}_i - a_i^{\hat{v},\hat{l}} \leq
        M \cdot g_{i,(\hat{v},\hat{l},\hat{s})}^{(v,l,s)},
        \quad
        \forall i \in S,
        \forall h \in H,
        \forall k \in K_h \cap A_{B_i} \cap K_i,
        \forall u \in G_{h,k}, v=(h,k,u),
        \forall \hat{h} \in H,
        \forall \hat{k} \in K_{\hat{h}} \cap A_{B_i} \cap K_i,\\
        \forall \hat{u} \in G_{\hat{h},\hat{k}},
        \hat{v}=(\hat{h},\hat{k},\hat{u}),
        \forall l \in L_v,
        \forall \hat{l} \in L_{\hat{v}},
        \forall s \in \{Arr\},
        \forall \hat{s} \in \{Arr\},
        \end{multline}
    \end{subequations}

% Populating the intermediate binary variable $g$ when $s=Arr$ and $\hat{s}=Dep$:
    \begin{subequations} \label{eq:43}
        \begin{multline} \label{eq:43a}
        d_i^{\hat{v},\hat{l}} - a^{v,l}_i \leq
        M \cdot \Bigl( 1 - g_{i,(\hat{v},\hat{l},\hat{s})}^{(v,l,s)} \Bigr),
        \quad
        \forall i \in S,
        \forall h \in H,
        \forall k \in K_h \cap A_{B_i} \cap K_i,
        \forall u \in G_{h,k}, v=(h,k,u),
        \forall \hat{h} \in H,\\
        \forall \hat{k} \in K_{\hat{h}} \cap A_{B_i} \cap K_i,
        \forall \hat{u} \in G_{\hat{h},\hat{k}},
        \hat{v}=(\hat{h},\hat{k},\hat{u}),
        \forall l \in L_v,
        \forall \hat{l} \in L_{\hat{v}},
        \forall s \in \{Arr\},
        \forall \hat{s} \in \{Dep\},
        \end{multline}
        \begin{multline} \label{eq:43b}
        a^{v,l}_i - d_i^{\hat{v},\hat{l}} \leq
        M \cdot g_{i,(\hat{v},\hat{l},\hat{s})}^{(v,l,s)},
        \quad
        \forall i \in S,
        \forall h \in H,
        \forall k \in K_h \cap A_{B_i} \cap K_i,
        \forall u \in G_{h,k}, v=(h,k,u),
        \forall \hat{h} \in H,
        \forall \hat{k} \in K_{\hat{h}} \cap A_{B_i} \cap K_i,\\
        \forall \hat{u} \in G_{\hat{h},\hat{k}},
        \hat{v}=(\hat{h},\hat{k},\hat{u}),
        \forall l \in L_v,
        \forall \hat{l} \in L_{\hat{v}},
        \forall s \in \{Arr\},
        \forall \hat{s} \in \{Dep\},
        \end{multline}
    \end{subequations}

% \subsubsection{Constraining the cell values of space-time matrices when the vehicles being compared don't overlap in spacetime} \label{sec:3.4.2}
% %\section{Example Section} \label{sec:example-section} \ref{sec:example-section}

% The cell values of the space-time matrices consist of the variables $o$ which store some proportion of the load-component specific residue $r$ available for transhipment during the arrival and departure time-stamps of all vehicles. This uses the comparison of a vehicle's visiting time with another at the same Transhipment Port for all its visiting possibilities individually, \ie across the various levels. The \autoref{sec:3.4.2} considers the cases when the visiting durations of the vehicles under comparison don't overlap, \ie they concerned vehicles don't meet during their corresponding layer-specific visit at the concerned Transhipment Port $S$. Due to this non-overlapping nature, it is easy to determine the available proportion of residue, the cell value $o$; which would either be the entire residue of the vehicle under comparison $v$ \wrt a reference vehicle $\hat{v}$ if $v$ has already departed the transhipment port with some residue value $r$ during its specific layer-wise visit $l$; otherwise its value would be $0$ if the concerned vehicle $v$ is yet to arrive in the transhipment port \wrt the reference vehicle $\hat{v}$'s visit duration.

% Populating variable $o$ as equal to the residue when $g=1$, for $s=Dep$ and $\hat{s}=Arr$:

{\noindent \textbf{Constraining the cell values of Space-Time Matrices when the Vehicles being compared don't overlap in SpaceTime (Temporal Residue Constraints)}:}
    \begin{subequations} \label{eq:44}
        \begin{multline} \label{eq:44a}
        o_{i,(\hat{v},\hat{l},\hat{s})}^{c,(v,l)} \leq
        r_i^{v,l,c} + M \cdot \Bigl( 1 - g_{i,(\hat{v},\hat{l},\hat{s})}^{(v,l,s)}\Bigr),
        \quad
        \forall i \in S,
        \forall c \in B_i,
        \forall h \in H,
        \forall k \in K_h \cap A_c \cap K_i,
        \forall u \in G_{h,k},
        v=(h,k,u),
        \forall \hat{h} \in H,\\
        \forall \hat{k} \in K_{\hat{h}} \cap A_c \cap K_i,
        \forall \hat{u} \in G_{\hat{h},\hat{k}},
        \hat{v}=(\hat{h},\hat{k},\hat{u}),
        \forall l \in L_v,
        \forall s \in \{Dep\},
        \forall \hat{l} \in L_{\hat{v}},
        \forall \hat{s} \in \{Arr\},
        \end{multline}
        
        \begin{multline} \label{eq:44b}
        o_{i,(\hat{v},\hat{l},\hat{s})}^{c,(v,l)} \geq
        r_i^{v,l,c} - M \cdot \Bigl( 1 - g_{i,(\hat{v},\hat{l},\hat{s})}^{(v,l,s)} \Bigr),
        \quad
        \forall i \in S,
        \forall c \in B_i,
        \forall h \in H,
        \forall k \in K_h \cap A_c \cap K_i,
        \forall u \in G_{h,k},
        v=(h,k,u),
        \forall \hat{h} \in H,\\
        \forall \hat{k} \in K_{\hat{h}} \cap A_c \cap K_i,
        \forall \hat{u} \in G_{\hat{h},\hat{k}},
        \hat{v}=(\hat{h},\hat{k},\hat{u}),
        \forall l \in L_v,
        \forall s \in \{Dep\},
        \forall \hat{l} \in L_{\hat{v}},
        \forall \hat{s} \in \{Arr\},
        \end{multline}
    \end{subequations}

% Populating variable $o$ as equal to the residue when $g=1$, for $s=Dep$ and $\hat{s}=Dep$:
    \begin{subequations} \label{eq:45}
        \begin{multline} \label{eq:45a}
        o_{i,(\hat{v},\hat{l},\hat{s})}^{c,(v,l)} \leq
        r_i^{v,l,c} + M \cdot \Bigl( 1 - g_{i,(\hat{v},\hat{l},\hat{s})}^{(v,l,s)} \Bigr),
        \quad
        \forall i \in S,
        \forall c \in B_i,
        \forall h \in H,
        \forall k \in K_h \cap A_c \cap K_i,
        \forall u \in G_{h,k},
        v=(h,k,u),
        \forall \hat{h} \in H,\\
        \forall \hat{k} \in K_{\hat{h}} \cap A_c \cap K_i,
        \forall \hat{u} \in G_{\hat{h},\hat{k}},
        \hat{v}=(\hat{h},\hat{k},\hat{u}),
        \forall l \in L_v,
        \forall s \in \{Dep\},
        \forall \hat{l} \in L_{\hat{v}},
        \forall \hat{s} \in \{Dep\},
        \end{multline}
        
        \begin{multline} \label{eq:45b}
        o_{i,(\hat{v},\hat{l},\hat{s})}^{c,(v,l)} \geq
        r_i^{v,l,c} - M \cdot \Bigl( 1 - g_{i,(\hat{v},\hat{l},\hat{s})}^{(v,l,s)} \Bigr),
        \quad
        \forall i \in S,
        \forall c \in B_i,
        \forall h \in H,
        \forall k \in K_h \cap A_c \cap K_i,
        \forall u \in G_{h,k},
        v=(h,k,u),
        \forall \hat{h} \in H,\\
        \forall \hat{k} \in K_{\hat{h}} \cap A_c \cap K_i,
        \forall \hat{u} \in G_{\hat{h},\hat{k}},
        \hat{v}=(\hat{h},\hat{k},\hat{u}),
        \forall l \in L_v,
        \forall s \in \{Dep\},
        \forall \hat{l} \in L_{\hat{v}},
        \forall \hat{s} \in \{Dep\},
        \end{multline}
    \end{subequations}

% Nullifying variable $o$ when $g=1$, for $s=Arr$ and $\hat{s}=Arr$:
    \begin{subequations} \label{eq:46}
        \begin{multline} \label{eq:46a}
        o_{i,(\hat{v},\hat{l},\hat{s})}^{c,(v,l)} \leq
        M \cdot \Bigl( 1 - g_{i,(\hat{v},\hat{l},\hat{s})}^{(v,l,s)} \Bigr),
        \quad
        \forall i \in S,
        \forall c \in B_i,
        \forall h \in H,
        \forall k \in K_h \cap A_c \cap K_i,
        \forall u \in G_{h,k},
        v=(h,k,u),
        \forall \hat{h} \in H,\\
        \forall \hat{k} \in K_{\hat{h}} \cap A_c \cap K_i,
        \forall \hat{u} \in G_{\hat{h},\hat{k}},
        \hat{v}=(\hat{h},\hat{k},\hat{u}),
        \forall l \in L_v,
        \forall s \in \{Arr\},
        \forall \hat{l} \in L_{\hat{v}},
        \forall \hat{s} \in \{Arr\},
        \end{multline}
        
        \begin{multline} \label{eq:46b}
        o_{i,(\hat{v},\hat{l},\hat{s})}^{c,(v,l)} \geq
        - M \cdot \Bigl( 1 - g_{i,(\hat{v},\hat{l},\hat{s})}^{(v,l,s)} \Bigr),
        \quad
        \forall i \in S,
        \forall c \in B_i,
        \forall h \in H,
        \forall k \in K_h \cap A_c \cap K_i,
        \forall u \in G_{h,k},
        v=(h,k,u),
        \forall \hat{h} \in H,\\
        \forall \hat{k} \in K_{\hat{h}} \cap A_c \cap K_i,
        \forall \hat{u} \in G_{\hat{h},\hat{k}},
        \hat{v}=(\hat{h},\hat{k},\hat{u}),
        \forall l \in L_v,
        \forall s \in \{Arr\},
        \forall \hat{l} \in L_{\hat{v}},
        \forall \hat{s} \in \{Arr\},
        \end{multline}
    \end{subequations}

% Nullifying variable $o$ when $g=1$, for $s=Arr$ and $\hat{s}=Dep$:
    \begin{subequations} \label{eq:47}
        \begin{multline} \label{eq:47a}
        o_{i,(\hat{v},\hat{l},\hat{s})}^{c,(v,l)} \leq
        M \cdot \Bigl( 1 - g_{i,(\hat{v},\hat{l},\hat{s})}^{(v,l,s)} \Bigr),
        \quad
        \forall i \in S,
        \forall c \in B_i,
        \forall h \in H,
        \forall k \in K_h \cap A_c \cap K_i,
        \forall u \in G_{h,k},
        v=(h,k,u),
        \forall \hat{h} \in H,\\
        \forall \hat{k} \in K_{\hat{h}} \cap A_c \cap K_i,
        \forall \hat{u} \in G_{\hat{h},\hat{k}},
        \hat{v}=(\hat{h},\hat{k},\hat{u}),
        \forall l \in L_v,
        \forall s \in \{Arr\},
        \forall \hat{l} \in L_{\hat{v}},
        \forall \hat{s} \in \{Dep\},
        \end{multline}
        
        \begin{multline} \label{eq:47b}
        o_{i,(\hat{v},\hat{l},\hat{s})}^{c,(v,l)} \geq
        - M \cdot \Bigl( 1 -g_{i,(\hat{v},\hat{l},\hat{s})}^{(v,l,s)}\Bigr),
        \quad
        \forall i \in S,
        \forall c \in B_i,
        \forall h \in H,
        \forall k \in K_h \cap A_c \cap K_i,
        \forall u \in G_{h,k},
        v=(h,k,u),
        \forall \hat{h} \in H,\\
        \forall \hat{k} \in K_{\hat{h}} \cap A_c \cap K_i,
        \forall \hat{u} \in G_{\hat{h},\hat{k}},
        \hat{v}=(\hat{h},\hat{k},\hat{u}),
        \forall l \in L_v,
        \forall s \in \{Arr\},
        \forall \hat{l} \in L_{\hat{v}},
        \forall \hat{s} \in \{Dep\},
        \end{multline}
    \end{subequations}

% \subsubsection{Constraining the cell values of the space-time matrices when the visits of the vehicles being compared overlap in spacetime, residues being positive} \label{sec:3.4.3}

% Constraints in \autoref{sec:3.4.2} allow an exact value of the variable $o$ when $g=1$. However, when $g=0$ for both $s=Arr$ and $s=Dep$, \ie when the vehicles under comparison overlap during their visiting times at the same transhipment port, those constraints become redundant and the below constraints are developed to populate the $o$ variables of the space-time matrix within certain bounds. The $o$ variable is allowed to take any value within these bounds, but as per the optimization criteria, it should take the upper bound due to minimization and when the residues are positive.

% Constraining the cell values within their limits of $0$ to the maximum value being the residue $r$ itself (this is irrespective of whether the vehicles overlap in spacetime):

{\noindent \textbf{Constraining the cell values of Space-Time Matrices when the visits of the Vehicles overlap in SpaceTime (Residues being positive)}:}
    \begin{multline} \label{eq:48}
        0 \leq o_{i,(\hat{v},\hat{l},\hat{s})}^{c,(v,l)}
        + M \cdot \Bigl( 1-n_i^{v,l,c} \Bigr),
        \quad
        \forall i \in S,
        \forall c \in B_i,
        \forall h \in H,
        \forall k \in K_h \cap A_c \cap K_i,
        \forall u \in G_{h,k},
        v=(h,k,u),
        \forall \hat{h} \in H,\\
        \forall \hat{k} \in K_{\hat{h}} \cap A_c \cap K_i,
        \forall \hat{u} \in G_{\hat{h},\hat{k}},
        \hat{v}=(\hat{h},\hat{k},\hat{u}),
        \forall l \in L_v,
        \forall \hat{l} \in L_{\hat{v}},
        \forall \hat{s} \in \{Arr\},
    \end{multline}
    
    \begin{multline} \label{eq:49}
        o_{i,(\hat{v},\hat{l},Arr)}^{c,(v,l)} \leq o_{i,(\hat{v},\hat{l},Dep)}^{c,(v,l)}
        + M \cdot \Bigl( 1-n_i^{v,l,c} \Bigr),
        \quad
        \forall i \in S,
        \forall c \in B_i,
        \forall h \in H,
        \forall k \in K_h \cap A_c \cap K_i,
        \forall u \in G_{h,k},
        v=(h,k,u),\\
        \forall \hat{h} \in H,
        \forall \hat{k} \in K_{\hat{h}} \cap A_c \cap K_i,
        \forall \hat{u} \in G_{\hat{h},\hat{k}},
        \hat{v}=(\hat{h},\hat{k},\hat{u}),
        \forall l \in L_v,
        \forall \hat{l} \in L_{\hat{v}},
    \end{multline}
    
    \begin{multline} \label{eq:50}
        o_{i,(\hat{v},\hat{l},\hat{s})}^{c,(v,l)} \leq r_i^{v,l,c}
        + M \cdot \Bigl( 1-n_i^{v,l,c} \Bigr),
        \quad
        \forall i \in S,
        \forall c \in B_i,
        \forall h \in H,
        \forall k \in K_h \cap A_c \cap K_i,
        \forall u \in G_{h,k},
        v=(h,k,u),
        \forall \hat{h} \in H,\\
        \forall \hat{k} \in K_{\hat{h}} \cap A_c \cap K_i,
        \forall \hat{u} \in G_{\hat{h},\hat{k}},
        \hat{v}=(\hat{h},\hat{k},\hat{u}),
        \forall l \in L_v,
        \forall \hat{l} \in L_{\hat{v}},
        \forall \hat{s} \in \{Dep\},
    \end{multline}

% Upper bounds:
    \begin{multline} \label{eq:51}
        o_{i,(\hat{v},\hat{l},Arr)}^{c,(v,l)} \leq
        \dfrac{a_i^{\hat{v},\hat{l}}-a_i^{v,l}}{U^{k,c}}
        + M \cdot \Bigl(
        g_{i,(\hat{v},\hat{l},Arr)}^{(v,l,Arr)}
        + g_{i,(\hat{v},\hat{l},Arr)}^{(v,l,Dep)}
        \Bigr)
        + M \cdot \Bigl( 1-n_i^{v,l,c} \Bigr),
        \quad
        \forall i \in S,
        \forall c \in B_i,
        \forall h,\hat{h} \in H,
        \forall k \in K_h \cap A_c \cap K_i,\\
        \forall u \in G_{h,k},
        v=(h,k,u),
        \forall l \in L_v,
        \forall \hat{k} \in K_{\hat{h}} \cap A_c \cap K_i,
        \forall \hat{u} \in G_{\hat{h},\hat{k}},
        \hat{v}=(\hat{h},\hat{k},\hat{u}),
        \forall \hat{l} \in L_{\hat{v}},
        if~U^{k,c}>0,
    \end{multline}
    
    \begin{multline} \label{eq:52}
        o_{i,(\hat{v},\hat{l},Dep)}^{c,(v,l)} \leq
        \dfrac{d_i^{\hat{v},\hat{l}}-a_i^{v,l}}{U^{k,c}}
        + M \cdot \Bigl(
        g_{i,(\hat{v},\hat{l},Dep)}^{(v,l,Arr)}
        + g_{i,(\hat{v},\hat{l},Dep)}^{(v,l,Dep)}
        \Bigr)
        + M \cdot \Bigl( 1-n_i^{v,l,c} \Bigr),
        \quad
        \forall i \in S,
        \forall c \in B_i,
        \forall h,\hat{h} \in H,
        \forall k \in K_h \cap A_c \cap K_i,\\
        \forall u \in G_{h,k},
        v=(h,k,u),
        \forall l \in L_v,
        \forall \hat{k} \in K_{\hat{h}} \cap A_c \cap K_i,
        \forall \hat{u} \in G_{\hat{h},\hat{k}},
        \hat{v}=(\hat{h},\hat{k},\hat{u}),
        \forall \hat{l} \in L_{\hat{v}},
        if~U^{k,c}>0,
    \end{multline}

% Lower Bounds:
    \begin{multline} \label{eq:53}
        o_{i,(\hat{v},\hat{l},Arr)}^{c,(v,l)} \geq
        r_i^{v,l,c} - \dfrac{d_i^{v,l}-a_i^{\hat{v},\hat{l}}}{U^{k,c}}
        - M \cdot \Bigl(
        g_{i,(\hat{v},\hat{l},Arr)}^{(v,l,Arr)}
        + g_{i,(\hat{v},\hat{l},Arr)}^{(v,l,Dep)}
        \Bigr)
        - M \cdot \Bigl( 1-n_i^{v,l,c} \Bigr),
        \quad
        \forall i \in S,
        \forall c \in B_i,
        \forall h,\hat{h} \in H,
        \forall k \in K_h \cap A_c \cap K_i,\\
        \forall u \in G_{h,k},
        v=(h,k,u),
        \forall l \in L_v,
        \forall \hat{k} \in K_{\hat{h}} \cap A_c \cap K_i,
        \forall \hat{u} \in G_{\hat{h},\hat{k}},
        \hat{v}=(\hat{h},\hat{k},\hat{u}),
        \forall \hat{l} \in L_{\hat{v}},
        if~U^{k,c}>0,
    \end{multline}
    
    \begin{multline} \label{eq:54}
        o_{i,(\hat{v},\hat{l},Dep)}^{c,(v,l)} \geq
        r_i^{v,l,c} - \dfrac{d_i^{v,l}-d_i^{\hat{v},\hat{l}}}{U^{k,c}}
        - M \cdot \Bigl(
        g_{i,(\hat{v},\hat{l},Dep)}^{(v,l,Arr)}
        + g_{i,(\hat{v},\hat{l},Dep)}^{(v,l,Dep)}
        \Bigr)
        - M \cdot \Bigl( 1-n_i^{v,l,c} \Bigr),
        \quad
        \forall i \in S,
        \forall c \in B_i,
        \forall h,\hat{h} \in H,
        \forall k \in K_h \cap A_c \cap K_i,\\
        \forall u \in G_{h,k},
        v=(h,k,u),
        \forall l \in L_v,
        \forall \hat{k} \in K_{\hat{h}} \cap A_c \cap K_i,
        \forall \hat{u} \in G_{\hat{h},\hat{k}},
        \hat{v}=(\hat{h},\hat{k},\hat{u}),
        \forall \hat{l} \in L_{\hat{v}},
        if~U^{k,c}>0,
    \end{multline}

% \subsubsection{Constraining the cell values of the space-time matrices when the visits of the vehicles being compared overlap in spacetime, residues being negative} \label{sec:3.4.4}

% Constraining the cell values within their limits of $0$ to the minimum value being the residue $r$ itself (this is irrespective of whether the vehicles overlap in spacetime):

{\noindent \textbf{Constraining the cell values of Space-Time Matrices when the visits of the vehicles being compared overlap in SpaceTime (Residues being negative)}:}
    \begin{multline} \label{eq:55}
        0 \geq o_{i,(\hat{v},\hat{l},\hat{s})}^{c,(v,l)}
        - M \cdot n_i^{v,l,c},
        \quad
        \forall i \in S,
        \forall c \in B_i,
        \forall h \in H,
        \forall k \in K_h \cap A_c \cap K_i,
        \forall u \in G_{h,k},
        v=(h,k,u),
        \forall \hat{h} \in H,\\
        \forall \hat{k} \in K_{\hat{h}} \cap A_c \cap K_i,
        \forall \hat{u} \in G_{\hat{h},\hat{k}},
        \hat{v}=(\hat{h},\hat{k},\hat{u}),
        \forall l \in L_v,
        \forall \hat{l} \in L_{\hat{v}},
        \forall \hat{s} \in \{Arr\},
    \end{multline}
    \begin{multline} \label{eq:56}
        o_{i,(\hat{v},\hat{l},Arr)}^{c,(v,l)} \geq o_{i,(\hat{v},\hat{l},Dep)}^{c,(v,l)}
        - M \cdot n_i^{v,l,c},
        \quad
        \forall i \in S,
        \forall c \in B_i,
        \forall h \in H,
        \forall k \in K_h \cap A_c \cap K_i,
        \forall u \in G_{h,k},
        v=(h,k,u),
        \forall \hat{h} \in H,\\
        \forall \hat{k} \in K_{\hat{h}} \cap A_c \cap K_i,
        \forall \hat{u} \in G_{\hat{h},\hat{k}},
        \hat{v}=(\hat{h},\hat{k},\hat{u}),
        \forall l \in L_v,
        \forall \hat{l} \in L_{\hat{v}},
    \end{multline}
    \begin{multline} \label{eq:57}
        o_{i,(\hat{v},\hat{l},\hat{s})}^{c,(v,l)} \geq r_i^{v,l,c}
        - M \cdot n_i^{v,l,c},
        \quad
        \forall i \in S,
        \forall c \in B_i,
        \forall h \in H,
        \forall k \in K_h \cap A_c \cap K_i,
        \forall u \in G_{h,k},
        v=(h,k,u),
        \forall \hat{h} \in H,\\
        \forall \hat{k} \in K_{\hat{h}} \cap A_c \cap K_i,
        \forall \hat{u} \in G_{\hat{h},\hat{k}},
        \hat{v}=(\hat{h},\hat{k},\hat{u}),
        \forall l \in L_v,
        \forall \hat{l} \in L_{\hat{v}},
        \forall \hat{s} \in \{Dep\},
    \end{multline}

% Lower bounds:
    \begin{multline} \label{eq:58}
        o_{i,(\hat{v},\hat{l},\hat{s})}^{c,(v,l)} \geq
        \dfrac{a_i^{\hat{v},\hat{l}}-a_i^{v,l}}{-U^{k,c}}
        - M \cdot \Bigl(
        g_{i,(\hat{v},\hat{l},\hat{s})}^{(v,l,Arr)}
        + g_{i,(\hat{v},\hat{l},\hat{s})}^{(v,l,Dep)}
        \Bigr)
        - M \cdot n_i^{v,l,c},
        \quad
        \forall i \in S,
        \forall c \in B_i,
        \forall h,\hat{h} \in H,
        \forall k \in K_h \cap A_c \cap K_i,\\
        \forall u \in G_{h,k},
        v=(h,k,u),
        \forall l \in L_v,
        \forall \hat{k} \in K_{\hat{h}} \cap A_c \cap K_i,
        \forall \hat{u} \in G_{\hat{h},\hat{k}},
        \hat{v}=(\hat{h},\hat{k},\hat{u}),
        \forall \hat{l} \in L_{\hat{v}},
        \forall \hat{s} \in \{Arr\},
        if~U^{k,c}>0,
    \end{multline}
    
    \begin{multline} \label{eq:59}
        o_{i,(\hat{v},\hat{l},\hat{s})}^{c,(v,l)} \geq
        \dfrac{d_i^{\hat{v},\hat{l}}-a_i^{v,l}}{-U^{k,c}}
        - M \cdot \Bigl(
        g_{i,(\hat{v},\hat{l},\hat{s})}^{(v,l,Arr)}
        + g_{i,(\hat{v},\hat{l},\hat{s})}^{(v,l,Dep)}
        \Bigr)
        - M \cdot n_i^{v,l,c},
        \quad
        \forall i \in S,
        \forall c \in B_i,
        \forall h,\hat{h} \in H,
        \forall k \in K_h \cap A_c \cap K_i,\\
        \forall u \in G_{h,k},
        v=(h,k,u),
        \forall l \in L_v,
        \forall \hat{k} \in K_{\hat{h}} \cap A_c \cap K_i,
        \forall \hat{u} \in G_{\hat{h},\hat{k}},
        \hat{v}=(\hat{h},\hat{k},\hat{u}),
        \forall \hat{l} \in L_{\hat{v}},
        \forall \hat{s} \in \{Dep\},
        if~U^{k,c}>0,
    \end{multline}

% Upper Bounds:
    \begin{multline} \label{eq:60}
        o_{i,(\hat{v},\hat{l},\hat{s})}^{c,(v,l)} \leq
        r_i^{v,l,c} + \dfrac{d_i^{v,l}-a_i^{\hat{v},\hat{l}}}{U^{k,c}}
        + M \cdot \Bigl(
        g_{i,(\hat{v},\hat{l},\hat{s})}^{(v,l,Arr)}
        + g_{i,(\hat{v},\hat{l},\hat{s})}^{(v,l,Dep)}
        \Bigr)
        + M \cdot n_i^{v,l,c},
        \quad
        \forall i \in S,
        \forall c \in B_i,
        \forall h,\hat{h} \in H,
        \forall k \in K_h \cap A_c \cap K_i,\\
        \forall u \in G_{h,k},
        v=(h,k,u),
        \forall l \in L_v,
        \forall \hat{k} \in K_{\hat{h}} \cap A_c \cap K_i,
        \forall \hat{u} \in G_{\hat{h},\hat{k}},
        \hat{v}=(\hat{h},\hat{k},\hat{u}),
        \forall \hat{l} \in L_{\hat{v}},
        \forall \hat{s} \in \{Arr\},
        if~U^{k,c}>0,
    \end{multline}   
    \begin{multline} \label{eq:61}
        o_{i,(\hat{v},\hat{l},\hat{s})}^{c,(v,l)} \leq
        r_i^{v,l,c} + \dfrac{d_i^{v,l}-d_i^{\hat{v},\hat{l}}}{U^{k,c}}
        + M \cdot \Bigl(
        g_{i,(\hat{v},\hat{l},\hat{s})}^{(v,l,Arr)}
        + g_{i,(\hat{v},\hat{l},\hat{s})}^{(v,l,Dep)}
        \Bigr)
        + M \cdot n_i^{v,l,c},
        \quad
        \forall i \in S,
        \forall c \in B_i,
        \forall h,\hat{h} \in H,
        \forall k \in K_h \cap A_c \cap K_i,\\
        \forall u \in G_{h,k},
        v=(h,k,u),
        \forall l \in L_v,
        \forall \hat{k} \in K_{\hat{h}} \cap A_c \cap K_i,
        \forall \hat{u} \in G_{\hat{h},\hat{k}},
        \hat{v}=(\hat{h},\hat{k},\hat{u}),
        \forall \hat{l} \in L_{\hat{v}},
        \forall \hat{s} \in \{Dep\},
        if~U^{k,c}>0,
    \end{multline}

% \subsubsection{Constraining any collection to occur at transhipment ports only when the concerned material is available, \ie previously deposited by some other vehicle at the same $S$} \label{sec:3.4.5}

% The \autoref{sec:3.4} is used to develop the space-time matrix, for each load-type at each Transhipment Port $S$, considering compatibilities for allowable transhipments, populated by respective $o$ variables.
% The typical space-time matrix may be imagined as a matrix with rows as various space-stamps, each compatible vehicles' all levels in our case; and columns as various time-stamps. For each vehicle level in the rows, there are two time-stamps namely arrival and departure. Therefore all space-time matrices constructed by us are rectangular in nature with columns being twice that of rows, \ie the time-stamps correspond to the arrival and departure times of the compatible vehicles' levels.
% To ensure any collection happens only with available resources at transhipment ports, we constrain the sum of the cells across each time axis (each column) for every space-time matrix to be non-negative. This ensures that at every time-stamp the total component-specific residue at any $TP$ never becomes negative, \ie no collection of loads takes place before the load was deposited.

% Constraining the cascaded-temporal sum of residues to be non-negative at each time-stamp, for every load-type at every $S$:

{\noindent \textbf{Constraining any collection to occur at Transhipment Ports only when the concerned material is available (i.e. previously deposited by some other vehicle at the same TP)}:}
    \begin{multline} \label{eq:62}
        \sum_{\substack{h \in H,\\
        k \in K_h \cap A_c \cap K_i,\\
        u \in G_{h,k},\\
        v=(h,k,u)}}
        ~\sum_{l \in L_v}
        o_{i,(\hat{v},\hat{l},\hat{s})}^{c,(v,l)}
        \geq 0,
        \quad
        \forall i \in S,
        \forall c \in B_i,
        \forall \hat{h} \in H,
        \forall \hat{k} \in K_{\hat{h}} \cap A_c \cap K_i,
        \forall \hat{u} \in G_{\hat{h},\hat{k}},
        \hat{v}=(\hat{h},\hat{k},\hat{u}),
        \forall \hat{l} \in L_{\hat{v}},\\
        \forall \hat{s} \in \{Arr,Dep\},
    \end{multline}

% \subsubsection{Constraining final residues at Transhipment Ports to be null} \label{sec:4.5.6}

% Since Transhipment Ports are modelled to allow \ul{the transfer of (compatible) cargos during the operation period}, they should not contain any resource (pickup or delivery load type) finally at the end of the emergency relief-and-rescue operation. Therefore, we need to provide an additional constraint which is set so as to nullify the sum of all residues (of each resource type separately) at each Transhipment Port. Due to the nature of the problem, this is necessary to be constrained for only the PickUp Load Types, since \autoref{eq:62} is enough for the Delivery Load Types. If this constraint is not provided, PickUp quantities after transhipment initiation will remain at the respective Transhipment Ports instead of finally being deposited at Relief Centres.

% All residues at Transhipment Ports are summed up and equalled to zero:

{\noindent \textbf{Constraining final Residues at Transhipment Ports to be null}:}
    \begin{equation} \label{eq:62.5}
        \sum_{\substack{h \in H,\\
        k \in K_h \cap A_c \cap K_i,\\
        u \in G_{h,k},\\
        v=(h,k,u)}}
        ~\sum_{l \in L_v}
        r^{v,l,c}_i
        \leq 0,
        \quad
        \forall i \in S,
        \forall c \in B_i \cap P,
    \end{equation}

% In some instances, instead of using $\forall c \in B_i \cap P$ in \autoref{eq:62.5}, we do use $\forall c \in B_i$; this is mainly done to check the formulation run times for different combinations of constraint usages.

{\noindent \textbf{Vehicle Capacity (Volume) Constraints}:}

% \subsection{Vehicle Capacity Constraints}  \label{sec:3.5}
    
% Constraining all flow variables associated with an edge \wrt the total available vehicle-volume:
    \begin{subequations} \label{eq:63}
    \begin{multline} \label{eq:63a}
        \sum_{c \in C_k}
        \Bigl( E_c \cdot y^{v,l,c}_{i,j} \Bigr)
        \leq E^k \cdot x_{i,j}^{v,l},
        \quad
        \forall h \in H, \forall k \in K_h, \forall u \in G_{h,k}, v=(h,k,u),
        l=1,
        \forall i,j \in V_k \cup h,
        i \neq j, (if~i \in W \Rightarrow j \neq h),\\
        (if~i=h \Rightarrow j \notin R),
    \end{multline}
    
    \begin{multline} \label{eq:63b}
        \sum_{c \in C_k}
        \Bigl( E_c \cdot y^{v,l,c}_{i,j} \Bigr)
        \leq E^k \cdot x_{i,j}^{v,l},
        \quad
       \forall h \in H, \forall k \in K_h, \forall u \in G_{h,k}, v=(h,k,u),
       \forall l \in L_v \smallsetminus \{{1}\}, \forall i \in V_k,
       \forall j \in V_k \cup h,
       i \neq j,\\
       (if~i \in W \Rightarrow j \neq h),
    \end{multline}
    
    \begin{equation} \label{eq:63c}
        \sum_{c \in C_k}
        \Bigl( E_c \cdot y^{v,l,m,c}_i \Bigr)
        \leq E^k \cdot x^{v,l,m}_i,
        \quad
       \forall h \in H, \forall k \in K_h, \forall u \in G_{h,k}, v=(h,k,u),
       \forall l,m \in L_v,|l-m|=1,
       \forall i \in V_k,
    \end{equation}
    \end{subequations}

{\noindent \textbf{Vehicle Capacity (Weight) Constraints}:}
% Constraining all flow variables associated with an edge \wrt the total allowable vehicle-weight:
    \begin{subequations} \label{eq:64}
    \begin{multline} \label{eq:64a}
        \sum_{c \in C_k}
        \Bigl( F_c \cdot y^{v,l,c}_{i,j} \Bigr)
        \leq F^k \cdot x_{i,j}^{v,l},
        \quad
        \forall h \in H, \forall k \in K_h, \forall u \in G_{h,k}, v=(h,k,u),
        l=1,
        \forall i,j \in V_k \cup h,
        i \neq j,
        (if~i \in W \Rightarrow j \neq h),\\
        (if~i=h \Rightarrow j \notin R),
    \end{multline}
    
    \begin{multline} \label{eq:64b}
        \sum_{c \in C_k}
        \Bigl( F_c \cdot y^{v,l,c}_{i,j} \Bigr)
        \leq F^k \cdot x_{i,j}^{v,l},
        \quad
       \forall h \in H, \forall k \in K_h, \forall u \in G_{h,k}, v=(h,k,u),
       \forall l \in L_v \smallsetminus \{{1}\}, \forall i \in V_k,
       \forall j \in V_k \cup h,
       i \neq j,\\
       (if~i \in W \Rightarrow j \neq h),
    \end{multline}
    
    \begin{equation} \label{eq:64c}
        \sum_{c \in C_k}
        \Bigl( F_c \cdot y^{v,l,m,c}_i \Bigr)
        \leq F^k \cdot x^{v,l,m}_i,
        \quad
       \forall h \in H, \forall k \in K_h, \forall u \in G_{h,k}, v=(h,k,u),
       \forall l,m \in L_v,|l-m|=1,
       \forall i \in V_k,
    \end{equation}
    \end{subequations}

{\noindent \textbf{Populating variables used in the objective function}:}
% \subsection{Populating variables for the objective function}\label{sec:3.6}
% Eq. \ref{eq:65} assigns the MakeSpan variable of $z_x$, such that it can take the maximum among all vehicle's route durations.

\begin{subequations} \label{eq:65}
% For the VTs that follow open-tour, the last route back towards the depot is subtracted from the total time:
\begin{equation} \label{eq:65a}
z_x \geq
a_h^{v,l} - \sum_{i \in S_k,N_k,R_k} T^k_{i,h} \cdot x^{v,l}_{i,h},
\quad
\forall h \in H, \forall k \in K_h \cap O, \forall u \in G_{h,k},
v=(h,k,u),
\forall l \in L_v,
\end{equation}
% For the vehicle types that follow closed-tour, the entire tour length till arrival at the respective depot is considered:
\begin{equation} \label{eq:65b}
z_x \geq a_h^{v,l},
\quad
\forall h \in H, \forall k \in K_h \setminus O, \forall u \in G_{h,k},
v=(h,k,u),
\forall l \in L_v,
\end{equation}
\end{subequations}

% Similarly, equation \ref{eq:66} is used to populate the variable of $z_s$ such that it is equal to the sum of all vehicle-return-times at the Vehicle Depot. For the vehicles which don't need to return back at the starting depot, the time taken to return from the last vertex is subtracted as in Eq. \ref{eq:66a}. The $z_s$ variable represents the sum of all route durations and can be used as an alternate objective function.

\begin{subequations} \label{eq:66}
\begin{equation} \label{eq:66a}
z_{s_{OT}} =
\sum_{\substack{h \in H,\\
        k \in K_h \cap O,\\
        u \in G_{h,k},\\
        v=(h,k,u)}}~
\sum_{l \in L_v}~
\Bigl(
a_h^{v,l} - \sum_{i \in S_k,N_k,R_k} T^k_{i,h} \cdot x^{v,l}_{i,h}
\Bigr),
\end{equation}

\begin{multicols}{2}
\begin{equation} \label{eq:66b}
z_{s_{CT}} =
\sum_{\substack{h \in H,\\
        k \in K_h \setminus O,\\
        u \in G_{h,k},\\
        v=(h,k,u)}}~
\sum_{l \in L_v}~
a_h^{v,l},
\end{equation}

\begin{equation} \label{eq:66c}
z_s = z_{s_{CT}} + z_{s_{OT}}
\end{equation}
\end{multicols}
\end{subequations}

{\noindent \textbf{Variable Definitions}:}
% \subsection{Variable definitions}  \label{sec:3.7}

\begin{subequations} \label{eq:0}
\begin{equation} \label{eq:0a}
x^{v,l,m}_i \in \{0,1\}, \quad
\forall h \in H, \forall k \in K_h, \forall u \in G_{h,k}, v=(h,k,u),
\forall l,m \in L_v,
|l-m|=1,
\forall i \in V_k,
\end{equation}

\begin{multline} \label{eq:0b}
x_{i,j}^{v,1} \in \{0,1\}, \quad
\forall h \in H, \forall k \in K_h, \forall u \in G_{h,k}, v=(h,k,u),
\forall i,j \in V_k \cup h,
i \neq j,
if~i \in W \Rightarrow j \neq h,
if~i=h \Rightarrow j \notin R,
\end{multline}

\begin{equation} \label{eq:0c}
x_{i,j}^{v,l} \in \{0,1\}, \quad
\forall h \in H, \forall k \in K_h, \forall u \in G_{h,k}, v=(h,k,u),
\forall l \in L_v \smallsetminus \{{1}\},
\forall i \in V_k,
\forall j \in V_k \cup h,
i \neq j,
if~i \in W \Rightarrow j \neq h,
\end{equation}

\begin{equation} \label{eq:0d}
y^{v,l,m,c}_i \geq 0, \quad
\forall h \in H, \forall k \in K_h, \forall u \in G_{h,k}, v=(h,k,u),
\forall l,m \in L_v,
|l-m|=1,
\forall c \in C_k, \forall i \in V_k,
\end{equation}

\begin{multline} \label{eq:0e}
y_{i,j}^{v,1,c} \geq 0, \quad
\forall h \in H, \forall k \in K_h, \forall u \in G_{h,k}, v=(h,k,u),
\forall c \in C_k,
\forall i,j \in V_k \cup h,
i \neq j,
if~i \in W \Rightarrow j \neq h,
if~i=h \Rightarrow j \notin R,
\end{multline}

\begin{multline} \label{eq:0f}
y_{i,j}^{v,l,c} \geq 0, \quad
\forall h \in H, \forall k \in K_h, \forall u \in G_{h,k}, v=(h,k,u),
\forall l \in L_v \smallsetminus \{{1}\},
\forall c \in C_k,
\forall i \in V_k,
\forall j \in V_k \cup h,
i \neq j,
if~i \in W \Rightarrow j \neq h,
\end{multline}

\begin{equation} \label{eq:0g}
r_i^{v,l,c} \in \mathbb{R}, \quad
\forall h \in H, \forall k \in K_h, \forall u \in G_{h,k}, v=(h,k,u),
\forall l \in L_v,
\forall i \in S_k,
\forall c \in B_i \cap C_k,
\end{equation}

\begin{equation} \label{eq:0h}
n_i^{v,l,c} \in \{0,1\}, \quad
\forall h \in H, \forall k \in K_h, \forall u \in G_{h,k}, v=(h,k,u),
\forall l \in L_v,
\forall i \in S_k,
\forall c \in B_i \cap C_k,
\end{equation}

\begin{equation} \label{eq:0i}
e_i^{v,l,c} \geq 0, \quad
\forall h \in H, \forall k \in K_h, \forall u \in G_{h,k}, v=(h,k,u),
\forall l \in L_v,
\forall i \in S_k,
\forall c \in B_i \cap C_k,
\end{equation}

\begin{equation} \label{eq:0j}
a_i^{v,l} \geq 0,
\forall h \in H, \forall k \in K_h, \forall u \in G_{h,k}, v=(h,k,u),
\forall l \in L_v,
\forall i \in V_k \cup h,
\end{equation}

\begin{equation} \label{eq:0k}
d_h^{v,1} \geq 0,
\forall h \in H, \forall k \in K_h, \forall u \in G_{h,k}, v=(h,k,u),
\end{equation}

\begin{equation} \label{eq:0l}
d_i^{v,l} \geq 0,
\forall h \in H, \forall k \in K_h, \forall u \in G_{h,k}, v=(h,k,u),
\forall l \in L_v,
\forall i \in V_k,
\end{equation}

\begin{equation} \label{eq:0m}
t^{v,l,m}_i \geq 0, \quad
\forall h \in H, \forall k \in K_h, \forall u \in G_{h,k}, v=(h,k,u),
\forall l,m \in L_v,
|l-m|=1,
\forall i \in V_k,
\end{equation}

\begin{equation} \label{eq:0n}
t_{i,j}^{v,1} \geq 0, \quad
\forall h \in H, \forall k \in K_h, \forall u \in G_{h,k}, v=(h,k,u),
\forall i,j \in V_k \cup h, 
i \neq j,
if~i \in W \Rightarrow j \neq h,
if~i=h \Rightarrow j \notin R,
\end{equation}

\begin{equation} \label{eq:0o}
t_{i,j}^{v,l} \geq 0, \quad
\forall h \in H, \forall k \in K_h, \forall u \in G_{h,k}, v=(h,k,u),
\forall l \in L_v \smallsetminus \{{1}\},
\forall i \in V_k,
\forall j \in V_k \cup h,
i \neq j,
if~i \in W \Rightarrow j \neq h,
\end{equation}

\begin{equation} \label{eq:0p}
w_i^{v,l} \geq 0,
\forall h \in H, \forall k \in K_h, \forall u \in G_{h,k}, v=(h,k,u),
\forall l \in L_v,
\forall i \in S_k,
\end{equation}

\begin{multline} \label{eq:0q}
o_{i,(\hat{v},\hat{l},\hat{s})}^{c,(v,l)} \in \mathbb{R}, \quad
\forall i \in S,
\forall c \in B_i,
\forall h \in H, \forall k \in K_h \cap A_c \cap K_i, \forall u \in G_{h,k},
v=(h,k,u),
\forall \hat{h} \in H,
\forall \hat{k} \in K_{\hat{h}} \cap A_c \cap K_i,\\
\forall \hat{u} \in G_{\hat{h},\hat{k}},
\hat{v}=(\hat{h},\hat{k},\hat{u}),
\forall l \in L_v,
\forall \hat{l} \in L_{\hat{v}},
\forall \hat{s} \in \{Arr,Dep\},
\end{multline}

\begin{multline} \label{eq:0r}
g_{i,(\hat{v},\hat{l},\hat{s})}^{(v,l,s)} \in \{0,1\}, \quad
\forall i \in S,
\forall h \in H, \forall k \in K_h \cap A_{B_i} \cap K_i, \forall u \in G_{h,k},
v=(h,k,u),
\forall \hat{h} \in H,
\forall \hat{k} \in K_{\hat{h}} \cap A_{B_i} \cap K_i,\\
\forall \hat{u} \in G_{\hat{h},\hat{k}},
\hat{v}=(\hat{h},\hat{k},\hat{u}),
\forall l \in L_v,
\forall \hat{l} \in L_{\hat{v}},
\forall s,\hat{s} \in \{Arr,Dep\},
\end{multline}

\begin{multicols}{2}
\begin{equation} \label{eq:0s}
b_i^c \geq 0,
\forall i \in N,
\forall c \in D \cup P
\end{equation}
\begin{equation} \label{eq:0t}
b_i^c \geq 0,
\forall i \in W,
\forall c \in D
\end{equation}
\begin{equation} \label{eq:0u}
q_i^c \geq 0,
\forall i \in N,
\forall c \in D \cup P
\end{equation}
\begin{equation} \label{eq:0v}
q_i^c \geq 0,
\forall i \in W,
\forall c \in D
\end{equation}
\begin{equation} \label{eq:0w}
z_x \geq 0.
\end{equation}
\begin{equation} \label{eq:0x}
z_s \geq 0.
\end{equation}
\begin{equation} \label{eq:0y}
z_{s_{OT}} \geq 0.
\end{equation}
\begin{equation} \label{eq:0z}
z_{s_{CT}} \geq 0.
\end{equation}
\end{multicols}

\end{subequations}

{\noindent \textbf{Detailed Explanation of the Constraints used in this MILP Formulation}:}

Explaining the above constraints in details, Eq. \ref{eq:2} allows at most a single path to start from the respective depot for each vehicle, to any of the allowed elements (which are in the same SMTS) in the sets $W,N$ and $S$ for the base level. Eq. \ref{eq:3} constrains the number of incoming paths/edges into the vehicle's respective depot, across all it's levels, to be limited to $1$, \ie only when that vehicle is used. Eq. \ref{eq:4} conserves the route flows for each vertex for all layers constraining the sum of all incoming edges to be equal to the sum of all outgoing edges for that vertex; eq. \ref{eq:4a} is specific for the first layer, while eq. \ref{eq:4c} is specific for the last layer. Only the first layer has all the vertices in sets $W,S$ and $N$ being connected with edges (decision variables) from the vehicle's depot. All layers have connections from all vertices in the sets $S,N$ and $R$ to the vehicle's depot. For Simultaneous Delivery and PickUp Nodes, only one visit in total is allowed by any of vehicles; which is constrained in Eq. \ref{eq:7}.

We use single flow variables for each of the load-types and as per compatibility with the vehicle type; any outflow from the depots is made zero by Eq. \ref{eq:8}. All inflows into the depots are made zero by Eq. \ref{eq:9}; these nullifications will help us for the remainder of the flow formulations as we may not consider any of the flow variables into or out of any $VD$. For Warehouses, there will always be an increase from incoming links to outgoing links at each level of vehicle for each compatible delivery load type (Eq. \ref{eq:10}). (In case of non-existing variables during the equation framing, \eg when $l$ takes values of $1$ or $l_v^{max}$, the variables must be ignored. Since the constraints differ slightly when ignoring these non-existent variables, constraints dealing with the inter-layer variables are bifurcated into 3 constraints, 1 for the lowest layer, 1 for the in-between layers, and 1 for the top-most layer, similar to the bifurcations in Eq. \ref{eq:4}). There would be no change in the amounts of pickup loads at Warehouses (Eq. \ref{eq:11}). Each Warehouse capacity is constrained through Eq. \ref{eq:12}. For Relief Centres, there will always be a decrease in the outgoing pickup resource flow values \wrt the incoming pickups as the evacuees (picked-up from Nodes) deboard at each RC (Eq. \ref{eq:13}); no change in the delivery values would occur during any vehicle's visit to the RCs (Eq. \ref{eq:14}). Capacity of each RC for each type of the pickup load type is constrained through Eq. \ref{eq:15}. During any vehicle's visit at Nodes, the outgoing pickup value for each Node is constrained to be greater than the incoming pickup values through Eq. \ref{eq:16}. For the delivery variables, the sum of all incoming values at each level for each Node should be lesser than the sum of the outgoing values for each Cargo Type (CT) individually, during any vehicle's visit to Nodes (Eq. \ref{eq:17}). Minimal resource commitment at each Node is ensured by Eqns. \ref{eq:18} and \ref{eq:19}. At Transhipment Port, the residues of each type of $TP$-compatible load during each visit of any vehicle is calculated using Eq. \ref{eq:20}. For the loads which are not $TP$-compatible, no changes in the flow quantity is allowed at the respective $TP$ (Eq. \ref{eq:21}). 

For developing the Time Constraints, we visualize the following process-flow during the operation at any typical TP during any vehicle visit \ie at every layer of every vehicular trip.
\begin{enumerate}
  \item Vehicle reaches $TP$ at an assumed time $T_0$
  \item Loading/Unloading of Vehicle-\&-$TP$ compatible cargo is done. For each of the compatible load categories, the time considered is the sum of the waiting time $T_W$ loading/unloading time $T_{U/L}$.
  Here we assume this process to be done without any parallelization approach. As a future research area; or problem variant, this could be edited for task parallelization which could allow the transhipment operations to become organized. As a suggestion, simultaneous cargo-specific queues could also be pondered upon, depending on manpower available at the site and the number of extra vehicle-crew.
  \item Vehicle departs the TP during this same visit at a time $T_z$ where,
  $T_Z = T_0 + \sum_{\substack{Compatible \\ Cargos}} \Bigl( T_W+T_{U/L} \Bigr) $
\end{enumerate}
For a VD all departure times are considered $0$ (Eq. \ref{eq:22}), assuming that it takes no time to prepare the vehicles (fuelling, maintenance, etc.) or the crew (assigning the route instructions, providing a dry run of the tasks involved during their specific stops, etc.). This constraint may be modified to allow the condition when certain vehicles are available at a later time, allowing other available vehicles to start the relief operation immediately. In this work, we consider this modification only during for the literature instances (modelling vehicular fixed cost as the start time from depot).

Travel time calculation (on the same layer) is done by Eqns. \ref{eq:23}, \ref{eq:24} and \ref{eq:25}; the travel time calculation between layers is done by Eqns. \ref{eq:26}, \ref{eq:27} and \ref{eq:28}. Arrival times are calculated using Eqns. \ref{eq:29} and \ref{eq:29.5}.
Departure times are calculated in Eqns. \ref{eq:36} (for Transhipment Ports), \ref{eq:37} (for Nodes), \ref{eq:38} (for Warehouses), and \ref{eq:39} (for Relief Centres). The residues need to be arranged \wrt time so that we may ensure that at any point in time, the sum of residues of a specific load-type at a specific $TP$ is never negative; Eqns. \ref{eq:30}, \ref{eq:31}, \ref{eq:32}, \ref{eq:33}, \ref{eq:34} and \ref{eq:35} are used to calculate the modulus of the residues. By comparing Arrival and Departure times of vehicles at Transhipment Ports, Eq. \ref{eq:40} populates the intermediate binary variable $g$ when $s=Dep$ and $\hat{s}=Arr$, Eq. \ref{eq:41} populates variable $g$ when $s=Dep$ and $\hat{s}=Dep$, Eq. \ref{eq:42} populates variable $g$ when $s=Arr$ and $\hat{s}=Arr$ and Eq. \ref{eq:43} populates $g$ for when $s=Arr$ and $\hat{s}=Dep$.

The cell values of the space-time matrices consist of the variables $o$ which store some proportion of the (load-specific) residue $r$ available for transhipment during the arrival and departure times of all vehicles. This uses the comparison of a vehicle's visiting time with another at the same Transhipment Port for all its visiting possibilities individually, \ie across the various levels. We first consider the case when the visiting durations of the vehicles under comparison don't overlap, \ie they concerned vehicles don't meet during their corresponding layer-specific visit at the concerned Transhipment Port $S$. Due to this non-overlapping nature, it is easy to determine the available proportion of residue (\ie the cell value $o$) which would either be the entire residue of the vehicle under consideration $v$ \wrt a reference vehicle $\hat{v}$ if $v$ has already departed the transhipment port with some residue value $r$ during its specific layer-wise visit $l$; otherwise its value would be $0$ if the concerned vehicle $v$ is yet to arrive in the transhipment port \wrt the reference vehicle $\hat{v}$'s visit duration. Constraining the cell values of space-time matrices when the vehicles being compared don't overlap in spacetime is ensured through the Eqns. \ref{eq:44}, \ref{eq:45}, \ref{eq:46} and \ref{eq:47}, which populate the variable $o$ when $g=1$.

However, when $g=0$ for both $s=Arr$ and $s=Dep$, \ie when the vehicles under comparison overlap during their visiting times at the same transhipment port, those constraints become redundant and the new constraints are developed to populate the $o$ variables of the space-time matrix within certain bounds such that the rate of loading/un-loading of the resources is not breached (due to minimization, the $o$ variable should take the upper bound when the residues are positive). When the residues are positive (and the visits of the vehicles being compared overlap in spacetime), we constrain the cell values of the space-time matrices within their limits of $0$ to the maximum value being the residue $r$ itself, through Eqns. \ref{eq:48}, \ref{eq:49}, \ref{eq:50}, \ref{eq:51}, \ref{eq:52} \ref{eq:53} and \ref{eq:54} (the Eqns. \ref{eq:48}, \ref{eq:49} and \ref{eq:50} are applicable irrespective of whether the vehicles overlap in spacetime or dont). When the residues are negative (and vehicles being compared overlap in spacetime), we constrain the cell values of the space-time matrices within their limits of $0$ to the minimum value being the residue $r$ itself, through Eqns. \ref{eq:55}, \ref{eq:56}, \ref{eq:57}, \ref{eq:58}, \ref{eq:59} \ref{eq:60} and \ref{eq:61} (the Eqns. \ref{eq:55}, \ref{eq:56} and \ref{eq:57} are applicable irrespective of whether the vehicles overlap in spacetime or don't).

The typical space-time matrix may be imagined as a matrix with rows as various space-stamps (\ie each compatible vehicle's levels) and columns as various time-stamps (refer to Figure S1 and Figure S4 in the Supplementary document). For each vehicle level in the rows, there are two time-stamps namely arrival and departure. Therefore all space-time matrices constructed by us are rectangular in nature with columns being twice that of rows, \ie the time-stamps correspond to the arrival and departure times of the compatible vehicles' levels. To ensure any resource collection (\ie loading into vehicles) happens only when there are available resources at transhipment ports, we constrain the sum of the cells across each time axis (each column) for every space-time matrix to be non-negative. This ensures that at every time-stamp the total resource-specific residue at any $TP$ never becomes negative, \ie no collection of loads takes place before the load was deposited, and only the loading of materials previously deposited by some other vehicles is allowed into the vehicles that arrive later (thereby eliminating pseudo-loading, and ensuring any transfer of pickup or delivery types from the specific $TP$ happens only happen when enough material is available at that $TP$). Having developed the space-time matrix (\ie calculated each cell of the matrix, depicted by the variable $o$), for each load-type at each Transhipment Port $S$ (considering compatibilities for allowable transhipments) populated by respective $o$ variables, we constrain the cascaded-temporal sum of residues to be non-negative at each time-stamp, for every load-type at every $S$ (Eq. \ref{eq:62}).

Since Transhipment Ports are modelled to allow the transfer of (compatible) cargos during the operation period, they should not contain any resource (pickup or delivery load type) finally at the end of the emergency relief-and-rescue operation. Therefore, we need to provide an additional constraint which is set so as to nullify the sum of all residues (of each resource type separately) at each Transhipment Port. Due to the nature of the problem, this is necessary to be constrained for only the PickUp Load Types, since Eq. \ref{eq:62} is enough for the Delivery Load Types. If this constraint is not provided, PickUp quantities after transhipment initiation will remain at the respective Transhipment Ports instead of finally being deposited at Relief Centres. thus, all residues at Transhipment Ports are summed up and equalled to zero (Eq. \ref{eq:62.5}) to constrain the final residues at each Transhipment Ports to be null for each resource type. (In some instances, instead of using $\forall c \in B_i \cap P$ in Eq. \ref{eq:62.5}, we do use $\forall c \in B_i$; this is mainly done to check the formulation run times for different combinations of constraint usages.)

To ensure vehicle capacities are never breached, all flow variables associated with an edge are constrained \wrt the total available vehicle-volumes (Eq. \ref{eq:63}), and \wrt the total allowable vehicle-weights (Eq. \ref{eq:64}). Eq. \ref{eq:65} assigns the MakeSpan objective variable of $z_x$, such that it can take the maximum among all vehicle's route durations; for the $VTs$ that follow open-tours the last route back towards the depot is subtracted from the total time (Eq. \ref{eq:65a}), and for the  $VTs$ that follow closed-tours the entire tour length till arrival at the respective depot is considered (Eq. \ref{eq:65b}). Similarly, equation \ref{eq:66} is used to populate the variable of $z_s$ such that it is equal to the sum of all vehicle-return-times at $VD$. For the vehicles which don't need to return back at the starting depot, the time taken to return from the last vertex is subtracted as in Eq. \ref{eq:66a}. The $z_s$ variable represents the sum of all route durations and may be used as an alternate objective function. Finally, all variables are defined in the Eqns. \ref{eq:0}, with some additional constraints available in \ref{Extra Optional Constraints to the MILP} (these additional constraints are considered during the computational analyses).

\subsection{The Study-Example Result using our MILP Formulation} \label{Case Study MILP Result}

The optimization process, as followed in this paper, minimizes the highest route duration among all vehicles at the first step (cascade) and places this value as a constraint on the route duration of each vehicle. Then the vehicle with the highest route duration is removed from being further compared and the same optimization is done again. This time (in a subsequent cascade) another vehicle is found to be having the highest route duration (this route duration value will always be less than the values found in the previous cascades, given that the previous cascade-steps were run till optimality); this route duration is set as the upper bound of the vehicle route durations of all vehicles which were compared in this cascade-step (\ie excluding vehicles having the highest route duration in any previous cascade-step), and then the optimization step is re-run (after discarding the vehicle with the highest route duration in the latest cascade, from being compared any further).

\begin{figure}
    \centering
    \includegraphics[width=1.05\linewidth]{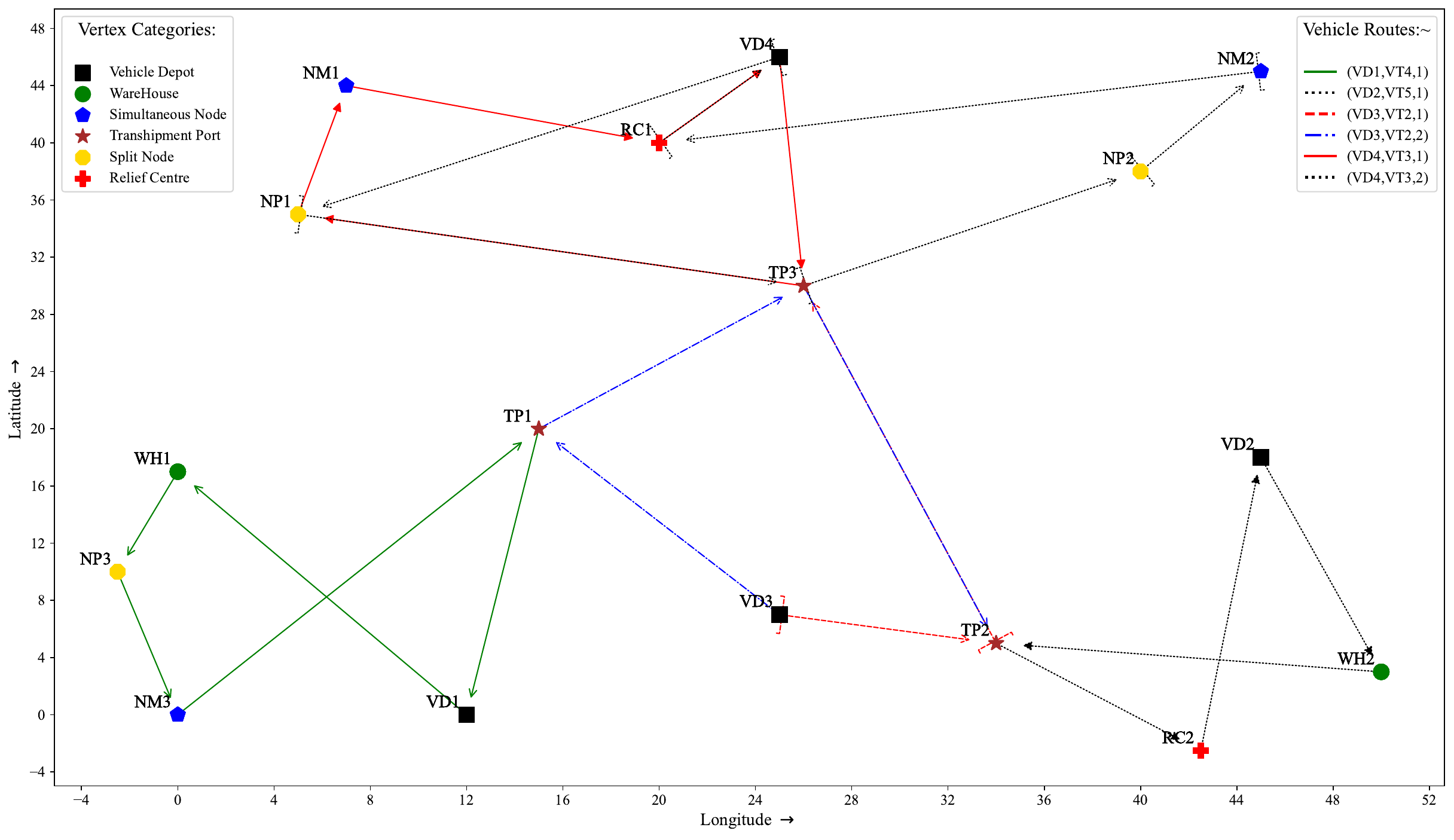} % 0.975 fits
    \caption{\textbf{MILP result for the illustrative example (for the Integer flow problem)}}
    \label{fig: casestudy solution milp}
\end{figure}

The result of the cascaded optimization for the example instance in \autoref{fig:Picture_Instance} (and considering all the resource flow variables $y$ to be integers) is shown in \autoref{fig: casestudy solution milp}. The optimal (average optimality gap across all cascades is exactly 0) routes in the order of their cascades are (with the route durations indicated at the end in brackets):

$\Downarrow$ \textbf{Cascade Numbers} \hfill \textbf{Route Durations} $\Downarrow$
\begin{enumerate}[label=\textbf{\scriptsize\arabic*}:] % Custom label format
    \item Vehicle (VD1, VT4, 1): VD1 $\rightarrow$ WH1 $\rightarrow$ NP3 $\rightarrow$ NM3 $\rightarrow$ TP1 $\rightarrow$ VD1 \hfill (115.083)
    \item Vehicle (VD2, VT5, 1): VD2 $\rightarrow$ WH2 $\rightarrow$ TP2 $\rightarrow$ RC2 $\rightarrow$ VD2 \hfill (99.522)
    \item Vehicle (VD3, VT2, 1): VD3 $\rightarrow$ TP2 $\rightarrow$ TP3 \hfill (99.074)
    \item Vehicle (VD3, VT2, 2): VD3 $\rightarrow$ TP1 $\rightarrow$ TP3 $\rightarrow$ TP2 \hfill (79.459)
    \item Vehicle (VD4, VT3, 1): VD4 $\rightarrow$ TP3 $\rightarrow$ NP1 $\rightarrow$ NM1 $\rightarrow$ RC1 $\rightarrow$ VD4 \hfill (75.502)
    \item Vehicle (VD4, VT3, 2): VD4 $\rightarrow$ NP1 $\rightarrow$ TP3 $\rightarrow$ NP2 $\rightarrow$ NM2 $\rightarrow$ RC1 $\rightarrow$ VD4 \hfill (57.39)
\end{enumerate}
Throughout this paper, we round all route durations (or their sums, while reporting the sum of all routes) to 3 decimal places, and all times to 2 decimals. The IntegralityFocus parameter in Gurobi for all instance runs (throughout this study) is always maintained at 1 (to allow rare usage of Gurobi's integer tolerances). For this run, when the Big M value was set at 3E3, and Gurobi's FeasibilityTol (feasibility tolerance for constraints) was set at 1E-5 (reduced from the default value of 1E-6), the total optimization time (for all the cascades) was 21.88 minutes; when the FeasibilityTol was changed to 1E-7 with the Big M value being 8E3, the full cascaded optimization happened in 31.05 minutes. 

The sum of all routes after the first cascade was 542.04 while after the last cascade was 526.03 (for the detailed results, \ie, the resource transfer amounts, waiting times \etc, please refer to the online repository mentioned after \autoref{Dataset Hosting}). This showcases the fact that conventional min-max optimization (\ie the single step Makespan minimization), does not fully solve the problem since this single step is unable to reach into smaller vehicular routes for internal optimizations. To find better solutions, multi-step optimization is necessary; we perform this multi-step minimization terming them cascades and as we progress through each cascade, the internal (smaller) vehicular routes further improve thus strengthening the solution progressively. This improvement through Cascaded Makespan Minimization is further highlighted (as the encoding $\textit{I}$) through the computational results in \autoref{Tab: Small Integer Results}. The Cascaded Makespan Minimization is a lexicographic/hierarchical optimization approach used for the multi-step optimization problem.

When we allow the resource flow variables to take continuous values, we get a different (fully optimal across all cascades) result where the value of the first cascade changes to 114.363 (this is lower than the integer problem's largest vehicle route duration of 115.083) and the sum of all route durations at the end of all cascades becomes 524.97; however, continuous problems are not realistic and don't have real-world applications as non-integer resources do not have any meaning in reality.

\section{Heuristic Development based on Preferential Selection of Routes (PSR) from Node-variant-specific Decision Trees (DTs): Route Generation, Integration and Perturbation (GIP)} \label{Heuristic Development}

\begin{figure}
    \centering
    \includegraphics[width=1\linewidth]{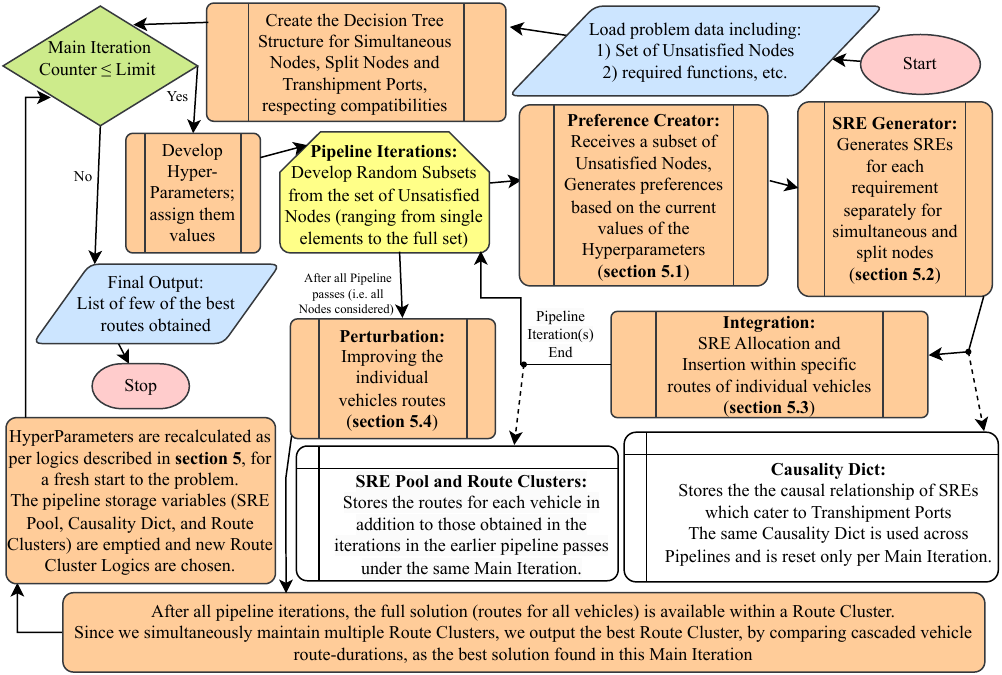} %1.04064
    \caption{\textbf{Overview of the proposed heuristic}}
    \label{fig:MainHeuristicFlowchartOverview}
\end{figure}

The main components of the DT-based PSR-GIP Heuristic are briefed below, along with the integration methodology of routes, and utilization of multiple passes of the heuristic's Pipeline  (full details is elaborated separately as Heuristic Methods in the supplementary section S3). The flowchart in \autoref{fig:MainHeuristicFlowchartOverview} describes the overview of the Heuristic pipeline.

\subsection{Preference Creation and Decision Trees} \label{Preference Creation and Decision Trees}
\begin{figure}
    \centering
    \includegraphics[width=1.05\linewidth]{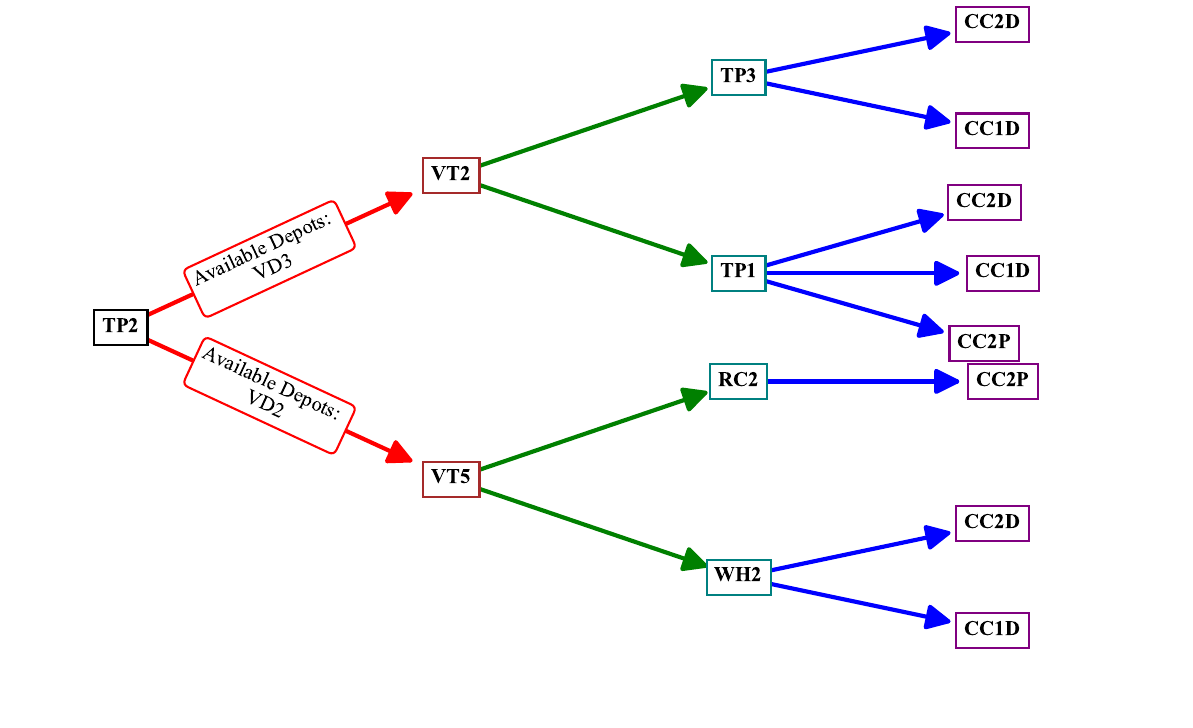} % Use the value of 0.9 for fitting in a page with other things
    \caption{\footnotesize\textbf{Decision Tree Structure of TP2 from the illustrative example. The Trunk in this case is the TP2; the two Branches are VT2 and VT5 (with their respective set of possible VDs); the 4 Twigs (2 on each Branch) are TP3, TP1, RC2 and WH2; finally, the 8 Leafs of the DTS are all the CTs indicated on the right-most side.}}
    \label{fig:Imagine_TP2}
\end{figure}

The preferences of vertices to satisfy a node (or a transhipment port) are stored within decision trees. Each Multi-step Decision Tree consists of Node $\rightarrow$ Vehicle Types $\rightarrow$ Vertices $\rightarrow$ Cargo Types, as Trunk $\rightarrow$ Branches $\rightarrow$ Twigs $\rightarrow$ Leaves, and this is deemed as the Decision Structure (DTS) of the problem. Apart from the exact sequence of routes or multi-trips of vehicles, this Decision Structure is able to capture most of the essential information in the considered problem and especially all the compatibility considerations. Elaborating on the Decision Structure (a sample is provided in \autoref{fig:Imagine_TP2}), from a single Trunk, multiple Branches emerge; from each Branch, multiple Twigs emerge (Twigs from one Vehicle Type-Branch can represent similar Vertices as represented by Twigs from another Vehicle Type-Branch); and ultimately from each Twig, multiple Leaves emerge (Leaves from one Vertex-Twig can represent similar CTs as represented by Leaves from another Vertex-Twig).

\small

Numerical weights are placed at the ends of these trees, which are used to calculate preference scores for route creation in order to satisfy requirements (detailed methodology for calculation of these preferences for nodes is provided in the supplementary sections S3.1.3 and S3.1.4; for TPs the preference scoring methods are provided in section S3.1.5 which leverages a complex recursive function developed for estimating and fathoming deep transhipments). The \autoref{fig:Degree Depiction in Supplementary} in \ref{Visualizing the Transhipment Degree} further highlights the stages/degrees of transhipment, with deep transhipments indicating higher degrees.

\subsection{Smallest Route Element (SRE) development} \label{Smallest Route Element (SRE) development}
\begin{figure}
    \centering
    \includegraphics[width=1\linewidth]{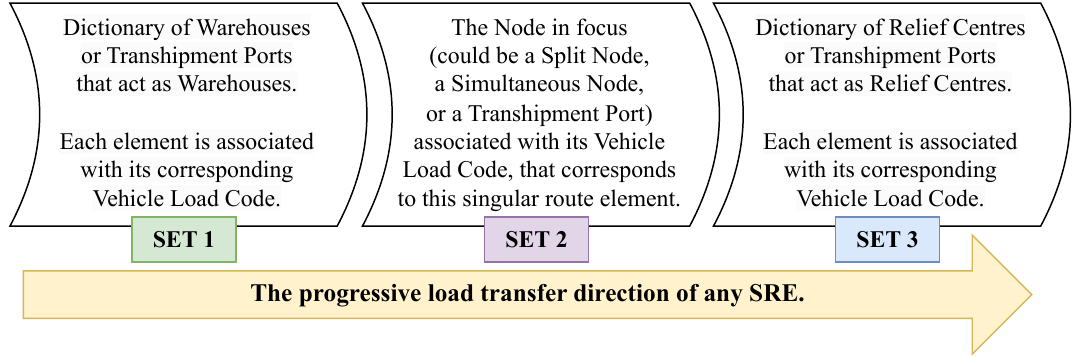}
    \caption{\textbf{Smallest Route Element. Each SRE is specific to a Vehicle Type and is associated with a list of possible Vehicle Depots; this is used during the assigning of an SRE to any vehicle of that specific Vehicle Type originating from any of the associated VDs.}}
    \label{fig:Smallest Route Element depiction}
\end{figure}

A Smallest Route Element is an (almost) independent portion within the route of a vehicle. It may be only connected with another SRE due to the causal relationship during transhipments. During any transhipment, the resources which are to be delivered from any TP can happen only after those specific amount (or higher) of each resource reaches that specific Port in $S$; this creates a necessary causality constraint for each transhipment. This causal relation comes into effect whenever any resource needs to be transhipped and these relationships are stored within a Causality Dict (CD). There are three segments to an SRE as shown in \autoref{fig:Smallest Route Element depiction}, termed as SETs 1, 2 and 3.

During each Pipeline pass within a Main Iteration, requirements of
A set of unsatisfied nodes (Split and/or Simultaneous) is created from among all the unsatisfied nodes (without replacement) during each Pipeline pass, and each of these nodes are satisfied by creating SREs to fulfill their requirements (the full process of developing SREs for each of the node types is explained in the supplementary section S3.2, along with the details for generating SREs during transhipment processes while also logging each transfer into the Causality Dict in section S3.2.1 to maintain temporal causality during each unique transhipment).

\subsection{Integration of SREs within allocated Vehicle's Route6} \label{Integration of SREs within allocated Vehicle's Route6}

During every Pipeline pass (\ie Pipeline Iteration), we obtain the newly created SREs developed in the current pipeline, these are also stored within a Global SRE Pool (emptied after every Main Iteration), as well as a Pipeline SRE Pool (emptied after each Pipeline Iteration). For the Allocation and Insertion of SREs (from the Pipeline SRE Pool only), we develop each SRE to take a form similar to the structure of a Route, naming this metamorphosed SRE as a Route Portion (RoPr).

\subsubsection{Description of the Route6 structure of all vehicular routes} \label{Description of the Route6 structure of all vehicular routes}

Each vehicle's route is stored within six different arrays:
\begin{enumerate}[label=\textbf{\scriptsize List \arabic*}:]
    \item Stores the vertices orderly as visited by a vehicle. TPs, if any visited by the vehicle, retain their altered new name in this list.
    \item Maintains the count of the number of Perturbations of the vertex at this same position in \textbf{List 1}.
    \item Stores the Vehicle Load Code (VLC) for the vertex in the corresponding position in \textbf{List 1}.
    \item Stores StatusCodes of the vehicle after each visit to vertices. StatusCodes indicate the current amount of each Cargo Type present within a vehicle.
    \item Stores the TimeTuple for this specific visit. A TimeTuple has the form:
    \begin{itemize}
        \item (Vertex Entering Time, Loading Unloading Time, Vertex Leaving Time) if it is not a for Transhipment Port
        \item If the concerned vertex being visited (at the same position in \textbf{List 1}) is a Transhipment Port, an additional information of Waiting Time is necessary and the TimeTuples for Transhipment Port visits are as: (Vertex Entering Time, Loading Unloading Time, Waiting Time, Vertex Leaving Time)
    \end{itemize}
    \item To appropriately implement the entire Heuristic, we need to distinguish between two identical SREs (which may be formed due to the requirement of multiple trips, possibly due to capacity constraints of vehicles). For this we maintain a reference ID of the SRE (during our code implementation, the reference ID was set as the position of the SRE within the Global SRE Pool).
\end{enumerate}
A VLC represents the actual transfer of CTs \wrt the vehicle (it is can be imagined by considering the vehicle as an enclosed Gaussian surface, and any entry of a resource into the surface is attached with a positive number, and vice versa).

\begin{table}[H]
\centering
\caption{A. \textbf{A sample Compacted\_Route6 example of a complete VT3 trip from \autoref{fig:Picture_Instance}}\label{longTab: Sample Route6 Examples}}
% \noindent
\resizebox{1\textwidth}{!}{ % Resize to fit the width of the page
\begin{tabular}{|*{11}{m{0.095\textwidth}|}} % 3 columns with a repeating pattern
% \hline
\hline	\textbf{List 1}: &	VD4 &	NP1	&	RC1	&	TP3\_G	& NP2 &	NM2	& RC1 & NP1 & NM1 & VD4\\
\hline	\textbf{List 2}: &	N/A &	0	&	0	&	2	& 0 &	0	& 0 & 0 & 0 & N/A\\
\hline	\textbf{List 3}: &	\{\} &	\{1P:25\}	&	\{1P:-25\}	&	\{1D:21, 2D:23\}	& \{1D:-4, 2D:-7\} &	\{1D:-3, 2D:-2, 1P:20\}	& \{1P:-20\} & \{1D:-5, 2D:-6\} & \{1D:-9, 2D:-8\} & \{\}\\
\hline	\textbf{List 4}: &	\{\} &	\{1P:25\}	&	\{\}	&	\{1D:21, 2D:23\}	& \{1D:17, 2D:16\} &	\{1D:14, 2D:14, 1P:20\}	& \{1D:14, 2D:14\} & \{1D:9, 2D:8\} & \{1D:-9, 2D:-8\} & \{\}\\
\hline	\textbf{List 5}: &	( , , 0) &	(13.98, 3.75, 17.73) &	(27.23, 3.75, 30.98) &	(38.14, 16.35, $\omega$, 54.49+$\omega$) &	(64.38+$\omega$, 3.9, 68.28+$\omega$) &	( 73.59+$\omega$, 4.95, 78.54+$\omega$) &	(93.61+$\omega$, 3, 96.61+$\omega$) &	(106.11+$\omega$, 4.05, 110.16+$\omega$) &	(115.62+$\omega$, 6.45, 122.07+$\omega$) &	(132.66+$\omega$, , ) \\
\hline	\textbf{List 6}: & & \textit{SRE 30} & \textit{SRE 30} & \textit{SRE 23, SRE 54, SRE 89, SRE 17} & \textit{SRE 23} & \textit{SRE 54} & \textit{SRE 54} & \textit{SRE 89} & \textit{SRE 17} & \\
\hline
\end{tabular}
} % For Resize Box Closure
\end{table}

\addtocounter{table}{-1}

\tiny
\begin{table}[H]
\centering
\caption{B. \textbf{Expansion of the compacted element from \autoref{longTab: Sample Route6 Examples}A}\label{longTab: Compact Route6 Element Expanded}}
% \noindent
\resizebox{0.9\textwidth}{!}{ % Resize to fit the width of the page
% \begin{tabular}{|*{6}{m{0.2\textwidth}|}} % 3 columns with a repeating pattern

\begin{tabular}{|m{0.175\textwidth}|m{0.075\textwidth}|*{4}{m{0.175\textwidth}|}} % 3 columns with a repeating pattern
% \begin{tabular}{|*{6}{c|}}

% \hline

\cline{1-1} \cline{3-6}
\textbf{Compact Element} & & \multicolumn{4}{|c|}{\textbf{Original Elements in expanded form}} \\
\cline{1-1} \cline{3-6}

TP3\_G & & TP3\_G & TP3\_G & TP3\_G & TP3\_G \\
\cline{1-1} \cline{3-6}
0 & & 0 & 1 & 1 & 0 \\
\cline{1-1} \cline{3-6}
\{1D:21, 2D:23\} & & \{1D:9, 2D:8\} & \{1D:3, 2D:2\} & \{1D:5, 2D:6\} & \{1D:4, 2D:7\} \\
\cline{1-1} \cline{3-6}
\{1D:21, 2D:23\} & & \{1D:9, 2D:8\} & \{1D:12, 2D:10\} & \{1D:17, 2D:16\} & \{1D:21, 2D:23\} \\
\cline{1-1} \cline{3-6}
(38.14, 16.35, $\omega$, 54.49+$\omega$) & & (38.14, 6.45, $\omega_1$, 44.59+$\omega_1$) & (44.59+$\omega_1$, 1.95, $\omega_2$, 46.54+$\omega_1+\omega_2$) & (46.54+$\omega_1+\omega_2$, 4.05, $\omega_3$, 50.59+$\omega_1+\omega_2+\omega_3$) & (50.59+$\omega_1+\omega_2+\omega_3$, 3.9, $\omega_4$, 54.49+$\omega_1+\omega_2+\omega_3+\omega_4$) \\
\cline{1-1} \cline{3-6}
\textit{SRE 23, SRE 54, SRE 89, SRE 17} & & \textit{SRE 17} & \textit{SRE 54} & \textit{SRE 89} & \textit{SRE 23} \\
\cline{1-1} \cline{3-6}

\end{tabular}
} % For Resize Box Closure
\end{table}
\small

If multiple same vertices from different SREs occur sequentially in a Route6, we do not combine them together. The combined Route6 as shown in \autoref{longTab: Sample Route6 Examples}A may only be developed after the Perturbation step (section \ref{Perturbations to find the best Route Cluster}), when the final VDs are assigned at the end of each vehicle's route (depending on if closed trip is considered); till before that step, we maintain the vertices separately as shown for one such vertex in \autoref{longTab: Compact Route6 Element Expanded}B.\newline

Each SRE is moulded into a structural form same as a Route6, termed a RoPr (see supplementary section S3.3.2); depending on the number of permutations possible within each SET of any SRE, multiple RoPr can be formed with are together termed as RoPr-Combos. A complete solution consisting of all vehicle routes is termed as a Route Cluster (RoCu). We develop multiple allocation logics to integrate RoPrs within RoCus; initially each RoCu is labelled with a specific allocation-integration logic and as we introduce dynamic shuffling among these logics, the solution is found to improve. Detailed descriptions of these processes and algorithms can be found in the supplementary section S3.3. While a RoPr is being integrated into a Route6 within a RoCu, similarity of its individual elements are checked with the existing elements of the Route6 allowing to match similar vertices (supplementary section S3.3.6). A novel Waiting Time calculation method based on geometrical logics is devised (supplementary section S3.3.7) which is necessary to maintain temporal causality during every individual transhipment.

\subsection{Finding neighbouring route structures through Perturbations of every Route6 within a Route Cluster} \label{Perturbations to find the best Route Cluster}

This function is used to find the neighbouring solutions of a Route6. It takes a RoCu (termed a\_RoCu in the pseudocode) and appends the final VDs at the ends of those VTs that have closed trip requirements (thus adding an additional TimeTuple to the routes with a closed trip). Next, iteration over each of its Route6 takes place, to retain successful perturbations as detailed in the below pseudocode. For every Main Iteration, we use a user-defined limit on the number of perturbations allowed to a RoCu (termed P in the below function definition); this number is spread across the different Routes6 in that RoCu proportional to their vertex lengths.

\noindent\hrulefill

\noindent\textbf{Pseudocode of Function to find neighbouring route structures of a Route6}

\noindent\hrulefill
\begin{algorithmic}[1]

\Function{neighbouring\_route\_structure\_generation}{a\_RoCu, P}
        \State append the return VDs (as applicable) across all routes in a\_RoCu
	\For{each Route6 within a\_RoCu}
		\For{no. of iterations equal to $ P \times \frac{length~of~this~Route6}{sum~of~lengths~of~all~Route6~within~RoCu}$:}
			\State randomly selects a vertex to be perturbed (within this Route6) \Comment{The randomly selected vertex element (\ie the whole column of List 1 to 6 together) will be shifted within that route}
			\State upper and lower positional bounds are obtained for a feasible shift
			\State shift the full vertex column element to any position within the identified bound
			\State after the shift, recomputation of the StatusCodes for this Route6 is done
			\State all TimeTuples are updated \Comment{During the TimeTuple updation, it may be necessary to trigger the Waiting Time calculation function multiple times; also other route durations will be affected due to this shift, as they may have connected waiting times to this route}
			\If{all newly computed StatusCodes are within the vehicle's weight and volume limitations (\ie no capacity breach detected), \textbf{and} all new Waiting Time computations are successful, \textbf{and} all TimeTuples across that RoCu have been updated}
				\If{there is any improvement seen in the overall cascaded vehicle route durations of this RoCu}
					\State we retain this shifted RoCu (as an update to the original RoCu) \Comment{Subsequent perturbations are performed on the updated RoCu}
					\State increment the List 2, of the specific vertex which is moved, by 1
				\Else
					\State we revert back to the previous RoCu before this single-vertex shift
				\EndIf
			\EndIf
		\EndFor
	\EndFor
\EndFunction
\noindent\hrulefill
\end{algorithmic}

We find the feasible left and right Position bounds (for step 6 in the above pseudocode) till where it may be shifted without hampering the relative position of SET 1, 2 and 3 of its respective SRE (this is done by searching across Route6 and marking the Positions of elements from the same SRE). For every SRE, we need to always ensure that all SET 1 elements are positioned before all SET 2 elements, both of which are positioned before SET 3 elements. If any of the SETs 1 or 3 is empty, then the bound of the SRE 2 vertex will be unlimited in that direction (and limited by the first occurrence of the other SET 3 or 1 element in the opposite direction of Route6's columns).

As a future area of research, incorporating a threshold during perturbations to allow finding better solutions via worse solutions could be considered; adaptive thresholds \cite{10406762} could also be integrated within simulated annealing approaches.

\subsection{Overall Heuristic Overview} \label{Overall Heuristic Overview}

Each Main Iteration of the Heuristic is run as multiple Pipeline passes (depicted in \autoref{fig:MainHeuristicFlowchartOverview}). We start by sending each unsatisfied nodes one at a time through this pipeline, generate SREs for the full satisfaction of this node, and integrate RoPr of the SREs into all available Route Clusters; in the next pass we send another unsatisfied node and continue likewise till all nodes have been passed through the pipeline. This creates the final solution across all RoCus which is perturbed. The best cascaded solution among all clusters is deemed the best solution found in this Main Iteration. 

When we send all the nodes together in a  single pipeline pass, all the SREs pertaining to every node in the problem is developed. All these SREs are then integrated into route clusters and then perturbation takes places (after adding the final return vehicle depot for closed route-vehicle types). 

The final best solution from the Heuristic is output as the final best cascaded solution from among the best solutions from each Main Iteration. In our implementation, all Main Iterations are independent; the no. of nodes being passed in any Main Iteration can vary from a single node up to the full problem requirement being considered at once in a single pipeline. We believe this pipelining approach also plays a crucial role in generating diverse quality of solutions.

Our Heuristic does not have any memory-element being passed from one Main Iteration to another (considering this may be in the future research scope), since we already develop multiple diverse solutions through the different clusters with their respective logics labels (shuffled for all route clusters independently after an SRE has been integrated into them).

\subsection{Result of the illustrative example problem-instance using our heuristic} \label{Case Study Heuristic Result}
\begin{figure}
    \centering
    \includegraphics[width=1.01\linewidth]{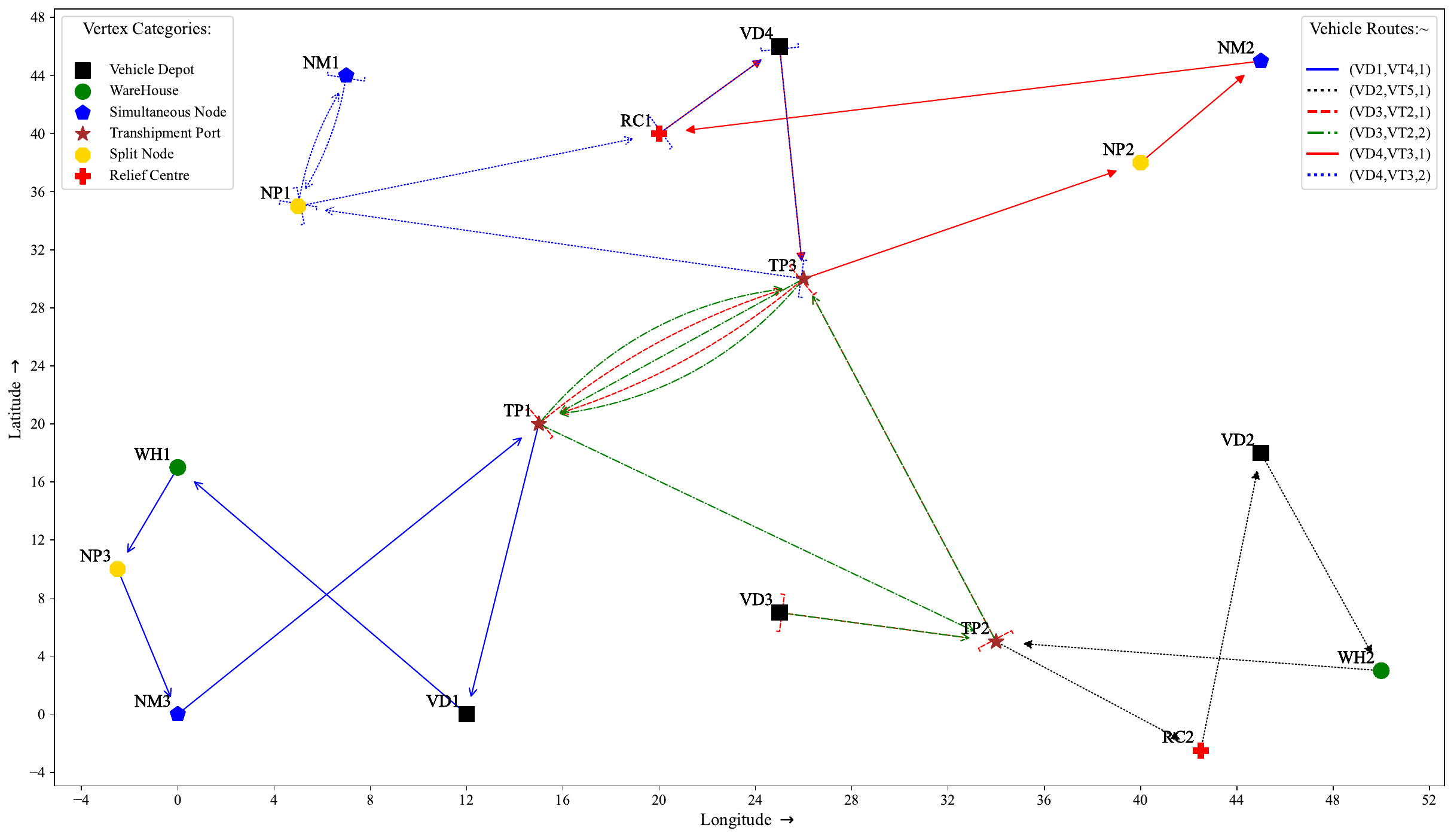} % 0.975 fits
    \caption{\textbf{A Heuristic result for our Study-Example Instance (Integer flow problem)}}
    \label{fig: casestudy solution heuristic}
\end{figure}

A good heuristic result of the cascaded optimization for the example instance in \autoref{fig:Picture_Instance} (and considering integral resource transfer) is shown in \autoref{fig: casestudy solution heuristic}, the routes being:

$\Downarrow$ \textbf{Cascade Number} \hfill \textbf{Route Durations} $\Downarrow$
\begin{enumerate}[label=\textbf{\scriptsize\arabic*}:] % Custom label format
    \item Vehicle (VD4, VT3, 1): VD4 $\rightarrow$ TP3 $\rightarrow$ NP2 $\rightarrow$ NM2 $\rightarrow$ RC1 $\rightarrow$ VD4 \hfill (116.87)    
    \item Vehicle (VD4, VT3, 2): VD4 $\rightarrow$ TP3 $\rightarrow$ NP1 $\rightarrow$ NM1 $\rightarrow$ NP1 $\rightarrow$ RC1 $\rightarrow$ VD4 \hfill (114.689)
    \item Vehicle (VD2, VT5, 1): VD2 $\rightarrow$ WH2 $\rightarrow$ TP2 $\rightarrow$ RC2 $\rightarrow$ VD2 \hfill (98.219)
    \item Vehicle (VD3, VT2, 2): VD3 $\rightarrow$ TP2 $\rightarrow$ TP3 $\rightarrow$ TP1 $\rightarrow$ TP3 $\rightarrow$ TP1 $\rightarrow$ TP2 \hfill (78.004)
    \item Vehicle (VD1, VT4, 1): VD1 $\rightarrow$ WH1 $\rightarrow$ NP3 $\rightarrow$ NM3 $\rightarrow$ TP1 $\rightarrow$ VD1 \hfill (74.002)
    \item Vehicle (VD3, VT2, 1): VD3 $\rightarrow$ TP2 $\rightarrow$ TP3 $\rightarrow$ TP1 $\rightarrow$ TP3 \hfill (64.398)
\end{enumerate}
The sum of all routes in this solution is 546.183 found in 4.21 minutes. The average of the highest route durations was 168.965 with an S.D. of 21.073; the best solution with the minimum sum of all routes was 511.47 (the average and S.D. of the solution with minimum sum of all routes being 678.94 and 90.303). We don't develop any separate algorithm for optimization based on the minimum sum of all routes, instead from the same pool of route clusters, we filter encountered solutions based on this criteria maintaining a separate list. These reportings had the average number of Main Iterations as 14.2 (S.D. 0.6) with an average time of 3.6 minutes (S.D. 1.98); for our study-example heuristic runs, we used the following bounds to obtain starting Heuristic input-parameters (from an uniform distribution):

\noindent{
\begin{tabular}{|m{0.7\textwidth}|m{0.1\textwidth}||m{0.1\textwidth}|}
\hline
\textbf{Heuristic Parameter} & \textbf{Lower Bound} & \textbf{Upper Bound} \\ \hline
No. of Internal Main Iterations & 10 & 20 \\ \hline
Max. No. of Perturbation in each Main Iteration: & 1729 & 8911 \\ \hline
Max. No. of Route Generation Logic Allowed from SRE: & 100 & 300 \\ \hline
\end{tabular}
\newline}

Detailed results are hosted online on GitHub (repository details is available after \autoref{Dataset Hosting}).

\section{Discussion and Computational Study} \label{Discussion and Computational Study}
The details of the supporting Compute Infrastructure are provided in Table S7 in the supplementary, along with the software versions used. We develop novel Instances for benchmarking the Heuristic \wrt MILP results.

\subsection{Small Instances} \label{Self Generated Small Instances}
The results from the small instances (dataset overview in supplementary Table S8) are presented in supplementary Table S11 and \autoref{Tab: Small Integer Results}. To effectively convey the comparison of the MILP solution with the Heuristic Algorithm (and to fit the results within a single screen-view), the column headings are encoded as described in \autoref{Tab: Encoded Column Headings}.

We generalize the makespan minimization through cascaded route optimization; we are able to show significant differences in the routes of the final cascade when compared to the first cascade's routes, as can be observed as $S_{FB}$ in Tables S11 and \ref{Tab: Small Integer Results}. This happens because there remains a slack even in the optimum solution after minimizing the makespan allowing the smaller routes to be further optimized, a single makespan minimization will never be able to achieve this.

As an alternative to the Cascaded Time Minimization, Makespan minimization followed by minimizing the sum of all the remaining route durations (apart from the highest route which forms the Makespan) could also be thought; the utility of these different optimization approaches must be researched \wrt application areas as we believe this approach of cascaded time minimization is best suited for our consideration within disaster management scenarios.

\begin{table}[ht]
\centering
\caption{\textbf{Encoded column headings for representing Computational Results of Small Instances from Supplementary Table S8, in the Table S11 for Continuous problems and the \autoref{Tab: Small Integer Results} for the Integer problems}\label{Tab: Encoded Column Headings}}
% \noindent
\resizebox{1\textwidth}{!}{ % Resize to fit the width of the page
% \begin{tabular}{|m{0.6875\textwidth}|m{0.75\textwidth}|m{0.055\textwidth}|}
\begin{tabular}{|m{0.5\textwidth}|m{0.65\textwidth}|m{0.066\textwidth}|} % 3 columns with a repeating pattern

\multicolumn{3}{m{0.225\textwidth}}{\textbf{\(\nearrow\)ENCODINGs\(\searrow\)}}\\
\cline{2-3}
\multicolumn{1}{c|}{} & \textbf{Explanation of the Encoded Symbols} & \textbf{\footnotesize Encoded Symbol}\\
\hline
\multirow{9}{=}{\textbf{Cascaded Makespan Minimization using our developed exact MILP Formulation}:\newline
\newline Here we use Gurobi's Warm Start to progress towards higher cascades; this happens by default when limited changes in an existing optimization model is performed allowing Gurobi to restart from the previous state(and solution). \newline\newline
At the end of every cascade, we add/update constraints limiting the route durations of all vehicles (except those vehicles which had the maximum route duration in any previous cascade, and whose route durations have been removed from constraining the objective function through the \autoref{eq:65}). This addition of upper limits to vehicle route durations (of unfixed vehicles at the end of every cascade) allows the route duration of any vehicle to flexibly change (\ie reduce further) in the subsequent cascades.
}	&	The highest Vehicle Route Duration (after completion of all Cascades) (i.e. Makespan)	&	$	O	$	\\
\cline{2-3}
	&	The sum of all vehicle's Route Durations (after all Cascades are complete)	&	$	S_F	$	\\
    \cline{2-3}
	&	Total Number of Cascades (equal to the Number of Vehicles)	&	$	\mathbb{C}	$	\\
    \cline{2-3}
	&	Time Taken for all Cascades of the specific Level to complete. Each Cascade was allowed to run for:
    \begin{itemize}
        \item $\frac{180}{\mathbb{C}}$ mins for the Continuous problems
        \item $\frac{500}{2^{\mathbb{C}-1}}$ mins for the Integer problems
    \end{itemize}
	&	$	T	$	\\
    \cline{2-3}
	&	Gap \% as reported by Gurobi after first Cascade	&	$	G_O	$	\\
    \cline{2-3}
	&	Average Gap \% across all Cascades	&	$	G_A	$	\\
    \cline{2-3}
	&	Levels used in this solution: \newline
    Levels provided to each vehicle were increased from 1 in case the model was found infeasible, or, no feasible solution was found in any of the previous cascades.	&	$	L	$	\\
    \cline{2-3}
	&	The sum of all vehicle's Route Durations, after the 1st Cascade	&	$	S_I	$	\\
    \cline{2-3}
	&	Improvement in the Sum of all Vehicle's Route Durations, w.r.t. final Cascade \newline
    Calculated as: \quad $\frac{S_I-S_F}{S_F}*100$ (\%) \newline
    &	$	I	$	\\
\hline
\multirow{9}{=}{\textbf{PSR-GIP Heuristic}:
Cascaded Makespan Minimization sense of Optimization
\begin{itemize}
    \item For each Continuous Instance, 250 independent runs of the Heuristic Algorithm were performed (and reported in Table S11)
    \item For each Integer Instance, 175 independent runs of the Heuristic Algorithm were performed  (and reported in \autoref{Tab: Small Integer Results})
\end{itemize}}	&	The Largest Route Duration among all vehicles of a solution (Best solution among All runs)	&	$	L_B	$	\\
\cline{2-3}
	&		Average of the Highest Vehiclular Route Duration (across All runs)	&	$	L_A	$	\\
    \cline{2-3}
	&		S.D. of the Highest Vehiclular Route Duration	&	$	L_{SD}	$	\\
    \cline{2-3}
	&		The Aggregate of all Vehicular Route Durations of the Best solution	&	$	S_B	$	\\
    \cline{2-3}
	&		Average (across All runs) of the Aggregated Vehicular Route Durations	&	$	S_A	$	\\
    \cline{2-3}
	&		S.D. (across All runs) of the Aggregated Vehicular Route Durations	&	$	S_{SD}	$	\\
    \cline{2-3}
	&		No. of Vehicles used in the Best Solution / Out of Total Available Vehicles	&	$	V	$	\\
    \cline{2-3}
	&	Average Heuristic RunTime (minutes)	&	$	R_A	$	\\
    \cline{2-3}
	&	S.D. Of Heuristic RunTimes (minutes)	&	$	R_{SD}	$	\\
\hline
\multirow{2}{=}{\textbf{Comparative Gap \%} of Heuristic Best solution \wrt MILP Formulation}	&	Gap \% \wrt the Largest Vehicle's Route Duration\newline
Calculated as: \quad $\frac{L_B-O}{L_B}*100$\newline	&	\textbf{$	OL_B	$}	\\
\cline{2-3}
	&	Gap \% \wrt the Aggregate of all Vehicle's Route Durations\newline
Calculated as: \quad $\frac{S_B-S_F}{S_B}*100$\newline	&	$	S_{FB}	$	\\
\hline
\end{tabular}
} % For Resize Box Closure
\end{table}

The computational results considering the flow variables as Continuous are presented in Table S11 under supplementary section S4.1, and the results when considering the flow variables as Integers are presented in \autoref{Tab: Small Integer Results} The parameters used during these MILP runs (for both the Continuous and Integer problem) along with the (non-default) Gurobi settings are:%\newline

\footnotesize
\noindent{
% \begin{tabular}{|c|c|}
\begin{tabular}{|>{\centering\arraybackslash}m{0.4\textwidth}|>{\centering\arraybackslash}m{0.55\textwidth}|}
\hline
Big M (used in Formulation) & 3000 \\ \hline
Constraint Feasibility Tolerance (FeasibilityTol) & $1E-5$ (default is $1E-6$) \\ \hline
IntegralityFocus & 1 (default is 0, setting it to 1 helps allow better satisfaction of integer variables by preventing over-exploitation of integer tolerances) \\ \hline
\end{tabular}
\newline}\small

The Heuristic's Parameter values used were:%\newline

\footnotesize
\noindent{
\begin{tabular}{|m{0.5\textwidth}|m{0.2\textwidth}|m{0.2\textwidth}|}
\hline
\textbf{Heuristic User-Inputs} & \textbf{Continuous Problem} & \textbf{Integer Problem} \\ \hline
No. of Internal Main Iterations & 10 & 15 \\ \hline
Max. No. of Perturbation in each Main Iteration: & 1729 & 4104 \\ \hline
Max. No. of Route Generation Logic Allowed from SRE: & 100 & 200 \\ \hline
\end{tabular}
\newline}\small

We clearly observe from the column \textbf{$	OL_B$} in \autoref{Tab: Small Integer Results}, that the Relative Gap \% for the Cascaded Makespan Optimization results as found by our Heuristic (best solution) is within 10\% of the solution found by the MILP Formulation (the Average Relative Gap \% across all Integer Instances of \autoref{Tab: Small Integer Results} is \textbf{3.1\%}). Some points to observe for these Integer problems are:
\begin{enumerate}
    \item For the instance \instbox{\footnotesize S6}, the formulation was infeasible with Levels=1; no solution was found with Levels=2, 3, or 4.
    
    \item For the instance \instbox{\footnotesize S7}, the optimal solution using only 1 level proved to have an incomplete solution space, as the heuristic solution was found to be better. The MILP for this case was also run using levels equal to 2 (the relative gaps of the Heuristic’s best-found solution w.r.t. both the MILP runs is reported).
    
    \item For the instance \instbox{\footnotesize S9} the vehicle with the largest trip length in the first cascade also had the largest trip length in the final cascade (just like in most of the other instances), but this route duration value got decreased in the final cascade (compared to what had been the cascade 1 objective function, and subsequently used for constraining all vehicle trip durations at the end of Cascade 1). The Gap \% in this case should therefore be slightly lower, but it has not been adjusted.
    Further, no solution was found during the second Cascade with Levels=1; (as Gurobi Warm Start did not seem to continue (even after repeated attempts)).
    
\end{enumerate}

\begin{landscape}
\tiny
\vspace*{\fill}  % Push down
\begin{table}[ht]
% \vfill
\centering
\caption{\textbf{Computational Results of Small Instances from Supplementary Table S8 with resource flow variables constrained to being Integers}\label{Tab: Small Integer Results}}
% \noindent
\resizebox{1.55\textwidth}{!}{ % Resize to fit the width of the page
\begin{tabular}{*{21}{|c}|}
\cline{2-21}
\multicolumn{1}{c|}{} & \multicolumn{9}{|
c|}{\makecell{\textbf{Cascaded Makespan Minimization} \\ \textbf{using our developed exact MILP Formulation}}} & \multicolumn{9}{|c|}{\textbf{PSR-GIP Heuristic}} & \multicolumn{2}{|c|}{\makecell{\textbf{Comparative} \\ \textbf{Gap \%}}} \\
\hline
\textbf{ID}	&	\textbf{$O$}	&	\textbf{$S_F$}	&	\textbf{$\mathbb{C}$}	&	\textbf{$T$}	&	\textbf{$G_O$}	&	\textbf{$G_A$}	&	\textbf{$L$}	&	\textbf{$S_I$}	&	\textbf{$I$}	&	\textbf{$L_B$}	&	\textbf{$L_A$}	&	\textbf{$L_{SD}$}	&	\textbf{$S_B$}	&	\textbf{$S_A$}	&	\textbf{$S_{SD}$}	&	\textbf{$V$}	&	\textbf{$R_A$}	&	\textbf{$R_{SD}$}	&	\textbf{$OL_B$}	&	\textbf{$S_{FB}$}	\\
\hline
\hline
\instbox{S1}	&	64.516	&	79.629	&	2	&	0.00	&	0.00	&	0.00	&	1	&	116.307	&	46.06	&	64.516	&	65.456	&	2.679	&	79.629	&	80.721	&	2.674	&	2/2	&	0.61	&	0.03	&	0.00	&	0.00	\\
\instbox{S2}	&	64.278	&	152.772	&	4	&	1.23	&	0.00	&	0.00	&	1	&	172.665	&	13.02	&	64.278	&	65.159	&	1.54	&	157.264	&	175.42	&	13.161	&	4/4	&	1.35	&	0.11	&	0.00	&	2.86	\\
\instbox{S3}	&	68.122	&	251.274	&	4	&	133.39	&	0.00	&	13.45	&	1	&	256.254	&	1.98	&	68.326	&	88.089	&	10.249	&	194.985	&	241.438	&	34.554	&	4/4	&	2.34	&	0.2	&	0.30	&	-28.87	\\
\instbox{S4}	&	88.614	&	195.078	&	3	&	750.45	&	48.82	&	27.22	&	3	&	195.173	&	0.05	&	88.614	&	109.421	&	7.866	&	195.108	&	220.429	&	17.067	&	3/3	&	2.67	&	0.19	&	0.00	&	0.02	\\
\instbox{S5}	&	50.569	&	262.56	&	7	&	992.23	&	98.33	&	99.76	&	2	&	262.56	&	0.00	&	53.828	&	57.707	&	2.172	&	268.134	&	279.325	&	7.604	&	7/7	&	3.29	&	0.37	&	6.05	&	2.08	\\
\instbox{S6}	&	No Sol.	&	-	&	-	&	-	&	-	&	-	&	4	&	-	&	-	&	102.277	&	104.979	&	3.29	&	441.597	&	462.145	&	29.464	&	7/7	&	5.11	&	0.36	&	-	&	-	\\\hline
\multirow{2}{*}{\instbox{S7}}	&	151.937	&	388.611	&	5	&	3.24	&	0.00	&	0.00	&	1	&	479.326	&	23.34	&	\multirow{2}{*}{141.194}	&	\multirow{2}{*}{141.861}	&	\multirow{2}{*}{0.267}	&	\multirow{2}{*}{450.8}	&	\multirow{2}{*}{386.727}	&	\multirow{2}{*}{31.213}	&	\multirow{2}{*}{5/5}	&	\multirow{2}{*}{1.54}	&	\multirow{2}{*}{0.2}	&	-7.61	&	13.80	\\
	&	141.809	&	364.835	&	5	&	968.77	&	26.13	&	85.23	&	2	&	456.791	&	25.20	&		&		&		&		&		&		&		&		&		&	-0.44	&	19.07	\\\hline
\instbox{S8}	&	136.139	&	421.667	&	5	&	968.77	&	26.99	&	85.40	&	2	&	569.818	&	35.13	&	152.2	&	153.058	&	0.703	&	376.127	&	390.237	&	19.929	&	5/5	&	3.45	&	0.23	&	10.55	&	-12.11	\\
\instbox{S9}	&	101.8	&	448.559	&	7	&	992.28	&	90.67	&	98.67	&	2	&	680.747	&	51.76	&	111.582	&	125.755	&	9.022	&	467.843	&	517.12	&	45.458	&	6/7	&	5.42	&	0.56	&	8.77	&	4.12	\\
\instbox{S10}	&	104.776	&	506.891	&	7	&	992.23	&	74.87	&	87.48	&	1	&	620.241	&	22.36	&	113.16	&	129.072	&	8.011	&	497.505	&	546.797	&	45.472	&	7/7	&	8.9	&	0.71	&	7.41	&	-1.89	\\
\instbox{S11}	&	84.667	&	354.589	&	6	&	455.05	&	0.00	&	18.29	&	1	&	386.301	&	8.94	&	84.931	&	85.467	&	0.43	&	368.968	&	384.883	&	20.224	&	6/6	&	3.06	&	0.18	&	0.31	&	3.90	\\
\instbox{S12}	&	125.107	&	604.593	&	7	&	992.33	&	30.94	&	88.49	&	1	&	679.992	&	12.47	&	136.38	&	152.257	&	8.448	&	610.26	&	632.575	&	43.788	&	7/7	&	8.23	&	1.12	&	8.27	&	0.93	\\
\instbox{M13}	&	84.418	&	407.068	&	6	&	52.70	&	0.00	&	17.15	&	1	&	501.11	&	23.10	&	84.738	&	86.012	&	1.022	&	392.075	&	435.048	&	30.336	&	6/6	&	4	&	0.19	&	0.38	&	-3.82	\\
\instbox{M14}	&	87.179	&	367.202	&	6	&	419.69	&	0.00	&	26.46	&	1	&	470.793	&	28.21	&	89.928	&	98.814	&	2.127	&	436.253	&	464.425	&	30.888	&	6/6	&	5.05	&	0.27	&	3.06	&	15.83	\\
\instbox{M15}	&	68.48	&	253.605	&	4	&	190.39	&	0.00	&	45.15	&	1	&	261.564	&	3.14	&	74.547	&	109.044	&	11.183	&	205.648	&	310.256	&	51.143	&	4/4	&	3.09	&	0.36	&	8.14	&	-23.32	\\
\instbox{M16}	&	119.277	&	303.135	&	4	&	694.25	&	0.00	&	0.00	&	1	&	337.684	&	11.40	&	119.277	&	119.277	&	0	&	305.026	&	306.5	&	5.875	&	4/4	&	2.57	&	0.28	&	0.00	&	0.62	\\
\instbox{M17}	&	102.921	&	395.028	&	4	&	937.52	&	57.81	&	73.49	&	1	&	410.966	&	4.03	&	102.921	&	103.614	&	0.536	&	398.577	&	397.922	&	6.311	&	4/4	&	2.27	&	0.19	&	0.00	&	0.89	\\
\instbox{M18}	&	49.185	&	142.648	&	4	&	183.42	&	0.00	&	0.00	&	1	&	160.395	&	12.44	&	49.185	&	49.448	&	0.439	&	112.618	&	135.41	&	16.943	&	3/4	&	1.55	&	0.21	&	0.00	&	-26.67	\\
\instbox{M19}	&	88.789	&	359.811	&	6	&	244.79	&	0.00	&	16.82	&	1	&	397.975	&	10.61	&	96.23	&	103.382	&	2.641	&	377.329	&	398.492	&	10.852	&	6/6	&	7.02	&	0.53	&	7.73	&	4.64	\\
\instbox{M20}	&	136.831	&	525.789	&	6	&	478.88	&	0.00	&	3.87	&	1	&	556.593	&	5.86	&	138.032	&	147.49	&	3.028	&	582.489	&	584.866	&	12.852	&	6/6	&	5.91	&	0.57	&	0.87	&	9.73	\\
\instbox{M21}	&	41.799	&	233.326	&	6	&	1.90	&	0.00	&	0.00	&	1	&	221.536	&	-5.05	&	42.731	&	58.902	&	8.879	&	172.905	&	214.22	&	23.809	&	5/6	&	1.69	&	0.22	&	2.18	&	-34.94	\\
\instbox{M22}	&	34.485	&	159.547	&	5	&	968.78	&	5.79	&	51.90	&	1	&	159.567	&	0.01	&	35.196	&	38.115	&	1.05	&	160.112	&	165.9	&	4.989	&	5/5	&	4.9	&	0.36	&	2.02	&	0.35	\\
\instbox{M23}	&	98.765	&	461.996	&	6	&	219.32	&	0.00	&	6.98	&	1	&	537.231	&	16.28	&	99.609	&	107.032	&	4.454	&	504.868	&	525.437	&	29.017	&	6/6	&	6.82	&	0.26	&	0.85	&	8.49	\\
\instbox{M24}	&	32.71	&	147.564	&	6	&	839.42	&	25.92	&	28.86	&	1	&	186.007	&	26.05	&	33.739	&	38.214	&	1.619	&	160.435	&	165.559	&	6.598	&	6/6	&	8.27	&	0.4	&	3.05	&	8.02	\\
\instbox{M25}	&	88.097	&	383.595	&	6	&	241.99	&	0.00	&	33.33	&	1	&	435.233	&	13.46	&	92.666	&	104.848	&	4.649	&	372.48	&	403.182	&	18.486	&	6/6	&	12.15	&	1.24	&	4.93	&	-2.98	\\
\hline
\end{tabular}
} % For Resize Box Closure
\end{table}
{\small The prefix in the ID is a rough estimation of the problem size for each instance (applicable for all the Tables S8, S9, \ref{Tab: Small Integer Results}, S11, \ref{longTab: Large Instance Results Integer} and S12):

$\instbox{S}\rightarrow Small \quad\quad \instbox{M}\rightarrow Medium \quad\quad \instbox{L}\rightarrow Large$}
% \# Infeasible: 1; No Solution: 2, 3, 4
\vspace*{\fill}  % Push down
\end{landscape}

\begin{landscape}
\begin{multicols}{2}
\subsection{Large Instances} \label{Self Generated Large Instances}
The overview of the dataset details for the large instances we developed can be found in Supplementary Table S9.
These Instances are run as Continuous and Integer problems, affecting the variable $y$ in the MILP formulation, and the Heuristic is also tuned to match the same. The computational results of the problems with the resource flow variables as continuous are presented in Supplementary Table S12, while the integer consideration of the same problems is presented in \autoref{longTab: Large Instance Results Integer}. During both the computational cases (both Integer and Continuous Problems), the Gurobi (non-default) settings for the MILP runs were:\newline

\scriptsize
\noindent{
% \begin{tabular}{|c|>{\centering\arraybackslash}m{0.4\textwidth}|}
\begin{tabular}{|m{0.2\textwidth}|m{0.475\textwidth}|}
\hline
Big M (used in Formulation) & 1E5 \\ \hline
MIP Focus & 1 (to focus on finding feasible solutions) \\ \hline
Integer Feasibility Tolerance (IntFeasTol) & $1E-7$ (default is $1E-5$) \\ \hline
IntegralityFocus & 1 (default is 0, setting it to 1 helps allow better satisfaction of integer variables by preventing over-exploitation of integer tolerances) \\ \hline
\end{tabular}
\newline}
\small

and the Heuristic Parameter values used were:\newline

\scriptsize
\noindent{
\begin{tabular}{|c|c|}
\hline
No. of Internal Main Iterations & 10 \\ \hline
Max. No. of Perturbation in each Main Iteration: & 1729 \\ \hline
Max. No. of Route Generation Logic Allowed from SRE: & 100 \\ \hline
\end{tabular}
\newline}
\small

We always compare only the fully solved problems for all instances considered across all datasets in \autoref{Discussion and Computational Study}, as we don't consider partial satisfaction anywhere in this study. For every Main Iteration within every Instance-runs for the results reported in 
Tables S11, \ref{Tab: Small Integer Results}, S13 and
\ref{longTab:Literature_Integer_Instances} there is no Unsatisfiable Portion ever encountered by our Heuristic (showing the strength of our designed algorithm). Only for few of the complex large instances in Tables S12 and \ref{longTab: Large Instance Results Integer} do we encounter some of the Instance-runs within which every Main Iteration failed to find any solution.

For the Integer Problems in \autoref{longTab: Large Instance Results Integer}:
\begin{enumerate}
    \item In the instance \instbox{\footnotesize L31}: 35 out of 44 Instance runs were successful 
    \item In the instance \instbox{\footnotesize L32}: 20 out of 30 Instance runs were successful; the remaining 10 runs had no Main Iteration that was able to satisfy all the resource requirements of the problem to their fullest
    \item In the instance \instbox{\footnotesize L34}: 25 out of 30 Instance runs were successful
\end{enumerate}

These indicate that for very complex large instances, with multiple intricate compatibility considerations within the problems, the Heuristic is able to generate good solutions at least \~70\% of the time or more.
\end{multicols}

% \footnotesize
\small
\begin{longtable}{|m{0.025\linewidth}|m{0.085\linewidth}|m{0.07\linewidth}|m{0.095\linewidth}|m{0.05\linewidth}|m{0.08\linewidth}|m{0.12\linewidth}|m{0.11\linewidth}|m{0.065\linewidth}|m{0.065\linewidth}|m{0.075\linewidth}|}
\caption{\textbf{Computational Results for the Large Instances (instance details in Supplementary Table S9) considering Integer resource flow}\label{longTab: Large Instance Results Integer}}\\
\hline
\multirow{2}{=}{\newline\textbf{ID}} & \multicolumn{4}{|c|}{\textbf{MILP Formulation}} & \multicolumn{5}{|c|}{\textbf{PSR-GIP Heuristic}} & \multirow{2}{=}{\newline\textbf{Relative Comparative Gap \% \newline\newline $\frac{b-a}{b}*100$}} \\
\cline{2-10}
& Objective Value (Cascade 1 only, if any) \newline\newline $(a)$	& Gap \% (as reported by Gurobi) & Levels Used (I$\rightarrow$Infeasible; \newline N$\rightarrow$NFSF; \newline F$\rightarrow$Feasible) \#& Time (mins) & Best Objective (among at least 30 runs) \newline\newline $(b)$ & Average (of the largest route duration among all Vehicles across successful Heuristic runs) & S.D. (of the vehicle with maximum route duration across successful Heuristic runs) & Average Heuristic RunTime (mins) & S.D. of Heuristic Runtimes (mins) & \\
\hline
\endfirsthead
\hline
\endhead
\hline
\instbox{\footnotesize S26}	&	468.118	&	0	&	F:5	&	0.07	&	468.118	&	468.118	&	0	&	0.3	&	0	&	0.00	\\
\instbox{\footnotesize S27}	&	645.164	&	0	&	F:3	&	53.55	&	645.164	&	657.601	&	30.382	&	0.53	&	0.04	&	0.00	\\
\instbox{\footnotesize M28}	&	1007.834	&	92.569	&	N:1,2; F:5	&	1440	&	994.522	&	1121.995	&	57.748	&	1.85	&	0.23	&	-1.34	\\
\instbox{\footnotesize M29}	&	NFSF	&	-	&	N:1	&	1440	&	14417.25	&	16847.385	&	971.673	&	31.69	&	8.83	&	-	\\
\instbox{\footnotesize L30}	&	NFSF	&	-	&	I:5; N:6	&	OEAM	&	220.544	&	384.252	&	71.529	&	71.62	&	10.22	&	-	\\
\instbox{\footnotesize L31}	&	NFSF	&	-	&	N:1	&	1440	&	809.61	&	1222.469	&	247.362	&	21.94	&	8.08	&	-	\\
\instbox{\footnotesize L32}	&	NFSF	&	-	&	N:5,3	&	1440	&	994.712	&	1514.093	&	243.862	&	14.09	&	9.27	&	-	\\
\instbox{\footnotesize L33}	&	NFSF	&	-	&	I:3; N:7,5	&	1440	&	470.141	&	819.59	&	272.271	&	11.16	&	4.19	&	-	\\
\instbox{\footnotesize L34}	&	NFSF	&	-	&	N:5	&	1440	&	10550.721	&	14602.567	&	2868.328	&	20.38	&	10.87	&	-	\\
\instbox{\footnotesize L35}	&	NFSF	&	-	&	N:5,2	&	1440	&	1828.904	&	2749.606	&	693.81	&	56.65	&	41.45	&	-	\\
\instbox{\footnotesize L36}	&	NFSF	&	-	&	N:7,5,4	&	1440	&	340.644	&	518.671	&	82.805	&	83.1	&	13.67	&	-	\\
\instbox{\footnotesize L37}	&	NFSF	&	-	&	I:2; N:3	&	1370	&	411.3	&	771.423	&	175.041	&	72.38	&	15.27	&	-	\\
\instbox{\footnotesize L38}	&	NFSF	&	-	&	I:1; N:2	&	OEAM	&	1208.275	&	1866.051	&	713.845	&	98.6	&	43.08	&	-	\\
\hline
\end{longtable}
\quad\quad NFSF $\Longrightarrow$ No Feasible Solution Found \quad\quad\quad\quad\quad\quad OEAM $\Longrightarrow$ Optimization Exhausted Available Memory
\newline
% \newline
\indent\scriptsize\# All Vehicles were provided with the same number of levels, and the result of those instances is briefed. For each of the values mentioned, separate 24 hour MILP runs were performed. The instance results in this table correspond to the last value in the Levels Used column. Further, for levels indicating Infeasibility (due to lack of solution space restricting multi-trips), all lower-Level values were also infeasible (since the solution space is even smaller then with much lesser number of multi-visits possible at any vertex).

\end{landscape}

\begin{landscape}
\small
% \footnotesize
% \scriptsize
% \tiny

\subsection{Datasets from Literature: Considering reduced version of our problem} \label{Benchmarking using previously considered reduced versions of our problem}

We take instances from the Table 5 of \cite{NUCAMENDIGUILLEN2021}; the Instance Name $I_1$ corresponds to their case study instance \footnote{The corrected datasets as sent by the authors of \cite{NUCAMENDIGUILLEN2021} are available along with all other datasets (see after \autoref{Dataset Hosting})}. The overview of the moulded dataset to our problem's specifications can be found in Supplementary Table S10. We use these instances to generate Continuous and Integer Problems, depending on the allowed variable values of the resource flow variables. Continuous variables would mean that fractional pickup and delivery is allowed at any vertex, whereas restricting them to Integers corresponds more towards realistic scenarios (as in real-world, Cargo Type units would not be sub-divided).
The table \autoref{longTab:Literature_Integer_Instances} corresponds to case where the problems are run using integer variables (for both the MILP and the Heuristic); the MILP computations were performed on the Computer C\_1, and the Heuristic Algorithms were run on the Computer C\_4. It is observed that here too, the relative gap is significantly low, and is generally lower than the gaps found for the continuous case; this indicates that the Heuristic is suitable for real-world optimization problems.

\begin{longtable}{|m{0.02\linewidth}|m{0.105\linewidth}|m{0.085\linewidth}|m{0.09\linewidth}|m{0.15\linewidth}|m{0.085\linewidth}|m{0.075\linewidth}|m{0.07\linewidth}|m{0.05\linewidth}|m{0.05\linewidth}|m{0.06\linewidth}|}
% Paper Details along with spatiotemporal details and Objectives of Papers reviewed
\caption{\textbf{Computational Results of the Instances (instances described in Supplementary Table S10) considering Integer Resource Flows}\label{longTab:Literature_Integer_Instances}}\\
% \hline
% \multicolumn{2}{| c |}{Begin of Table}\\
\hline
\multirow{3}{=}{\rotatebox{75}{\textbf{Instance}}} & \multicolumn{3}{|c|}{\textbf{Cascaded Makespan Minimization}} & \multicolumn{6}{|c|}{\textbf{PSR-GIP Heuristic \#}}& \multirow{3}{=}{\newline\textbf{Relative Gap \% \newline\newline$\frac{(b-a)*100}{b}$}}\\

\cline{2-10}
& \multicolumn{3}{|c|}{\textbf{Our Exact MILP Formulation*}} & \multicolumn{4}{|c|}{\textbf{Cascaded Makespan Minimization Optimization}}& \multicolumn{2}{|c|}{\textbf{Runtime}}&\\

\cline{2-10}
& \textbf{Objective Function Value of Cascade 1 \newline ($a$)} & \textbf{No. of Vehicles used} &\textbf{MILP Gap \% (reported by Gurobi)}
& \textbf{Best Solution of Largest Route Duration among all Vehicles \newline ($b$)} & \textbf{No. of Vehicles used in Best Solutions} & \textbf{Average of Best Solutions} & \textbf{S.D. of Best Solutions} & \textbf{Avg. (mins)} & \textbf{S.D. (mins)} & \\
\hline

\hline
\endfirsthead

$	I_1	$	&	40.706	&	9	&	60.46	&	45.706	&	10	&	73.074	&	20.267	&	3.79	&	0.45	&	10.94	\\
$	I_2	$	&	38.77	&	9	&	58.486	&	38.77	&	11	&	51.104	&	15.029	&	3.27	&	0.41	&	0.00	\\
$	I_3	$	&	40.706	&	10	&	60.46	&	42.48	&	10	&	53.46	&	13.662	&	2.95	&	0.41	&	4.18	\\
$	I_4	$	&	40.706	&	9	&	60.46	&	40.943	&	8	&	48.967	&	9.138	&	2.49	&	0.32	&	0.58	\\
$	I_5	$	&	38.77	&	8	&	58.481	&	39.135	&	10	&	50.309	&	11.61	&	2.6	&	0.39	&	0.93	\\
$	I_6	$	&	39.793	&	7	&	59.554	&	39.793	&	11	&	51.301	&	11.699	&	2.97	&	0.37	&	0.00	\\
$	I_7	$	&	40.706	&	7	&	60.447	&	40.706	&	10	&	49.133	&	11.256	&	2.57	&	0.32	&	0.00	\\
$	I_8	$	&	40.706	&	8	&	60.458	&	40.713	&	8	&	50.167	&	11.03	&	2.56	&	0.28	&	0.02	\\
$	I_9	$	&	33.787	&	6	&	0	&	34.787	&	7	&	40.808	&	5.341	&	1.62	&	0.23	&	2.87	\\
$	I_{10}	$	&	40.706	&	8	&	60.458	&	40.943	&	9	&	54.043	&	13.705	&	2.7	&	0.39	&	0.58	\\
$	I_{11}	$	&	40.706	&	11	&	60.46	&	42.274	&	11	&	59.865	&	22.614	&	3.84	&	0.51	&	3.71	\\
$	I_{12}	$	&	40.706	&	9	&	59.064	&	40.706	&	9	&	46.674	&	8.355	&	2.16	&	0.24	&	0.00	\\
$	I_{13}	$	&	38.77	&	9	&	58.486	&	38.77	&	9	&	44.285	&	6.551	&	2.83	&	0.3	&	0.00	\\
$	I_{14}	$	&	40.706	&	9	&	60.46	&	40.706	&	11	&	52.394	&	8.634	&	3.44	&	0.5	&	0.00	\\
$	I_{15}	$	&	40.706	&	10	&	60.46	&	41.405	&	9	&	50.961	&	7.086	&	3.04	&	0.37	&	1.69	\\
\hline
\end{longtable}
Gurobi (non-default) settings for the MILP runs:
\noindent{
\begin{tabular}{|c|c|}
\hline
Big M (used in Formulation) & 3000 \\ \hline
Constraint Feasibility Tolerance (FeasibilityTol) & $1E-5$ (default is $1E-6$) \\ \hline
\end{tabular}
\newline
\newline
% \par
}

Heuristic Parameters used:
\noindent{
\begin{tabular}{|c|c|}
\hline
No. of Internal Main Iterations & 20 \\ \hline
Max. No. of Perturbation in each Main Iteration: & 8911 \\ \hline
Max. No. of Route Generation Logic Allowed from SRE: & 300 \\ \hline
\end{tabular}
\newline
% \newline
\par
}
* All MILPs were run for 24 hours with the Instance $I_9$ completing in 14.6 hrs (only first Cascades were run)
\newline
% \newline
\indent\# The Best, Average and S.D. are reported comparing 99 Heuristic runs for each Instance
\end{landscape}
\small

\section{Conclusion} \label{Conclusion}
Disaster Management is an encompassing topic and in this study, we cover one important aspect of optimizing the time for emergency operations. We discuss a rich vehicle routing problem specific to disaster scenarios, defining a novel problem. For the first time, we separate the idea of depots from being packed with the all functionalities (generally depots combine the functionality of our warehouses, relief centres and vehicle depots), and allow the consideration of different functional nodes (multiple in number); we believe this is also an appropriate way of distributing disaster response resources and facilities across geographies, instead of maintaining all of them in a single space which may become critical in the case of man-made emergencies. We consider multiple cargo types having their separate (un)loading times for each vehicle type; vehicle-to-road (or vehicle-to-mode) compatibilities are similarly considered as in \cite{10406762}. This problem conception focuses on a decentralized scheme where all vehicles, resources and pickups need not to be in the same place. In case it is necessary to combine any of the functional vertices (VDs, WHs or RCs), for the same geographical location, they may be created separately (at the same latitude, longitude) in the same spot.

% Summary
We propose a complex MILP formulation for the conceived problem and develop an encompassing Heuristic, ensuring it to be robust across diverse instance categories. Our formulation provides the solution space for multiple trips by introducing new levels to the vehicles. Initially, since we do not know how many trips would be required, we provide all vehicles with level 1 (single trip solution space), and try to optimize within this solution space. If the problem is found to be infeasible, we increase this solution space by incrementing the value of levels provided to each vehicle by one. The formulation developed is powerful enough to find optimal solutions to many instances, in spite of having intricate and complex compatibility constraints implemented within it.

To help solve this problem within real-world timeframes, we design a novel heuristic approach. The heuristic starts by creating decision trees to understand the problem feasibility and compatible cargo transfers from a PRV to a Node using heterogeneous vehicles stationed on the same SMTS. These CTs while effectively representing the feasible problem space (inherently capturing vehicle-load, and TP-cargo compatibility) also store preferences which allow preferential SRE generation. The smallest route elements may be connected to each other (causality connections as represented through the CD) if they are part of the same transhipment trail. All generated SREs are allocated to individual vehicles and are inserted (\ie integrated via matching) within the existing routes; the purpose of matching vertices here is to prevent a vertex from being visited multiple times in the same vehicle's route and try to merge all the tasks that need to be performed at that vertex in a single visit. Feasible perturbations (ensuring vehicles' capacities are not violated) are also considered to allow the finding of better neighbourhood solutions.

We develop novel datasets (small, medium and large instances) and compare the results of our formulation and heuristic; the results obtained are within acceptable relative gaps. The complexity of our novel feature-rich problem arises due to the consideration of transhipment at Ports which take up most of the space both in the formulation and heuristic design. Future research directions for the Heuristic are highlighted at their appropriate contextual places; mentioning some of the most important future research scopes briefly here: the computation of waiting times during every unique transhipment port visit of any vehicle would need an advanced accurate approach which could allow better solutions.

The problem complexity grows due to necessity of minimizing time, while also ensuring transhipments. We firmly believe these aspects both together have never been considered before; considering these aspects together mandates temporal causality to come into the picture (had this been distance/ cost minimization with transhipments, the formulation and heuristic development would have been simpler). Our main contributions lie in the novel problem conception highlighting the importance of the cascaded (multi-step) optimization approach, as well as development of a formulation and heuristics that are able to capture the aspect of causality; as rightly highlighted in \cite{Kern2025}, the need for such causal modelling arises from the fact that complex and accurate decision-making mandates analyzing the problem through the temporal picture. We further stress on the fact that removing waiting times should not make much complexity difference (as causality still remains), but could be drastic in terms of solution quality; the objective value could degrade dramatically. Waiting time considerations should always be included when solving such problems, allowing a wider optimality search where sometimes waiting (a vehicle staying idle) instead of routing (the vehicle mandatorily moving, if the problem-modelling does not allow it to wait) could prove to be better overall.

We have also developed a GIS plugin, utilizing the open-sourced QGIS software \cite{QGIS_software}, as part of a project to design a spatial decision support system for usage across disasters/emergencies (the details of the plugin and its developmental steps is discussed in an upcoming paper). Our Heuristic Algorithm has been integrated into this QGIS plugin and is being thoroughly tested for robustness across many disaster scenarios.

The developed formulation and heuristic may be also used for slightly modified problems than the one considered by us. For example, this heuristic can also take care of last-mile pickup-delivery problems by considering each specific shipment as a unique cargo type. 
Furthermore, as suggested in \cite{ALLAHYARI2015756}, combination of covering and routing can be considered in disaster-specific scenarios.

% Future Works
One major feature addition that may be considered in the context of much more extensive disaster rescue problems could be the introduction of another new vehicle type that can send out smaller vehicles (rescuer craft) and receive them (the rescuer craft with the rescued). This is practically done in the case of flood situations where a land-based vehicle carries uninflated lifeboats, which are then inflated to reach out to stranded individuals. Generally, during this process, the land-vehicle remains immobile, but a future research consideration of allowing the mother-vehicle to remain mobile while the sidekick-vehicles perform their tasks and intercept the mother vehicles at an appropriate location (to tag along/ get absorbed). Such problems would also prove to be fruitful when sending out rescue drones from moving vehicles; future research in this direction seems necessary and the call of the hour given the rising tides of both natural and man-made catastrophes.

\footnotesize
% \scriptsize

\section*{Online Repository for Datasets and our Study-Example Results}  \label{Dataset Hosting}
GitHub repository containing all problem instances and our study-example results: \hyperlink{https://github.com/SanTanBan/rich-MultiTrip-MultiModal-VRP}{https://github.com/SanTanBan/rich-MultiTrip-MultiModal-VRP}.%\newline
This repository also contains videos of a few solutions of the MILP for showcasing the layered modelling approach. Elaborate in-depth description of the formulation development steps (splitting the basic equations from the transhipment specific constraints) is available in \cite{my_MS_Thesis} including details on the PSR-GIP Heuristic and major pseudocodes. For further details, please contact the first/corresponding author.

\section*{CRediT authorship contribution statement} \label{CRediT authorship contribution statement}
% \noindent \textbf{Santanu Banerjee}: Conceptualization, Data curation, Software, Supervision, Validation, Visualization, Problem Conception, Literature Review, Methodology and Investigation, Formulation Development, Heuristic Construction, Coding Formulation and Heuristic in Python, Dataset Construction for Computational Result Comparison, Writing Original Draft, Creating Flowchart Visualization and Pictures, Formal Analysis, Editing and Reviewing the Manuscript, Integration of Heuristic Algorithm with GIS Plugin.
% \textbf{Goutam Sen}: Project Administration, Funding Acquisition, Supervision, Resources, Review.
% \textbf{Siddhartha Mukhopadhyay}: Development of GIS Plugin, Software, Integration of Heuristic Algorithm within GIS Plugin.
\noindent \textbf{Santanu Banerjee}: Conceptualization, Methodology, Investigation, Data curation, Formal analysis, Validation, Software, Visualization, Writing – original draft, Writing – review \& editing.
\textbf{Goutam Sen}: Project administration, Funding acquisition, Supervision, Resources, Writing – review \& editing.
\textbf{Siddhartha Mukhopadhyay}: Software.

\section*{Acknowledgments} \label{Acknowledgments}
This research was conducted through the project entitled “Development of SDSS Tools for Addressing Emergency / Disaster Management”, (Code: IIT/SRIC/R/AEM/2021/101), sponsored through the Kalpana Chawla Space Technology Cell (KCSTC) in IIT Kharagpur, funded by the National Remote Sensing Centre (NRSC) under the Indian Space Research Organization (ISRO).

\section*{Declaration of Competing Interest} \label{Declaration of Competing Interest}
The authors declare that they have no known competing financial interests or personal relationships that could have appeared to influence the work reported in this paper in any way.

% \tiny
\scriptsize

\bibliographystyle{elsarticle-num}
\bibliography{references}
% \bigskip


%% \begin{thebibliography}{00}


%% \bibitem[ ()]{}

%% \end{thebibliography}

% \small

\clearpage
% \vspace*{\fill}  % Push down
\appendix

\footnotesize
\section{List of Abbreviations} \label{List of Abbreviations}
\begin{longtable}{lllll}
TP	&	Transhipment Port	&	&	PRV	&	Primary Resource Vertex	\\
VRP	&	Vehicle Routing Problem	&	&	PGF	&	Preference Generator Function	\\
VD	&	Vehicle Depot	&	&	SRE	&	Smallest Route Element	\\
MILP	&	Mixed Integer Linear Programming	&	&	CD	&	Causality Dict	\\
WH	&	Warehouse	&	&	VLC	&	Vehicle Load Code	\\
$N$	&	Node	&	&	TRS	&	Transhipment Trail Stemming	\\
RC	&	Relief Centre	&	&	RoPr	&	Route Portion	\\
DCT	&	Delivery Cargo Type	&	&	RoPr-Combos	&	Route Portion Combinations	\\
PCT	&	PickUp Cargo Type	&	&	RoCu	&	Route Cluster	\\
VT	&	Vehicle Type	&	&	S.D.	&	Standard Deviation	\\
CT	&	Cargo Type	&	&	CPU	&	Central Processing Unit	\\
SMTS	&	Single-Mode Transportation Segment	&	&	RAM	&	Random Access Memory	\\
MMTN	&	Multi-Modal Transportation Network	&	&	OR	&	Operations Research	\\
PSR	&	Preferential Selection of Routes	&	&	CVRP	&	Capacitated Vehicle Routing Problem	\\
DT	&	Decision Tree	&	&	NFSF	&	No Feasible Solution Found	\\
GIP	&	Generation, Integration and Perturbation	&	&	OEAM	&	Optimization Exhausted Available Memory	\\
DTS	&	Decision Tree Structure	&	&	\wrt	&	with respect to	\\
\end{longtable}

\vspace*{\fill}  % Push down

\small
\section{Visualizing the Transhipment Degree} \label{Visualizing the Transhipment Degree}
In our DT-based PSR-GIP Heuristic, transhipments are generally allowed when the combined resources within a potential/degree is insufficient to satisfy the entire requirement.

\autoref{fig:Degree Depiction in Supplementary} depicts stages/degrees of transhipment (which we may imagine as potentials); the resources available directly in the same SMTS (\ie degree 0) of a requesting Node try to satisfy the demands of that Node as a priority, when the resource availability is insufficient (insufficieny may be due to non availability of Warehouses or Relief Centres in that SMTS, lack of resources available at PRVs which may be unable to satify all the demands of all Nodes in that SMTS, or the requirement of a special resource available only at a distant Warehouse which is not in that SMTS) the next higher degrees are searched for the necessary resources to satisfy the Node's requirements. The degrees in this regard may be compared with energy potentials, and a higher energy potential os accessed only when there is dearth of resources in the PRVs available in lesser transhipment degrees.

The \autoref{fig:Degree Depiction in Supplementary} highlights an example scenario of this problem; when a non-zero transhipment code is generated indicating that access to a higher degree is necessary through a TP, the concerned Trunk-TP is referred as the DTS for the next stage of resource availability check in the higher degree(s). Notice that the $W1$ accessible via $VT8$ also shows the cargo-type $2D$ indicating that the $VT9$ is incompatible with this DCT $2D$ (and therefore this CT is not a leaf from the $W1$ twig from the branch $VT9$).
\begin{landscape}
\begin{figure}
    \centering
    \includegraphics[width=1.035\linewidth]{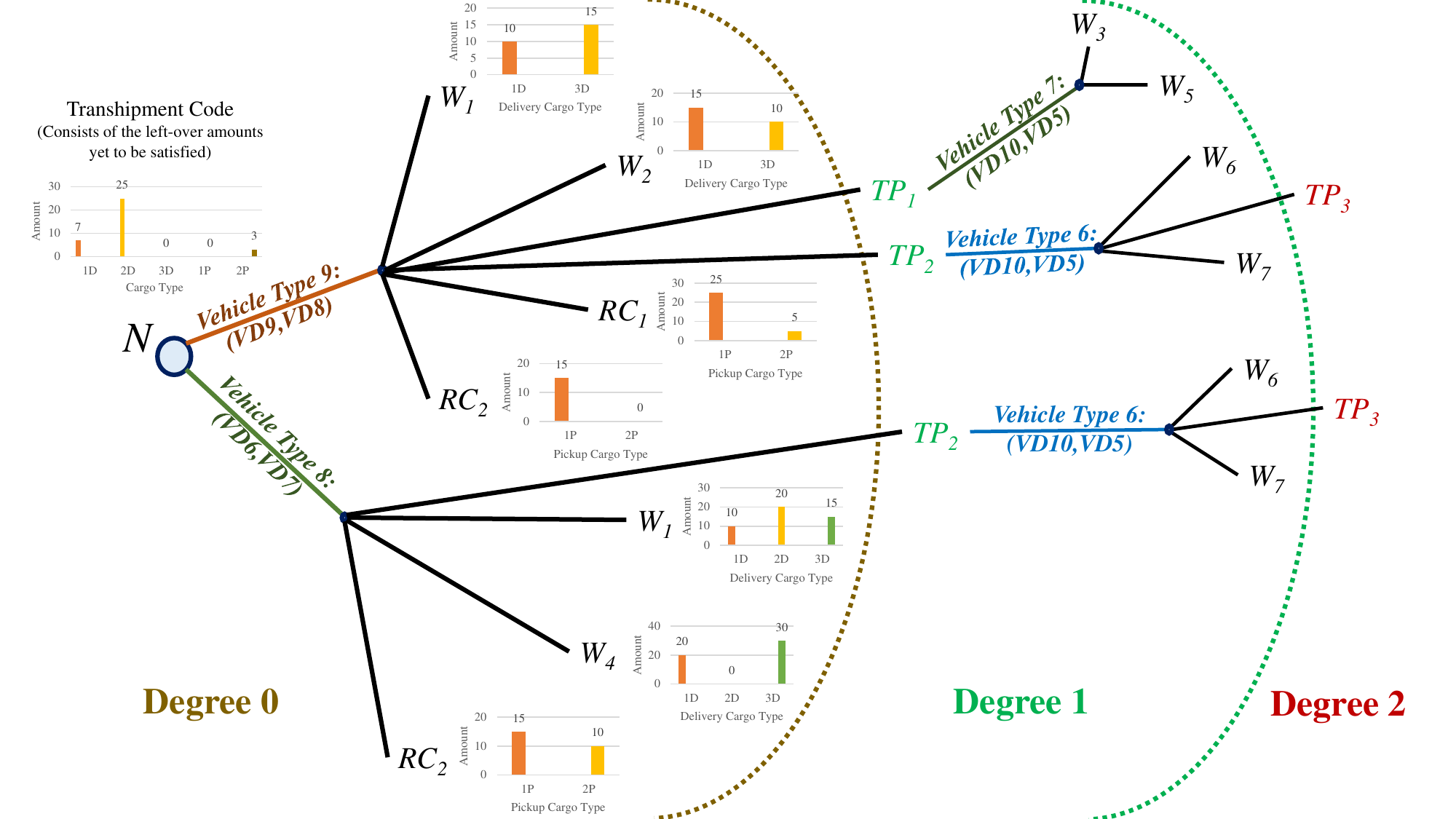}
    \caption{\textbf{A hypothetical problem-scenario illustrating the increasing transhipment degrees}}
    \label{fig:Degree Depiction in Supplementary}
\end{figure}
\end{landscape}

% \appendix

\section{Additional Constraints to the MILP} \label{Extra Optional Constraints to the MILP}

\subsection{Extra Optional Constraints} \label{Extra Optional Constraints}

Eq. \ref{eq:5} allows a new level to be used only if the level below cannot be used. The usage of these constraints have shown to reduce the computation time of the formulation.

\begin{subequations} \label{eq:5}
    \begin{multline} \label{eq:5a}
    \sum_{\substack{z \in V_k,h\\
                    i \neq z \\
                    if~i \in W \Rightarrow z \neq h}}
    x_{i,z}^{v,l} +
    \sum_{\substack{z \in V_k,h \\
                    z \neq j \\
                    if~j \in R \Rightarrow z \neq h}}
    x_{z,j}^{v,l}
    \geq
    x^{v,l+1}_{i,j}, \quad
    \forall h \in H,
    \forall k \in K_h,
    \forall u \in G_{h,k},
    v=(h,k,u),
    if~l_v^{max}>1,
    l=1,
    \forall i,j \in V_k,\\
    i \neq j,
    \end{multline}
    
    \begin{multline} \label{eq:5b}
    x^{v,l,l-1}_i +
    x^{v,l-1,l}_j +
    \sum_{\substack{z \in V_k,h \\
                    i \neq z \\
                    if~i \in W \Rightarrow z \neq h}}
    x_{i,z}^{v,l} +
    \sum_{\substack{z \in V_k \\
                    z \neq j}}
    x_{z,j}^{v,l}
    \geq
    x^{v,l+1}_{i,j}, \quad
    \forall h \in H,
    \forall k \in K_h,
    \forall u \in G_{h,k},
    v=(h,k,u),
    if~l_v^{max}>1,\\
    \forall l \in L_v \smallsetminus \{{1,l_v^{max}}\},
    \forall i,j \in V_k,
    i \neq j,
    \end{multline}
\end{subequations}

Eq. \ref{eq:6} reduces extra level usages for the final connections back to the depots.
\begin{subequations} \label{eq:6}
    \begin{equation} \label{eq:6a}
    \sum_{\substack{z \in V_k,h \\
                    i \neq z}}
    x_{i,z}^{v,l}
    \geq
    x^{v,l+1}_{i,h}, \quad
    \forall h \in H,
    \forall k \in K_h,
    \forall u \in G_{h,k},
    v=(h,k,u),
    if~l_v^{max}>1,
    l=1,
    \forall i \in N_k \cup S_k \cup R_k,
    \end{equation}
    
    \begin{multline} \label{eq:6b}
    x^{v,l,l-1}_i +
    \sum_{\substack{z \in V_k,h \\
                    i \neq z}}
    x_{i,z}^{v,l} \geq
    x^{v,l+1}_{i,h}, \quad
    \forall h \in H,
    \forall k \in K_h,
    \forall u \in G_{h,k},
    v=(h,k,u),
    if~l_v^{max}>1,
    \forall l \in L_v \smallsetminus \{{1,l_v^{max}}\},\\
    \forall i \in N_k \cup S_k \cup R_k,
    \end{multline}
\end{subequations}

% ////////////////////Commented Equations \newline
Eq. \ref{eq:6.25} prevents the inter-level connections connecting the same two vertices between two consecutive levels (in opposite directions), to not be used together.
    \begin{equation} \label{eq:6.25}
    x^{v,l-1,l}_i + x^{v,l,l-1}_i \leq 1, \quad
    \forall h \in H,
    \forall k \in K_h,
    \forall u \in G_{h,k},
    v=(h,k,u),
    %if~l_v^{max}>1,
    \forall l \in L_v \smallsetminus \{{1}\},
    \forall i \in V_k,
    \end{equation}

Eq. \ref{eq:6.75} ensures the sum of incoming edges into any vertex to be at most $1$ only if that vehicle is used.
\begin{subequations} \label{eq:6.75}
    \begin{equation} \label{eq:6.75a}
    x^{v,l+1,l}_i +
    \sum_{\substack{j \in V_k,h \\ i \neq j \\ if~i \in R \Rightarrow j\neq h}}
    x_{j,i}^{v,l}
    \leq
    \sum_{j \in W_k,N_k,S_k}
    x_{h,j}^{v,1}, \quad
    \forall h \in H,
    \forall k \in K_h,
    \forall u \in G_{h,k},
    v=(h,k,u),
    l=1,
    \forall i \in V_k,
    \end{equation}
    
    \begin{multline} \label{eq:6.75b}
    x^{v,l-1,l}_i +
    x^{v,l+1,l}_i +
    \sum_{\substack{j \in V_k \\ i \neq j}}
    x_{j,i}^{v,l}
    \leq
    \sum_{j \in W_k,N_k,S_k}
    x_{h,j}^{v,1}, \quad
    \forall h \in H,
    \forall k \in K_h,
    \forall u \in G_{h,k},
    v=(h,k,u),
    \forall i \in V_k,\\
    \forall l \in L_v \smallsetminus{\{1, l_v^{max}\}},
    \end{multline}
    
    \begin{equation} \label{eq:6.75c}
    x^{v,l-1,l}_i +
    \sum_{\substack{j \in V_k \\ i \neq j}}
    x_{j,i}^{v,l}
    \leq
    \sum_{j \in W_k,N_k,S_k}
    x_{h,j}^{v,1}, \quad
    \forall h \in H,
    \forall k \in K_h,
    \forall u \in G_{h,k},
    v=(h,k,u),
    \forall i \in V_k,
    if~l_v^{max}>1,
    l=l_v^{max},
    \end{equation}
\end{subequations}
% ////////////////////Commented Equations \newline \newline \newline \newline \newline

\subsection{Additional Capacity Constraints}  \label{Additional Capacity Constraints}

The optional Eq. \ref{oeq:63} is only required for specific problem cases, when some Cargo Type needs to have volume and weight as null; which can happen when a dummy load to created to visit a vertex for some activity. An example of such a task would be the clearance of fallen trees which could be modelled as a cargo-associated-task, and the location where this activity needs to take place would be modelled as a Node location, since this is an activity/task there would be no weight/volume utilization of the concerned vehicle; however the loading/unloading time at this node would represent the tree clearance time. In case such a task is modelled as a PCT, then a dummy RC could be created at/near the same dummy Node location to dispose off the fallen trees; otherwise if the task if modelled as a DCT, then a dummy $WH$ would indicate the point of the fallen tree location (basically in the middle of a road) and the disposal location of the cut-and-removed tree trunk would be represented as a dummy Node (very near to the dummy $WH$).

    \begin{subequations} \label{oeq:63}
    \begin{multline} \label{oeq:63de}
        y^{v,l,c}_{i,j}
        \leq M \cdot x_{i,j}^{v,l},
        \quad
       \forall h \in H, \forall k \in K_h, \forall u \in G_{h,k}, v=(h,k,u),
       \forall l \in L_v,
       \forall i,j \in V_k \cup h,
       (if~l = 1 \Rightarrow i \neq h),
       i \neq j,\\
       (if~i \in W \Rightarrow j \neq h),
       \forall c \in C_k,
    \end{multline}
    \begin{equation} \label{oeq:63f}
        y^{v,l,m,c}_i
        \leq M \cdot x^{v,l,m}_i,
        \quad
       \forall h \in H, \forall k \in K_h, \forall u \in G_{h,k}, v=(h,k,u),
       \forall l,m \in L_v,|l-m|=1,
       \forall i \in V_k,
       \forall c \in C_k,
    \end{equation}
    \end{subequations}

\end{document}

\endinput

\begin{frontmatter}



\title{\small SUPPLEMENTARY MATERIALS for the MANUSCRIPT TITLED: ``Rich Vehicle Routing Problem in Disaster Management enabling Temporally-causal Transhipments across Multi-Modal Transportation Networks''}


\author[1]{\href{https://www.linkedin.com/in/santanu-banerjee-093929150/}{Santanu Banerjee}~\orcidlink{0000-0001-9861-7030}\corref{cor1}}\ead{santanu@kgpian.iitkgp.ac.in}
\author[1]{Goutam Sen}\ead{gsen@iem.iitkgp.ac.in}
\author[1]{Siddhartha Mukhopadhyay~\orcidlink{0009-0005-5764-7299}}\ead{sid.mpadhyay@kgpian.iitkgp.ac.in}

\cortext[cor1]{Corresponding author}

\affiliation[1]{organization={Department of Industrial and Systems Engineering (ISE), Indian Institute of Technology (IIT) Kharagpur},%Department and Organization
            %addressline={}, 
            city={Kharagpur},
            postcode={721302}, 
            state={West Bengal},
            country={India}}

\begin{abstract}
This supplementary material provides additional information and results to complement the main manuscript. It includes detailed explanation of how some important constraints in the formulation were developed, extended Heuristic Description with implementation details of the algorithm, along with the computational results of the continuous instances.
\end{abstract}

% \tableofcontents
\end{frontmatter}
\tableofcontents
\clearpage

\begin{landscape}
\section{Literature Review of Individual Papers} \label{Literature Review of Individual Papers}

\footnotesize

\vspace*{\fill}  % Push down

\begin{longtable}{|m{0.06\linewidth}|m{0.055\linewidth}|m{0.25\linewidth}|m{0.17\linewidth}|m{0.22\linewidth}|m{0.075\linewidth}|m{0.07\linewidth}|}
% Paper Details along with spatiotemporal details and Objectives of Papers reviewed
\caption{A. \textbf{Discussed Papers and Solution Methodologies used}\label{longTab: LitRev-1}}\\
% \hline
% \multicolumn{2}{| c |}{Begin of Table}\\
\hline
\multirow{2}{=}{\textbf{Reference Paper}} & \multirow{2}{=}{\textbf{Keyword used for search}} & \multicolumn{3}{|c|}{\textbf{Problem Details}}& \multicolumn{2}{|m{0.15\linewidth}|}{\textbf{InterConnections and IntraDependencies}}\\
\cline{3-7}
& & \textbf{Description, Problem Class, and Optimization Approach (of Objective)} & \textbf{Types of Formulations provided and/or discussed} & \textbf{Heuristic used and/or Novelty} & \textbf{Load (PCT or DCT) to Vehicle Type Compatibility} & \textbf{Road to Vehicle Type Compatibility}\\
\hline\hline
\endfirsthead

\hline
\multicolumn{7}{|c|}{Continuation of \autoref{longTab: LitRev-1}. A.}\\
\hline
\multirow{2}{=}{\textbf{Reference Paper}} & \multirow{2}{=}{\textbf{Keyword used for search}} & \multicolumn{3}{|c|}{\textbf{Problem Details}}& \multicolumn{2}{|m{0.15\linewidth}|}{\textbf{InterConnections and IntraDependencies}}\\
\cline{3-7}
& & \textbf{Description, Problem Class, and Optimization Approach (of Objective)} & \textbf{Types of Formulations provided and/or discussed} & \textbf{Heuristic used and/or Novelty} & \textbf{Load (PCT or DCT) to Vehicle Type Compatibility} & \textbf{Road to Vehicle Type Compatibility}\\
\hline
\endhead

% \hline
% \endfoot

% \hline
% \multicolumn{2}{| c |}{End of Table}\\
% \hline\hline
% \endlastfoot

This study \cite{my_MS_Thesis}	&	-	&	(Open or Closed) VRP cascaded Makespan Minimization	&	MILP	&	Time based transhipment at Transhipment Ports; Cascaded Time Minimization	&	Yes	&	\vspace*{0.75cm} Yes \vspace*{0.75cm}	\\	\hline
\hline
\cite{doi:10.1287/trsc.37.2.153.15243}	&	Capaci-tated VRP	&	Minimize total distance travelled by all vehicles	&	MILP	&	Branch-and-cut algorithm with new Cutting Planes	&	No	&	No	\\
\hline
\hline
\cite{Koch2018}	&	\multirow{6}{=}{\rotatebox{90}{Simultaneous Delivery and Pickup}}	&	3L-VRPSDPTW, Minimizes the Total travel Distance	&	-	&	Hybrid algorithm of Adaptive Large Neighbourhood Search is used. Authors consider loading from the back side, as well as loading from the long side (in case of 2 compartments) which have not been studied in Literature before	&	No	&	No	\\	\cline{1-1} \cline{3-7}
\cite{MIN1989377}	&		&	VRPSDP, Minimizes the total travel time of the route	&	MILP	&	First paper on multiple vehile routing problem considering simultaneous delivery and pickup; uses three-phase sequential procedure, analogous to “cluster-first route-second”	&	No	&	No	\\	\cline{1-1} \cline{3-7}
\cite{OZTAS2022117401}	&		&	VRPSDP	&	-	&	Hybrid Meta-Heuristic Algorithm (ILS-RVND-TA)	&	No	&	No	\\	\cline{1-1} \cline{3-7}
\cite{Rieck2013}	&		&	Symmetric and Assymetric VRPSDP, Minimize transportation cost of the distances travelled	&	 Two MILP Formulations, one as a Vehicle Flow and another as a Commodity Flow, are proposed.	&	Authors expect future work to include multiple objective functions. Pre-processing techniques for reducing domain, and effective cutting planes are outlined.	&	No	&	No	\\	\cline{1-1} \cline{3-7}
\cite{S0217595905000522}	&		&	Minimize Total Route Cost	&	MILP requires minimum number of trucks needed as an input; unrestricted return quantity is allowed; minimizes total cost of all routes	&	Route construction Heuristic	&	No	&	No	\\	\cline{1-1} \cline{3-7}
\cite{10406762}	&		&	HFVRPSDP, Minimization of the sum of all heterogeneous vehicle route used for minimizing fixed plus variable cost of routing	&	A more compact MILP (than existing previously in the literature) is provided highlighting improvement in solution time	&	A novel hybrid heuristic is designed ehanced with function-flavours for MachineTuning, and a hyerparameter of JellyFishing	&	No	&	Yes	\\	\hline
\hline
\cite{ALLAHYARI2015756}	&	\multirow{5}{=}{\rotatebox{90}{Multiple Depots}}	&	Multi-Depot Covering Tour Vehicle Routing Problem; minimizes total routing cost and allocation cost of customers to nearby nodes	&	Two MILP formulations developes: Flow and Node based	&	Hybrid solution developed using GRASP; step-wise clustering and routing; iterated local search and simulated annealing	&	No	&	No	\\	\cline{1-1} \cline{3-7}
\cite{NUCAMENDIGUILLEN2021113846}	&		&	Minimization of Fixed plus Variable cost of all vehicle routes, Multi-Depot Open Location Routing Problem	&	MILP	&	Metaheuristic developed performs supplier clustering during construction phase, and has an additional improvement phase	&	No	&	No	\\	\cline{1-1} \cline{3-7}
\cite{YoshiakiSHIMIZU20162016jamdsm0004}	&		&	Handling cost at depot, Routing transportation cost and Fixed charges of vehicles and Opening depots	&	Combinatorial optimization problem	&	Modified tabu search	&	No	&	No	\\	\cline{1-1} \cline{3-7}
\cite{Barma_Dutta_Mukherjee_2019}	&		&	Minimize total routing cost or distance of an Multiple Depot Vehicle Routing Problem	&	MILP	&	Bio-inspired   meta-heuristic  method  named  Discrete  Antlion  Optimization  algorithm  (DALO) followed  by  the  2-opt  algorithm  for  local  searching  is  used  to  	&	No	&	No	\\	\cline{1-1} \cline{3-7}
\cite{VENKATANARASIMHA201363}	&		&	Minimize the tour-length of the vehicle traveling the longest distance	&	-	&	Ant colony optimization; decomposition into different Single Depot VRPs	&	No	&	No	\\
\hline
\hline
\cite{KRAMER2019162}	&	\multirow{3}{=}{\rotatebox{90}{rich-VRP}}	&	Periodic Demands; Functional Nodes are considered here; Minimize vehicle type and traveling distance dependant cost as well as the cost to supply the auxiliary depots	&	No Exact formulation	&	Multi Start Iterated Local Search	&	No	&	Yes (by Customer-Node compatibility)	\\	\cline{1-1} \cline{3-7}
\cite{Rabbouch2021}	&		&	Distribution problem with time constraints including ready and due dates for customers already assigned to depots; Minimizes the total distance	&	MILP	&	Genetic Algorithms	&	No	&	No	\\	\cline{1-1} \cline{3-7}
\cite{RIECK2014863}	&		&	Location-routing problem (involves locating hub facilities as a mandatory transfer point between supply and demand nodes); Objective function selects hub locations and vehicle routes, minimizing overall operational and transportation cost	&	MILP	&	Heuristic uses multi-start procedure based on fix+optimize, and a genetic algorithm	&	No	&	No	\\	\hline
\end{longtable}

\vspace*{\fill}  % Push down

\clearpage
\small

\vspace*{\fill}  % Push down

\addtocounter{table}{-1}

\begin{longtable}{|m{0.08\linewidth}|m{0.09\linewidth}|m{0.0875\linewidth}|m{0.087\linewidth}|m{0.0925\linewidth}|m{0.0825\linewidth}|m{0.11\linewidth}|m{0.12\linewidth}|m{0.1\linewidth}|}

\caption{B. \textbf{Network Features in discussed literature}\label{longTab: LitRev-2}}\\
% \hline
% \multicolumn{2}{| c |}{Begin of Table}\\
\hline
\multirow{2}{=}{\textbf{Ref. Paper}} & \multicolumn{8}{|c|}{\textbf{Network Features Considered}}\\
\cline{2-9}
& \textbf{Closed VRP / Open VRP} & \textbf{No. of Trips allowed for each vehicle} & \textbf{Time Windows considered} & \textbf{Maximum Route Duration constrained?} & \textbf{Waiting Time consideration} & \textbf{Loading / Unloading / Service Time considered} & \textbf{Max. No. of Node Visits per route (provided as constraint)} & \textbf{Multi-Modal Transportation Network} \\
\hline
\endfirsthead

\hline
\multicolumn{9}{|c|}{Continuation of \autoref{longTab: LitRev-2}. B.}\\
\hline
\multirow{2}{=}{\textbf{Ref. Paper}} & \multicolumn{8}{|c|}{\textbf{Network Features Considered}}\\
\cline{2-9}
& \textbf{Closed VRP / Open VRP} & \textbf{No. of Trips allowed for each vehicle} & \textbf{Time Windows considered} & \textbf{Maximum Route Duration constrained?} & \textbf{Waiting Time consideration} & \textbf{Loading / Unloading / Service Time considered} & \textbf{Max. No. of Node Visits per route (provided as constraint)} & \textbf{Multi-Modal Transportation Network} \\
\hline
\endhead

% \hline
% \endfoot

% \hline
% \multicolumn{2}{| c |}{End of Table}\\
% \hline\hline
% \endlastfoot
\hline
Our study	&	Both	&	Multiple	&	No (but can be very easily implemented within the current formulation)	&	No (can be easily extended)	&	Yes	&	Yes	&	Unconstrained, we consider this answered as time minimization would be better than customer connstraining	&	Yes	\\	\hline
\cite{doi:10.1287/trsc.37.2.153.15243}	&	Closed	&	1	&	No	&	No	&	No	&	No	&	No	&	No	\\	\hline
\cite{Koch2018}	&	Closed	&	1	&	Yes	&	No	&	Yes	&	Yes	&	No	&	No	\\	\hline
\cite{MIN1989377}	&	Closed	&	1	&	No	&	No	&	No	&	No	&	No	&	No	\\	\hline
\cite{OZTAS2022117401}	&	Closed	&	1	&	No	&	No	&	No	&	No	&	No	&	No	\\	\hline
\cite{Rieck2013}	&	Closed	&	1	&	No	&	No	&	No	&	No	&	No	&	No	\\	\hline
\cite{S0217595905000522}	&	Closed	&	Multiple	&	No	&	No	&	No	&	No	&	No	&	No	\\	\hline
\cite{10406762}	&	Closed	&	1	&	No	&	No	&	No	&	No	&	No	&	No	\\	\hline
\cite{ALLAHYARI2015756}	&	Closed	&	1	&	No	&	No	&	No	&	No	&	No	&	No	\\	\hline
\cite{NUCAMENDIGUILLEN2021113846}	&	Open	&	1	&	No	&	No	&	No	&	No	&	No	&	No	\\	\hline
\cite{YoshiakiSHIMIZU20162016jamdsm0004}	&	Closed	&	1	&	No	&	No	&	No	&	No	&	No	&	No	\\	\hline
\cite{Barma_Dutta_Mukherjee_2019}	&	Closed	&	1	&	No	&	No	&	No	&	No	&	Yes	&	No	\\	\hline
\cite{VENKATANARASIMHA201363}	&	Closed	&	1	&	No	&	Yes	&	No	&	No	&	No	&	No	\\	\hline
\cite{KRAMER2019162}	&	Closed	&	Multiple (across periods)	&	Yes (Flexible)	&	Yes	&	Yes	&	Yes	&	Yes	&	No	\\	\hline
\cite{Rabbouch2021}	&	Closed	&	1	&	Yes (Hard)	&	Yes	&	No	&	Yes	&	No	&	No	\\	\hline
\cite{RIECK2014863}	&	Closed	&	$\leq$ 3	&	No	&	No	&	No	&	No	&	No	&	No	\\	\hline
\end{longtable}
% \begin{itemize}
%     \item[$^*$] Due to Disaster, Congestion/ Traffic or damaged road
%     \item[$^{**}$] Influence of altered road condition, and/ or speed considerations need to be considered within our initial static Travel-Time Matrix.
% \end{itemize}

\vspace*{\fill}  % Push down

\clearpage
\small

\vspace*{\fill}  % Push down

\addtocounter{table}{-1}

\begin{longtable}{|m{0.08\linewidth}|m{0.11\linewidth}|m{0.09\linewidth}|m{0.09\linewidth}|m{0.165\linewidth}|m{0.09\linewidth}|m{0.11\linewidth}|m{0.15\linewidth}|}

\caption{C. \textbf{Vehicle and Cargo features in the discussed literature}\label{longTab: LitRev-3}}\\
% \hline
% \multicolumn{2}{| c |}{Begin of Table}\\
\hline
\multirow{2}{=}{\textbf{Ref. Paper}} & \multicolumn{4}{|c|}{\textbf{Vehicle Features Considered}} & \multicolumn{3}{|c|}{\textbf{Cargo Features Considered}}\\
\cline{2-8}
& \textbf{Homogeneous / Heterogeneous Fleet} & \textbf{Weight restricted Vehicle} & \textbf{Volume restricted Vehicle} & \textbf{Common Vehicle Capacity for all types of Load commodity (\ie no compartmentalization)} & \textbf{No. of PickUp Types considered (N/A if not allowed)} & \textbf{No. of Delivery Types considered (N/A if not allowed)} & \textbf{Split / Simultaneous resource transfer (Delivery or PickUp)}\\
\hline
\endfirsthead

\hline
\multicolumn{8}{|c|}{Continuation of \autoref{longTab: LitRev-3}. C.}\\
\hline
\multirow{2}{=}{\textbf{Ref. Paper}} & \multicolumn{4}{|c|}{\textbf{Vehicle Features Considered}} & \multicolumn{3}{|c|}{\textbf{Cargo Features Considered}}\\
\cline{2-8}
& \textbf{Homogeneous / Heterogeneous Fleet} & \textbf{Weight restricted Vehicle} & \textbf{Volume restricted vehicle} & \textbf{Common Vehicle Capacity for all types of Load commodity (\ie no compartmentalization)} & \textbf{No. of PickUp Types considered (N/A if not allowed)} & \textbf{No. of Delivery Types considered (N/A if not allowed)} & \textbf{Split / Simultaneous resource transfer (Delivery or PickUp)}\\
\hline
\endhead

% \hline
% \endfoot

% \hline
% \multicolumn{2}{| c |}{End of Table}\\
% \hline\hline
% \endlastfoot
\hline
Our study	&	Heterogeneous	&	Yes	&	Yes	&	Yes	&	$\geq$ 1	&	$\geq$ 1	&	Simultaneous and/or Split	\\	\hline
\cite{doi:10.1287/trsc.37.2.153.15243}	&	Homogeneous	&	\multicolumn{2}{|c|}{Single Capacity considered}			&	N/A	&	N/A	&	1	&	None	\\	\hline
\cite{Koch2018}	&	Homogeneous	&	Yes	&	Yes (3D packing considered)	&	Both considered	&	\multicolumn{2}{|m{0.19\linewidth}|}{Each item is unique in identity and dimensions and multiple such items in the backhaul can be returned by the customer}			&	Simultaneous Delivery and Pickup	\\	\hline
\cite{MIN1989377}	&	Homogeneous	&	-	&	Yes	&	Yes	&	1	&	1	&	Simultaneous Delivery and Pickup	\\	\hline
\cite{OZTAS2022117401}	&	Homogeneous	&	\multicolumn{2}{|c|}{Single Capacity considered}			&	Yes	&	1	&	1	&	Simultaneous Delivery and Pickup	\\	\hline
\cite{Rieck2013}	&	Homogeneous	&	\multicolumn{2}{|c|}{Single Capacity considered}			&	Yes	&	1	&	1	&	Simultaneous	\\	\hline
\cite{S0217595905000522}	&	Homogeneous	&	\multicolumn{2}{|c|}{Single Capacity considered}			&	Yes	&	1	&	1	&	Simulteneous Delivery and Pickup	\\	\hline
\cite{10406762}	&	Heterogeneous	&	\multicolumn{2}{|c|}{Single Capacity considered}			&	Yes	&	1	&	1	&	Simulteneous Delivery and Pickup	\\	\hline
\cite{ALLAHYARI2015756}	&	Homogeneous	&	\multicolumn{2}{|c|}{Single Capacity considered}			&	Yes	&	N/A	&	1	&	Single Instance Delivery	\\	\hline
\cite{NUCAMENDIGUILLEN2021113846}	&	Heterogeneous	&	\multicolumn{2}{|c|}{Single Capacity considered}			&	Yes	&	N/A	&	1	&	Split Delivery	\\	\hline
\cite{YoshiakiSHIMIZU20162016jamdsm0004}	&	Any (flexible model to consider heterogeneous fleet)	&	\multicolumn{2}{|c|}{Single Capacity considered}			&	Yes	&	1	&	1	&	Simultaneous	\\	\hline
\cite{Barma_Dutta_Mukherjee_2019}	&	Homogeneous	&	\multicolumn{2}{|c|}{Single Capacity considered}			&	Yes	&	N/A	&	1	&	Single Instance Delivery	\\	\hline
\cite{VENKATANARASIMHA201363}	&	Homogeneous	&	\multicolumn{2}{|c|}{Single Capacity considered}			&	Yes	&	N/A	&	1	&	Single Instance Delivery	\\	\hline
\cite{KRAMER2019162}	&	Heterogeneous	&	\multicolumn{2}{|c|}{Single Capacity considered}			&	Yes	&	N/A	&	1	&	Split across multiple-periods	\\	\hline
\cite{Rabbouch2021}	&	Heterogeneous	&	\multicolumn{2}{|c|}{Single Capacity considered}			&	Yes	&	N/A	&	1	&	Single Instance Delivery	\\	\hline
\cite{RIECK2014863}	&	Homogeneous	&	\multicolumn{2}{|c|}{Single Capacity considered}			&	Yes	&	$\geq$ 1	&	$\geq$ 1	&	Delivery and/or Pickup	\\	\hline
\end{longtable}

\vspace*{\fill}  % Push down

\clearpage
\small
\vspace*{\fill}  % Push down

\addtocounter{table}{-1}

\begin{longtable}{|m{0.08\linewidth}|m{0.14\linewidth}|m{0.1\linewidth}|m{0.11\linewidth}|m{0.085\linewidth}|m{0.11\linewidth}|m{0.1\linewidth}|m{0.07\linewidth}|m{0.07\linewidth}|}

\caption{D. \textbf{Node Features in discussed literature}\label{longTab: LitRev-4}}\\
% \hline
% \multicolumn{2}{| c |}{Begin of Table}\\
\hline
\multirow{2}{=}{\textbf{Ref. Paper}} & \multicolumn{8}{|c|}{\textbf{Node Features Considered}}\\
\cline{2-9}
& \textbf{Separate Vehicle and Resource parking Locations (N/A if no resource transfer is considered)} & \textbf{Number of Vehicle Depots (or Resource Locations) allowed} & \textbf{Capacity constrained Depot (and/or Limited number of Vehicles)} & \textbf{Capacity constrained Warehouse(s) or equivalent} & \textbf{Relief Points (separate final location to drop all PickUps) considered} & \textbf{Transhipment Ports (or similar features) are available?} & \textbf{Many-to-Many Transport} & \textbf{Multi-visit allowed at Vertices / Nodes}\\
\hline
\endfirsthead

\hline
\multicolumn{9}{|c|}{Continuation of \autoref{longTab: LitRev-4}. D.}\\
\hline
\multirow{2}{=}{\textbf{Ref. Paper}} & \multicolumn{8}{|c|}{\textbf{Node Features Considered}}\\
\cline{2-9}
& \textbf{Separate Vehicle and Resource parking Locations (N/A if no resource transfer is considered)} & \textbf{Number of Vehicle Depots (or Resource Locations) allowed} & \textbf{Capacity constrained Depot (and/or Limited number of Vehicles)} & \textbf{Capacity constrained Warehouse(s) or equivalent} & \textbf{Relief Points (separate final location to drop all PickUps) considered} & \textbf{Transhipment Ports (or similar features) are available?} & \textbf{Many-to-Many Transport} & \textbf{Multi-visit allowed at Vertices / Nodes}\\
\hline
\endhead

% \hline
% \endfoot

% \hline
% \multicolumn{2}{| c |}{End of Table}\\
% \hline\hline
% \endlastfoot
\hline
Our study	&	Yes	&	>0	&	Yes	&	Yes	&	Yes	&	Yes	&	Yes	&	Yes	\\	\hline
\cite{doi:10.1287/trsc.37.2.153.15243}	&	N/A	&	1	&	Yes	&	N/A	&	N/A	&	No	&	No	&	No	\\	\hline
\cite{Koch2018}	&	No	&	1	&	Yes	&	N/A	&	N/A	&	No	&	No	&	No	\\	\hline
\cite{MIN1989377}	&	No	&	-	&	Yes	&	N/A	&	N/A	&	No	&	No	&	No	\\	\hline
\cite{OZTAS2022117401}	&	No	&	1	&	Yes	&	N/A	&	N/A	&	No	&	No	&	No	\\	\hline
\cite{Rieck2013}	&	No	&	1	&	Yes	&	N/A	&	No	&	No	&	No	&	No	\\	\hline
\cite{S0217595905000522}	&	No	&	1	&	Yes	&	N/A	&	N/A	&	No	&	No	&	Yes	\\	\hline
\cite{10406762}	&	No	&	1	&	Yes	&	N/A	&	N/A	&	No	&	No	&	No	\\	\hline
\cite{ALLAHYARI2015756}	&	No	&	Multiple	&	Yes	&	No	&	N/A	&	No	&	Yes	&	No	\\	\hline
\cite{NUCAMENDIGUILLEN2021113846}	&	Yes	&	Multiple	&	Yes	&	Yes	&	N/A	&	No	&	Yes	&	Only For the single Split Delivery Node	\\	\hline
\cite{YoshiakiSHIMIZU20162016jamdsm0004}	&	No	&	Multiple	&	Yes	&	N/A	&	N/A	&	No	&	Yes	&	No	\\	\hline
\cite{Barma_Dutta_Mukherjee_2019}	&	No	&	Multiple	&	Yes	&	N/A	&	N/A	&	No	&	Yes	&	No	\\	\hline
\cite{VENKATANARASIMHA201363}	&	No	&	Multiple	&	Yes	&	No	&	No	&	No	&	Yes	&	No	\\	\hline
\cite{KRAMER2019162}	&	No	&	Multiple	&	Yes	&	Yes	&	No	&	Yes	&	Yes	&	Yes	\\	\hline
\cite{Rabbouch2021}	&	No	&	Multiple	&	Yes	&	Yes	&	No	&	No	&	Yes	&	No	\\	\hline
\cite{RIECK2014863}	&	No	&	Multiple	&	Yes	&	Yes	&	No	&	Yes	&	Yes	&	Yes	\\	\hline
\end{longtable}

\vspace*{\fill}  % Push down

\end{landscape}
\small

\clearpage

\section{Additional comments regarding our Exact Formulation} \label{MILP appendix}

\begin{longtable}{|m{0.4\linewidth}|m{0.6\linewidth}|}
\caption{\textbf{Multi-functionality of Graph Vertices}\label{longTab: Multi-functionality of Graph Vertices}}\\
 \hline
 \multicolumn{2}{|c|}{\textbf{Types of Vertices based on functionality}} \\
 
 \hline
 \textbf{Mathematical View Point} & \textbf{Nomenclature Description}\\
 \hline
 \endfirsthead

 \hline
 \multicolumn{2}{|c|}{Continuation of \autoref{longTab: Multi-functionality of Graph Vertices} \textbf{Multi-functionality of Graph Vertices}}\\
 \hline
 \textbf{Mathematical View Point} & \textbf{Nomenclature Description}\\
 \hline
 \endhead

\hline Source / Origin only & \multirow{2}{=}{Not considered in our Problem or Formulation, applicable only for open VRPs}\\
\cline{1-1} Sink / Destination only & \\

%% \hline Source / Origin only & Not considered in our Problem or Formulation\\
%% \hline Sink / Destination only & Not considered in our Problem or Formulation\\
\hline Both Source \& Sink \ie Vertex acts as both Origin \& Destination & Vehicle Depot ($VD$) [these vertices are also the time-origin for corresponding vehicles]\\
\hline PickUp only & Warehouses ($W$)\\
\hline	Delivery only	&	Relief Centers ($RC$)	\\
\hline	PickUp \& Delivery (Simultaneous)	&	Simultaneous Nodes ($NM$) [for single vehicle relief-rescue]	\\
\hline	PickUp \& Delivery (Split)	&	Split Nodes ($NP$) [allows multiple vehicles to perform pre-defined operations over the entire course of the mission]	\\
\hline	Time-Dependent Cargo-Compatible Echelon Points	&	Transhipment Ports ($TP$) [also acts as Multimodal Junctions]	\\
\hline	Route Split	&	For deploying Water Boats from a Land Vehicle for urban-flood relief [not considered in current formulation or problem]	\\
\hline	Route Capture	&	For accepting boats with rescued victims back into the Mother Vehicle [not considered in our problem or formulation]	\\
\hline	Allows both Route Splitting and Capture	&	Deploy Points [not considered in current formulation and suggested as a possible future extension]	\\
\hline

\end{longtable}

\subsection{Considerations while running the MILP} \label{Optimization MethodologyConsiderations while running the MILP}

The equations \ref{eq:4}, \ref{eq:6.75}, \ref{eq:10}, \ref{eq:11}, \ref{eq:13}, \ref{eq:14}, \ref{eq:16}, \ref{eq:17}, \ref{eq:20}, \ref{eq:21}, \ref{eq:29}, \ref{eq:37}, \ref{eq:38}, \ref{eq:39} get slightly altered while taking the extreme values of $l_v$, since when $l=1$ the incoming variables $x$ or $y$ connecting the bottom-most layer $1$ from any below layer is non-existent; and similarly when $l=l_v^{max}$ outgoing variables from the top-most layer $l_v^{max}$  to any other higher layer is non-existent. Further if $l_v^{max}=1$, then both the above non-existent variables would further reduce such constraints, and these equations respectively reduce to reduced versions of equations \ref{eq:4a}, \ref{eq:6.75a}, \ref{eq:10a}, \ref{eq:11a}, \ref{eq:13a}, \ref{eq:14a}, \ref{eq:16a}, \ref{eq:17a}, \ref{eq:20a}, \ref{eq:21a}, \ref{eq:29a}, \ref{eq:37a}, \ref{eq:38a}, \ref{eq:39a}; the alteration process of one equation (Eq. \ref{eq:4}) is shown below into Eq. \ref{eq:68} (which is a reduced form of Eq. \ref{eq:4a} without any of the inter-layer variables as they are non-existent variables).
\begin{equation} \label{eq:68}
\sum_{\substack{j \in V_k,h \\ i \neq j \\ if~i \in R \Rightarrow j\neq h}}
x_{j,i}^{v,l}
=
\sum_{\substack{j \in V_k,h \\ i \neq j \\ if~i \in W \Rightarrow j \neq h}}
x_{i,j}^{v,l},
\quad
\forall h \in H,
\forall k \in K_h,
\forall u \in G_{h,k},
v=(h,k,u),
l=1,
\forall i \in V_k,
\end{equation}

Increased values of $l_v^{max}$ expand the solution space, and the optimum solution of the problem would depend on these values. Therefore, since it is hard to estimate mathematically the required vehicular levels necessary for each vehicle so that the solution space creation (\ie variable generation for the exact formulation) is apt for finding the optimal solution, we ideate determining upper bounds on (vehicle type specific) chunks of these variables in section \ref{Maximum Vehicle-Level Requirement Estimation}. However, creating the massive solution space based on the upper bound becomes computationally inefficient and we develop an approach to progressively increase these $l_v^{max}$ values (starting all from 1) till no major solution improvement is seen. We may also alternatively assign a good number of levels (\ie enough vehicle layers) to each individual vehicle to solve the problem once (with any among the minimization approaches of Mi-Max, Makespan reduction, or Cascaded Time Minimization) without going through the process of increasing the layers progressively as highlighted in section \ref{Optimization MethodologyConsiderations while running the MILP}.

\subsubsection{Maximum Vehicle-Level Requirement Estimation} \label{Maximum Vehicle-Level Requirement Estimation}
Our proposed exact formulation is dependent on the available levels $l \in L_v$ to each individual vehicle $v$. A way to calculate the possible upper limit of the vehicle-levels would be to solve the formulation separately many times (with a sufficient number of available levels, to maybe find only the first few feasible solutions) while considering only a single type of vehicle available each time (homogeneous fleet across all depots), and collecting the maximum level usage for each of these vehicle types; these computed values of the levels could then be used in the original bigger problem considering all vehicles in the exact formulation. However, this would be time-taking, cumbersome and infeasible problems could arise, hindering the collection process of this rough upper bound of the levels.

Instead, we solve the original problem considering all vehicles having $l_v^{max}$ as $1$ initially. Next, we progress further by  a combination of the logics highlighted below:
\begin{enumerate}
    \item For vehicles which are unused, we may choose any of the following strategies:
        \begin{enumerate}
    \item Not consider them further, \ie remove them from the next stages.
    \item Progressively change the level value as stages progress, by changing the existing level value by some (preferably negative fraction, say $-0.25$) number so that the changes seem gradual. The actually assigned number of vehicle layers in the MILP needs to be rounded off values in this case.
        \end{enumerate}    
    \item For the used vehicles, we may progressively change the level value as the stages progress, by changing the existing level value by some (preferably positive integer, say $1$) number.
\end{enumerate}
In both cases, the changes may be done on the number of levels assigned in the previous stage, or the number of levels used (in the previous stage).

In our case, we consider that every next stage should have one more level for each vehicle than what was used in the previous stage by that same vehicle (so an unused vehicle will have $1$ level to ply on in the subsequent stage). This would progressively decrease the search space by decreasing the number of vehicles, and on the other hand increase the search space by providing additional levels to the used vehicles. A variation of this logic could be using the number of vehicle-specific-levels being one more than what was used by that vehicle in some of the last few stages (by developing a stage-wide memory element).
 
 We consider that through this process, we would be able to progress towards the optimal solution (the optimal solution could have also been found by providing a very large number of levels to each vehicle in the first place, however at the cost of eons of computational time); and consider it over when any additional stage of change in the available vehicles and their levels does not yield any improvement. However, since this approach does not guarantee the exploration of the entire solution space, sub-optimal solutions may be found.

It may be possible that a vehicle, that is being unused at lower values of its available levels, starts getting utilized suddenly from a higher value of the levels available to it (also depending on the availability of other vehicles and their respective available levels). We may try to find even better solutions and refute (or confirm) the above possibility, with all vehicles having the same available levels equal to one more than the maximum of all levels provided to any vehicle in the final stage of the solution above; however, we find this computationally extremely inefficient.

% \subsubsection{Optimization Methodology} \label{Optimization Methodology}
% The objective function proposed focuses on the Min-Max approach, and we may also alter the objective function slightly to minimize the makespan, \ie use equations this instead of as the objective function. In our problem, this Min-Max approach minimizes the maximum route duration

\subsection{Explaining the constraint development process through Space-Time matrices used in the development of our MILP formulation} \label{Explaining the constraint development process for the Space-Time matrix used in the MILP formulation}

To ensure temporal-causality during transhipments is maintained, we develop constraints; one important component during the development of these constraints is the space-time matrix, which is elaborately explained in this section.

\begin{figure}
    \centering
    \includegraphics[width=1\linewidth]{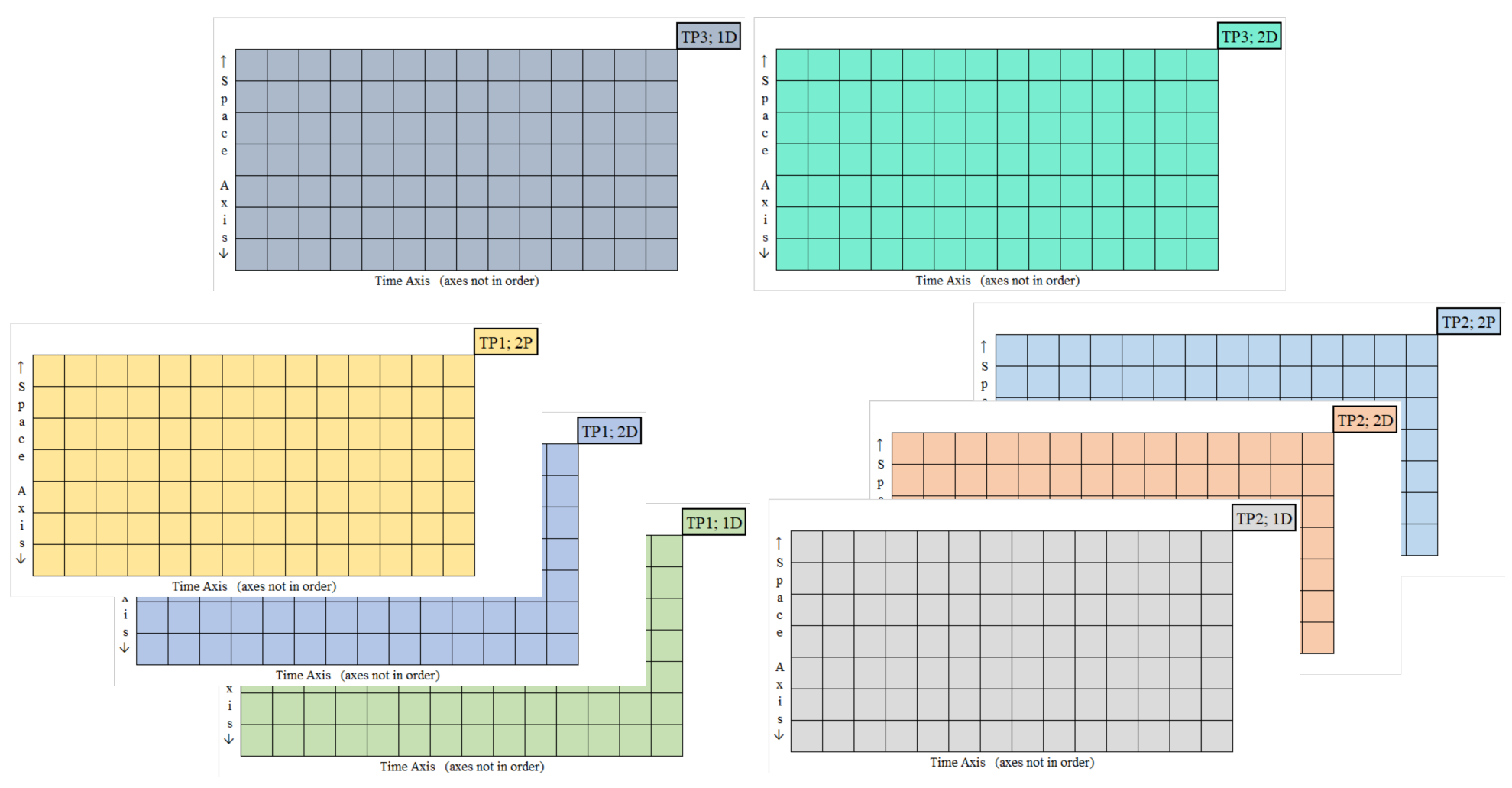} %1.04064
    \caption{\textbf{Space-Time Matrices for our study-example instance. The axes are not in any order, indicating that the Time value may not be steadily increasing rightwards, unlike as is generally the case.}}
    \label{fig:Explanation_Space_Time}
\end{figure}

For each TP, we develop individual space-time matrices for each commodity allowed to be transhipped through that TP. For our study-example instance:
\begin{enumerate}
    \item For TP1 and TP2: there will be three separate space-time matrices, each for the commodities 1D, 2D, and 2P
    \item For TP3: there will be two separate space-time matrices, each for the DCTs 1D and 2D (see \autoref{fig:Picture_Instance} for the compatibilities, and the \autoref{fig:Explanation_Space_Time} to understand the Space-Time Matrices for our Study-Example instance)
\end{enumerate}

\begin{figure}
    \centering
    \includegraphics[width=1\linewidth]{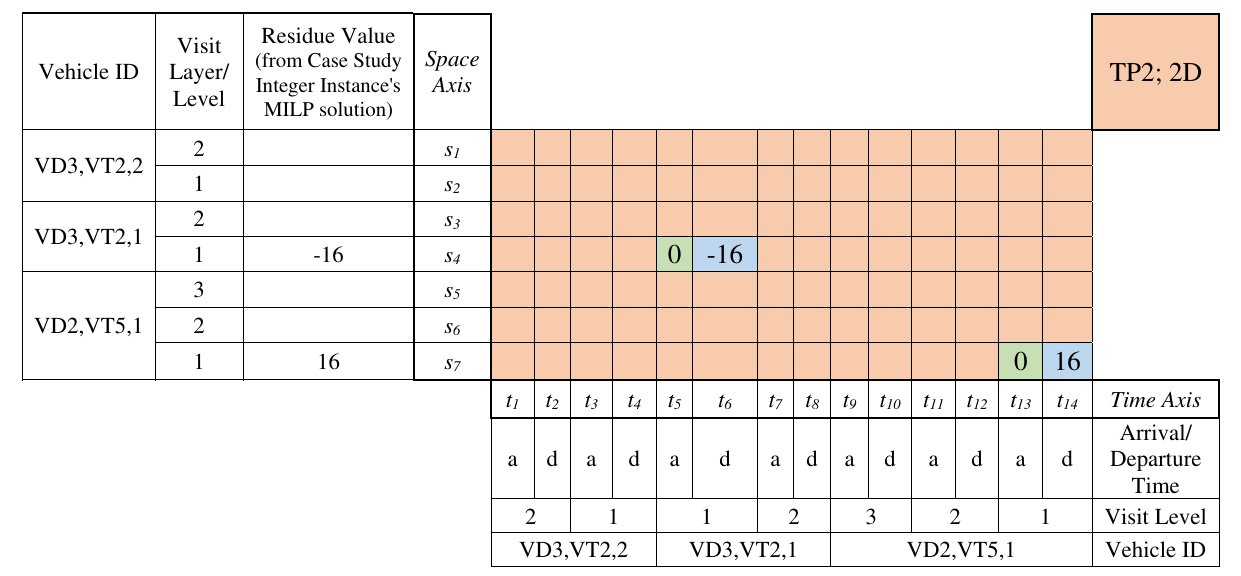}
    \caption{\textbf{One space-time matrix is elaborated}}
    \label{fig:ST Matrix Explanation 1}
\end{figure}
\begin{figure}
    \centering
    \includegraphics[width=1\linewidth]{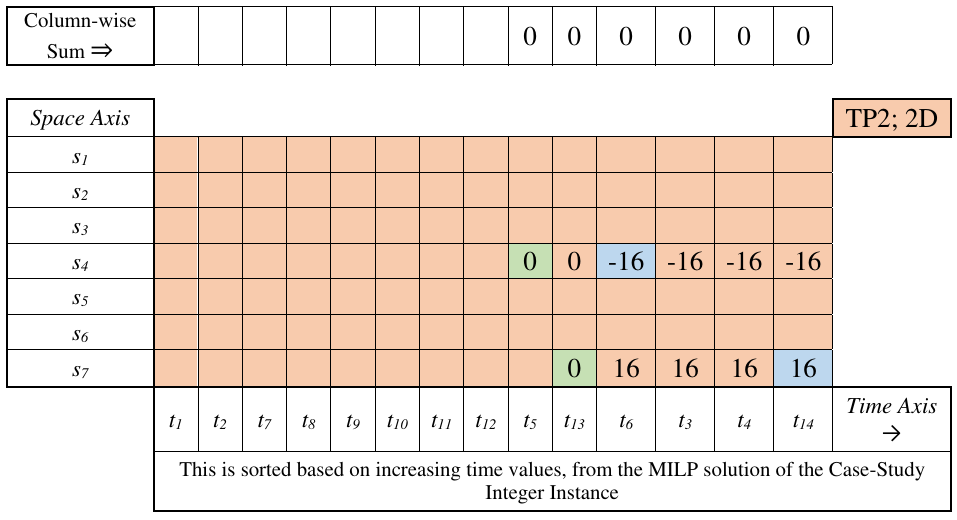}
    \caption{\textbf{The necessary values in the ST Matrix are populated (as per our study-example Integer solution from MILP); all the other unpopulated cells will have null values.}}
    \label{fig:ST Matrix Explanation 2}
\end{figure}

Each such matrix has its Y-axis (the vertical Space Axis) indicating the various vehicles and their levels, which may be used to visit the concerned TP. The X-axis (the horizontal Time Axis) indicates the possible time instances associated with vehicle visits. In case each visit is associated with a single time, then the number of rows and columns in the matrix would have been equal; however in our case, we have two time-stamps associated with any vehicle's visit at any vertex, namely the Arrival Time and the Departure Time, therefore for each row heading, we will have two associated columns (one for the arrival and one for the departure timings), thus these matrices have double the number of columns \wrt rows.
Each cell of the matrix contains the value of the quantity of the resource that was transhipped at that row-level till the time indicated in the column-bottom times (\autoref{fig:ST Matrix Explanation 1}); if there has been no transhipment in that row, then the values will be 0. Each cell is essentially the variable $o$ of the MILP.

In \autoref{fig:ST Matrix Explanation 1}, the residue values correspond to the variable $r$ in the MILP. Two residue values are shown here which indicate the transfer of this cargo-type 2D from the vehicle (VD2,VT5,1) to the TP2 (positive value), and from the TP2 to the vehicle (VD3,VT2,1) which is indicated with the opposite sign. We see that only the first levels are used in this solution. The residue value is populated along the entire corresponding row only when the times are higher than the departure time. Similarly for all times less than the arrival time, the cells in a row will contain null value; this can be easily understood by the fact that if there has been an amount of the CT transfer at this TP, the transfer details should be maintained subsequently as well allowing other vehicles to take up the CT and further tranship it. The arrival time of the vehicle for which there is a non-zero residue is depicted as green cells in \autoref{fig:ST Matrix Explanation 1}, with the departure time-corresponding cells depicted in blue.
The intermediate variable $g$ in the MILP allows the comparison between times of any two cells in this matrix, to ultimately allocate the residue value or a null value in the cell. This can be further clearly understood in \autoref{fig:ST Matrix Explanation 2} where for the reader's clarify and ease, the time axis is sorted in increasing rightward value; observe here that (for each row) the cells on the left of the green cells are null, and all the cells on the right of the blue-departure cell will continue to have the full residue values. In the case between the green and blue, \ie during loading/unloading, the cells are constrained to be able to take up and value less than the full residue amount without breaching the loading-unloading time constraints.
Ultimately, to ensure temporal causality is maintained, each time-column must not have any net-negative residue; a net-negative residue would indicate that some pseudo-amount of a CT has been picked up from the TP, which had not been earlier sent to that TP. This is the first time in literature that such kind of temporal causality is implemented directly through an MILP; we use the Eq.\ref{eq:62} to implement this.
\begin{figure}
    \centering
    \includegraphics[width=1\linewidth]{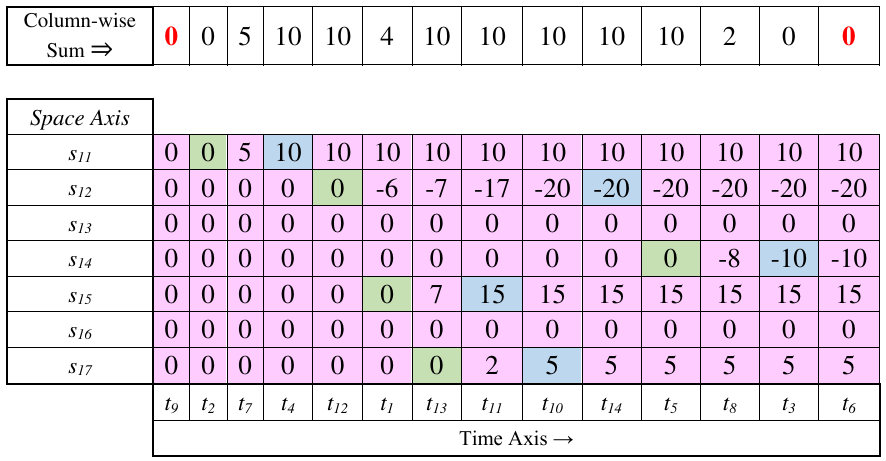}
    \caption{\textbf{A hypothetical Space-Time matrix for some arbitrary TP-CT combination}}
    \label{fig:ST Matrix Explanation 3}
\end{figure}

A hypothetical ST matrix is shown in \autoref{fig:ST Matrix Explanation 3} for better clarity. Notice that all individual column sums are either zero or positive; the first column-sum value (i.e. for the lowest time-column) must be zero. The final column-sum value (i.e. for the highest value of visit time) should also be zero to ensure that no extra resource is left at the TP after the entire emergency operation is completed, this is done using the constraint in Eq.\ref{eq:62.5}.

\clearpage

\section{Heuristic Methodology and Components} \label{Heuristic Methodology and Components}
Each component of our DT-based PSR-GIP Heuristic is expanded and explained elaborately in details here.

\subsection{Preference Generation and storing HyperParameterized Scores within Decision Tree Structures} \label{Preference Generation and storing within Decision Trees' Branches and Leaves}

We create and maintain a decision tree-like structure of the different levels of possible decisions needed to satisfy any Node. The different stages of the Decision Tree Structure (DTS) are as follows:

\begin{enumerate}
  \item Trunk: The first element of the tree is the Trunk which refers to any of the Nodes (for each Simultaneous or Split Node as the Trunk, we develop separate Decision Trees) or Transhipment Ports (which are treated similar to Simultaneous Nodes while considering a resource under transhipment).
  \item Branch: The next decision stage involves choosing Branches originating from the single Trunk, where each Branch represents a unique Vehicle Type (along with the information about the set of possible Vehicle Depots that are in the same SMTS as the Trunk).
  \item Twig: Further, Twigs emerging from each Vehicle Type-Branch represent the Vertices accessible by the specific Vehicle Type; these Vertices include those Warehouses, Transhipment Ports and Relief Centres that are located on the same SMTS as the Node-Trunk of the Decision Tree is located upon.
  \item Leaf: Finally, Leaves emerge from each Vertex-Twig of the Decision Tree, where a Leaf represents a Cargo-Type; the Twigs which represent Warehouses have only Delivery Cargo Types as Leaves whereas the Twigs representing Relief Centres have only PickUp Cargo Types as Leaves; the Leaves (\ie Delivery or PickUp Cargo Types) emerging from Transhipment Port-Twigs must be compatible to be transhipped at the respective Port. Additionally, all Cargo Type-Leaves must be compatible to be carried by the respective Vehicle Type-Branch from which those Leaves have emerged, through Twigs.
\end{enumerate}

Only those Vehicle Type-Branches are considered within any Decision Tree, for whom there are a non-zero number of vehicles available at any of the corresponding Vehicle Depots (the Vehicle Depots necessarily being present on the same SMTS); all Branches within a Decision Tree represent unique Vehicle Types, in case different Vehicle Depots on the same SMTS have the same Vehicle Types in diverse quantities having a single Branch of that specific Vehicle Type suffices (the information of all Vehicle Depots, on the same SMTS, which contains the particular Vehicle Type is maintained for each respective Vehicle Type-Branch separately) similar to the Vehicle Depot details for each Vehicle Type (as mentioned in \autoref{fig:Imagine_TP2}).

Additionally, when considering the Vehicle Type-Branches for the Decision Trees with Simultaneous Nodes as Trunks, only those Vehicle Types are allowed to form Branches that can satisfy the entire resource demand of the Node-Trunk in one go; this means that the Vehicle Capacity (both Volume and Weight, separately) must be greater than or at least equal to the combined amount of PickUp requested at the Simultaneous Node-Trunk, and the same Vehicle Capacities must also not be less than the combined amount (Volume and Weight) of Delivery requested by the same Simultaneous Node-Trunk. Thus Branches within Decision Trees of different Simultaneous Nodes (even when located in the same SMTS) can be very different depending on the actual combined amount (Volume and Weight) of the PickUp or Delivery Cargo quantities requested. One of the many decision trees for our illustrative example Instance (section \ref{An Example Case Study}) is depicted in \autoref{fig:Imagine_TP2} with the Trunk as "TP2", the Vehicle Types compatible to approach the Trunk Nodes as branches ("VT2", and "VT5" in the figure) alongwith their respective depots, the subsequent Twigs inidicating the vertices that are connected to the Trunk directly through the respective Branches, and finally the Leaves representing the compatible cargo types that are transferable (across Transhipment Ports) or available to be satisfied (for Warehouses and Relief Centres).

At each Leaf within DTSs, we store HyperParameters which help in the calculation of a preference $Score$ for satisfaction of a Node-Trunk by a Vertex-Twig in the same Decision Tree; this is detailed in sections \ref{Preference Generator Function for Split Node} and \ref{Preference Generator Function for Simultaneous Node}.

\subsubsection{Requirement Code, Compare Code and Score Calculator} \label{Compare Code Calculator}
Nodes indicated as Trunks in the Decision Trees have requirements for some/all of the Cargo types which is stored as a Requirement Code (this has been implemented in code as a Pythonic dictionary with the keys being the Cargo names and the respective values indicating the amount of Delivery and PickUp requirements for each Cargo Type).

A Compare Code is calculated for a Trunk \wrt any of its Twigs, \ie a Node's Requirement Code is compared with the available resources (for Delivery Loads) or resource acceptance capability (for PickUp Loads) of a Vertex Twig (this Vertex Twig can be either a Warehouse or a Relief Centre) to generate a Compare Code (in the code implementation, a Python dictionary is generated for each Compare Code calculation, where the keys remain the Cargo names and the values indicate the maximum amount of that specific resource satisfaction possible by that Vertex Twig). For any Cargo Type, the value in the Compare Code is the minimum of: the necessity of that Cargo Type at the Node-Trunk, and the availability of that same Cargo Type at the Vertex Twig.

Having computed the Compare Code for a Node-Trunk---Vertex-Twig pair, a score is developed to signify the amount of resource satisfaction the Vertex Twig can provide to the Node-Trunk, and this score essentially translates to a preference as described subsequently.

For our study-example instance (\autoref{fig:Picture_Instance}), the Requirement Code for the Split Node ``NP3" is \{``2P":6, ``1D":11, ``2D":13\} (as the other compatible pickup cargo of ``1P" is null for this node; the Cargo Type names ending in ``P" indicate a PickUp type and those ending in ``D" indicate a Delivery type). When comparing this node with the Warehouse ``WH1", the Compare Code is \{``2P":0, ``1D":11, ``2D":13\}. If the Requirement Code had been \{``2P":6, ``1D":41, ``2D":33\}, then the Compare Code \wrt the same Warehouse would have been \{``2P":0, ``1D":20, ``2D":20\} as the resource limits of the concerned Warehouse would have been reached. Similarly, for the Requirement Code of the Simultaneous Node ``NM2" being \{``1P":20, ``1D":3, ``2D":2\}, its Compare Code \wrt the Relief Centre ``RC1" is \{``1P":20, ``1D":0, ``2D":0\}; had the Requirement Code been \{``1P":100, ``2P":100, ``2D":13\}, the Compare Code would have been \{``1P":50, ``2P":0, ``2D":0\} as ``RC1" has a maximum of 50 units of capacity availability to accept the PickUp type ``1P".

% We also develop a possible availability estimate for each Transhipment Port which is discussed in \ref{Recursive function to enable Transhipment of a single load-type}.

Before describing the preference score generation, it would be apt to mention that any symbol used in section \ref{Compare Code Calculator} is only for demonstration purposes for the score calculation and does not carry any relation/significance with previously defined symbols in \autoref{longTab: Sets and Parameters} or \autoref{longTab: Variable Definitions} unless specifically mentioned. For each Cargo Type in the Compare Code, we develop a fraction $f$ indicating the ratio of its value in the Compare Code \wrt the Requirement Code (only if there is a non-zero value in the Requirement Code); if this fraction $f$ is $0$, it would mean that the Twig Vertex was unable to satisfy any amount of the specific  Cargo Type for this demand, while $f$ equal to $1$ indicates that the entire amount of this requirement was satisfied for the specific Cargo Type. For a Requirement Code of \{``1P":100, ``2P":100, ``2D":13\} and a Compare Code being \{``1P":50, ``2P":0, ``2D":0\}, we have the Fractions as $\{f_{1P}:0.5, f_{2P}:0, f_{2D}:0\}$. 

For effective representation, each fraction $f$ is subscripted by its respective Cargo Type of its Leaf. Each individual Leaf is associated with two HyperParameters one to act as a Multiplier (designated as $m$ in the Eqns. \ref{Eq.70}, \ref{Eq.71} and \ref{Eq.72}) and another to act as an Exponent (designated as $p$ in Eqns. \ref{Eq.70} and \ref{Eq.72}) to the fraction $f$ for the preference score generation. These HyperParameters are random values stored in each of the Cargo-Leaves. All HyperParameter values remain constant during every Main Iteration and are changed before the next Main iteration. The HyperParameter values are taken from Normal distribution ($\mu$ and $\sigma$ chosen randomly between 0 and 1) or Uniform distribution (upper and lower bounds chosen randomly between 0 and 1), allowing different HyperParameters to be derived from different distributions (or be equal to unity as a rare probability). During the score calculation, a $Numerator$ and a $Denominator$ are computed separately for each Twig of the Decision Tree using the Compare Codes generated by comparing the Requirement Code of the Trunk-Node and the resource satisfiability at the Twig.

For the $Numerator$ calculation, each Cargo Type specific satisfiability fraction is raised to the power of a unique Exponent-HyperParameter value (indicated in $p$ in Eq. \ref{Eq.70}) and then multiplied by a unique Multiplier-HyperParameter (indicated in $m$ in Eq. \ref{Eq.70}). Each of these HyperParameters are Leaf specific and are represented with the appropriate subscripts of the specific Cargo Type represented by the Leaf. Assuming a sample Requirement Code to be $\{C_1:100, C_2:0, C_3:125, C_4:10, C_5:15, C_6:0, C_7:5\}$, a sample Compare Code could be $\{C_1:80, C_2:0, C_3:125, C_4:0\}$ or $\{C_5:15, C_6:0, C_7:0\}$, when comparing the Node requirements with a Warehouse or a Relief Centre Respectively. Here we assume that the first 3 cargo types are Delivery with $C_4$ and $C_5$ being PickUp load types; thus the Compare Codes are different when calculated with Vertex Twigs which are Warehouses (the first sample Compare Code above), versus when the Vertex Twigs are Relief Centres (the second sample Compare Code above). The Compare Codes contain possible satisfaction values for all Cargo types for which the Vertex Twig has Leaves. The final score is calculated as the ratio of the $Numerator$ and the $Denominator$, as shown for the first sample Compare Code in Eqns. \ref{Eq.70} and \ref{Eq.71} respectively.
\begin{equation} \label{Eq.70}
    Numerator =  
    m_{C_1} * (f_{C_1}) ^ { p_{C_1} } +
    m_{C_3} * (f_{C_3}) ^ { p_{C_3} } +
    m_{C_4} * (f_{C_4}) ^ { p_{C_4} }
\end{equation}
\begin{equation} \label{Eq.71}
    Denominator =  
    m_{C_1} +
    m_{C_2} +
    m_{C_3} +
    m_{C_4}
\end{equation}

It can be observed that the calculation of the $Numerator$ may have some fraction $f$ values missing (due to $0$ value in the Requirement Code), for which the corresponding Multiplier-HyperParameter values are still added up in the $Denominator$, this is done so that the scores for the same Node-Trunk across different Vertex Twigs (of a similar kind, \ie Warehouses or Relief Centres) may be comparable. A sample $Score$ calculation is also shown for the second sample Compare Code in Eq. \ref{Eq.72}.
\begin{equation} \label{Eq.72}
    Score = 
    \frac{Numerator}{Denominator}
    =
    \frac
    {m_{C_5} * (f_{C_5}) ^ { p_{C_5} } +
    m_{C_7} * (f_{C_7}) ^ { p_{C_7} }}
    {m_{C_5} +
    m_{C_6} +
    m_{C_7}}
\end{equation}

A Requirement Code is compared with each Vertex Twig in its DTS to generate such Compare Codes and Scores. Two other HyperParameters (these are single Hyperparameters for the entire problem) toggle the method of forming the $Score$:
\begin{itemize}
  \item Time consideration during Scoring: If this HyperParameter is ON/True, the Travel Time from the Node-Trunk to the Vertex-Twig is multiplied to the whole Denominator.
  \item Denominator consideration during scoring: If this HyperParameter is OFF/False, then the entire Denominator during all $Score$ calculations is made equal to unity.
\end{itemize}

\subsubsection{Recursive function to enable Transhipment of individual Cargo Types} \label{Recursive function to enable Transhipment of a single load-type}

Here we develop degrees of transhipment according to the depth of the transhipment; the degrees can be imagined as potentials with possible transhipments in lower degrees (lower potentials) happening first before moving to higher potentials (deeper degrees).

This recursive function accepts the following six arguments:
%\begin{itemize}[label={\textbullet}]
\begin{enumerate}[label={(\Roman*)}] % If just use \begin{enumerate}, then labels will be wrong
    \item \label{item:One} Already considered Warehouses and Relief Centres 
    \item Vehicle Type which sent this requirement \label{item:Two}
    \item Cargo Type name and the minimum quantity that needs to be transhipped \label{item:Three}
    \item Name of the Transhipment Port which is being considered for this transhipment requirement \label{item:Four}
    \item Degree of Transhipment (set to the default value of 1) \label{item:Five}
    \item A Tabu List of Vehicle Type and Transhipment Port paired combination which have already been considered to satisfy this specific transhipment requirement (the initial default value is set to an empty list) \label{item:Six}
    \item A List of SMTS which have been previously considered for Transhipment \label{item:Seven}
\end{enumerate}
%\end{itemize}

When this function is initially called from other functions, only the first four arguments are passed, the last three arguments are used only during recursive calls of this same function from within itself. This current SMTS can be inferred from the Branch-Twig information in \autoref{item:Two} and \ref{item:Four} (since one Vehicle Type is allowed to ply on a unique SMTS in our considered problem), for updating the list in \autoref{item:Six} during recursions.

The primary objective of this function is to generate preference scores for satisfying a Transhipment Port similar to that of a Split Node (the case of satisfying all the requirements of any single Transhipment Port together, similar to a Simultaneous Node, is not dealt with in our problem and could be a future research direction); for this, the DTSs with each individual Transhipment Ports as Trunks are considered. 

To generate the preference score (similar to as described in section \ref{Compare Code Calculator}, we need to have the Compare Code first, however since Transhipment Ports inherently don't contain static resources, we use this function to find out how much of a specific resource can be obtained (for Delivery types, from subsequent Transhipment Ports or Warehouses) or distributed (for PickUp types, from subsequent Transhipment Ports or Relief Centres), alongwith the depth of Transhipment necessary for the same. This depth of Transhipment, which we have termed as degree, may not necessarily be a single number which can be directly inferred.
For example, in the study-example instance described in \autoref{fig:Picture_Instance}, there are no Warehouses in the SMTS ``Road\_A" and therefore, all Deliveries must be resourced using the Transhipment Port ``TP3". In this case, when this function is called using the first four argument values as below:
\begin{enumerate}[label={(\arabic*)}]
    \item ``RC1", \label{eg:point1}
    \item ``VT3",
    \item ``1D": 5 (assuming that the request to this function was sent from ``NP1" as the Trunk),
    \item ``TP3" (as one of the Twigs in the DTS with Trunk as ``NP1"), \label{eg:point4}
\end{enumerate}
 there would be two types of function returns, the first Component contains only the vertices containing Static resources that may be able to (fully or partially) satisfy the requirement of the Trunk Node specifically for the requested type of Cargo (``1D" in this case) alongwith the information of their depth/degree and the quantity of the resource that they have, and the second Component containing the same information clubbed together such that only vertices that are approachable directly via a single new journey (i.e. upto one higher degree) is mentioned (with a combined degree for Transhipment ports requiring multiple subsequent transhipments) alongwith quantity of resources (available at downstream degrees for Transhipment Ports, or present in the vertices containing static resources like Warehouses or Relief Centres).

Whenever this function is queried, the Decision Tree with the Trunk being the Transhipment Port in \autoref{item:Four} above is considered; with its Branches being the Vehicle Types that have access to this Trunk through different SMTSs; the Twigs representing the different Warehouses, Transhipment Ports and Relief Centres which are accessible by the respective Vehicle Type-Branch; and the final Leaves representing those Cargo Types which are compatible to be carried on the respective Vehicle Type-Branch as well as allowed to be transhipped through the Transhipment Port-Trunk. Since this function takes care of only one resource type at a time, all leaves of only this Cargo Type is considered during a specific query and throughout during recursions of that same query.

The functional approach is described in pseudocode:

% Pseudocode Heading
\noindent\hrulefill

\noindent\textbf{Recursive function to satisfy transhipment requirement of a single CT}

\noindent\hrulefill
\begin{algorithmic}[1]

\Function{Transhipment\_of\_a\_CT\_from\_a\_TP}{considered\_PRVs, previous\_VT, (CT\_name,Q), TP\_name, transhipment\_Degree=1, Tabu\_List, used\_SMTS=[ ]}

  CT\_Satisfaction\_Counter = 0
  \For{each VT-Branch in the DTS with TP\_name as Trunk (except when the VT-Branch is same as previous\_VT):}
    \For{each Twig emerging from the outer iterator:}
      \If{the iterator Twig is a TP, \textbf{or} the iterator Twig is among considered\_PRVs}
        \State \Continue \Comment{iterate over to the next Twig, till a suitable PRV is encountered}
      \EndIf

      \For{each CT-Leaf associated with the iterator Twig, \textbf{and} if the Leaf represents CT\_name:}
        \State Satisfiable\_Leaf\_Amount = min(available quantity of CT\_name at iterated Twig, Q)

	\If{the iterator Twig is not within considered\_PRVs} \Comment{Updates are done only once for each unique vertex, since a single PRV can be encountered in multiple Twigs on different Branches in this single Decision Tree and considering them multiple times would be a faulty approach}
          \State Satisfiable\_Leaf\_Amount is stored in Temp\_DS and associated with its PRV (\ie iterator Twig)
	  \State CT\_Satisfaction\_Counter $\gets$ CT\_Satisfaction\_Counter + Satisfiable\_Leaf\_Amount
	\EndIf

      \EndFor

    \EndFor
  \EndFor

  \State considered\_PRVs is extended with the information from Temp\_DS
  \If{CT\_Satisfaction\_Counter $\geq$ Q}
    \State This indicates all requested resources can be transhipped from PRVs that are just one degree away
    \If{Global\_Transhipment\_Fathom is \textbf{False}}
      \State \textbf{return} considered\_PRVs, considered\_PRVs \Comment{Since there is no deep transhipment yet, both return Components are same}
    \EndIf
  \EndIf

  \State used\_SMTS is updated with the information of this current network-segment, \ie from the previous\_VT and TP\_name information

  \State marge\_$\mathbb{C}_1$ = [ ]
  \State marge\_$\mathbb{C}_2$ = [ ]

  \For{a\_VT\_Branch among all Branches in the DTS with TP\_name as Trunk, except when a\_VT\_Branch represents previous\_VT:}
    \For{a\_TP\_Twig among all Twigs associated with a\_VT\_Branch:}

      \If{the Branch-Twig combination pair (a\_VT\_Branch, a\_TP\_Twig) is present in Tabu\_List}
        \State \Continue \Comment{We skip for the combinations already considered and thus present in the Tabu\_List}
      \EndIf

      \State Dynamic\_Tabu\_List is generated by updating Tabu\_List with this new Branch-Twig (i.e., (a\_VT\_Branch, a\_TP\_Twig)) combination pair

      \For{each CT-Leaf associated with the iterator TP-Twig, \textbf{and} if the Leaf represents CT\_name:} \Comment{Some CT-Leafs could be missing due to incompatibility issues with VT-Branch or TP-Trunk}

        \State new\_Q = 0 \Comment{new\_Q is the subsequent transhipment request amount}
	\If{FullFathom is \textbf{True}}
	  \State new\_Q = Q
	\Else
	  \State new\_Q = Q - CT\_Satisfaction\_Counter
	\EndIf

	\State r$\mathbb{C}_1$,r$\mathbb{C}_2$ = \textsc{Transhipment\_of\_a\_CT\_from\_a\_TP}(considered\_PRVs, a\_VT\_Branch, (CT\_name,new\_Q), a\_TP\_Twig, transhipment\_Degree+1, Dynamic\_Tabu\_List, used\_SMTS)\Comment{The amount of resources transhipable at a\_TP\_Twig is found through this recursive functional return data of return Components r$\mathbb{C}$}

        \State r$\mathbb{C}_2$ is stored in the respective CT-Leaf in the DTS with the Trunk being a\_TP\_Twig \Comment{The value of relative degree is best assumed 0 during this storage as the same function-calculation can also happen independent of previous-degrees; therefore the value of transhipment\_Degree is substracted from each of the internal degree reportings while storing}

        \State marge\_$\mathbb{C}_1$ is enhanced with r$\mathbb{C}_1$ data
        \State marge\_$\mathbb{C}_2$ is enhanced with r$\mathbb{C}_1$ data, in the format of Component\_2

      \EndFor
    \EndFor
  \EndFor

  \State \textbf{return} marge\_$\mathbb{C}_1$, marge\_$\mathbb{C}_2$

\EndFunction
\noindent\hrulefill
\end{algorithmic}

Elaborating on the first set of return from step 20 above, for each PRV in Component\_1 \ie in considered\_PRVs, there are three values associated: the first value indicates the Degree of Transhipment, the second value indicates the respective Satisfiable\_Leaf\_Amount for when the specific PRV-Twig was being considered for the satisfaction of the Trunk's requirements, and the third value indicates the maximum amount of CT\_name present at the concerned PRV-Twig.

A unique SMTS can be directly identified with the information of a VT and a TP. Expanding on the step 28, deeper recursions happening from different Branches/Twigs can only consider combination-pairs not encountered previously, but parallel recursions may consider the same combination-pairs. Similarly, PRVs already used in a Potential/Degree cannot be encountered in any other Potential/Degree within the same recursion (such pseudo transhipments are rejected through the usage of the Tabu\_List, they may however be encountered in the separate parallel recursions).

Higher Degrees of Transhipments indicate complex problems needing multiple stages of resource transhipment for proper allocation and transfer; these Degrees could be compared with Potentials where a lower Potential for Transhipment is preferred and the algorithm moves to an higher potential for new transhipment search only when all possibilities within a lower Potential value are unable to satisfy the resource transhipment requirement completely.

Each Dynamic\_Tabu\_List is an update over the function’s original Tabu\_List passed in its argument. To be more clear, Tabu Lists increase through progressively deeper recursive function calls, but change across recursive function calls made during other iteration-scans while traversing the DTS. Detailing on the recursive functional \textit{return} component pair r$\mathbb{C}$ from step 39:
\begin{enumerate} [label={(\Alph*)}]
    \item \textbf{Component\_1} consists of only those PRVs which can satisfy the concerned transhipment amount; since the PRVs can be on different SMTSs, there will be intermediate transhipments necessary. Each PRV information in Component\_1 is associated with the maximum quantity of the CT that is available and satisfiable at the PRV, and the degree of transhipment. If the transhipment degrees of some of the PRVs are higher than the rest, it would mean that access to that PRV is possible only through more subsequent transhipments across different SMTSs.
    \item \textbf{Component\_2} represents the satisfaction of the concerned transhipment amount through usage of Vertices (could be WHs, RCs, and TPs) in the just next transhipment degree, \ie on the same SMTS as that of  a\_TP\_Twig. Thus, the Component\_2 can be derived from the Component\_1. The PRVs used through subsequent deep-transhipments are compacted within Component\_2's TP data; this is represented as the Downstream\_Information where the PRVs that are accessed through subsequent transhipment are associated within their specific TP-vertex inside Component\_2 along with their full Transhipment Degree values. This in addition to the information types of Component\_1, \ie the transhipment degree, and amount of CT satisfiable by a vertex, each TP-vertex in Component\_2 also has the additional Downstream\_Information of the subsequent PRVs which can be encountered through that TP-vertex. The satisfiable/available resource quantities (\ie the second element here) are used during the Compare Coding of TPs.
\end{enumerate}

We now elaborately explain the processing of the Component Pairs received from recursive function calls, \ie the steps 41 and 42. The different Component\_Pairs captured (during the iteration-scan) at each Leaf are aggregated and moulded into a final \textit{return} component pair in step 46:
\begin{enumerate}
    \item To obtain the aggregated Component\_1, all the first components caught from the different recursive function calls are merged into marge\_$\mathbb{C}_1$ (step 41).
    \item marge\_$\mathbb{C}_2$: To obtain the aggregated Component\_2 (step 42), all the first components from different recursions are used to infer a single degree representation. This Adjusted Degree of Transhipment is calculated for each encountered TP, and is a value between the maximum and minimum of all deeper Transhipment Degrees encountered (depending on the value of a HyperParameter Transhipment\_Degree\_Slider, which is toggled per Main Iteration). The Adjusted Degree of a TP may best be weighted with the resource satisfiability, across all PRV degrees under its TP; this Adjusted Degree closely represents the depth of transhipments necessary. The calculation of the weighted average-Adjusted Degrees is shown through an independent example (not related to our study-example).
    
    Assume the following Component Pairs as the functional \textit{return}s (in step 39) caught during different recursions (\ie these component pairs are stored in r$\mathbb{C}_1$ and r$\mathbb{C}_2$ across parallel recursions):
    \begin{enumerate}[label=\textbf{Eg. \arabic*}:] % Custom label format
        \item Component\_1: ``WH1'': (1, 20, 20), ``WH2'': (2, 15, 25), 
        % The inferred degree from this Component could be:
        % \begin{itemize}[label=\#]  %[\textbullet]
        %     \item 2, if maximum among the degrees is considered
        %     \item 1, if minimum among all degrees is considered
        %     \item anything between 1 and 2, depending on some conceived HyperParameter
        %     \item 
        % \end{itemize}
        We use weighted sum to calculate the inferred degree for a TP in Component 2 as $\frac{1*20 + 2*15}{20 + 15} = \frac{10}{7}$, with the weights being the satisfiable resource quantities as mentioned in the second tuple elements.
    
        Component\_2: ``TP1'': ($\frac{10}{7}$, 35, \{1:\{``WH1'':20\}, 2:\{``WH2'':15\}\})

        This is done assuming that the ''TP1'' is used for a deeper transhipment to access ''WH1'' and ``WH2''.
        
        \item Component\_1: ``WH3'': (1, 15, 15), ``WH4'': (2, 10, 10), ``WH5'': (3, 10, 10)
        
        Inferred degree: $\frac{1*15 + 2*10 + 3*10}{15 + 10 + 10} = \frac{13}{7}$
        
        Component\_2: ``TP2'': ($\frac{13}{7}$, 35, \{1:\{``WH3'':15\}, 2:\{``WH4'':10\}, 3:\{``WH5'':10\}\})

        \item Component\_1: ``WH5'': (2, 10, 10), ``WH2'': (3, 25, 25), ``WH4'': (3, 10, 10)
        
        Inferred degree: $\frac{2*10 + 3*25 + 3*10}{10 + 15 + 10} = \frac{25}{7}$
        
        Component\_2: ``TP3'': ($\frac{25}{7}$, 35, \{2:\{``WH5'':10\}, 3:\{``WH2'':25, ``WH4'':10\}\})
    \end{enumerate}

    Now, for aggregating the above assumed recursion results, we will develop the Component Pair:
    \begin{enumerate}
        \item Aggregated Component\_1 (step 41): Here we group all the first components from the above three results to obtain: ``WH1'': (1, 20, 20), ``WH2'': (2, 25, 25), ``WH3'': (1, 15, 15), ``WH4'': (2, 10, 10), ``WH5'': (2, 10, 10).
        
        This indicates the maximum amount of resources that can be leveraged from unique vertices and through the lowest potential degree. Observe that all the vertices in the third result are already available in either of the previous two results, thus we always include the lower degrees for vertices with multiple occurrences and consider the higher of the resource satisfiable quantity (\ie the second elements of their tuples).
         
        \item Aggregated Component\_2 (step 42): If all the three Component\_2s had been for the same TP (in the examples \textbf{1, 2} and \textbf{3} above), say ``TP1'', then we needed to obtain a common degree; we could then calculate a weighted degree from the Aggregated Component\_1 as well by $\frac{1*20 + 2*25 + 1*15 + 2*10 + 2*10}{20 + 25 + 15 + 10 + 10} = \frac{125}{80}$. For the other elements like resource satisfiability, we take the best values; and take union of the downstream information set also with the respective best values and degrees (similar to as is done in the Aggregated Component\_1).
        
        However, all the Component\_2s have unique TPs (the same TPs will be encountered if there are more than one VTs plying on a SMTS) in our example, and therefore we can directly group them up: ``TP1": ($\frac{10}{7}$, 35, \{1:\{``WH1'':20\}, 2:\{``WH2'':15\}\}), ``TP2'': ($\frac{13}{7}$, 35, \{1:\{``WH3'':15\}, 2:\{``WH4'':10\}, 3:\{``WH5'':10\}\}), ``TP3'': ($\frac{25}{7}$, 35, \{2:\{``WH5'':10\}, 3:\{``WH2'':25, ``WH4'':10\}\}).
    \end{enumerate}
    The HyperParameters Global\_Transhipment\_Fathom, and FullFathom, toggle between True and False values every Main Iteration; for the above example they may just be assumed False.
\end{enumerate}

% =====--------------------------------------------

\noindent \textbf{Continuing our study-example discussions:}

When this function is queried with the inputs \autoref{eg:point1}-\ref{eg:point4}, since there is no PRV in any of the higher degree SMTS (in this case there is only a single SMTS of degree one, which is SMTS ``Air\_A''). Thus the first functional \textit{return} (of step 20) is not executed and the CT\_Satisfaction\_Counter remains 0. Thus the next portion of scanning the TP-Twigs in the DTS starts (on the single available Branch of ``VT2'', hosted at ``VD3''). Now there are two ports encountered during the iteration across the Twigs, namely ``TP1'' and ``TP2'' for which this same function is called independently (as two separate parallel recursion calls).

For the first recursion call, the functional arguments passed are:
\begin{enumerate}[label={(\arabic*)}]
    \item ``RC1'', (this is the list of all PRVs encountered in the previous SMTS ``Road\_A'', and in the quering SMTS ``Air\_A'')
    \item ``VT2'',
    \item ``1D'': 5, (as there was no reduction)
    \item ``TP1'',
    \item degree = 2
    \item (``VT3'', ``TP3''),
    \item ``Road\_A'' (this was inferred using the Branch-Twig information, and the current SMTS being used will be inferred at the beginning of the next recursive function summon).
\end{enumerate}
The functional \textit{return} obtained being ``WH1'': (2, 5, 20) both for the components 1 and 2; the first degree refers to the relative degree depth where the search was done, the next value of 5 is essentially the minimum among the resource requested and the resource available, and the next 20 represents the maximum resource availability at ``WH1''.

For the second recursion call, the functional arguments passed are same as above with the only difference being in the fourth argument being ``TP2''; the functional \textit{return} obtained is ``WH2'': (2, 5, 20) both for the components 1 and 2. The original function instance that receives these recursion \text{return}s in step 39, combines them depending on the component:
\begin{enumerate} [label={(\Alph*)}]
 \item For Component\_1: The final aggregation is obtained as \{``WH1'': (2, 5, 20), ``WH2'':(2, 5, 20)\}
 \item For Component\_2: The information is represented as \{``TP1'': (2, 5, \{2:\{``WH1'':20\}\}), ``TP2'': (2, 5, \{2:\{``WH2'':20\}\})\}. Here the first value of 2 for each TP represents the inferred (aggregated) degree from within the downstream information; the downstream information is the third element representing the currently available resources and the respective relative degree of that vertex \wrt the PRV-Twigs, arranged as per degree-wise dictionary-keys).
\end{enumerate}

\subsubsection{Preference Generator Function for Split Node} \label{Preference Generator Function for Split Node}

This function develops preferences for satisfying the Split Node passed to it as an argument. The DTS with the unsatisfied Split Node as its Trunk is already available before these Preference Generator Functions (PGFs) are called. The task of PGFs (both the Split and Simultaneous PGFs) is to assign preferences to the Twigs of the specific Decision Tree such that higher values indicate more affinity of the Node-Trunk to be satisfied by a Twig.

The passed Split Node-Trunk has its own Requirement Code and the function discussed in section \ref{Compare Code Calculator} is used to generate the Compare Codes along with the preference $Scores$ initially only for all Warehouse and Relief Centre Twigs in the Decision Tree. It is then checked if the Warehouse and the Relief Centres directly accessible by the Split Node (them being on the same SMTS) are able to fully satisfy the requirements of the Split Node. This is done by adding all the generated Compare Codes (Cargo type-wise separately) and comparing if each cargo amount obtained from this addition is at least equal to or exceeds the Requirement Code amounts.

If the satisfaction requirements cannot be fulfilled by the $W's$ or $R's$ on the same SMTS, a Transhipment Code is generated which contains only the extra amounts for each Cargo Type, these extra amounts being the difference of the specific Cargo Type in the Requirement Code \wrt the summed up Compare Codes; the Transhipment Code is populated with those corresponding Cargo Types whose excess amounts are positive, (as after the difference, negative values indicate that they may be satisfied through the $W's$ and $R's$). During the actual satisfaction and route element generation, this consideration will not hold every time as the resources in the $W's$ and $R's$ will continually decrease and new Transhipment requirements would need to be catered to then.

Transhipment Codes can also takeup the entire initial value of the Requirement Code, depending on whether the single HyperParameter of Full\_Demand\_Passed\_for\_Enhanced\_Variability is set as On/True; this allows the function to try and find Transhipment possibilities for the entire amount and is essentially a stress testing of the Heuristic, since in reality, this sort of a situation would be the absolute worst scenario.

After the Transhipment Code is finalized (\ie with amounts within the Transhipment Code being greater than zero), the maximum transhipment possible for each Cargo Type within the Transhipment Code is found using the function described in section \ref{Recursive function to enable Transhipment of a single load-type}; the Recursive function for Transhipment calculation is called for each Transhipment Port-Twig in the DTS, with the Trunk being the Split Node. The function \textit{return}s the information of maximum possible resource satisfaction (similar to a Compare Code, and termed as the Satisfaction Code) for each Cargo Type within the Transhipment Code.

The maximum and minimum Transhipment Degree, which is encountered while satisfying each Cargo type, is also inferred from this functional \textit{return}. A Transhipment\_Degree\_Slider\_for\_$N^P$ HyperParameter (common for the entire problem) is used to find a value between these maximum and minimum Transhipment Degree bounds; this Adjusted Transhipment Degree is used to find the preference $Score$ for the Transhipment Port using the function described in section \ref{Preference Scoring based on Satisfaction Codes of Transhipment Ports considering Adjusted Transhipment Degrees}.

It should be mentioned that DTS of the TPs are populated with the necessary Component Pairs which are generated during the function-usages of sections \ref{Preference Generator Function for Split Node} and \ref{Preference Generator Function for Simultaneous Node} (the internal recursions of the Transhipment function in \ref{Recursive function to enable Transhipment of a single load-type} help in populating the possible transhipment resource quantities along with their respective degree estimates).

The Split Nodes which require transhipment to be satisfied are identified and stored within a Mandatory\_$N^P$\_Transhipment\_Requirement dictionary along with the Cargo Type and the respective quantity which is necessary to be transhipped (since the specific Split Node can be partially satisfied without transhipments).

\subsubsection{Preference Generator Function for Simultaneous Node} \label{Preference Generator Function for Simultaneous Node}

This function is necessarily slightly different from the PGF for Split Nodes, since the PGF for Split Nodes considers partial satisfaction, therefore utilizing Leaves across different Branches is allowed there. However, in this case, to satisfy Simultaneous Nodes, all the Leaves utilized within any Decision Tree must be from one single Vehicle Type-Branch (this simple solution of developing routes from the same branch for a Simultaneous Node has been another major advantage of using the DTS; other advantages include inherent filtering of incompatibilities to create the necessary feasible solution space). During the actual route generation, multiple trips of that specific Vehicle (of the same type as the Branch) may be required.

A parameter within this function is the Consideration\_with\_Transhipment, which remains set at a default False state initially. If the HyperParameter of Full\_Demand\_Passed\_for\_Enhanced\_Variability is set as On/True, then the functional parameter of Consideration\_with\_Transhipment is changed to True at this function start. Here too, similar to the PGF for Split Nodes, we develop a Requirement Code for the passed Simultaneous Node. 

We iterate across Vehicle Type-Branches, and across the Warehouse- and Relief Centre-Twigs within each Branch to generate Compare Codes and corresponding $Scores$ for each Twig using the function described in section \ref{Compare Code Calculator}. During the iteration, only those Vehicle Type-Twigs are considered which are able to carry the full resource requirement of the Simultaneous Node (both in terms of Volume and Weight). The Requirement Code of the Simultaneous Node-Trunk is compared with the aggregate of all the Compare Codes \wrt all Vertex-Twigs in the same Vehicle Type-Branch and any remaining resource satisfaction requirement is checked. If all resources are able to be satisfied, then we proceed to check the next Vehicle Type-Branch, otherwise, depending on if the Consideration\_with\_Transhipment HyperParameter is True, we further check the same Vehicle Type-Branch with the Vertex-Twigs now being the Transhipment Ports. The recursive function for Transhipment estimation (\ref{Recursive function to enable Transhipment of a single load-type}) is used to find the estimated resource quantity for the specific resources which were not satisfiable previously without transhipment. After all the Transhipment Ports accessible by the Vehicle Type-Branch is iterated, it is now checked if this Branch (along with all its Twigs of Warehouses, Relief Centres, and Transhipment Ports) is able to satisfy the Requirement Code completely. If the Branch is able to fully satisfy the entire Requirement Code, then an additional $Score\_Sum$ is stored for each Branch; otherwise this Branch is rejected. The $Score\_Sum$ for a Branch is the sum of $Scores$ of all Twigs originating from that Branch.

Here too, all Simultaneous Nodes which were found to contain all the Branches with Transhipment requirements are stored within a Mandatory\_$N^M$\_Transhipment\_Requirement, and the function calls to the recursive Transhipment estimation function (\ref{Recursive function to enable Transhipment of a single load-type}) helps in populating the resource leverage estimates for each Transhipment Port.

\subsubsection{Preference Scoring based on Satisfaction Codes of Transhipment Ports considering Adjusted Transhipment Degrees} \label{Preference Scoring based on Satisfaction Codes of Transhipment Ports considering Adjusted Transhipment Degrees}
This function uses its inputs of the Satisfaction Code, Transhipment Code, and an additional Degree Code (which contains Cargo Type specific Adjusted Transhipment Degrees), to \textit{return} a score. This score is meant to be a preference for the transhipment of resources through a specific Transhipment Port-Twig to satisfy a Node-Trunk via a specific Vehicle Type-Branch; this function is similar to the function for generating Compare Codes (\ref{Compare Code Calculator}) with the primary difference being that the function in \ref{Compare Code Calculator} is specifically meant to generate preference scores for those Twigs which represent Warehouse and Relief Centres; for scoring the Transhipment Ports represented as the Twigs, we use this function (\ref{Preference Scoring based on Satisfaction Codes of Transhipment Ports considering Adjusted Transhipment Degrees}).

Information of the Trunk, Branch and Twig is also passed to this function so that it is able to develop the score using the appropriate HyperParameters at each Cargo Type-Leaf similar to as described in the Eqns. \ref{Eq.70}-\ref{Eq.72}. In this case, the fraction $f$, as is calculated for each Cargo Type separately, is the ratio of amounts of the same Cargo Type in the Satisfaction Code \wrt the Transhipment Code.

Apart from the two HyperParameters affecting the $Denominator$ during scoring (Travel Time multiplication in the $Denominator$, or making the $Denominator$ unity), we also consider the following two additional HyperParameters in this function to utilize the Adjusted Degree:
\begin{itemize}
  \item Degree\_affects\_Multiplier: If this HyperParameter is ON/True, the Multiplier-HyperParameters are divided by the Adjusted Degree obtained for the specific Cargo-type in the Degree Code.
  \item Degree\_affects\_Exponent: If this HyperParameter is ON/True, the Exponent-HyperParameters are multiplied by the Adjusted Degree found in the Degree Code for the specific Cargo-type.
\end{itemize}

These are two separate Hyperparameters for the entire Problem, which focus on toggling the scoring mechanisms during different Main Iterations. The logic behind using the Adjusted Degree to alter the preference (through alterations in the Exponent- and/or Multiplier-HyperParameters) is based on the natural sense that if some Cargo is able to be obtained using lesser number of Transhipments, then we should actually prefer that, instead of going for higher number of Transhipments to obtain that Cargo Type (in full or in part). For very complex problems with reduced feasible regions, this may not always be the case and could allow overlooking the optimal solution which may only be obtained using some counter-intuitive logic; however, for practical real scenarios, we firmly believe that lesser number of Transhipments would be better, and thus our Heuristic moves to a higher Potential/Degree for the Transhipment only when no more Cargo Load may be satisfied in the lower Potentials (keeping the relative Transhipment Degree low); therefore no Transhipment (\ie Degree $0$) is the most preferred.

% *********************************************************************
% *********************************************************************

\subsection{Creating Smallest Route Elements (SREs)  while maintaining the Causality Dict (CD)} \label{Creating Smallest Route Elements (SREs) while maintaining a Causality Dict}

\begin{figure}
    \centering
    \includegraphics[width=1\linewidth]{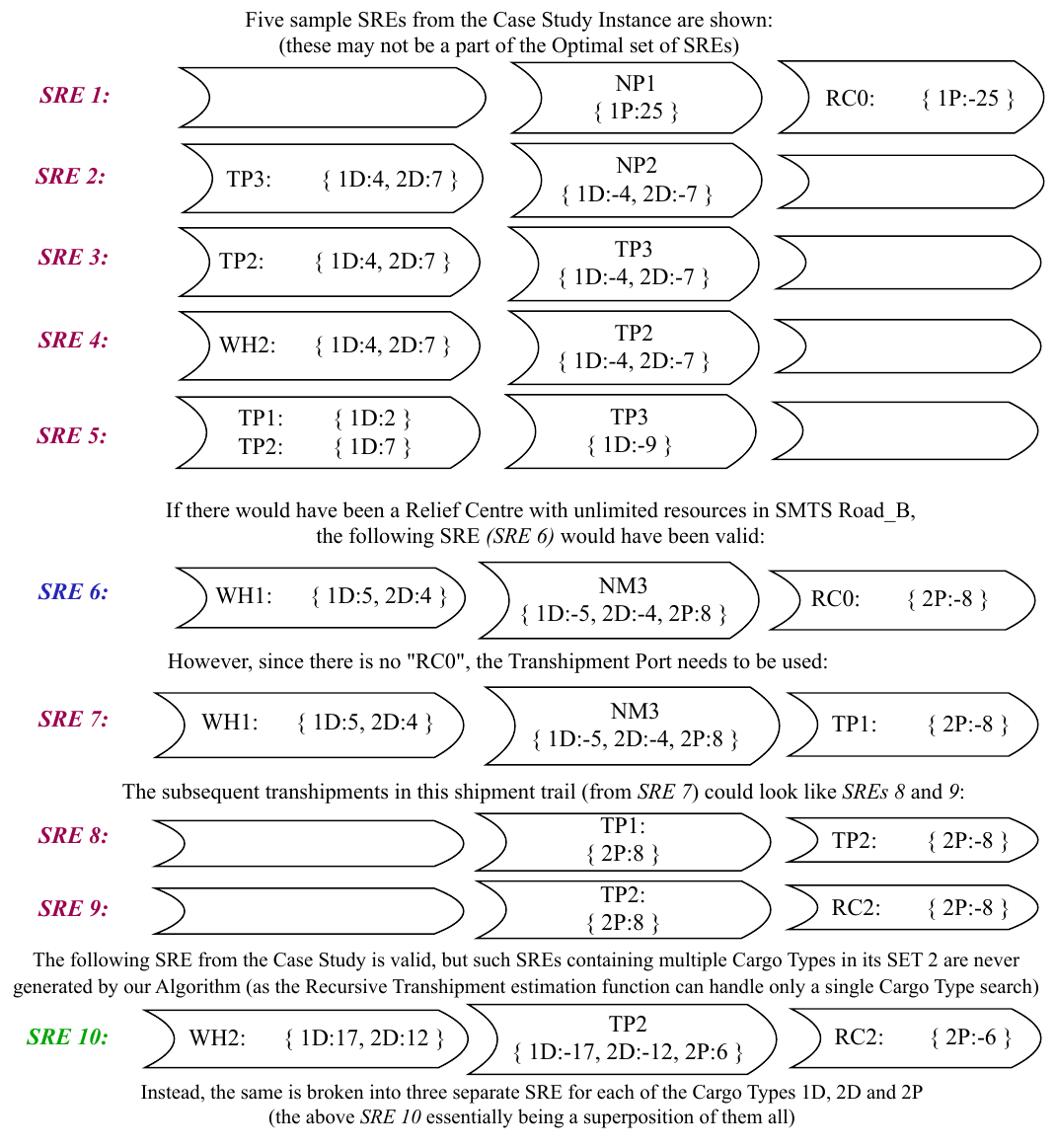}
    \caption{\textbf{Sample SREs (numbered examples) pertaining to our study-example instance in \autoref{An Example Case Study}}}
    \label{fig:SRE Examples}
\end{figure}

The term Vehicle Load Code (VLC) in \autoref{fig:Smallest Route Element depiction} refers to the amount of cargo (of each type) that was transferred from a vertex into the vehicle, in case some quantity has been transferred from inside the vehicle to the vertex (as is applicable at Relief Centres, or during some transhipments) the quantity would be negative. In simple terms, the Vehicle outline is imagined to have a Gaussian surface (closed smooth topology) with any resource quantity entering inside the Gaussian volume to be positive and negative values indicating exit of some resources. It may be observed, from the examples in \autoref{fig:SRE Examples}, that at no instant in time will any Cargo Quantity be present as a negative number within the vehicle, since that would indicate an exit of some non-existent resources from within the Gaussian enclosure.

The SRE is prepared for a specific Vehicle Type respecting the constraint of not exceeding that Vehicle Type's capacity. All SREs are subsequently allocated to a specific Vehicle (a Vehicle Depot and the unique vehicle type number in that Depot) in section \ref{Integration of SREs into vehicle routes}. A Split Node (or a Transhipment Port) can have multiple SREs to satisfy whereas a Simultaneous Node (as per the definition of Simultaneous Nodes, they must be satisfied fully using a single vehicular trip) can have only a single SRE. This can be seen in the example with the dummy Relief Centre ``RC0" in \autoref{fig:SRE Examples} (and the next corrected example with ``TP1" in place of the ``RC0"). In other words the SET 2 here (and for any Simultaneous Node for that matter) contains the full requirement code of the specific node (``NM3" in our example in the $SRE~6$) but in the form of the sign-corrected VLC.

The SET 2 in \autoref{fig:Smallest Route Element depiction} will contain only a single element (i.e. the prime vertex which is the focus of the route element) along with its VLC. The SET 1 and SET 3 may contain multiple elements; either of SET 1 or SET 3 may be empty but not both. The addition of all VLCs (i.e. separately adding the amount of transfers for each load type across all the vertices across an SRE's SETs) within an SRE must be 0, indicating that the load transferred from some vertex was delivered to another (in the direction of the arrow in \autoref{fig:Smallest Route Element depiction}). To be more clear: SET 1 transfers only Delivery load types to SET 2; and SET 2 transfers only Pickup load types to SET 3. 

In the subsequent sections, we discuss how to create SREs for the different Types of Nodes and TPs. The SETs 1 and 3 of an SRE can contain multiple vertex elements, whose relative positioning is perturbed using methods described in section \ref{Perturbations to find the best Route Cluster}, while maintaining causality constraints (developed during transhipments). Perturbations allow deeper exploration of neighboring solutions, and are made to respect the Causality constraints through the CD, and also the vehicle capacity constraints.

\subsubsection{Explaining the Causality Dict} \label{Explaining a Causality Dict}
All causality constraints are stored in the Causality Dict; whenever a TP is used within any of the SETs of an SRE, the details of this allocation is recorded in the CD. For example, if we assume all the (ten) numbered SREs from \autoref{fig:SRE Examples} are actually being generated sequentially one after another in some way during the program execution, the CD will be updated right after the SREs 2, 3, 4, 5, 7, 8, 9 and 10 (since the other SREs numbered 1 and 6 do not contain any TP anywhere within any SET) are generated.

The CD is structured as a Table with mainly two columns (apart from the first indexing column), the left one for representing superior causal events (\ie events which should take place before in time), and the right one for representing inferior causal events (\ie events which should take place after all the respective superior causal events have taken place). It is depicted having the following three column headings:\newline

\scriptsize
\noindent{
\begin{tabular}{|m{0.3\textwidth}|m{0.3\textwidth}|m{0.3\textwidth}|}
\hline
A unique ID corresponding to the TP-in-concern at which the Causal events are being reported & \textbf{LEFT Column: \newline Superior Causal Events} & \textbf{RIGHT Column: \newline Inferior Causal Events} \\ \hline
\end{tabular}
\newline}
\small

Two types of $SREs$ may be observed in \autoref{fig:SRE Examples}, SREs not containing any TP in their SET 2, and SREs that contain a TP in their SET 2. Whenever $SREs$ get generated, the TPs concerned within any of the SETs in that $SRE$ is searched, and the TPs found in the SETs 1 and 3 are put in separate rows of the CD as per the superior or inferior causal logic. All TPs encountered in the SET 1 are causally inferior, since to satisfy these allocations, further transhipments of Delivery resources will be necessary. However, if any TP is encountered in the SET 3, it is deemed to be causally superior as subsequent pickups from within the resources in the respective VLC associated with the encountered TP will be able to take place only afterwards (after the resources get deposited in the TP in the first place). We only look at the SETs 1 and 3 when any new SRE is generated as the SET 2 is developed internally (using the logic mentioned in section \ref{SRE creation for Transhipment Ports}) to continue a transhipment trail.

The usage of the CD with examples from \autoref{fig:SRE Examples} is explained step-by-step:
\begin{itemize}
    \item When the $SRE~2$ from \autoref{fig:SRE Examples} gets generated (using the function described in section \ref{SRE Generation for Split Nodes}), all TPs concerned within any of the SETs in that $SRE$ are searched for, in this case a single TP is present in its SET 1, which is observed to be causally inferior. Therefore, after the generation of $SRE~2$, the CD is updated as below (along with the respective $SRE$ source information):

\scriptsize
\noindent{
\begin{tabular}{|m{0.15\textwidth}|m{0.3\textwidth}|m{0.3\textwidth}|}
\hline
TP Name with a Unique ID & \textbf{LEFT Column: \newline Superior Causal Events} & \textbf{RIGHT Column: \newline Inferior Causal Events} \\ \hline
TP3\_A & & $SRE~2$\\ \hline
\end{tabular}
\newline}
\small

    If in SET 1 of the $SRE~2$, an additional TP would have been present (with its VLC) as say TP0: \{3D:3, 4D:1\}, this would have been assigned another fresh row in the CD as:

\scriptsize
\noindent{
\begin{tabular}{|m{0.15\textwidth}|m{0.3\textwidth}|m{0.3\textwidth}|}
\hline
TP Name with a Unique ID & \textbf{LEFT Column: \newline Superior Causal Events} & \textbf{RIGHT Column: \newline Inferior Causal Events} \\ \hline
TP3\_A & & $SRE~2$\\ \hline
TP0\_B & & $SRE~2$\\ \hline
\end{tabular}
\newline}
\small
    
    \item After every new inferior causal event addition (as unique rows) in the CD, its superior causal events are generated using the function described in section \ref{SRE creation for Transhipment Ports}, and inserted in the LEFT column of the corresponding row. For the $SRE~2$, the next continuation in the transhipment trail is shown as $SRE~3$ with the CD now being:

\scriptsize
\noindent{
\begin{tabular}{|m{0.15\textwidth}|m{0.3\textwidth}|m{0.3\textwidth}|}
\hline
TP Name with a Unique ID & \textbf{LEFT Column: \newline Superior Causal Events} & \textbf{RIGHT Column: \newline Inferior Causal Events} \\ \hline
TP3\_A & $SRE~3$ & $SRE~2$\\ \hline
\end{tabular}
\newline}
\small
    
    The meaning that may be inferred from this row of the CD is that at the TP3 (as seen from the TP Name in the index), the $SREs~2$ and $3$ are causally linked such that the superior causal event of $SRE 3$ should end before (in time) the inferior causal event of $SRE~2$.
    
    \item Right after the generation of $SRE~3$, similar search to detect the presence of any TPs in its SET 1 and 2 is conducted; the usage of TP2 is observed in its SET 1. The CD therefore gets updated to:

\scriptsize
\noindent{
\begin{tabular}{|m{0.15\textwidth}|m{0.3\textwidth}|m{0.3\textwidth}|}
\hline
TP Name with a Unique ID & \textbf{LEFT Column: \newline Superior Causal Events} & \textbf{RIGHT Column: \newline Inferior Causal Events} \\ \hline
TP3\_A & $SRE~3$ & $SRE~2$\\ \hline
TP2\_C & & $SRE~3$\\ \hline
\end{tabular}
\newline}
\small

    \item The function in section \ref{SRE creation for Transhipment Ports} gets triggered to satisfy the further superior causal elements in this transhipment trail. This results in the generation of the $SRE~4$, and further updating the CD as below:
    
\scriptsize
\noindent{
\begin{tabular}{|m{0.15\textwidth}|m{0.3\textwidth}|m{0.3\textwidth}|}
\hline
TP Name with a Unique ID & \textbf{LEFT Column: \newline Superior Causal Events} & \textbf{RIGHT Column: \newline Inferior Causal Events} \\ \hline
TP3\_A & $SRE~3$ & $SRE~2$\\ \hline
TP2\_C & $SRE~4$ & $SRE~3$\\ \hline
\end{tabular}
\newline}
\small

    \item This ends this transhipment trail, since $SRE~4$ contains only the PRVs in its SET 1 and 3. The final meaning of this CD can be fathomed by subtly imagining the relative causal relations obtained between the $SREs~4, 3$ and $2$; the causal superiority diminishes across them as they transport the same resources from a PRV (in $SRE~4$) to a node (in $SRE~2$). These SREs may (and will) be assigned to different vehicles (also of different types, in this case), and keeping these records using the novel devised method of Causality Dict ensures that the temporal causality is maintained during all transhipments. The CD is also used to calculate the waiting time (if any) of vehicles at TPs (separately during each of their unique visits, incase the same TP is visited by a vehicle multiple times), since it may be better for a vehicle to wait for a short amount of time to pickup an incoming resource.
\end{itemize}

In a very similar fashion, when TPs are encountered in the SET 3 of an SRE, the SRE is added to the RIGHT causally superior column in the CD, and the same function described in section \ref{SRE creation for Transhipment Ports} is used to develop the subsequent SREs extending the transhipment trail till the resources reach a PRV. The evolution of the CD in this case can be seen by taking the example of $SRE~7$ from \autoref{fig:SRE Examples}:

\begin{itemize}
    \item After the $SRE~7$ is generated (using functional methods from section \ref{SRE Generation for Simultaneous Nodes with or without transhipments}), it is scanned to search for the existence of any TP in its SETs 1 and 3; after the TP1 is found in SET 3, the CD is updated to (only the rows relevant to this transhipment trail are shown):

\scriptsize
\noindent{
\begin{tabular}{|m{0.15\textwidth}|m{0.3\textwidth}|m{0.3\textwidth}|}
\hline
TP Name with a Unique ID & \textbf{LEFT Column: \newline Superior Causal Events} & \textbf{RIGHT Column: \newline Inferior Causal Events} \\ \hline
TP1\_D & $SRE~7$ & \\ \hline
\end{tabular}
\newline}
\small

    \item Next, using the function from section \ref{SRE creation for Transhipment Ports}, the subsequent SRE is generated to resolve the resource commitments done to the TP1 in $SRE~7$. This is shown as the $SRE~8$; the CD now gets updated to:

\scriptsize
\noindent{
\begin{tabular}{|m{0.15\textwidth}|m{0.3\textwidth}|m{0.3\textwidth}|}
\hline
TP Name with a Unique ID & \textbf{LEFT Column: \newline Superior Causal Events} & \textbf{RIGHT Column: \newline Inferior Causal Events} \\ \hline
TP1\_D & $SRE~7$ & $SRE~8$ \\ \hline
\end{tabular}
\newline}
\small
     \item The newly generated SRE is scanned and here too, a TP is encountered in the SET 3; this results in the CD to get updated to:

\scriptsize
\noindent{
\begin{tabular}{|m{0.15\textwidth}|m{0.3\textwidth}|m{0.3\textwidth}|}
\hline
TP Name with a Unique ID & \textbf{LEFT Column: \newline Superior Causal Events} & \textbf{RIGHT Column: \newline Inferior Causal Events} \\ \hline
TP1\_D & $SRE~7$ & $SRE~8$ \\ \hline
TP2\_E & $SRE~8$ & \\ \hline
\end{tabular}
\newline}
\small

    Observe that the index of the second row above contains the same TP which is encountered in the SET 3 (with an additional alphabet-tag to keep all row IDs unique).
     
     \item The creation of a new row in the CD triggers the function in section \ref{SRE creation for Transhipment Ports} to resolve the latest resource commitments. It generates $SRE~9$, updating the CD to:
     
\scriptsize
\noindent{
\begin{tabular}{|m{0.15\textwidth}|m{0.3\textwidth}|m{0.3\textwidth}|}
\hline
TP Name with a Unique ID & \textbf{LEFT Column: \newline Superior Causal Events} & \textbf{RIGHT Column: \newline Inferior Causal Events} \\ \hline
TP1\_D & $SRE~7$ & $SRE~8$ \\ \hline
TP2\_E & $SRE~8$ & $SRE~9$ \\ \hline
\end{tabular}
\newline}
\small

    This marks the end of this transhipment trail, as all vertices in the SETs 1 and 3 of latest SRE generated (\ie $SRE~9$) are PRVs. This CD helps imagine the resource flow from the relatively temporally-supreme causal SRE of $SRE~7$ downstream along the transhipment trail, through the relatively causally-poorer SRE of $SRE~8$, to final resource destination (a PRV).
\end{itemize}

In all the above examples, each individual cell contains a single SRE, however, it may so happen during the SRE generation process that the SETs 1 or 3 (or both) may contain multiple TPs (interestingly the same TP could also be in both SETs, and so we uniquely number-tag all TPs used during SRE creations to enable the concrete creation of corresponding rows in the CD). To handle all these cases, the following sequence of steps is followed:
\noindent\hrulefill

\noindent\textbf{Pseudocode to fathom each Transhipment Trail}

\noindent\hrulefill
\begin{algorithmic}[1]
% \begin{algorithmic}
% \Require $n \geq 0$
% \Ensure $y = x^n$
\Search {for any TP across the SETs 1 and 3 in every newly created SRE}
    \For{each TP usage in these SETs:}
        \State A unique row is created in the CD (if the same TP is in SET 1 and 3, two different CD rows are created having different Unique ID tags). Each row is named starting with the concerned TP found.
        \For{each generated CD row:}
        \State The SRE ID is placed in the appropriate CD column (RIGHT or LEFT depending on the temporal causality, as discussed previously).
        \State The function in section \ref{SRE creation for Transhipment Ports} is triggered to continue the respective transhipment trail. This function generates new SREs to satisfy assigned resource commitments to TPs, which may be multiple in number.
        
        Each newly created SRE is placed in the opposite column (of the respective row) from which its generation was triggered.
        
        For each of these newly created SRE, this pseudo-code is recursively called from Step 1 (\textbf{Search}).
        \EndFor
    \EndFor
\EndSearch
\noindent\hrulefill
\end{algorithmic}

Elaborating on Step 6 above, all cells in the CD either have none or a single SRE ID assigned till before the first time Step 6 is executed. Through the Step 6, multiple SREs may be generated during multi-trips as well as normally when different vehicles are used to satisfy one segment of the transhipment trail. This can be imagined similar to the route structure of a tree (or in the delta of rivers), where single elements converge (or diverge). A Transhipment Trail dealing in PickUp Cargo Types start from the Node similar to the main river element, and stem out via different SREs (river outlets) distributing this pickup to PRVs or TPs, the pickups from TPs in this transhipment trail are further sent via further (increasingly causally inferior) SREs which may be compared to the final river-branchings, before it reaches the sea (\ie any of the PRVs, some outlets may reach oceans faster which other may branch multiple times more and reach different oceans). Similarly, the aggregation of delivery resources from PRVs (nutrients from the soil) occurs through the formation of different SREs (smallest roots), which transport these to some TPs after which they may be transported further upstream together (roots joining to form tertiary, secondary and primary structures) untill the delivery reaches the nodes (stems of plants). We term this phenomenon as the Transhipment Trail Stemming (TRS). We observe that each unique transhipment trail generated by our algorithm may be either comparable to the pickup distribution (river delta stemming outwards always) or delivery aggregation (plant root network stemming inwards always), but not both. This is because of the way the Transhipment Satisfaction function has been designed to handle one Cargo Type at a time. The only place where the TRSs of different delivery and pickup types of transhipment trails may meet are at the Nodes. The TP aggregation points in the Delivery TRSs, and the TP distribution points in the PickUp TRSs may be handling multiple types of resources, but of the same Delivery or PickUp type.
At the end of Step 6, after the placement of the newly generated SREs in the opposite columns (of the respective rows) some cells in the CD may contain multiple elements, as in the hypothetical case mentioned in \autoref{tab: Hypothetical Case for Waiting Time Image}.

\begin{table}[htb] % Preferable here, or at the top of this page, or at the bottom
\centering
\caption{A few rows of a hypothetical Causality Dict\label{tab: Hypothetical Case for Waiting Time Image}}
% \noindent
\resizebox{0.75\textwidth}{!}{ % Resize to fit the width of the page
\begin{tabular}{|m{0.15\textwidth}|m{0.3\textwidth}|m{0.3\textwidth}|}
\hline
TP Name with a Unique ID & \textbf{LEFT Column: \newline Superior Causal Events} & \textbf{RIGHT Column: \newline Inferior Causal Events} \\ \hline
TP19\_ZA & $SRE~45, SRE~51, SRE~31$ & $SRE~74$ \\ \hline
TP25\_ZB & $SRE~94$ & $SRE~62, SRE~88, SRE~18$ \\ \hline
\multicolumn{3}{c}{where some assumed details of each SRE could be:} \\ % Row without vertical lines
\hline
SRE ID & \textbf{Vehicle Type} & \textbf{Possible set of Vehicle Depots to be allocated from} \\ \hline
\hline
$SRE~45$ & VT7 & VD16, VD18 \\ \hline
$SRE~51$ & VT8 & VD15 \\ \hline
$SRE~31$ & VT8 & VD15 \\ \hline
$SRE~74$ & VT9 & VD17, VD19 \\ \hline
\hline
$SRE~94$ & VT9 & VD17, VD19 \\ \hline
$SRE~62$ & VT8 & VD15 \\ \hline
$SRE~88$ & VT7 & VD16, VD18 \\ \hline
$SRE~18$ & VT8 & VD15 \\ \hline
\end{tabular}
} % For Resize Box Closure
\end{table}

As described previously through the TRS, our algorithm design prevents any row in the CD to have more than one SRE in both columns simultaneously, as this would mean aggregation and branching happening in the same transhipment trail. If such a requirement is necessary, separate transhipment requests are generated forming different trails. For example, consider a hypothetical case where two nodes NP8 and NP9 have delivery demands of water bottle packs in the respective quantities of 18 and 19, and there is no PRV present in their SMTS named SMTS\_X; only a single transhipment port named TP72 is situated in their SMTS. This TP72 also has access to another SMTS named SMTS\_Y containing only two warehouses, W6 and W8, both containing 20 units of water. 
The considered MMTN consists only of two SMTS having different modes of transport, and each SMTS has 2 vehicles of the same type stationed at a single Vehicle Depot in both the SMTSs.

Two vehicles are available in SMTS\_Y, and so it would (generally) be pertinent to send each one from the single depot and fetch resources from either of the Warehouses to dump them at the TP72. Meanwhile, the two vehicles in SMTS\_Y would be journeying to TP72 and wait untill the resources arrive, only to leave (one after the other) for either of the Nodes to complete the final delivery. In this case, if we focus on the happenings at TP72, we find that resources aggregated from two warehouses get distributed to the nodes; however in our implementation, each request from an individual node will constitute as a separate independant transhipment request initiating a transhipment trail. So our algorithm is able to completely handle all such cases, and much more complex ones (as may be seen in the datasets); by imagining the superposition of the individual TRSs will one be able to picture the aggregation and the distribution together. Further, incase the nodes had requests of other types of delivery, these would constitute as other independant transhipment trails for subsequent SRE generations.

\subsubsection{SRE creation for Split Nodes} \label{SRE Generation for Split Nodes}

This function accepts a subset of unsatisfied Split Nodes, and generated SREs to satisfy them (to the fullest extent possible). In case for some nodes, transhipment is deemed necessary, the function described in section \ref{SRE Generation for Split Nodes with Mandatory Transhipment} is triggered. The Algorithm proceeds through the following steps:

\noindent\hrulefill

\noindent\textbf{Pseudocode to generate SREs for Split Nodes}

\noindent\hrulefill
\begin{algorithmic}[1]
\If{any remaining Unsatisfied Split Nodes, within the set of Nodes passed to the function, are left to be satisfied:}
    \State $send\_all\_remaining\_for\_Mandatory\_Transhipment \gets True$
        \For{each of the Twigs in every DTS, with the Split Nodes as Trunks, iterate in the decreasing order of preference scores, irrespective of the Trunk the Twigs belong to (scores already generated using section \ref{Preference Generator Function for Split Node}):}	\Comment{This means that the first iteration could be pointing to a Twig in some DTS and the other to another DTS originating from a different source Split Node as the Trunk, and potentially visiting newer Twigs from previously iterated DTSs.}
            \If{the encountered Twig's original Trunk is within the current set of Unsatisfied Split Nodes, and if the Twig itself refers to a PRV (\ie is a Warehouse or a Relief Centre):}
                \State $Trips \gets 0$ \Comment{Zero ensures atleast 1 trip with the encountered Twig}
                \While{$e^{-Trips \times Multi-Tripping~HyperParameter~value~from~this~Twig} > rand_0^1$} \Comment{All Twigs of DTSs is associated with a Multi-Trip-HyperParameter varying between 0 and 1; $rand_a^b$ gives a random number from a uniform distribution ranging from $a$ to $b$}
                    
                    \State The actual SRE generation process starts \Comment{Detailed below}

                    \State Current resource availability at the PRV visited in the SRE is adjusted (decremented)
                    \State $send\_all\_remaining\_for\_Mandatory\_Transhipment \gets False$ \Comment {Since there has been an SRE creation in this iteration across the Unsatisfied Split Nodes}
                    
                    \If{the Twig still has resources left to be leveraged, and the respective Trunk (being already in the list of Unsatisfied Split Nodes) still has unsatisfied requirements, that may be satisfied by the PRV Twig through subsequent transfers:}
                        \State $Trips \gets Trips+1$
                        \State \Continue  \Comment{starts next iteration after checking the condition in step 6}

                        \ElsIf{if the Node is fully satisfied}
                            \State remove it from the list of Unsatisfied Nodes
                            \State \Break  \Comment{Exit this loop to proceed with the iteration of Step 3}
                        
                    \EndIf

		    \EndWhile
		\EndIf

            \If{$send\_all\_remaining\_for\_Mandatory\_Transhipment=True$}
                \State all the leftover (partial or full) unsatisfied requirements are set to be satisfied using the functional methods described in section \ref{SRE Generation for Split Nodes with Mandatory Transhipment} through transhipments
            \Else
                \State Execution starts from step 1 \Comment{Process flow goes back to step 1}
            \EndIf
            
	\EndFor

\ElsIf{no remaining Split Nodes are left to be satisfied}
	\State return the execution flow ending this function
\EndIf
\noindent\hrulefill
\end{algorithmic}

\textbf{Detailing the SRE generation process from Step 7}:

The details of the resources available at the PRV Twig is matched with the resources requested by the Split Node (Unsatisfied Trunk). The Vehicle Type to be considered for this SRE creation is the same as the DTS Branch connecting the Trunk (Unsatisfied Node) with the Twig (PRV). Respecting both the Vehicle Volume and Weight constraints, as much possible of the requested requirement is tried to be satisfied. This already accounts for the Vehicle-Cargo compatibility, since any incompatible combinations would not contain the subsequent Leaves of that Cargo Type in this DTS.

The resource allocation will therefore hault, only if:
    \begin{itemize}
        \item There is no more requirement left to be satisfied at the Node, or,
        \item Any leftover requirement needs to be allocated to another Cargo-compatible VT, or,
        \item Some of the requested CTs are unavailable/partially available at the PRV, or,
        \item The vehicle capacity has been reached.
    \end{itemize}
In any case, provided that some resources were allocated for a trip, the actual SRE set creation would be done as follows:
    \begin{enumerate}
        \item An empty SRE with three SETs is created  
        \item VLC creation: The allocated Cargo Types along with their respective quantities with the appropriate signage (indicating inflow or outflow from the Vehicular Gaussian Surface) is prepared.
        \item The (Unsatisfied) Split Node (\ie the Trunk of this DTS) is assigned to the VLC and stored in the SET 2.
        \item Similarly, the VLC for the Twig is prepared. It may be noted that since a PRV may contain either Delivery Load Types, or have space to accommodate Pickup Load Types, only one among the SETs 1 or 3 will be populated with the VLC. Appropriate signage in the VLC is ensured (essentially opposite to that used in SET 2 creation).
        \item If the PRV Twig is a Warehouse, the PRV and VLC combo is assigned within the SET 1; if the PRV Twig is a Relief Centre, the PRV and VLC combo is assigned within the SET 3.
        \item Any SRE is incomplete without its associated VT and VD information. In this case, the Vehicle Type information, along with the possible set of VDs, is inferred from the DTS Branch connecting the Twig and the Trunk. These respective information is tagged along with every SRE.
        \item The SRE is provided with a Unique ID, and appended into the list of all SREs.
    \end{enumerate}

\subsubsection{SRE creation to satisfy Mandatory Transhipment for some Split Nodes, as necessary} \label{SRE Generation for Split Nodes with Mandatory Transhipment}

This function satisfies the Mandatory Transhipment for some of the Split Nodes which were already identified during preference generation, as well as other Split Nodes which may have been later identified to require transhipment (through the processes in section \ref{SRE Generation for Split Nodes}). Whenever this function is triggered, it tries to fulfill all mandatory transhipments necessary, to or from Split Nodes; the algorithm is designed as follows:

\noindent\hrulefill

\noindent\textbf{Pseudocode to generate SREs for Split Nodes requiring transhipment}

\noindent\hrulefill
\begin{algorithmic}[1]
\Ensure HyperParameter of Transhipment\_Trip\_Setting takes the True value with an 80\% probability before every Main Iteration gets initiated
\For{each request present within the list of all currently generated Mandatory Transhipment requirements for Split Nodes:}
    \State Transhipment requests of different CTs to/from the same Node is spread separately within the Mandatory Transhipment requirement List. A single request consists of three different information:
    \begin{enumerate}
        \item the Split Node which requires transhipment,
        \item a single Cargo Type which is to be transhipped, and
        \item the Quantity of that Cargo Type requested.
    \end{enumerate}
    From the DTS containing the Split Node at the Trunk, all non-PRV Twigs (\ie Twigs which correspond to TPs) are sorted in the decreasing order of their scores; this is termed as Filtered DTS-TP
    \While{the current remaining amount of resources requested is non-zero}
        \For{each TP within Filtered DTS-TP}
            \If{the requested Cargo is obtainable via the concerned TP}
                \If{Current\_Estimate = 0} \Comment{computation described below}
                    \State \Continue \Comment{proceeds to check the next lower-scored candidate within the Filtered DTS-TP}
    		\Else
    		      \State $Trips \gets 0$
                    \While{$e^{-Trips \times Multi-Tripping~HyperParameter~value~of~Twig} \geq rand_0^1$}
                        \State The maximum amount of resource transfer possible is calculated by finding the minimum among:
    \begin{itemize}
        \item Current\_Estimate,
        \item Maximum amount of the Cargo Type (from step 2.ii.) that can be loaded into the Vehicle Type (that connects the Twig with the Trunk), respecting the Volume and Weight capacities of the vehicle
        \item Quantity of the resource requested (for transhipment, in step 2.iii.)
    \end{itemize}
    and is dubbed Amt\_Sat (the amount satisfiable through this SRE).
                        \State An empty SRE is created; the VD information associated with the SRE is taken from the same VT-Branch connecting the Trunk and the Twig.
                        \State An appropriate signage-adjusted VLC is generated (using the Amt\_Sat value) to indicate the resource transfer at the Split Node. This is stored in the SET 2 of the fresh SRE.
                        \State A similar VLC (opposite in its signage \wrt the previous) is developed and assigned to the TP-Twig (\ie the current iterator of Iteration 4); this is stored appropriately in either the SETs 1 or 3 (depending on if the Cargo Type is of Delivery or PickUp respectively).
                        \State Update (decrease) the resource requirements at the Split Node for the concerned CT by the Amt\_Sat value, as some (or all, depending on the Amt\_Sat) of it is now shifted to the TP through this SRE creation.
                        \State A new row in the CD is generated with this SRE ID placed in the appropriate column (RIGHT in case the CT from Step 2.ii. is of Delivery Type, and LEFT in case it is a PCT).
                        \State To satisfy the resources diverted to/from the TP-Twig, the function SRE creation for TPs (section \ref{SRE creation for Transhipment Ports} is triggered; this automatically takes care of all subsequent transhipment requirements generated from this transhipment trail. \Comment{This function is automatically triggered whenever a new row in the CD is added}
                        \State Decrease both the transhipment requirement amount from step 2.iii., and the Current\_Estimate by value of Amt\_Sat, and find their minimum (say this minimum value is $m$)
                        \If{$m \geq 0$ \textbf{and} Transhipment\_Trip\_Setting\_HyperParameter=True}
                            \State $Trips \gets Trips+1$
                            \State \Continue \Comment{rechecks the \textbf{while} condition in step 10}
                        \ElsIf{$m=0$}
                            \State \Break \Comment{starts the next iteration of step 4}
                        \EndIf
                    \EndWhile
    
    		\EndIf
            \ElsIf{the concerned TP is unable to leverage any amount of the resource type (or if it is found incompatible to transfer the cargo type)}
                \State \Continue \Comment{proceeds to check the next lower-scored candidate within the Filtered DTS-TP}
            
            \EndIf
        \EndFor

        \If{no SRE was created in the latest full run of the \textbf{for} loop in step 4}
            \State this transhipment requirement is deemed unsatisfiable: the details are stored within Unsatisfiable\_Portion, and the next Main Iteration is initiated
        \EndIf

        \If{Transhipment\_Trip\_Setting HyperParameter=False}
            \State \Break \Comment{flow continues with iterations of the \textbf{for} loop from step 1}
        \EndIf
        
    \EndWhile	
\EndFor
\noindent\hrulefill
\end{algorithmic}

To ensure a TP can satisfy some amount of the requested CT (elaborating the step 5 above) the Estimated Compare Code of the iterated TP-Twig in its DTS found during the preference generation of Split Nodes in section \ref{Preference Generator Function for Split Node} is looked into. If this Estimated Compare Code contains the requested CT in some non-zero quantity, it would indicate that if this TP-Twig is committed with some amount of that resource, it will be able to leverage the resource upto the quantity in the Estimated Compare Code. 

As the Heuristic progresses and transhipment requests are generated, a lot of the Estimated Compare Codes describing the resource leveraging capability gets redundant. This happens because some of the resources (for delivery, or their accepting capability for pickup) in some PRV may have decreased, as they may have been sent away (for delivery), or accepted (for pickups, thus decreasing the accepting capacity of the RC) by the PRV. We don't create fresh scores for the Twig after every such allocations (as this would be more computationally intensive) because then all scores across all DTSs involving that specific PRV  directly or indirectly (through TPs) would need to be updated. Instead we focus on getting a new, accurate and current estimate of the resource leveraging capability; for PRVs this is directly inferred from the current amount of resource they have, but for the TPs, this is made possible by referring to the Component\_2 (developed earlier using the function in section \ref{Recursive function to enable Transhipment of a single load-type}) stored in the DTS; this Component\_2 contains the details of the PRVs from/to which resources may be leveraged, and so we just find the current available quantity of resource (of the Cargo Type under consideration from Iteration 1) at those PRVs and sum them together creating a Current\_Estimate.

Expanding on the step 32 above, this scenario never arises unless the problem itself is infeasible, essentially requiring transhipment when no TP is available in that SMTS. Complex cases containing disconnected SMTSs forming multiple independent MMTNs are deemed to be separate problems, and may also be solved together using our Heuristic or MILP Formulation as a single problem.

\subsubsection{SRE creation for Simultaneous Nodes (with or without transhipments)} \label{SRE Generation for Simultaneous Nodes with or without transhipments}

This function accepts a set of Simultaneous Nodes that are yet to be satisfied, and follows the following logical steps:

\noindent\hrulefill

\noindent\textbf{Pseudocode to generate SREs for Simultaneous Nodes}

\noindent\hrulefill
\begin{algorithmic}[1]

\Function{SRE\_creation\_for\_Simultaneous\_Nodes}{Unsatisfied\_$N^M$}
	\For{each DTS-Trunk, where the Simultaneous Node-Trunk is present within Unsatisfied\_$N^M$:}
		\For{each VT-Branch (of the DTS considered in the outer loop)}
			\State compare the Requirement Code of the Simultaneous Node-Trunk, with the aggregate of all the resources available across all PRV-Twigs connected to this VT-Branch \Comment{This is done by considering only the CTs within their respective Leaves since the DTS preserves the compatibilities (thus ensuring we are not considering any resource transfer onboard an incompatible vehicle)}
			\If{it is found that all the combined PRV-Twigs are unable to satisfy the requirement fully}
				\State TP-Twigs are also considered to satisfy the requested requirement \Comment{Elaborated below}
				\If{it is found that this Branch will be able to satisfy the Simultaneous Node fully, either using only PRV-Twings or also using TP-Twigs}
					\State Create an empty SRE
					\State Modify the Requirement Code to develop a VLC with the appropriate signage to be assigned to the Simultaneous Node; this Simultaneous Node-VLC combo is stored in the SET 2 of the SRE
					\State Relieve resource requirements across all CTs \Comment{Detailed below}
				\Else
					\State this portion of the requested requirement is deemed unsatisfiable (details are placed in Unsatisfiable\_Portion): As per our algorithm design, this could mean infeasibility in the problem; in rare cases it could also indicate bad HyperParameter combination of the current Main Iteration, which are revised while starting the next Main Iteration immediately.
				\EndIf
			\EndIf
		\EndFor
	\EndFor
\EndFunction
\noindent\hrulefill
\end{algorithmic}

Elaborating the step 6 above, for each TP-Twig, we have the details of the PRVs that it can leverage resources to/from through its stored Component\_2 (developed earlier using section \ref{Recursive function to enable Transhipment of a single load-type}). We aggregate all the Component\_2s (from the existing CT Leaves for a TP-Twig), across all the TP-Twigs connected with the iterator 2's VT-Branch; then we sum up the resources available to be leveraged (for each CT individually) from within each Component\_2's PRV set, across the aggregated Component\_2s. During this process we cautiously ensure that no PRV-resource of the same CT gets summed up multiple times (as the same CT-Leafs from different Twigs can share some of the same PRVs within their Component\_2). This gives us the CT-specific leveraging potential of the entire VT-Branch and allows us to compare the Requirement Code to check if it can be fully satisfied.

To maximally relieve resource requirements across all the CTs (detailing the step 10 above), we iterate through the PRV-Twigs of the VT-Branch from Iteration 2 in descending preferential order (section \ref{Preference Generator Function for Simultaneous Node}). We don't need to consider the VT capacity issues in this case as they have already been considered when this VT-Branch was deemed viable to satisfy this full requirement simultaneously. All these assignments are separately encoded to develop appropriately signed VLCs and stored in the appropriate SETs (SET 1 or 3 of the SRE) associated with their concerned PRV. If it was earlier found (in step 5) that the requirement satisfaction would necessitate transhipment(s) as well, we iterate through the TP-Twigs (in the decreasing preferential order). For each TP-Twig:
    \begin{itemize}
    \item We calculate the maximum amount of all resources that may be leveraged utilizing the Component\_2 PRV information from across each CT-Leaf associated with this specific TP-Twig.
    \item We develop separate VLC's for the PickUp and Delivery CTs (this was automatically ensured when iterating through the PRV-Twigs as each PRV is specific to either PCTs or DTCs, but for TPs, we need to maintain then separately).
    \item The VLC(s) developed are assigned (with appropriate signs representing load transfer into/out of a vehicle) to the appropriate SETs (SET 1 for VLC with only DCTs, and SET 3 for the VLC with only PCTs) bearing association to the TP-Twig. This is the only place in our entire algorithm design where both the SETs 1 and 3 may be filled up (tweaking this aspect is another broad area with future research scope).
    \item Once the SRE development leveraging a TP-Twig to satisfy some/all of the requirement is completed, the SRE is assigned a unique ID, associated with the VD information obtained from the corresponding VT-Branch.
    \item A new row in the CD is created, containing the TP-Twig name (tagged with an unique combination of alphabets) as the row ID; the generated SRE ID is place in an appropriate column in that row.
    \item The function to generate SREs, from TPs as DTS-Trunks (section \ref{SRE creation for Transhipment Ports}), gets triggered with the information of the newly generated CD-row, to satisfy the subsequent transhipments from the TRS originating from this new CD-row.
    \end{itemize}

A future investigation area could be the utilization of Branch scores (developed as an aggregation of the preferences of the connecting Twigs) for choosing the VT. Our design considers the first VT-Branch that can fully satisfy the Simultaneous Node-Trunk as the Branch scoring approach could be computationally intensive, requiring one pass of this entire function to calculate the current leveraging potential for each Branch, score then, and then perform the exact task done by this function (section \ref{SRE Generation for Simultaneous Nodes with or without transhipments}) in a second pass.

% \subsubsection{SRE creation for Simultaneous Nodes that cannot be satisfied fully using one Vehicular Branch of the Decision Tree} \label{SRE Generation for Split Nodes with Mandatory Transhipment}

\subsubsection{SRE creation for Transhipment Ports} \label{SRE creation for Transhipment Ports}

This recursive function is designed to be triggered anytime during a Main Iteration to satisfy any transhipment requirement generated through the addition of a new row in the CD.

This function accepts the relevant details of the TP name (from the passed row-Unique ID of the newly generated CD-row), the CT requested, and its amount (depicted as Q in the pseudo-code below). The row-Unique ID helps to query the SRE directly and get its vehicle information; this helps in reducing the subsequent computational search by removing the SMTS (\ie all VTs can ply on it) onto which the SRE-Vehicle Type plys on. Thus we get a Tabu List of VTs. The function proceeds accordingly:

\noindent\hrulefill

\noindent\textbf{Pseudocode of Function to generate SREs for Transhipment Ports}

\noindent\hrulefill
\begin{algorithmic}[1]

\Function{SRE\_creation\_for\_Transhipment\_Ports}{TP\_Name,CT,Q}
	\State $DeepTranshipment \gets False$
	\State $Unrestricted\_Amt\_of\_DeepTranshipment \gets False$
	\While{Q>0}
            \State \textbf{Initialize} $ScoreList \gets [~]$ \Comment{Empty list to store scores}
		\For{each Vehicle Type-Branch in the DTS with Trunk as TP\_Name; the Vehicle Type not being in the Tabu List}
			\For{each Twig connected to a Branch of the outer iterator}
				\If{$DeepTranshipment = False$}
					\If{the iterated Twig is a TP}
						\State \Continue \Comment{we leave the TPs for later and try to satisfy the current resource requirements only with PRVs (\ie without deeper transhipments right away)}
					\EndIf
				\EndIf
				\For{the specific Cargo Type-Leaf connected to the outer-iteration Twig which is the same as the CT in the function definition} \Comment{relevant compatibilities (Cargo Type -vs- TP) are taken care of by the DTS, allowing the algorithm to always search within the feasible region}
					\State Amt\_Sat is calculated
					\State $f=Amt\_Sat/Q$ \Comment{fraction of the amount satisfiable}
					\State $score = mf^p$ \Comment{Exponent $p$ and Multiplier $m$ are HyperParameters of the concerned Cargo Type-Leaf (discussed in section \ref{Compare Code Calculator})}
					\If{HyperParameter of Time consideration during Scoring is ON}
						\State $score = score/T$ \Comment{Here $T$ is the time taken to travel from the Trunk to the Twig (or vice versa) using the Branch-VT}
					\EndIf

					\If{the outer iterator Twig is a TP}
						\State{$score = score - Adjusted\_Degree \times DegreeScoreReduction$}
					\EndIf
					\State append $score$ to $ScoreList$
				\EndFor
			\EndFor
		\EndFor

		\If{if no positive Amt\_Sat was generated within any of the iterations within step 6}
			\If{$DeepTranshipment = False$}
				\State $DeepTranshipment \gets True$
				\State \Continue {execution initiates from step 4}
			\ElsIf{$DeepTranshipment=True$}
				\If{$Unrestricted\_Amt\_of\_DeepTranshipment = False$}
					\State $Unrestricted\_Amt\_of\_DeepTranshipment \gets True$
				\Else
					\State remaining transhipment amount yet to be satisfied is deemed to be an Unsatisfiable Portion within the entire problem: \Comment{We get Unsatisfiable Portions extremely rarely for very complex large instances, and that too only in very few of the Main Iterations (as this depends on the overall HyperParameter value combinations). Any addition in the Unsatisfiable Portion of the Problem indicates a need to shift to the next Main Iteration, as the current HyperParameter combinations are deemed unsuitable.}
				\EndIf
			\EndIf

		\EndIf \Comment{At this stage, we should have some positive Amt\_Sats}
		\State take the element with the highest score from the ScoreList and create a new SRE with it
		\State update the CD-row with this newly generated SRE appended in the opposite column
		\State update the requirement Q
	\EndWhile
\EndFunction
\noindent\hrulefill
\end{algorithmic}

Expanding on the calculation of the Amt\_Sat from step 14, for PRV-Twigs, the Amt\_Sat takes the minimum value from the following:
    \begin{itemize}
        \item requested amount of resource Q
        \item maximum quantity of that Cargo Type which can fit within the Vehicle's Volume, and Weight limitations,
        \item current amount of the resource available of specific Cargo Type in the PRV
    \end{itemize}
    For TP-Twigs, the Amt\_Sat takes the minimum value from among the following:
    \begin{itemize}
        \item requested transhipment amount Q
        \item maximum quantity of resource that can fit within the Vehicle Types, weight and volume restrictions
        \item current estimate of resource quantity possible to be levaraged by the TP-Twig (we consider this only if Unrestricted\_Amt\_of\_DeepTranshipment is set as False): We calculate this by summing the appropriate available resources at the PRVs associated from the Component\_2 stored in the DTS-Leaf.
    \end{itemize}

As previously mentioned, we consider the degrees as potentials which are incorporated in this scoring (if the Twig is a TP). The previously calculated Adjusted Degree (which is a cumulative measure of the degrees for the specific TP-Twig) is taken into consideration here and multiplied by a HyperParameter DegreeScoreReduction; this product is subtracted from the Score (see step 21 above). The HyperParameter of DegreeScoreReduction is varied after every Main Iteration and takes values between -2 and 8 from within a varying distribution (we allow negative values allowing the Heuristic to experiment with reversed potentials, \ie deeper transhipments then will be better than no or single-degree transhipment).

Elaborating on the SRE creation from step 39, the SET 2 of this SRE contains the VLC with the appropriate signage with the vertex assigned as the TP which was passed into the function definition along with its requirement. The VLC contains the maximum quantity of resource satisfiable which is the Amt\_Sat associated with the highest Score.
We create another VLC with the appropriate signage (opposite to the one assigned in SET 2) and assign this to the Vertex-Twig that had the highest score; this Vertex may be a PRV or a TP. This Vertex-VLC association is sent to either the SET 1 or 3 depending on the logical sequence of transhipment; as a short explanation/trick, whenever the Cargo is of Delivery Type, it is assigned to SET 1, otherwise it is assigned to SET 3.

The appropriate signage (as previously elaborated) always ensures resources exit from the vehicle only after they enter it at a previously visited vertex. So, the sign associated within the VLC in SET 2 is:
    \begin{enumerate}
        \item POSITIVE when Cargo is PickUp Type
        \item NEGATIVE when the Cargo is of Delivery Type
    \end{enumerate}
    and likewise opposite for the other VLC placed in either SET 1 or 3. The VT information about this SRE, and its possible choices of VDs is inferred from the VT-Branch step 6.

Elaborating on the step 41:
\begin{itemize}
    \item If the Vertex-Twig is a PRV:\newline
    We decrease (update) the resource availability of that PRV, since some of that reserve was used up as is represented in value by the Amt\_Sat of the highest Score (this may not necessarily be the highest among all Amt\_Sat due to the scoring process followed, an area of interest ripe with future research potential).
    \item If the Vertex-Twig is a TP:\newline
    We introduce a new row in the CD (\ie place the newly created SRE in the appropriate column of a new CD-row, with its row ID being the Vertex-TP tagged along with unique alphabet combination). We immediately call this same function of SRE creation for TPs (section \ref{SRE creation for Transhipment Ports}) starting a recursion chain to satisfy the TRS; as per our design, the TRS is satisfied similar to depth-first approach as it stems deeper. The details for the recursion function start are derived from the newly generated CD-row (our implementation of this function does not have any functional \textit{return}).
\end{itemize}
After step 41, in case when re-initiation of the first iteration happens, more than one SRE ID would be put in the same cell as shown in \autoref{tab: Hypothetical Case for Waiting Time Image}. It may be kept in mind that as per our algorithm design, multiple SRE insertions within the same cell in the CD can happen only through this function.

\subsection{Allocation of Vehicles enabling Integration of SREs within Vehicle routes} \label{Integration of SREs into vehicle routes}

Each SRE, as created and stored within the Global SRE Pool, also contains information on its vehicle type and additionally, the possible Vehicle Depots whose Vehicles of the concerned type can allow the allocation of this SRE (in its Route Portion form) within their routes. The additional information about the Vehicle Depots is necessary to ensure routes developed within an SMTS are not allocated to another SMTS.

\subsubsection{Selection Logic for the Allocation of individual SREs to unique vehicles based on the CD} \label{Selection Logic for the Allocation of individual SREs to unique vehicles}

Once we have generated SREs for this PASS within the Main Iteration (PASS indicates the Node Set sent for creation of SREs), we integrate them to vehicles of the same type as contained in the SRE information, and to a vehicle starting from any of the depots in the VD set of the SRE. However we do not send the SREs to be integrated in any sequence. This is because, SREs without any TP vertices within their SETs are essentially independant entities but those SREs containing some TP-Vertex within any of their SETs are in a way connected to the SRE with superior or inferior causality, this connection stems down as per the TRS structure. We therefore ensure that while selecting SREs for integration, we start with the absolute causally-superior SREs. The relative superiority is directly inferred from each row of the CD; to respect the absolute superiority while sending SREs for integration, we ensure that:
\begin{enumerate}[label=\textbf{\scriptsize Step \arabic*}:]
    \item If the SRE to be integrated is in the RIGHT Column of the CD (in any row), all its causally-superior SREs (on the corresponding rows LEFT Column) must have been already sent for integration.
    \item If the SRE to be integrated is in the LEFT Column of the CD, we check if it is anywhere else on the CD in the RIGHT Column. If it is found anywhere else on the RIGHT Column, we restart from \textbf{Step 1} above, otherwise we send it for integration.
\end{enumerate}

After sending out all the SREs from the CD for integration within vehicles, we send the remaining SREs at random since they don't contain any TP vertex within any of their SETs and therefore random integration will do no harm.

In case this is not followed, some Route Clusters could become infeasible later on due to complexity (impossibility) of Waiting Time calculation at TPs.

\subsubsection{Moulding an SRE into its Route Portion Combination(s)} \label{Moulding an SRE into a Route Portion (RoPr), and weighing logics for the different permutations in the RoPr of a single SRE}

An SRE in its own form cannot be directly integrated into a Route6. We create Route Portions from SREs as the first step towards integration into vehicle route, after each SRE is chosen for integration.

A possible RoPr of an SRE is the sequence of vertices of the SRE such that independant permutations of its SET 1 and SET 3 is allowed, while maintaining the relative sequence of vertices in the order SET 1 $\rightarrow$ SET 2 $\rightarrow$ SET 3. Alongwith the vertex list (as the List 1 in the RoPr), all the other lists similar to those in a Route6 are created (with some default values, like 0 for all List 2 elements, only populating the Loading+Unloading Time within all TimeTuples), and not initiating any StatusCode entries in the List 4 for now).

An example of such an RoPr is provided in \autoref{longTab: RoPr Combo}; for any SRE, the number of such combinations will be the value $min(1, |\text{SET}\_1|!\cdot|\text{SET}\_3|!)$. Reiterating the main points, all the SET 1 and SET 3 elements are allowed to be in any relative ordering in the RoPr, however, all SET 2 vertices must be placed after all SET 1 vertices, all of which must be placed before SET 3 vertices in the List 1 of any RoPr-Combo.

\begin{table}[ht]
\centering
\caption{\textbf{All RoPr-Combos of \textit{SRE 5} from the \autoref{fig:SRE Examples}}\label{longTab: RoPr Combo}}
% \noindent
\resizebox{1\textwidth}{!}{ % Resize to fit the width of the page
\begin{tabular}{|*{9}{m{0.08\textwidth}|}} % 3 columns with a repeating pattern
% \hline
\cline{2-4} \cline{7-9}
% & \multicolumn{3}{c|}{ASD} & & & \multicolumn{3}{c|}{ASD} \\ % Merging columns properly
\multicolumn{1}{m{0.125\textwidth}}{} & \multicolumn{3}{|c|}{\textbf{RoPr Combo 1}} & \multicolumn{1}{m{0.125\textwidth}}{} & \multicolumn{1}{m{0.125\textwidth}}{} & \multicolumn{3}{|c|}{\textbf{RoPr Combo 2}} \\
\cline{1-4} \cline{6-9}
\textbf{List 1}: &	TP1 &	TP2	&	TP3 &&	\textbf{List 1}: &	TP2 &	TP1	&	TP3 \\
\cline{1-4} \cline{6-9}
\textbf{List 2}: &	0 &	0	&	0	&&	\textbf{List 2}: &	0 &	0	&	0	\\
\cline{1-4} \cline{6-9}
\textbf{List 3}: &	\{1D:2\} & \{1D:7\} &	\{1D:-9\}	&& \textbf{List 3}: &	\{1D:7\} & \{1D:2\} &	\{1D:-9\}	\\
\cline{1-4} \cline{6-9}
\textbf{List 4}: &	\{\} &	\{\}	&	\{\} &&	\textbf{List 4}: &	\{\} &	\{\}	&	\{\} \\
\cline{1-4} \cline{6-9}
\textbf{List 5}: &	( , 0.22, ) &	(, 0.77, ) &	(, 0.99, ) &&	\textbf{List 5}: &	( , 0.77, ) &	(, 0.22, ) &	(, 0.99, ) \\
\cline{1-4} \cline{6-9}
\textbf{List 6}: & \textit{SRE 5} & \textit{SRE 5} & \textit{SRE 5} && \textbf{List 6}: & \textit{SRE 5} & \textit{SRE 5} & \textit{SRE 5} \\
\cline{1-4} \cline{6-9}
\end{tabular}
} % For Resize Box Closure
\end{table}

\subsubsection{Understanding the Integration Logics of Route Clusters} \label{Understanding the Integration Logics of Route Clusters}

An entire solution of the whole problem is stored as a Route Cluster (RoCu); a RoCu contains Route6 for all vehicles. We develop different logical steps during the process of Integration of an SRE into a Vehicle's Route6. These progress through:
\begin{enumerate}
    \item The selection of an SRE (section \ref{Selection Logic for the Allocation of individual SREs to unique vehicles}) initiating an Integration step
    \item Development of the RoPr-Combos from the SRE, such that each RoPr is similar in structure to a typical Route6 section (\ref{Moulding an SRE into a Route Portion (RoPr), and weighing logics for the different permutations in the RoPr of a single SRE})
    \item Selection of some of the Integration Logics, based on the combinations of RoPr-Combo, Allocation and Matching Logics (section \ref{Logics for selection of some RoPr Combos}, \ref{Allocation Logics of Route Portion within the allocated vehicle's route, across different Route Clusters}, and \ref{Matching of the same vertices of SREs within the route of a single vehicle, ensuring feasibility during the Insertion of an SRE}). In addition to these, we also consider comparing results with Waiting Time computation and without it (as indicated in the last column of \autoref{Tab: Some sample Integration Logics assigned to respective Route Clusters}). In total we can have a total of 554 (slightly above 4*7*9*2, as found by multiplying all the total no. of logical cases as shown in \autoref{Tab: Some sample Integration Logics assigned to respective Route Clusters} along with the inclusion of all the True rows for which Waiting Time consideration as False is available) unique Integration Logics.
    \item Using a user input, we take an upper bound on the no. of unique logics we will be considering for a Heuristic run, (with a lower bound set at half the user-input value). Some Integration Logics are considered and put into a set NonTabuLogic; the number of such logics lie between the bounds.
    \item We generate empty Route Clusters equal to the no. of initial elements in the NonTabuLogic set.
    \item We assign each RoCu with a unique logic label from the NonTabuLogic set.
    \item Whenever an SRE is chosen for Integration, it is integrated into each of the Route Clusters; the integration is performed as per the current logic label of the respective RoCu. Thus each Route Cluster keeps its own unique copy of the developing problem solution; each solution possibly being drastically different from the others.
\end{enumerate}

Even within each individual Integration Logic, there are multiple combinations possible between the selected Route Portions, and the selected Vehicles. Such a case can be seen for the currently assigned Route Cluster 3's Logic Label in the \autoref{Tab: Some sample Integration Logics assigned to respective Route Clusters}, where multiple Route Portions will be returned (also happens for the Route Cluster 4), as well as all the vehicles possible to be allocated to will be considered. During such cases, we create temporary twins of the RoCu and try out all these combinations; we retain the best RoCu-twin that has the minimum cascaded time of the route cluster, among all cluster-twins.

The cascaded vehicle route duration comparison, between two clusters, is done by:
\begin{enumerate}
    \item Updating all the TimeTuples of each Route6. This is done along with the Waiting Time computations if the corresponding logic of Waiting Time computation Consideration Logic is True, otherwise, all waiting times are considered 0
    \item Arranging the route durations in descending order, for each RoCu separately
    \item Compare the highest route duration between the clusters. If found the same we compare the next highest vehicle route, until we get a deterministic comparison.
    \item If we have compared all the vehicular routes across the Route Clusters in descending order and all are found to be the same in value, then we deem both solutions to be equally good (and possibly the same).
\end{enumerate}

During the algorithm development, shuffling of the logic labels proved extremely successful. We develop a HyperParameter to decide if such shufflings should be carried out; this HyperParameter To\_Shuffle\_or\_Not\_to\_Shuffle takes a True value 80\% of the time before every Main Iteration (but we recommend setting it to full, as the improvement has been found to be dramatic). If To\_Shuffle\_or\_Not\_to\_Shuffle is found to be True, the logic labels of the Route Clusters get shuffled after every SRE Integration. In some cases, when all the waiting times within every Route6 of a RoCu cannot be computed (discussed in section \ref{Waiting Time computation function for a Route Cluster}), we ditch that RoCu along with the current logic label associated with it (thus this ditched logic label is not considered during future Integration-Logic-Label shufflings within the same Main Iteration). A lot of scope remains in generating new Integration Logics, and comparing our shuffling approach with the concept of Crossover used in Genetic Algorithms.

\begin{table}[ht]
\centering
\caption{\textbf{Some sample Integration Logics assigned to respective Route Clusters (four example Logic-Labels are provided)}\label{Tab: Some sample Integration Logics assigned to respective Route Clusters}}
% \noindent
\resizebox{1\textwidth}{!}{ % Resize to fit the width of the page
\begin{tabular}{|m{0.2\textwidth}|m{0.15\textwidth}|m{0.2\textwidth}|m{0.3\textwidth}|m{0.2\textwidth}|}
% \hline
\hline
\textbf{Total no. of logical cases $\rightarrow$}: & 4 & 7 & 9 & 2\\
\hline
\multicolumn{1}{m{0.125\textwidth}|}{} & \textbf{RoPr selection Logic} & \textbf{Allocation Logic} & \textbf{Matching Logic\newline (Breach Consideration, LeftOver Mobility, TP Matching)} & \textbf{Waiting Time computation Consideration Logic \#}\\
\cline{2-5}
\hline
\textbf{Route Cluster 1}: &	MinTime &	VertexSimilarity\_1	&	OFF & True \\
\hline
\textbf{Route Cluster 2}: &	Random &	VertexSimilarity\_3	&	(True, False, True) & False\\
\hline
\textbf{Route Cluster 3}: &	All & All &	(False, True, True) & False\\
\hline
\textbf{Route Cluster 4}: &	Min2Time & Random &	(True, True, False) & True\\
\hline

\end{tabular}
} % For Resize Box Closure
\newline
\begin{flushright}
\scriptsize\# For every False row, its corresponding True row is always considered as a separate available logic, and assigned to some other Route Cluster not mentioned here
\end{flushright}
\end{table}
\small

\subsubsection{Logics for selection of some RoPr Combos} \label{Logics for selection of some RoPr Combos}

As previously mentioned, each Route Cluster is associated with a specific Integration\_Logic. The first component of a typical Integration\_Logic determines the selection logics for RoPrs from among the RoPr-Combos. We used four following different RoPr Selection Logics:
\begin{itemize}
    \item \textbf{Random}:
    Ramdomly selects one single RoPr from among the RoPr-Combos

    \item \textbf{MinTime}:
    Selects the RoPr with the minimum duration compared to the other RoPrs in the RoPr-Combo. This calculation assumes the arrival time of the vehicle at first vertex in any RoPr to be 0, for all the RoPrs within the RoPr-Combos; waiting times are also considered to be 0, incase any TP is encountered. The final Vertex leaving time (\ie the last element in the TimeTuple, of the last Vertex) is compared for comparing the minimum durations.

    \item \textbf{Min2Time}:
    Similar to MinTime, but selects two RoPrs, one with the lowest duration, and the other with the second lowest duration.

    \item \textbf{All}:
    Selects the entire RoPr-Combo set for the subseuqnet processes.
\end{itemize}

\subsubsection{Allocation Logics of Route Portion within the allocated vehicle's route, across different Route Clusters} \label{Allocation Logics of Route Portion within the allocated vehicle's route, across different Route Clusters}

This function \textit{return}s a set of suitable vehicles which can be used for integrating a RoPr. The SRE information, along with any one among the below allocation logics, are passed to this function. The exact Vehicle information is obtained from the VT and VD set associated with the original SRE from which the RoPr was generated; all the vehicles matching both the VT and VD description are considered while choosing. An additional input of the Vertices present in the SRE, is used in some of the logics. The allocation logics used by us are:
\begin{itemize}
    \item \textbf{Random}:
    A random vehicle from among the suitable set of vehicles is selected.
    \item \textbf{VertexSimilarity\_1}:
    In this case the vertices in the SRE are iterated, within which the current vertices in the concerned vehicle's Route6 is also iterated upon; their number of matches are counted. The vehicle with the maximum value of these matches is selected.
    \item \textbf{VertexSimilarity\_2 and 3}:
    Here we take a more practical approach by counting only the feasible (feasibility in this case is not \wrt the vehicle capacity constraints but\wrt the necessity of keeping all vertices of the SRE's SET 1 before its SET 2 which must be entirely before its SET 3) number of insertions possible. We do this by (setting a counter at zero, and) starting with the first vertex in the RoPr; then we proceed by:
    \begin{enumerate}[label=\textbf{\scriptsize Step \arabic*}:]
        \item Calculating the number of Vertex-matchings (\ie the same vertex currently being considered from the SRE being present in the RoPr), from after each individual Vertex-matching of the previous SRE-element as separate cases. The logic behind considering the matchings only from after the previous matchings (if any) lies in the fact that when actually integrating the RoPr into a Route6, the relative positioning of the RoPr elements within the Route6 will remain intact, even though its individual vertices can be spaced out to accomodate other (any number of) vertices between them.
        \item We add this number of feasible matching to a counter.
        \item If matchings are found for this RoPr vertex, we proceed to the next vertex and start from \textbf{Step 1}; otherwise we proceed to \textbf{Step 4}.
        \item We check if the program flow converges here due to no matching, or due to all RoPr vertices having been considered. In case of all the RoPr vertices were considered, it would mean that the final vertex in the RoPr found a match in the Route6 and all the other vertices were also able to be inserted in the same relative order; in this case the function \textit{return}s the counter as the value. If however, the process flow converged due to the case of no match found after some (or no) matchings, the function \textit{return}s a modified value of the counter, which:
        \begin{itemize}
            \item For \textbf{VertexSimilarity\_2} is kept at: $counter+10000*(i+1)$
            \item For \textbf{VertexSimilarity\_2} is used as: $counter*(i+1)$
        \end{itemize}
        where $i$ denotes the number of elements of the RoPr which were possible to be matched, before the next ($i+1^{th}$) vertex couldn't find any feasible match.
    \end{enumerate}
    The vehicle with the maximum value of these values is chosen.
    
    \item \textbf{MinTime}:
    For all the vehicles, we consider their time of leaving last vertex in their respective routes (\ie we compare the last element of the last TimeTuple, for every suitable vehicle for which this value exists) in case it had been already computed during its last RoPr insertion. We choose the vehicle which has this value as the lowest; in case none of the values exist we assign one at random.
    \item \textbf{Min2Time}:
    This is similar to the MinTime logic, but in this case we select both the minimum and the second minimum options. We assign one vehicle at random if the values don't exist yet, which will also happen when the vehicle routes are empty.
    \item \textbf{All}:
    The full set of suitable vehicles is returned.
\end{itemize}

\subsubsection{Matching of the same vertices of SREs within the route of a single vehicle, ensuring feasibility during the Insertion of an SRE} \label{Matching of the same vertices of SREs within the route of a single vehicle, ensuring feasibility during the Insertion of an SRE}
When we integrate a RoPr into a Route6, we are able to decide whether the integration will consider Matching during this Integration, In the case of Matching having being turned ON, it is done so by toggling the below mentioned MatchMaking Parameters, each of which can take True/False values independently. Thus this Matching function has 9 possible logical scenarios, the $2^3$ toggled combinations of the MatchMaking Parameters and the case without Matching. 

The MatchMaking Parameters used in this function are:
\begin{enumerate}
    \item \textbf{Breach Consideration}:
    This considers the Breach Point as LeftOver insertion upper bound
    \item \textbf{LeftOver Mobility}:
    Allows some insertion mobility to the leftover (or all) RoPr-Elements, leftover elements being those without any confirmed match.
    \item \textbf{TP Matching}:
    Keeping this True (recommended) allows the TP vertices to also be considering during the Matching process described in section \ref{Matching of the same vertices of SREs within the route of a single vehicle, ensuring feasibility during the Insertion of an SRE}.
\end{enumerate}
Additionally a functional variable representing the LeftOver-insertion-UpperBound is initially set to a Position value after the final vertex of Route6. \par

The main idea here is to try and place elements from the RoPr (starting from the first element, and always ensuring the relative position of all the elements from the RoPr are undisturbed during/after the insertion). The steps followed are shown:

\begin{enumerate}[label=\textbf{\scriptsize Step \arabic*}:]
    \item Iteration 1 initiation across the RoPr elements ((we start considering each of the RoPr elements from its extreme left)
     
    \item Logical Flow within Iteration 1: \newline
    We try to find a match (\ie same vertex) of the RoPr-Element (under consideration from \textbf{Step 1} or \textbf{3}) with any of the Route6 elements. For comparing with the Route6 elements, the search starts from after the vertex position where any previous RoPr-Element was inserted; as for the first RoPr-Element matching since will be no such position available, we start the matching search from the starting vertex (after the VD) of the Route6.

    If a match is found, we check if inserting the RoPr-Element just after this match Position results in any (weight or volume) capacity breach across any of the subsequent StatusCodes in a temporarily modified Route6; this is done by updating the StatusCodes starting from the inserted RoPr-Element (as a continuation of the StatusCode from the match position). Any StatusCode is calculated as the sum of the previous Element's StatusCode with the current Element's VLC; it represents the resource amount contained within the vehicle, and is therefore never negative for any CT. If no capacity breach is found, we confirm this Match and go to \textbf{Step 3}. Otherwise, when a Breach is reported:
    \begin{itemize}
        \item If Breach consideration is True, the LeftOver-insertion-UpperBound is set as the position available just before the Breach Position. If it is False, the LeftOver-insertion-UpperBound is unchanged.
        \item We go to \textbf{Step 4} reporting a capacity breach.
    \end{itemize}
    
    \item We insert the RoPr-Element just next (rightwards in positioning) to the matched Route6-Element; practically this indicates that the vertex instead of being visited multiple times may be visited a single time. After this we start the Matching consideration for the next RoPr-Element from a position in the Route6 after this Match Position (\ie go to \textbf{Step 1}).
    
    \item The RoPr elements which were not inserted yet (\ie starting from the last unmatched element position, and taking all the subsequent elements always in the same sequence) are termed LeftOver. We try to place the LeftOver RoPr elements together in a suitable Position of the Route6. This suitable position is searched:
    \begin{itemize}
        \item From within a range starting after the latest RoPr-Element insertion Position upto the LeftOver-insertion-UpperBound, if LeftOver Mobility is set as True.
        \item If it is False, the insertion happens at the LeftOver-insertion-UpperBound.
    \end{itemize}

\end{enumerate}

A possible matching from the example in \autoref{tab: Route6 showing insertion positions} could be:
VD1 $\rightarrow$ WH1 $\rightarrow$ WH1 $\rightarrow$ NP3 $\rightarrow$ NP3 $\rightarrow$ NP3 $\rightarrow$ TP1\_GA $\rightarrow$ TP1\_GP

When Matching is turned OFF, then the full RoPr is integrated at the very end of the Route6, as shown for the same example:
VD1 $\rightarrow$ WH1 $\rightarrow$ NP3 $\rightarrow$ NP3 $\rightarrow$ TP1\_GA $\rightarrow$ WH1 $\rightarrow$ NP3 $\rightarrow$ TP1\_GP

\begin{table}[ht]
\centering
\caption{A. \textbf{Sample Route6 from our Illustrative Example in section \ref{An Example Case Study} of a VT4 on SMTS Road\_B indicating the possible insertion positions for accomodating (none or multiple) elements of new incoming  RoPrs.}\label{tab: Route6 showing insertion positions}}
% \noindent
\resizebox{1\textwidth}{!}{ % Resize to fit the width of the page
\begin{tabular}{|*{11}{m{0.1175\textwidth}|}} % BEst fit at 0.1175
\multicolumn{1}{c}{} & \multicolumn{1}{c}{} & \multicolumn{1}{m{0.135\textwidth}}{\centering Position 1 \newline \textbf{$\downarrow$}} & \multicolumn{1}{c}{} & \multicolumn{1}{m{0.135\textwidth}}{\centering Position 2 \newline \textbf{$\downarrow$}} & \multicolumn{1}{c}{} & \multicolumn{1}{m{0.135\textwidth}}{\centering Position 3 \newline \textbf{$\downarrow$}} & \multicolumn{1}{c}{} & \multicolumn{1}{m{0.135\textwidth}}{\centering Position 4 \newline \textbf{$\downarrow$}} & \multicolumn{1}{c}{} & \multicolumn{1}{m{0.135\textwidth}}{\centering Position 5 \newline \textbf{$\downarrow$}}\\
\cline{1-2}	\cline{4-4} \cline{6-6} \cline{8-8} \cline{10-10}
\textbf{List 1}: &	VD1 &&	WH1	&&	NP3	&&	NP3	&& TP1\_GA & \multicolumn{1}{c}{} \\
\cline{1-2}	\cline{4-4} \cline{6-6} \cline{8-8} \cline{10-10}
\textbf{List 2}: &	N/A &&	0	&& 0 &&	0	&& 0 & \multicolumn{1}{c}{}\\
\cline{1-2}	\cline{4-4} \cline{6-6} \cline{8-8} \cline{10-10}
\textbf{List 3}: &	\{\} &&	\{1D:11, 2D:13\}	&&	\{1D:-11, 2D:-13\}	&&	\{2P:6\}	&& \{2P:-6\}	& \multicolumn{1}{c}{}\\
\cline{1-2}	\cline{4-4} \cline{6-6} \cline{8-8} \cline{10-10}
\textbf{List 4}: &	\{\} &&	\{1D:11, 2D:13\}	&&	\{\}	&&	\{2P:6\}	&& \{\}	& \multicolumn{1}{c}{}\\
\cline{1-2}	\cline{4-4} \cline{6-6} \cline{8-8} \cline{10-10}
\textbf{List 5}: &	( , , 0) &&	(12.88, 5.15, 18.03) &&	(22.54, 5.15, 27.69) &&	(27.69, 1.2, 28.89) &&	(41.3, 1.2, 0* 42.5) & \multicolumn{1}{c}{}\\
\cline{1-2}	\cline{4-4} \cline{6-6} \cline{8-8} \cline{10-10}
\textbf{List 6}: &&& \textit{SRE 59} && \textit{SRE 59} && \textit{SRE 72} && \textit{SRE 72} & \multicolumn{1}{c}{}\\
\cline{1-2}	\cline{4-4} \cline{6-6} \cline{8-8} \cline{10-10}
\end{tabular}
} % For Resize Box Closure
\end{table}

\addtocounter{table}{-1}

\begin{table}[H]
\centering
\caption{B. \textbf{RoPr of the \textit{SRE 7} from \autoref{fig:SRE Examples}; the \textit{SRE 7} has only this possible combination within its RoPr-Combos set. Each individual column element is also indicated.}\label{tab: ropr element wise gapped}}
% \noindent
\resizebox{0.75\textwidth}{!}{ % Resize to fit the width of the page
\begin{tabular}{|*{11}{m{0.1175\textwidth}|}} % BEst fit at 0.1175
\multicolumn{1}{c}{} & \multicolumn{1}{c}{} & \multicolumn{1}{m{0.135\textwidth}}{\centering Element 1 \newline \textbf{$\downarrow$}} & \multicolumn{1}{c}{} & \multicolumn{1}{m{0.135\textwidth}}{\centering Element 2 \newline \textbf{$\downarrow$}} & \multicolumn{1}{c}{} & \multicolumn{1}{m{0.135\textwidth}}{\centering Element 3 \newline \textbf{$\downarrow$}}\\
\cline{1-1}	\cline{3-3} \cline{5-5} \cline{7-7}
\textbf{List 1}:&&	WH1	&&	NP3	&&	TP1\_GP\\
\cline{1-1}	\cline{3-3} \cline{5-5} \cline{7-7}
\textbf{List 2}:&&	0	&& 0 &&	0\\
\cline{1-1}	\cline{3-3} \cline{5-5} \cline{7-7}
\textbf{List 3}:&&	\{1D:5, 2D:4\}	&&	\{1D:-5, 2D:-4, 2P:8\}	&&	\{2P:-8\}\\
\cline{1-1}	\cline{3-3} \cline{5-5} \cline{7-7}
\textbf{List 4}:&&		&&		&&	\\
\cline{1-1}	\cline{3-3} \cline{5-5} \cline{7-7}
\textbf{List 5}:&&	( , 2.15, ) &&	( , 3.75, ) &&	( , 1.6, 0*, )\\
\cline{1-1}	\cline{3-3} \cline{5-5} \cline{7-7}
\textbf{List 6}:&& \textit{SRE 7} && \textit{SRE 7} && \textit{SRE 7}\\
\cline{1-1}	\cline{3-3} \cline{5-5} \cline{7-7}
\end{tabular}
} % For Resize Box Closure
\newline
\begin{flushright}
\scriptsize\hfill * Due to superior causality, they naturally has zero waiting times
\end{flushright}
\end{table}
\small

\subsubsection{Waiting Time computation function for a Route Cluster} \label{Waiting Time computation function for a Route Cluster}
% Waiting Time adjustment function at TPs
 This function accepts a RoCu and updates it with the computed Waiting Times in the appropriate List of each vertex within all Route6s in the RoCu. It starts by taking a random vehicle's Route6 from the RoCu and proceeding as below:

\noindent\hrulefill

\noindent\textbf{Function to update all TimeTuples of a RoCu, and Compute Waiting times}

\noindent\hrulefill
\begin{algorithmic}[1]

\Function{Geometrically-Computing\_Waiting\_Time}{passed\_RoCu}
    \While{TimeTuple computation of some position in any Route6 is yet to be done}
	\For{each Route6 in the passed\_RoCu}
		\For{each vertex in the chosen Route6 from the outer iteration} \Comment{The vertices are taken either from the beginning, or from a TP-vertex whose waiting Time computation was not possible earlier (refer to step 9 below)}
			\State calculate TimeTuple in its List 5; if we encounter a TP-Vertex, we estimate the best waiting time (for that specific visit) as below:
			\If{the encountered TP represents an Inferior Causal event}
				\If{corresponding Superior Causal event(s) have already been encountered in this function previously}
                        \State compute waiting time using geometrical approach
				\Else
					\State stop this Route6's TimeTuple computation at this TP-vertex (to resume TimeTuple calculations from this vertex-position later)
                        \State \Break \Comment{start the next Route6's TimeTuple computations}
				\EndIf
			\Else
				\State waiting time is 0 \Comment{Superior Causal events have no waiting time}
			\EndIf
		\EndFor
	\EndFor
    \EndWhile

\EndFunction
\noindent\hrulefill
\end{algorithmic}

In step 5 above, the TimeTuple calculations from the beginning Vertex are started by putting the value of the first vertex-leaving-time as 0, that vertex being its VD (while comparing with the Literature Instances in \autoref{longTab: Literature Dataset Overview}, we had to give a VT-specific starting time delay, due to requirement of converting the literature problems' vehicle fixed cost into our problem specifications) and just adding up the subsequent travel time to the next vertex, adding the loading/unloading time at this next vertex, and keep doing so repeatedly till we encounter a TP.

At a TP, we need to estimate the best waiting time accurately; this not a static calculation. Each unique TP usage associates its TP Name with a unique alphabet tag.
% Further the same Route6 column (here we are always dealing with uncompacted Route6, the Compacted\_Route6 as depicted in \autoref{longTab: Sample Route6 Examples} is only developed at the end of a Main Iteration to output the final solution) will have the SRE ID in its List 6.
Using this unique TP name (of the encountered TP from step 5), we find the corresponding row of the CD, and check if the SRE ID is in the Superior or Inferior Causality column. If it is found to be of superior causality, we don't need to calculate the Waiting Time of this vehicle at this Vertex-TP as it will be 0; this may be easily understood by considering the superior causal events (\ie the unloading events depicted in blue) in the Figures \ref{fig: Waiting Time Calculation 3to1} and \ref{fig: Waiting Time Calculation 1to3}, as these unloading events shouldn't require any waiting time. In case the SRE ID is found to be of inferior causality (step 6), then we check if its corresponding superior causal event(s) have already been encountered in this function previously (step 7). If the time computation is not possible (through step 8) due to some superior causal event(s) not having been encountered previously in this function, we leave this Route6's TimeTuple computation here since we are unable to proceed with the calculation of the TimeTuples across this vertex and further vertices for this Route6. We move on to another Route6 and (re)start the process of TimeTuple computation (from where it was last left, if at all) re-initiating from step 3.

\begin{landscape}
\begin{figure}
    \centering
    \includegraphics[width=1.035\linewidth]{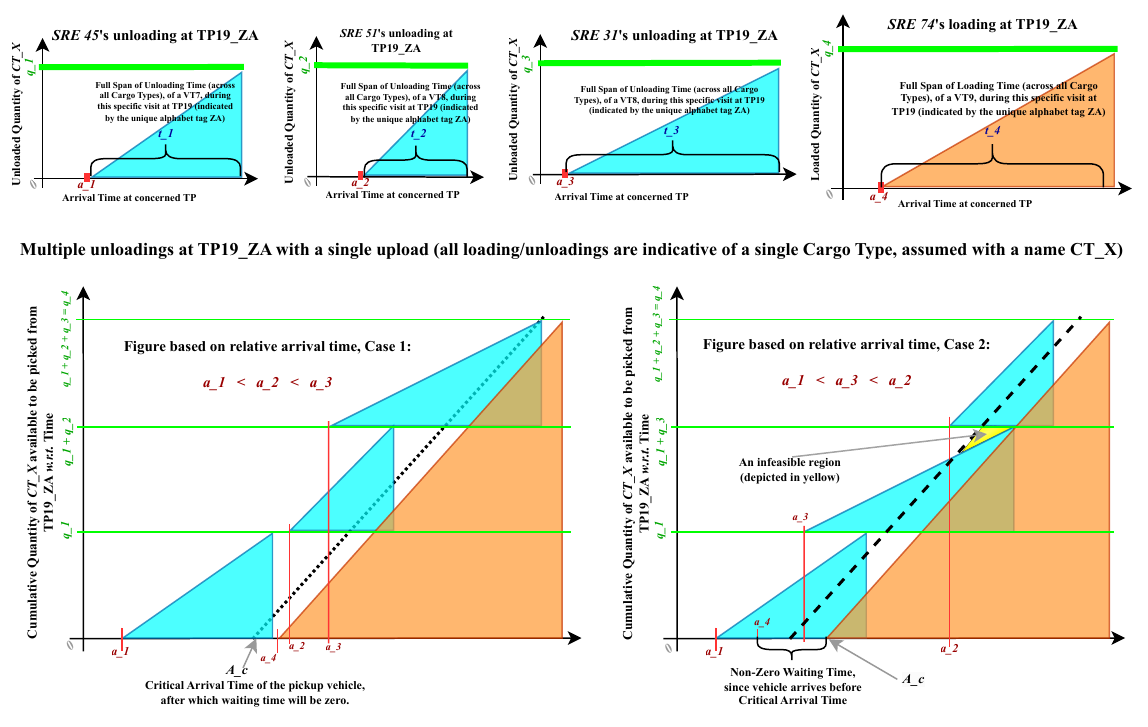} % 0.9125 fits
    \caption{\textbf{Visual representation of the Waiting Time calculation, for the first CD-row of \autoref{tab: Hypothetical Case for Waiting Time Image} (the loading or unloading is represented by colour and not by signage in this figure). Here two cases are considered to deepen the understanding of the approach used by us.}}
    \label{fig: Waiting Time Calculation 3to1}
\end{figure}
% \end{landscape}

% \begin{landscape}
\begin{figure}
    \centering
    \includegraphics[width=1.035\linewidth]{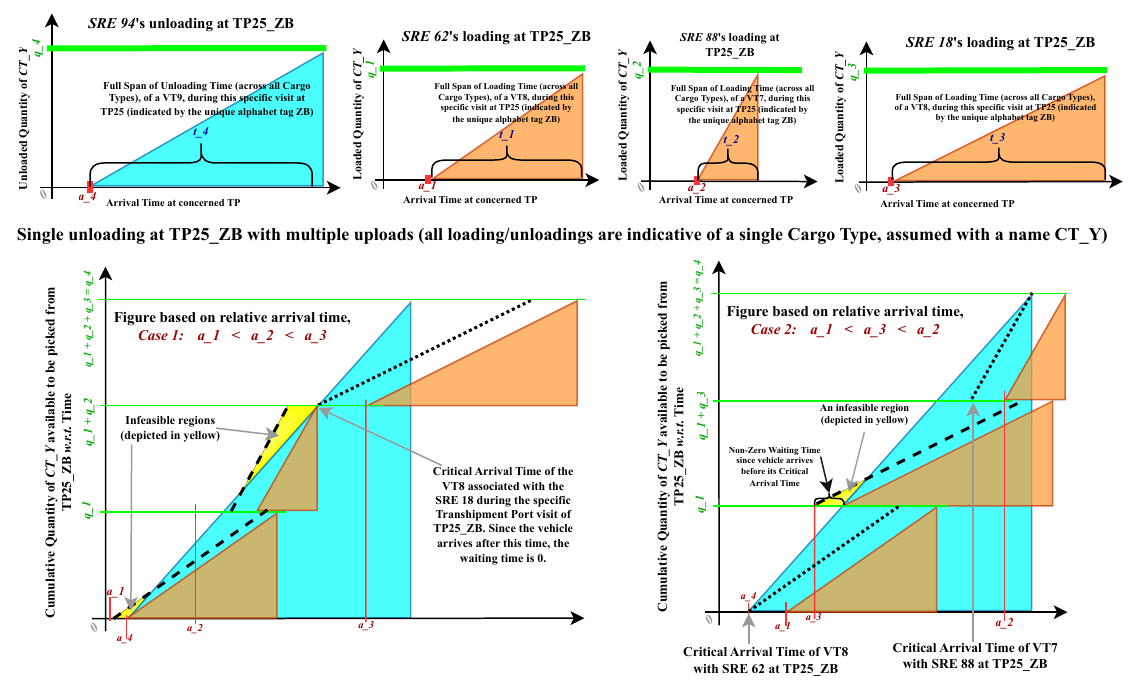} % 0.925 fits
    \caption{\textbf{Visual representation of the Waiting Time calculation, for the second CD-row of \autoref{tab: Hypothetical Case for Waiting Time Image}. Some of the infeasible regions (yellow regions near dashed lines) and Critical Arrival Times (calculated at the base of each dotted line) for each of the vehicles arriving at the TP with inferior causality (SREs depicted as orange loadings) are mentioned with arrows.}}
    \label{fig: Waiting Time Calculation 1to3}
\end{figure}
\end{landscape}

However, if all superior causal events are encountered, their current TimeTuples in their respective List 5 will be known and fully populated, and we can now proceed to compute waiting times using a geometrical approach (step 8). We use all these TimeTuples to construct the Waiting Time computation problem geometry (as shown in the first four small input figures in both the Figures \ref{fig: Waiting Time Calculation 3to1} and \ref{fig: Waiting Time Calculation 1to3}); the respective VLCs from each of their List 3 determines the heights of the individual triangular loading/unloading sections, and the full loading/unloading time (across all the CTs present in the respective VLC, and transferred in this specific TP-visit) determines the base of the right triangular region. We first place all the unloading (blue) triangles over each others' baseline, stacked from the 0 resource level in the order of their arrival times (refer to the stacking of the orange unloading triangles in \autoref{fig: Waiting Time Calculation 3to1}); their hypotenuses start from exactly their arrival times, since these superior causal events don't have any waiting time. Next we place all the inferior causal triangle geometries on top of their respective baselines starting from the zero resource level, stacked from the bottom up in the order of their arrival times (refer to the stacking of the blue loading triangle in \autoref{fig: Waiting Time Calculation 1to3}); for explanation purposes, we initially place them  at an imaginary infinite time at the very far right. Finally we try to push all the loading triangles towards the left , always keeping in mind that:
 
 \begin{itemize}
     \item The hypotenuse of the loading triangles (orange triangles in figures \autoref{fig: Waiting Time Calculation 3to1} and \autoref{fig: Waiting Time Calculation 1to3} do not cross the hypotenuse of any of the unloading triangles within their respective resource height
     \item The hypotenuse of the loading triangle does not intersect with its own arrival time
 \end{itemize}

    In \autoref{fig: Waiting Time Calculation 3to1} the resource height of the loading triangle ranges from level $0$ to $q\_4$. In the case 1, if the actual arrival time of the loading triangle $a_4$ would have been less than or equal to the Critical Arrival Time (calculated at base of the loading hypotenuse pushed left-wards till it reaches the edge of feasibility, and depicted as dotted lines in figures \ref{fig: Waiting Time Calculation 3to1} and \ref{fig: Waiting Time Calculation 1to3}) $A\_c$, the loading triangle could have then be pushed further left to reduce any waiting time. However, the (big, orange) loading triangle in case 2 cannot be shifted any further left as it then goes to the infeasible region (yellow region, depicted with an arrow in the figure); in this case the infeasible region corresponds to the loading vehicle having better loading speed (steeper hypotenuse) due to which it has been able to take up all resources delivered to the TP by the SRE 45 and SRE 31 quickly even though it started its loading process at the Critical Arrival Time.

    Finally, as per our design, we compute the waiting time (calculated for each loading triangle individually) as the difference between the actual arrival time of a loading triangle, and the time from where its hypotenuse starts from its triangle-base. However this is the waiting time estimated considering only the specific Cargo Type for which the resource levels were constructed. If the VLC has other cargo types which are being transferred, similar geometrical approach is taken to find the waiting times for them too; in terms of the triangle geometries, only their heights will change depending on the new resource levels of the other CT, but their base lengths will remain the same as we consider the entire loading/unloading time across all CTs (\ie for the full VLC). The final waiting time is the maximum among all these waiting times calculated for each CT-specific geometrical approach.
    
    This waiting time minimization can be vividly imagined by sliding the loading triangles leftwards from infinity, on their respective baselines, and stopping only if any part of their hypotenuse hits any part of the hypotenuses of any unloading triangle (within its resource levels on which it was slid); we also stop beforehand if the loading triangle's arrival time gets intersected by its hypotenuse. The intersection with the arrival times is prevented in the case of the SRE 18 triangle hitting $a\_3$ in Case 1 of \autoref{fig: Waiting Time Calculation 1to3}, and also for the SRE 88's triangle geometry hitting $a\_2$ in the Case 2; both of these triangles have their Critical Arrival Times shown as the base of their dotted hypotenuse pushed further leftwards, and as conveyed previously, they will have 0 waiting times.

    This is the exact logic we implement in code, and envision lot of future potential in utilizing and improving this idea of waiting time minimization. As an example, different stacking orders, and breaking triangle hypotenuse into multiple slide-able parts could be looked into. A caveat in our approach is that we consider the base width for each triangle to be the entire loading-unloading time for a VLC for simplifying our approach, which allows some sub-optimality in the waiting time minimization approach. (We are yet to devise a better strategy for the same!)

Resource Quantity (loading or unloading is represented by colour and not by signage in this figure). If it has not been made clear yet, the rate of loading and unloading activities get represented through the hypotenuse, across the total time during which all loading/unloading (for all cargo type) happens (this time being represented by a triangle's base-length).

Multiple passes of this process will allow all Waiting Times to be computed and fully update all the Time Tuples. In case, it so happens that no progress in the time calculation is observed across all the Route6 even after trying all of the Route6s once again (\ie every Route6 time computation is left at the same place, even after retrying from the stopped-Vertex calculation position of each respective Route6), we abandon this RoCu concluding the current Integration-Logic of the RoCu to be a TabuLogic, complications during Matchings section \ref{Matching of the same vertices of SREs within the route of a single vehicle, ensuring feasibility during the Insertion of an SRE}, or due to not following the sequential SRE Allocation process of section \ref{Selection Logic for the Allocation of individual SREs to unique vehicles}.

\clearpage

\section{Detailed Discussions of the Computational Results} \label{Detailed Discussions and Computational Results}
The specifications of the different computers used during the compilation of the computational results is provided in \autoref{longTab: Computers Used} (All computations were run on any of these computers).

\scriptsize
\begin{longtable}
{|m{0.1\linewidth}|m{0.2\linewidth}|m{0.2\linewidth}|m{0.2\linewidth}|m{0.2\linewidth}|}
\caption{\textbf{Computer Infrastructure Details}\label{longTab: Computers Used}}\\
 \hline
 & \textbf{C\_1} & \textbf{C\_2} & \textbf{C\_3} & \textbf{C\_4} \\
 \hline
 
 \endfirsthead

 % \hline
 % \multicolumn{2}{|c|}{Continuation of \autoref{longTab: Multi-functionality of Graph Vertices} \textbf{Multi-functionality of Graph Vertices}}\\
 % \hline
 % \textbf{Mathematical View Point} & \textbf{Nomenclature Description}\\
 % \hline
 \endhead

\hline	\textbf{Gurobi Optimizer Version}	&	11.0.0 build v11.0.0rc2 (win64 - Windows 11+.0 (22631.2))	&	11.0.1 build v11.0.1rc0 (win64 - Windows 11.0 (22621.2))	&	11.0.0 build v11.0.0rc2 (linux64 - "CentOS Linux 7 (Core)")	&	12.0.0 build v12.0.0rc1 (win64 - Windows 11.0 (22631.2))	\\
\hline	\textbf{Python Version}	&	3.11.5 | packaged by Anaconda, Inc. | (main, Sep 11 2023) [MSC v.1916 64 bit (AMD64)]	&	3.11.7 | packaged by Anaconda, Inc. | (main, Dec 15 2023) [MSC v.1916 64 bit (AMD64)]	&	3.9.19 (main, May 6 2024)	&	3.12.7 | packaged by Anaconda, Inc. | (main, Oct 4 2024) [MSC v.1929 64 bit (AMD64)]	\\
\hline	\textbf{CPU Model}	&	Intel(R) Xeon(R) W-2223 CPU @ 3.60GHz	&	Intel(R) Xeon(R) W-2133 CPU @ 3.60GHz	&	Intel(R) Xeon(R) Silver 4110 CPU @ 2.10GHz	&	13th Gen Intel(R) Core(TM) i7-13700	\\
\hline	\textbf{Thread Count}	&	4 physical cores, 8 logical processors, using up to 8 threads	&	6 physical cores, 12 logical processors, using up to 12 threads	&	16 physical cores, 32 logical processors, using up to 32 threads	&	16 physical cores, 24 logical processors, using up to 24 threads	\\
\hline	\textbf{Installed RAM}	&	16.0 GB (15.7 GB usable)	&	16.0 GB (15.7 GB usable)	&	62 GB	&	16.0 GB (15.8 GB usable)	\\
\hline
\end{longtable}
\small
We develop novel Instances for benchmarking the Heuristic \wrt MILP results. Overview of dataset details for the small and large instances developed by us can be found in \autoref{longTab: Small Dataset Overview} and \ref{longTab: Large Dataset Overview} respectively, with similar details of the literature instance provided in \autoref{longTab: Literature Dataset Overview}.

\subsection{Continuous Results of our Small Instances} \label{Continuous Results of our Small Instances}

Here we discuss about the small instances modelled as continuous problems. Some points to observe for the Continuous problems are:
\begin{enumerate}
    \item For the instance \instbox{S4}, the problem was found to be Infeasible with Levels 1 and 2. The Cascade 1 with Levels=3 gave a lower objective function value of 88.613 as the largest route duration among all vehicles, with a lower MILP Gap of 48.75\%, but its Cascade 2 could not find a feasible solution (Gurobi warm start failed to load).
    \item For the instance \instbox{S7}, the optimal solution using only 1 level proved to have incomplete solution space, as the heuristic solution was found to be better. Thus, the MILP for this case was also run using levels equal to 2 (the relative gaps of the Heuristic's best-found solution \wrt both the MILP runs is reported).
    \item For the instance \instbox{S8}, the Vehicle with the highest trip duration in the first Cascade was not the same Vehicle with the highest trip duration in the final Cascade.
    \item For the instance of \instbox{S11}, no solution was found during the fourth Cascade with Levels=1; Levels=2 bettered the solution.
\end{enumerate}

\begin{landscape}
\vspace*{\fill}  % Push down
\footnotesize
% \scriptsize

\begin{longtable}{|m{0.03\linewidth}|m{0.06\linewidth}|m{0.06\linewidth}|m{0.06\linewidth}|m{0.06\linewidth}|m{0.05\linewidth}|m{0.06\linewidth}|m{0.045\linewidth}|m{0.065\linewidth}|m{0.065\linewidth}|m{0.075\linewidth}|m{0.085\linewidth}|m{0.075\linewidth}|}
\caption{\textbf{Details of the Small Instances developed for our computational study}\label{longTab: Small Dataset Overview}}\\
\hline
\multirow{2}{=}{\textbf{ID}} & \textbf{No. of Vehicle Types} & \textbf{No. of PickUp Cargo Types} & \textbf{Delivery Cargo Types} & \textbf{Vehicle Depots} & \textbf{Ware Houses} & \textbf{Simul-taneous Nodes} & \textbf{Split Nodes} & \textbf{Trans-Shipment Ports} & \textbf{Relief Centres} & \textbf{No. of Modal Segments} & \multicolumn{2}{|c|}{\textbf{No. of Incompatibilities}} \\
\cline{2-13}
& \textbf{$|K|$} & \textbf{$|P|$} & \textbf{$|D|$} & \textbf{$|H|$} & \textbf{$|W|$} & \textbf{$|N^M|$} & \textbf{$|N^P|$} & \textbf{$|S|$} & \textbf{$|R|$} & \textbf{$|SMTS|$} &  \textbf{Vehicles -vs- Cargos} & \textbf{Cargos -vs- TPs} \\
\hline
\hline
\endfirsthead
\hline
\multirow{2}{=}{\textbf{ID}} & \textbf{No. of Vehicle Types} & \textbf{No. of PickUp Cargo Types} & \textbf{Delivery Cargo Types} & \textbf{Vehicle Depots} & \textbf{Ware Houses} & \textbf{Simul-taneous Nodes} & \textbf{Split Nodes} & \textbf{Trans-Shipment Ports} & \textbf{Relief Centres} & \textbf{No. of Modal Segments} & \multicolumn{2}{|c|}{\textbf{No. of Incompatibilities}} \\
\cline{2-13}
& \textbf{$|K|$} & \textbf{$|P|$} & \textbf{$|D|$} & \textbf{$|H|$} & \textbf{$|W|$} & \textbf{$|N^M|$} & \textbf{$|N^P|$} & \textbf{$|S|$} & \textbf{$|R|$} & \textbf{$|SMTS|$} &  \textbf{Vehicles -vs- Cargos} & \textbf{Cargos -vs- TPs} \\
\hline
\hline
\endhead
\instbox{S1}	&	2	&	1	&	1	&	2	&	1	&	2	&	2	&	1	&	1	&	2	&	0	&	0	\\
\instbox{S2}	&	2	&	1	&	1	&	2	&	1	&	3	&	3	&	1	&	1	&	2	&	0	&	0	\\
\instbox{S3}	&	2	&	1	&	1	&	2	&	1	&	4	&	4	&	1	&	1	&	2	&	0	&	0	\\
\instbox{S4}	&	3	&	0	&	2	&	2	&	1	&	0	&	6	&	1	&	0	&	3	&	0	&	0	\\
\instbox{S5}	&	5	&	0	&	1	&	5	&	2	&	0	&	3	&	5	&	0	&	5	&	0	&	0	\\
\instbox{S6}	&	5	&	0	&	1	&	5	&	2	&	0	&	3	&	5	&	0	&	5	&	0	&	0	\\
\instbox{S7}	&	5	&	0	&	1	&	5	&	2	&	0	&	2	&	8	&	0	&	5	&	0	&	0	\\
\instbox{S8}	&	5	&	0	&	1	&	5	&	2	&	1	&	1	&	8	&	0	&	5	&	0	&	0	\\
\instbox{S9}	&	5	&	0	&	1	&	5	&	2	&	2	&	2	&	8	&	0	&	5	&	0	&	0	\\
\instbox{S10}	&	5	&	0	&	1	&	5	&	2	&	2	&	2	&	8	&	0	&	5	&	0	&	0	\\
\instbox{S11}	&	3	&	1	&	1	&	2	&	1	&	4	&	5	&	2	&	1	&	3	&	0	&	0	\\
\instbox{S12}	&	5	&	0	&	1	&	5	&	2	&	2	&	2	&	10	&	0	&	5	&	0	&	0	\\
\instbox{M13}	&	3	&	1	&	1	&	2	&	1	&	5	&	5	&	2	&	1	&	3	&	0	&	0	\\
\instbox{M14}	&	3	&	1	&	1	&	2	&	1	&	5	&	6	&	2	&	1	&	3	&	0	&	0	\\
\instbox{M15}	&	2	&	1	&	1	&	2	&	1	&	5	&	5	&	1	&	1	&	2	&	0	&	0	\\
\instbox{M16}	&	3	&	1	&	2	&	4	&	2	&	3	&	1	&	3	&	1	&	3	&	2	&	3	\\
\instbox{M17}	&	3	&	1	&	2	&	4	&	2	&	1	&	3	&	3	&	1	&	3	&	2	&	3	\\
\instbox{M18}	&	3	&	1	&	2	&	4	&	2	&	2	&	2	&	3	&	1	&	3	&	2	&	3	\\
\instbox{M19}	&	4	&	1	&	2	&	4	&	2	&	1	&	2	&	2	&	1	&	4	&	3	&	3	\\
\instbox{M20}	&	4	&	1	&	2	&	4	&	2	&	2	&	1	&	2	&	1	&	4	&	3	&	3	\\
\instbox{M21}	&	6	&	2	&	0	&	3	&	0	&	3	&	3	&	2	&	2	&	3	&	4	&	0	\\
\instbox{M22}	&	3	&	1	&	2	&	4	&	2	&	3	&	3	&	3	&	1	&	3	&	2	&	3	\\
\instbox{M23}	&	3	&	1	&	2	&	2	&	1	&	5	&	5	&	2	&	1	&	3	&	0	&	0	\\
\instbox{M24}	&	3	&	1	&	2	&	2	&	1	&	6	&	4	&	2	&	1	&	3	&	0	&	0	\\
\instbox{M25}	&	4	&	1	&	2	&	4	&	2	&	2	&	2	&	2	&	1	&	4	&	3	&	3	\\
\hline
\end{longtable}
{\small The prefix in the ID indicates approximate (simplified) instance sizes (for all the Tables \ref{longTab: Small Dataset Overview}, \ref{longTab: Large Dataset Overview}, \ref{Tab: Small Integer Results}, \ref{Tab: Small Continuous Results}, \ref{longTab: Large Instance Results Integer} and \ref{longTab: Large Instance Results Continuous}):

$\instbox{S}\rightarrow Small \quad\quad \instbox{M}\rightarrow Medium \quad\quad \instbox{L}\rightarrow Large$}

\vspace*{\fill}  % Push down
\clearpage
\vspace*{\fill}  % Push down

\footnotesize

\begin{longtable}{|m{0.025\linewidth}|m{0.02\linewidth}|m{0.02\linewidth}|m{0.02\linewidth}|m{0.02\linewidth}|m{0.02\linewidth}|m{0.035\linewidth}|m{0.03\linewidth}|m{0.02\linewidth}|m{0.02\linewidth}|m{0.055\linewidth}|m{0.15\linewidth}|m{0.125\linewidth}|m{0.085\linewidth}|m{0.05\linewidth}|}
\caption{\textbf{Details of the Large Instances developed in this study for computational analysis}\label{longTab: Large Dataset Overview}}\\
\hline
\multirow{2}{=}{\textbf{ID}} & \multirow{2}{=}{\textbf{$|K|$}} & \multirow{2}{=}{\textbf{$|P|$}} & \multirow{2}{=}{\textbf{$|D|$}} & \multirow{2}{=}{\textbf{$|H|$}} & \multirow{2}{=}{\textbf{$|W|$}} & \multirow{2}{=}{\textbf{$|N^M|$}} & \multirow{2}{=}{\textbf{$|N^P|$}} & \multirow{2}{=}{\textbf{$|S|$}} & \multirow{2}{=}{\textbf{$|R|$}} & \multirow{2}{=}{\textbf{$|SMTS|$}} & \multicolumn{2}{|c|}{\textbf{No. of Incompatibilities}} & \multicolumn{2}{|c|}{\textbf{Computer Used \#}} \\
\cline{12-15}
& & & & & & & & & & & Cargos -vs- Vehicles & TPs -vs- Cargos & Continuous & Integer \\
\hline
\hline
\endfirsthead
\hline
\multirow{2}{=}{\textbf{ID}} & \multirow{2}{=}{\textbf{$|K|$}} & \multirow{2}{=}{\textbf{$|P|$}} & \multirow{2}{=}{\textbf{$|D|$}} & \multirow{2}{=}{\textbf{$|H|$}} & \multirow{2}{=}{\textbf{$|W|$}} & \multirow{2}{=}{\textbf{$|N^M|$}} & \multirow{2}{=}{\textbf{$|N^P|$}} & \multirow{2}{=}{\textbf{$|S|$}} & \multirow{2}{=}{\textbf{$|R|$}} & \multirow{2}{=}{\textbf{$|SMTS|$}} & \multicolumn{2}{|c|}{\textbf{No. of Incompatibilities}} & \multicolumn{2}{|c|}{\textbf{Computer Used \#}} \\
\cline{12-15}
& & & & & & & & & & & Cargos -vs- Vehicles & TPs -vs- Cargos & Continuous & Integer \\
\hline
\hline
\endhead
\instbox{S26}	&	3	&	0	&	1	&	2	&	1	&	0	&	2	&	1	&	0	&	3	&	0	&	0	&	C_1	&	C_1	\\
\instbox{S27}	&	2	&	1	&	1	&	2	&	1	&	2	&	2	&	1	&	1	&	2	&	0	&	0	&	C_3	&	C_3	\\
\instbox{M28}	&	5	&	0	&	1	&	5	&	2	&	2	&	2	&	5	&	0	&	5	&	0	&	0	&	C_1	&	C_3	\\
\instbox{M29}	&	3	&	3	&	5	&	4	&	2	&	7	&	9	&	3	&	2	&	3	&	2	&	3	&	C_2	&	C_1	\\
\instbox{L30}	&	8	&	2	&	4	&	4	&	3	&	20	&	20	&	6	&	3	&	4	&	8	&	1	&	C_1	&	C_2	\\
\instbox{L31}	&	10	&	2	&	3	&	7	&	4	&	10	&	10	&	12	&	4	&	6	&	2	&	26	&	C_1	&	C_3	\\
\instbox{L32}	&	11	&	2	&	3	&	7	&	4	&	10	&	10	&	7	&	4	&	7	&	3	&	11	&	C_3	&	C_1	\\
\instbox{L33}	&	11	&	2	&	3	&	7	&	4	&	8	&	8	&	12	&	4	&	7	&	3	&	26	&	C_2	&	C_1	\\
\instbox{L34}	&	11	&	2	&	4	&	7	&	4	&	8	&	8	&	6	&	4	&	7	&	9	&	0	&	C_1	&	C_3	\\
\instbox{L35}	&	10	&	3	&	3	&	7	&	4	&	10	&	10	&	12	&	4	&	6	&	6	&	16	&	C_2	&	C_1	\\
\instbox{L36}	&	11	&	2	&	3	&	7	&	4	&	30	&	30	&	10	&	4	&	6	&	8	&	1	&	C_3	&	C_3	\\
\instbox{L37}	&	11	&	2	&	4	&	7	&	6	&	40	&	50	&	10	&	5	&	6	&	9	&	1	&	C_2	&	C_3	\\
\instbox{L38}	&	11	&	2	&	4	&	7	&	7	&	60	&	50	&	10	&	5	&	6	&	9	&	1	&	C_3	&	C_2	\\
\hline
\end{longtable}
\# For both MILP runs using Gurobi, and Heuristic Algorithm implemented in Python

\vspace*{\fill}  % Push down
\footnotesize

\begin{longtable}{|m{0.025\linewidth}|m{0.02\linewidth}|m{0.02\linewidth}|m{0.02\linewidth}|m{0.02\linewidth}|m{0.02\linewidth}|m{0.03\linewidth}|m{0.025\linewidth}|m{0.02\linewidth}|m{0.02\linewidth}|m{0.055\linewidth}|m{0.15\linewidth}|m{0.125\linewidth}|}
\caption{\textbf{Dataset details of the Literature Instances discussed in the Table 5 in \cite{NUCAMENDIGUILLEN2021}}\label{longTab: Literature Dataset Overview}}\\
\hline
\multirow{2}{=}{\textbf{ID}} & \multirow{2}{=}{\textbf{$|K|$}} & \multirow{2}{=}{\textbf{$|P|$}} & \multirow{2}{=}{\textbf{$|D|$}} & \multirow{2}{=}{\textbf{$|H|$}} & \multirow{2}{=}{\textbf{$|W|$}} & \multirow{2}{=}{\textbf{$|N^M|$}} & \multirow{2}{=}{\textbf{$|N^P|$}} & \multirow{2}{=}{\textbf{$|S|$}} & \multirow{2}{=}{\textbf{$|R|$}} & \multirow{2}{=}{\textbf{$|SMTS|$}} & \multicolumn{2}{|c|}{\textbf{No. of Incompatibilities}} \\
\cline{12-13}
& & & & & & & & & & & Cargos -vs- Vehicles & TPs -vs- Cargos \\
\hline
\endfirsthead
\hline
\multirow{2}{=}{\textbf{ID}} & \multirow{2}{=}{\textbf{$|K|$}} & \multirow{2}{=}{\textbf{$|P|$}} & \multirow{2}{=}{\textbf{$|D|$}} & \multirow{2}{=}{\textbf{$|H|$}} & \multirow{2}{=}{\textbf{$|W|$}} & \multirow{2}{=}{\textbf{$|N^M|$}} & \multirow{2}{=}{\textbf{$|N^P|$}} & \multirow{2}{=}{\textbf{$|S|$}} & \multirow{2}{=}{\textbf{$|R|$}} & \multirow{2}{=}{\textbf{$|SMTS|$}} & \multicolumn{2}{|c|}{\textbf{No. of Incompatibilities}} \\
\cline{12-13}
& & & & & & & & & & & Cargos -vs- Vehicles & TPs -vs- Cargos \\
\hline
\endhead

$	I_1	$	&	8	&	1	&	0	&	13	&	0	&	13	&	0	&	0	&	1	&	1	&	0	&	0	\\
$	I_2	$	&	8	&	1	&	0	&	13	&	0	&	12	&	0	&	0	&	1	&	1	&	0	&	0	\\
$	I_3	$	&	8	&	1	&	0	&	13	&	0	&	11	&	0	&	0	&	1	&	1	&	0	&	0	\\
$	I_4	$	&	8	&	1	&	0	&	13	&	0	&	10	&	0	&	0	&	1	&	1	&	0	&	0	\\
$	I_5	$	&	8	&	1	&	0	&	13	&	0	&	10	&	0	&	0	&	1	&	1	&	0	&	0	\\
$	I_6	$	&	8	&	1	&	0	&	13	&	0	&	11	&	0	&	0	&	1	&	1	&	0	&	0	\\
$	I_7	$	&	8	&	1	&	0	&	13	&	0	&	10	&	0	&	0	&	1	&	1	&	0	&	0	\\
$	I_8	$	&	8	&	1	&	0	&	13	&	0	&	10	&	0	&	0	&	1	&	1	&	0	&	0	\\
$	I_9	$	&	8	&	1	&	0	&	13	&	0	&	7	&	0	&	0	&	1	&	1	&	0	&	0	\\
$	I_{10}	$	&	8	&	1	&	0	&	13	&	0	&	10	&	0	&	0	&	1	&	1	&	0	&	0	\\
$	I_{11}	$	&	8	&	1	&	0	&	13	&	0	&	13	&	0	&	0	&	1	&	1	&	0	&	0	\\
$	I_{12}	$	&	8	&	1	&	0	&	13	&	0	&	9	&	0	&	0	&	1	&	1	&	0	&	0	\\
$	I_{13}	$	&	8	&	1	&	0	&	13	&	0	&	11	&	0	&	0	&	1	&	1	&	0	&	0	\\
$	I_{14}	$	&	8	&	1	&	0	&	13	&	0	&	12	&	0	&	0	&	1	&	1	&	0	&	0	\\
$	I_{15}	$	&	8	&	1	&	0	&	13	&	0	&	11	&	0	&	0	&	1	&	1	&	0	&	0	\\
\hline
\end{longtable}
\vspace*{\fill}  % Push down
\end{landscape}

\begin{landscape}
\tiny
\vspace*{\fill}  % Push down
\begin{table}[ht]
\centering
\caption{\textbf{Computational Results of Small Instances from \autoref{longTab: Small Dataset Overview} with resource flow variables as Continuous}\label{Tab: Small Continuous Results}}
% \noindent
\resizebox{1.55\textwidth}{!}{ % Resize to fit the width of the page
\begin{tabular}{|*{21}{c|}} % 3 columns with a repeating pattern
% \hline
% \begin{tabular}{*{21}{|c}|}
\cline{2-21}
\multicolumn{1}{c|}{} & \multicolumn{9}{|
c|}{\makecell{\textbf{Cascaded Makespan Minimization} \\ \textbf{using our developed exact MILP Formulation}}} & \multicolumn{9}{|c|}{\textbf{PSR-GIP Heuristic}} & \multicolumn{2}{|c|}{\makecell{\textbf{Comparative} \\ \textbf{Gap \%}}} \\
\hline
\textbf{ID}	&	\textbf{$O$}	&	\textbf{$S_F$}	&	\textbf{$\mathbb{C}$}	&	\textbf{$T$}	&	\textbf{$G_O$}	&	\textbf{$G_A$}	&	\textbf{$L$}	&	\textbf{$S_I$}	&	\textbf{$I$}	&	\textbf{$L_B$}	&	\textbf{$L_A$}	&	\textbf{$L_{SD}$}	&	\textbf{$S_B$}	&	\textbf{$S_A$}	&	\textbf{$S_{SD}$}	&	\textbf{$V$}	&	\textbf{$R_A$}	&	\textbf{$R_{SD}$}	&	\textbf{$OL_B$}	&	\textbf{$S_{FB}$}	\\
\hline
\hline
\instbox{S1}	&	64.516	&	79.629	&	2	&	0.00	&	0.00	&	0.00	&	1	&	116.307	&	46.06	&	64.516	&	68.924	&	5.584	&	79.629	&	84.257	&	5.568	&	2/2	&	0.27	&	0.02	&	0.00	&	0.00	\\
\instbox{S2}	&	64.278	&	152.772	&	4	&	1.83	&	0.00	&	0.00	&	1	&	176.506	&	15.54	&	64.278	&	66.953	&	3.263	&	161.304	&	180.419	&	15.128	&	4/4	&	0.5	&	0.05	&	0.00	&	5.29	\\
\instbox{S3}	&	68.119	&	252.419	&	4	&	36.98	&	0.00	&	0.00	&	1	&	254.323	&	0.75	&	68.203	&	96.729	&	10.905	&	214.709	&	263.323	&	43.68	&	4/4	&	0.79	&	0.1	&	0.12	&	-17.56	\\
\instbox{S4}	&	88.614	&	195.078	&	3	&	270.13	&	56.21	&	32.95	&	4	&	195.078	&	0.00	&	88.964	&	115.921	&	7.279	&	195.458	&	237.882	&	23.354	&	3/3	&	0.89	&	0.06	&	0.39	&	0.19	\\
\instbox{S5}	&	52.875	&	275.042	&	7	&	466.72	&	98.40	&	99.77	&	2	&	302.602	&	10.02	&	54.152	&	59.225	&	2.82	&	268.003	&	281.271	&	9.298	&	7/7	&	0.97	&	0.14	&	2.36	&	-2.63	\\
\instbox{S6}	&	79.731	&	408.275	&	7	&	466.72	&	98.09	&	99.73	&	3	&	455.476	&	11.56	&	101.729	&	105.508	&	4.124	&	427.559	&	457.523	&	36.281	&	7/7	&	1.01	&	0.09	&	21.62	&	4.51	\\\hline
\multirow{2}{*}{\instbox{S7}}	&	151.936	&	411.598	&	5	&	93.39	&	0.00	&	2.19	&	1	&	529.535	&	28.65	&	\multirow{2}{*}{141.194}	&	\multirow{2}{*}{142.114}	&	\multirow{2}{*}{0.793}	&	\multirow{2}{*}{397.215}	&	\multirow{2}{*}{391.223}	&	\multirow{2}{*}{30.235}	&	\multirow{2}{*}{5/5}	&	\multirow{2}{*}{0.55}	&	\multirow{2}{*}{0.07}	&	-7.61	&	-3.62	\\
	&	141.756	&	371.003	&	5	&	411.00	&	24.13	&	84.83	&	2	&	517.192	&	39.40	&		&		&		&		&		&		&		&		&		&	-0.40	&	6.60	\\\hline
\instbox{S8}	&	135.093	&	462.105	&	5	&	411.00	&	23.98	&	54.39	&	2	&	637.363	&	37.93	&	152.2	&	153.34	&	1.097	&	376.127	&	390.819	&	20.245	&	5/5	&	1.06	&	0.11	&	11.24	&	-22.86	\\
\instbox{S9}	&	98.533	&	437.76	&	7	&	466.72	&	38.85	&	91.26	&	1	&	462.496	&	5.65	&	112.822	&	135.627	&	11.678	&	471.193	&	544.078	&	49.568	&	6/7	&	1.68	&	0.27	&	12.67	&	7.10	\\
\instbox{S10}	&	103.237	&	506.316	&	7	&	466.73	&	38.95	&	90.54	&	1	&	556.42	&	9.90	&	112.219	&	135.915	&	11.201	&	488.548	&	563.132	&	54.179	&	7/7	&	2.52	&	0.37	&	8.00	&	-3.64	\\
\instbox{S11}	&	69.263	&	342.179	&	6	&	441.00	&	11.60	&	67.27	&	2	&	390.568	&	14.14	&	84.931	&	85.962	&	0.697	&	370.338	&	386.504	&	24.154	&	6/6	&	0.99	&	0.09	&	18.45	&	7.60	\\
\instbox{S12}	&	128.183	&	588.457	&	7	&	466.73	&	39.31	&	91.32	&	1	&	622.294	&	5.75	&	134.628	&	162.127	&	13.757	&	602.437	&	654.272	&	52.118	&	7/7	&	2.2	&	0.45	&	4.79	&	2.32	\\
\instbox{M13}	&	84.412	&	366.032	&	6	&	158.28	&	0.00	&	12.21	&	1	&	428.96	&	17.19	&	84.738	&	87.592	&	3.199	&	398.048	&	438.734	&	30.055	&	6/6	&	1.22	&	0.1	&	0.38	&	8.04	\\
\instbox{M14}	&	87.179	&	371.954	&	6	&	182.38	&	0.00	&	31.27	&	1	&	508.285	&	36.65	&	89.631	&	101.862	&	4.489	&	486.238	&	469.631	&	32.32	&	6/6	&	1.5	&	0.14	&	2.74	&	23.50	\\
\instbox{M15}	&	68.467	&	254.914	&	4	&	61.46	&	0.00	&	25.27	&	1	&	257.391	&	0.97	&	69.116	&	121.11	&	13.901	&	200.21	&	330.504	&	52.413	&	4/4	&	1.13	&	0.19	&	0.94	&	-27.32	\\
\instbox{M16}	&	119.277	&	303.159	&	4	&	317.37	&	20.80	&	13.24	&	1	&	319.311	&	5.33	&	119.277	&	119.277	&	0	&	305.026	&	306.346	&	14.089	&	4/4	&	0.82	&	0.09	&	0.00	&	0.61	\\
\instbox{M17}	&	102.921	&	394.692	&	4	&	375.00	&	57.50	&	71.63	&	1	&	408.241	&	3.43	&	102.921	&	104.122	&	0.839	&	398.577	&	396.229	&	17.357	&	4/4	&	0.76	&	0.08	&	0.00	&	0.97	\\
\instbox{M18}	&	49.185	&	142.649	&	4	&	65.65	&	0.00	&	0.00	&	1	&	169.159	&	18.58	&	49.185	&	49.617	&	0.615	&	112.618	&	138.221	&	17.8	&	3/4	&	0.54	&	0.05	&	0.00	&	-26.67	\\
\instbox{M19}	&	88.749	&	357.989	&	6	&	96.74	&	0.00	&	0.00	&	1	&	362.771	&	1.34	&	96.711	&	105.443	&	3.777	&	379.025	&	402.532	&	15.582	&	6/6	&	1.96	&	0.2	&	8.23	&	5.55	\\
\instbox{M20}	&	136.408	&	548.69	&	6	&	365.45	&	0.01	&	17.06	&	1	&	549.547	&	0.16	&	136.459	&	150.172	&	4.194	&	555.35	&	588.425	&	17.661	&	6/6	&	1.65	&	0.18	&	0.04	&	1.20	\\
\instbox{M21}	&	40.954	&	234.258	&	6	&	2.49	&	0.00	&	0.00	&	1	&	234.258	&	0.00	&	44.579	&	64.411	&	9.103	&	185.333	&	227.243	&	23.471	&	6/6	&	0.6	&	0.07	&	8.13	&	-26.40	\\
\instbox{M22}	&	34.477	&	159.567	&	5	&	411.01	&	5.48	&	27.38	&	1	&	159.567	&	0.00	&	36.005	&	39.097	&	1.938	&	159.061	&	166.792	&	9.22	&	5/5	&	1.51	&	0.18	&	4.24	&	-0.32	\\
\instbox{M23}	&	98.71	&	428.758	&	6	&	265.37	&	0.00	&	34.83	&	1	&	428.758	&	0.00	&	100.327	&	110.737	&	6.608	&	521.344	&	537.761	&	35.245	&	6/6	&	2.08	&	0.18	&	1.61	&	17.76	\\
\instbox{M24}	&	32.71	&	145.474	&	6	&	267.36	&	0.00	&	29.46	&	1	&	184.859	&	27.07	&	34.47	&	39.575	&	1.744	&	158.88	&	168.197	&	6.51	&	6/6	&	2.35	&	0.19	&	5.11	&	8.44	\\
\instbox{M25}	&	88.017	&	357.433	&	6	&	264.32	&	0.00	&	33.33	&	1	&	359.41	&	0.55	&	93.039	&	108.06	&	6.502	&	378.379	&	410.372	&	19.501	&	6/6	&	3.2	&	0.27	&	5.40	&	5.54	\\
\hline
\end{tabular}
} % For Resize Box Closure
\end{table}
{\small The prefix in the ID indicates approximate (simplified) instance sizes (for all the Tables \ref{longTab: Small Dataset Overview}, \ref{longTab: Large Dataset Overview}, \ref{Tab: Small Integer Results}, \ref{Tab: Small Continuous Results}, \ref{longTab: Large Instance Results Integer} and \ref{longTab: Large Instance Results Continuous}):

$\instbox{S}\rightarrow Small \quad\quad \instbox{M}\rightarrow Medium \quad\quad \instbox{L}\rightarrow Large$}
\vspace*{\fill}  % Push down
\end{landscape}

\begin{landscape}
\vspace*{\fill}  % Push down
% \small
\subsection{Continuous Results of our Large Instances} \label{Continuous Results of our Large Instances}

Our developed large instances (overview of the same as may be found in \autoref{longTab: Large Dataset Overview}) were run both as continuous and integer problems. We discuss all integer results in the main paper, and the continuous problems in this supplementary paper.

For the Continuous Problems in \autoref{longTab: Large Instance Results Continuous}:
\begin{enumerate}
    \item In the instance \instbox{L32}: 20 out of 30 Instance runs were successful
    \item In the instance \instbox{L34}: 25 out of 30 Instance runs were successful
    \item In the instance \instbox{L38}: 30 out of 44 Instance runs were successful. By unsuccessful, we mean that the remaining 14 runs had every Main Iteration unable to satisfy the full resource requirement of the problem.
\end{enumerate}

% \small
% \footnotesize
% \scriptsize

\begin{longtable}{|m{0.02\linewidth}|m{0.085\linewidth}|m{0.07\linewidth}|m{0.095\linewidth}|m{0.05\linewidth}|m{0.08\linewidth}|m{0.12\linewidth}|m{0.1\linewidth}|m{0.065\linewidth}|m{0.06\linewidth}|m{0.085\linewidth}|}
\caption{\textbf{Computational Results for the Large Instances in \autoref{longTab: Large Dataset Overview} considering Continuous resource flow}\label{longTab: Large Instance Results Continuous}}\\
\hline
\multirow{2}{=}{\newline\textbf{ID}} & \multicolumn{4}{|c|}{\textbf{MILP Formulation}} & \multicolumn{5}{|c|}{\textbf{PSR-GIP Heuristic}} & \multirow{2}{=}{\newline\textbf{Relative Comparative Gap \% \newline\newline $\frac{b-a}{b}*100$}} \\
\cline{2-10}
& Objective Value (Cascade 1 only, if any) \newline\newline $(a)$	& Gap \% (as reported by Gurobi) & Levels Used (I$\rightarrow$Infeasible; \newline N$\rightarrow$NFSF; \newline F$\rightarrow$Feasible) \#& Time (mins) & Best Objective (among at least 30 runs) \newline\newline $(b)$ & Average (of the largest route duration among all Vehicles across successful Heuristic runs) & S.D. (of the vehicle with maximum route duration across successful Heuristic runs) & Average Heuristic RunTime (mins) & S.D. of Heuristic Runtimes (mins) & \\
\hline
\endfirsthead
\hline
\endhead
\hline
\instbox{S26}	&	468.118	&	0	&	F:3	&	0.04	&	468.118	&	468.118	&	0	&	0.31	&	0.01	&	0.00	\\
\instbox{S27}	&	645.164	&	0	&	F:3	&	63.06	&	645.164	&	657.154	&	29.116	&	0.49	&	0.05	&	0.00	\\
\instbox{M28}	&	NFSF	&	-	&	N:3,1	&	1440	&	993.032	&	1115.094	&	49.908	&	1.9	&	0.25	&	-	\\
\instbox{M29}	&	NFSF	&	-	&	N:1	&	1440	&	13985.943	&	16619.888	&	882.518	&	34.71	&	9.98	&	-	\\
\instbox{L30}	&	NFSF	&	-	&	I:6; N:8,7	&	1300	&	280.111	&	386.006	&	57.648	&	67.62	&	10.89	&	-	\\
\instbox{L31}	&	NFSF	&	-	&	N:1	&	1440	&	841.112	&	1250.869	&	294.932	&	21.51	&	9.61	&	-	\\
\instbox{L32}	&	NFSF	&	-	&	I:4; N:8,5	&	1440	&	989.239	&	1495.502	&	278.36	&	9.28	&	4.65	&	-	\\
\instbox{L33}	&	NFSF	&	-	&	I:4; N:5,10	&	1440	&	459.966	&	687.87	&	211.605	&	12.86	&	5.23	&	-	\\
\instbox{L34}	&	NFSF	&	-	&	I:4; N:5	&	1440	&	11814.237	&	16146.259	&	3345.709	&	12.29	&	5.39	&	-	\\
\instbox{L35}	&	NFSF	&	-	&	N:5,1	&	1440	&	1809.793	&	2952.262	&	1199.431	&	58.24	&	25.42	&	-	\\
\instbox{L36}	&	NFSF	&	-	&	I:3; N:4,5	&	1440	&	314.822	&	531.669	&	101.719	&	81.41	&	16.08	&	-	\\
\instbox{L37}	&	NFSF	&	-	&	I:1; N:2	&	OEAM	&	380.421	&	788.941	&	213.558	&	70.65	&	12.79	&	-	\\
\instbox{L38}	&	NFSF	&	-	&	I:2; N:3	&	OEAM	&	1203.092	&	1507.298	&	187.736	&	112.2	&	54.6	&	-	\\
\hline
% \newline
\end{longtable}
\quad NFSF $\Longrightarrow$ No Feasible Solution Found \quad\quad\quad OEAM $\Longrightarrow$ Optimization Exhausted Available Memory
\newline
\newline
\indent\# All Vehicles were provided with the same number of levels, and the result of those instances is briefed. For each of the values mentioned, separate 24 hour MILP runs were performed. The instance results in this table correspond to the last value in the Levels Used column. Further, for levels indicating Infeasibility (due to lack of solution space restricting multi-trips), all lower-Level values were also infeasible (since the solution space is even smaller then with much lesser number of multi-visits possible at any vertex).
\vspace*{\fill}  % Push down
\end{landscape}

\begin{landscape}
\vspace*{\fill}  % Push down
% \small
\subsection{Continuous Instance Runs of the Literature Dataset} \label{Continuous Instance Runs of the Literature Dataset}

This section contains the continuous problem results of instances from \cite{NUCAMENDIGUILLEN2021}. For both the cases of Continuous and Integer runs of this literature dataset (and essentially for all for MILP runs in our study), the Gurobi parameter of IntegralityFocus are always set to 1.
The table \autoref{longTab:Literature_Continuous_Instances} corresponds to case when the problems are run using continuous variables (for both the MILP and the Heuristic); all the Computations (for both the Exact Formulation and Heuristic Algorithm) were performed on the Computer C\_1. It is observed that here too, the relative gap is significantly low, and is generally lower than the gaps found for the continuous case; this indicates that the Heuristic is suitable for real-world optimization problems.

% \small
% \footnotesize
% \scriptsize
% \tiny

\begin{longtable}{|m{0.01\linewidth}|m{0.105\linewidth}|m{0.085\linewidth}|m{0.09\linewidth}|m{0.15\linewidth}|m{0.085\linewidth}|m{0.075\linewidth}|m{0.065\linewidth}|m{0.045\linewidth}|m{0.045\linewidth}|m{0.065\linewidth}|}
% Paper Details along with spatiotemporal details and Objectives of Papers reviewed
\caption{\textbf{Computational Results of the Instances in \autoref{longTab: Literature Dataset Overview} considering Continuous Resource Flows}\label{longTab:Literature_Continuous_Instances}}\\
\hline
\multirow{3}{=}{\rotatebox{85}{\textbf{Instance}}} & \multicolumn{3}{|c|}{\textbf{Cascaded Makespan Minimization}} & \multicolumn{6}{|c|}{\textbf{PSR-GIP Heuristic \#}}& \multirow{3}{=}{\newline\textbf{Relative Gap \% \newline\newline$\frac{(b-a)*100}{b}$}}\\
\cline{2-10}
& \multicolumn{3}{|c|}{\textbf{Our Exact MILP Formulation*}} & \multicolumn{4}{|c|}{\textbf{Cascaded Makespan Minimization Optimization}}& \multicolumn{2}{|c|}{\textbf{Runtime}}&\\

\cline{2-10}
& \textbf{Objective Function Value of Cascade 1 \newline ($a$)} & \textbf{No. of Vehicles used} &\textbf{MILP Gap \% (reported by Gurobi)}
& \textbf{Best Solution of Largest Route Duration among all Vehicles \newline ($b$)} & \textbf{No. of Vehicles used in Best Solutions} & \textbf{Average of Best Solutions} & \textbf{S.D. of Best Solutions} & \textbf{Avg. (mins)} & \textbf{S.D. (mins)} & \\
\hline

\hline
\endfirsthead

$	I_1	$	&	40.706	&	9	&	60.46	&	44.2	&	11	&	54.159	&	9.087	&	1.03	&	0.13	&	7.90	\\
$	I_2	$	&	38.77	&	10	&	58.49	&	38.773	&	11	&	48.297	&	5.678	&	0.92	&	0.1	&	0.01	\\
$	I_3	$	&	40.706	&	11	&	60.46	&	42.482	&	10	&	54.114	&	11.085	&	0.83	&	0.08	&	4.18	\\
$	I_4	$	&	40.706	&	10	&	60.46	&	40.713	&	8	&	53.959	&	10.428	&	0.73	&	0.06	&	0.02	\\
$	I_5	$	&	38.77	&	7	&	58.49	&	38.77	&	10	&	50.633	&	9.028	&	0.75	&	0.07	&	0.00	\\
$	I_6	$	&	39.793	&	9	&	59.55	&	42.48	&	11	&	53.718	&	11.445	&	0.81	&	0.08	&	6.33	\\
$	I_7	$	&	40.706	&	9	&	60.46	&	40.706	&	9	&	50.734	&	9.059	&	0.74	&	0.08	&	0.00	\\
$	I_8	$	&	40.706	&	7	&	60.46	&	41.235	&	10	&	55.669	&	10.712	&	0.77	&	0.08	&	1.28	\\
$	I_9	$	&	33.787	&	6	&	43.73	&	35.573	&	7	&	38.93	&	1.98	&	0.59	&	0.05	&	5.02	\\
$	I_{10}	$	&	40.706	&	7	&	60.46	&	40.934	&	9	&	51.604	&	7.046	&	0.78	&	0.06	&	0.56	\\
$	I_{11}	$	&	40.706	&	9	&	60.46	&	42.706	&	12	&	57.201	&	14.265	&	1	&	0.09	&	4.68	\\
$	I_{12}	$	&	40.706	&	6	&	60.46	&	40.706	&	9	&	49.679	&	6.883	&	0.72	&	0.06	&	0.00	\\
$	I_{13}	$	&	38.77	&	9	&	58.49	&	38.77	&	11	&	51.427	&	12.488	&	0.89	&	0.1	&	0.00	\\
$	I_{14}	$	&	40.706	&	8	&	60.46	&	40.713	&	10	&	56.295	&	13.578	&	0.96	&	0.09	&	0.02	\\
$	I_{15}	$	&	40.706	&	6	&	60.46	&	43.027	&	6	&	53.827	&	10.484	&	0.92	&	0.24	&	5.39	\\
\hline
\end{longtable}
Gurobi (non-default) settings for the MILP runs:
\noindent{
\begin{tabular}{|c|c|}
\hline
M & 1000 (Big M as used in Formulation) \\ \hline
MIPFocus & set at 2 (to focus more on proving optimality) \\ \hline
\end{tabular}
\newline
\par}

Heuristic Parameters used:
\noindent{
\begin{tabular}{|c|c|}
\hline
No. of Internal Main Iterations & 10 \\ \hline
Max. No. of Perturbation in each Main Iteration: & 1729 \\ \hline
Max. No. of Route Generation Logic Allowed from SRE: & 100 \\ \hline
\end{tabular}
\newline
\par}

* All MILP runs were provided with a runtime of 9999 seconds (only the first Cascades were run)

\# The Best, Average and S.D. are reported comparing 99 Heuristic runs for each Instance
\vspace*{\fill}  % Push down
\end{landscape}

\bibliographystyle{elsarticle-num}
\bibliography{references}